\pdfoutput=1
\PassOptionsToPackage{obeyspaces,spaces}{url}
\documentclass[a4paper,fleqn]{cas-sc}

\usepackage[authoryear]{natbib}
\usepackage{amsmath,amsfonts,amssymb}
\usepackage{bm}
\usepackage{xspace}

\usepackage{amsthm}
\makeatletter
\@ifpackageloaded{stix}{\let\mathds\mathbb}{\usepackage{dsfont}} 
\makeatother

\newtheorem{theorem}{Theorem}
\newtheorem{lemma}{Lemma}
\newtheorem{definition}{Definition}
\newtheorem{corollary}{Corollary}
\newtheorem{remark}{Remark}
\newtheorem{example}{Example}
\newtheorem{assumption}{Assumption}

\ExplSyntaxOn
\seq_new:N \g_paper_deferred_proofs_seq
\tl_new:N \g_paper_proof_placement_tl
\tl_new:N \l_paper_proof_placement_tl

\NewDocumentCommand{\setproofplacement}{m}
  {
    \str_case:nnF {#1}
      {
        {auto} { \tl_gclear:N \g_paper_proof_placement_tl }
        {inline} { \tl_gset:Nn \g_paper_proof_placement_tl {inline} }
        {appendix} { \tl_gset:Nn \g_paper_proof_placement_tl {appendix} }
      }
      {
        \PackageError{statementproof}
          {Unknown~global~proof~placement~`#1'}
          {Use~`auto',~`inline'~or~`appendix'.}
      }
  }

\NewDocumentEnvironment{statementproof}{ O{inline} m +b }
  {
    \tl_if_empty:NTF \g_paper_proof_placement_tl
      { \tl_set:Nn \l_paper_proof_placement_tl {#1} }
      { \tl_set_eq:NN \l_paper_proof_placement_tl \g_paper_proof_placement_tl }
    \str_case:VnF \l_paper_proof_placement_tl
      {
        {inline}
          { \begin{proof} #3 \end{proof} }
        {appendix}
          {
            \seq_gput_right:Nn \g_paper_deferred_proofs_seq
              { \begin{proof}[\bfseries Proof~of~\cref{#2}] #3 \end{proof} }
          }
      }
      {
        \PackageError{statementproof}
          {Unknown~proof~placement~`#1'}
          {Use~`inline'~or~`appendix'.}
      }
  }
  {}

\NewDocumentCommand{\printdeferredproofs}{}
  {
    \seq_if_empty:NF \g_paper_deferred_proofs_seq
      {
        \section{Proofs}\label{app:deferred_proofs}
        \seq_use:Nn \g_paper_deferred_proofs_seq { \par\medskip }
        \seq_gclear:N \g_paper_deferred_proofs_seq
      }
  }
\ExplSyntaxOff

\usepackage[obeyspaces,spaces]{url}
\makeatletter
\@ifpackageloaded{hyperref}{}{\usepackage[hidelinks]{hyperref}}
\makeatother
\usepackage{cleveref}
\crefname{assumption}{Assumption}{Assumptions}
\Crefname{assumption}{Assumption}{Assumptions}
\usepackage{autonum}
\numberwithin{equation}{section}

\makeatletter
\@ifclassloaded{cas-sc}{}{%
  \usepackage{titlesec}
  \titleformat{\paragraph}[runin]{\itshape\normalsize}{\theparagraph}{1em}{}}
\makeatother

\usepackage{booktabs}
\usepackage{colortbl}
\usepackage{multirow}

\usepackage{algorithm}
\usepackage{algpseudocode}

\usepackage{stackengine} 
\makeatletter
\@ifundefined{mathscr}{\usepackage{mathrsfs}}{} 
\makeatother
\newcommand\+[1]{\bm{#1}}
\makeatletter
\newcommand\doublehatsep{0.4}
\newsavebox\doublehat@inner
\newsavebox\doublehat@base
\newdimen\doublehat@raise
\newdimen\doublehat@inc
\def\doublehat@points#1{\strip@pt\dimexpr0.99626\dimexpr#1\relax\relax}
\newcommand\doublehat[1]{{\mathpalette\doublehat@{#1}}\widehat{#1}}
\newcommand\doublehat@[2]{%
  \sbox\doublehat@inner{$\m@th#1\widehat{#2}$}%
  \sbox\doublehat@base{$\m@th#1#2$}%
  \doublehat@inc=\dimexpr\ht\doublehat@inner-\ht\doublehat@base\relax
  \doublehat@raise=\doublehatsep\doublehat@inc
  \mathord{\rlap{\smash{%
      \pdfsave
      \pdfliteral{%
        \doublehat@points{-\wd\doublehat@inner}
        \doublehat@points{\ht\doublehat@base+\doublehat@raise+0.2\doublehat@inc}
        \doublehat@points{3\wd\doublehat@inner} \doublehat@points{\doublehat@inc+1pt} re W n}%
      \raise\doublehat@raise\copy\doublehat@inner
      \kern-\wd\doublehat@inner
      \pdfrestore}}}}
\makeatother
\newcommand\logical[1]{\doublehat{#1}}  
\newcommand\intermediate[1]{\widehat{#1}}      
\newcommand\intermediatecoord[1]{\intermediate{#1}} 
\newcommand\mat[1]{\+{\mathsf{#1}}}   
\newcommand\matcomp[1]{\mathsf{#1}}   
\newcommand\set[1]{\mathcal{#1}}      
\newcommand\funspace[1]{\set{#1}} 
\newcommand\bfunspace[1]{\set{#1}} 
\newcommand\fnal[1]{\mathscr{#1}}     

\newcommand\gdof{\fnal{S}} 

\newcommand\inprod[2]{\langle #1, #2 \rangle_{\funspace{L}_2}}

\newcommand\ellinprod[2]{\langle #1, #2 \rangle_{\ell_2}}

\newcommand\mappingc{\varphi}
\newcommand\mapping{\+\mappingc}
\newcommand\logicalmapping{\logical{\mapping}}           
\newcommand\intermediatemapping{\intermediate{\mapping}} 
\newcommand\pulledpolys[2]{\funspace{T}^{#1,#2}}
\newcommand\pushforward{\varPhi}

\newcommand\jacmat{\mat{J}}
\newcommand\jac{J}

\newcommand\indexone{\mathbb{I}}
\newcommand\indextwo{\mathbb{J}}
\newcommand\indexthree{\mathbb{K}}

\newcommand\bc{{0}}

\newcommand\polar[1]{\overline{#1}}
\newcommand\logicalpolar[1]{\vphantom{\logical{\polar{#1}}}\smash[t]{\logical{\polar{#1}}}\vphantom{\polar{#1}}}
\newcommand\intermediatepolar[1]{\vphantom{\intermediate{\polar{#1}}}\smash[t]{\intermediate{\polar{#1}}}\vphantom{\polar{#1}}}

\makeatletter
\@ifundefined{widecheck}{%
  \DeclareFontFamily{U}{mathx}{\hyphenchar\font45}
  \DeclareFontShape{U}{mathx}{m}{n}{<-> mathx10}{}
  \DeclareSymbolFont{mathx}{U}{mathx}{m}{n}
  \DeclareMathAccent{\widecheck}{0}{mathx}{"71}}{}
\makeatother
\newcommand\reduced[1]{\widetilde{#1}}
\newcommand\enriched[1]{\widecheck{#1}}

\newcommand{\pr}{{\mathfrak{r}}}

\newcommand\bsplineone{u}
\newcommand\Bsplineone{U}

\newcommand\bspline{v}
\newcommand\Bspline{V}

\newcommand\Bsplinethree{W}

\newcommand\defeq{\mathrel{\mathop:}=}
\newcommand\setdef[2]{\left\lbrace #1 \; : \; #2 \right\rbrace}

\newcommand\gradient{\nabla}

\DeclareMathOperator{\divergencet}{div}
\DeclareMathOperator{\gradientt}{grad}
\DeclareMathOperator{\curlt}{curl}

\newcommand\gradienth{\mathfrak{G}}
\newcommand\curlh{\mathfrak{R}}

\newcommand\internal\odot
\newcommand\external\otimes

\DeclareMathOperator\kernel{ker}
\DeclareMathOperator\image{im}

\newcommand\abs[1]{| #1 |}

\newcommand\rpar[1]{( #1 )}
\newcommand\srpar[1]{\left( #1 \right)}

\newcommand\sspar[1]{\left[ #1 \right]}

\newcommand\dd{\mathrm{d}}

\newcommand\xx{\mathrm{x}}
\newcommand\yy{\mathrm{y}}
\newcommand\zz{\mathrm{z}}
\newcommand\diff{\, \dd}

\newcommand\pard[2]{\frac{\partial #1}{\partial #2}}

\newcommand\totd[2]{\frac{\dd #1}{\dd #2}}

\newcommand\lhopital{L'Hôpital's rule\xspace}

\DeclareMathOperator{\spn}{span}

\makeatletter
\newcommand*{\transpose}{%
  {\mathpalette\@transpose{}}%
}
\newcommand*{\@transpose}[2]{%
  \raisebox{\depth}{$\m@th#1\intercal$}%
}
\makeatother

\usepackage{graphicx}
\usepackage{color}
\usepackage{tikz,pgfplots}
\pgfplotsset{compat=1.18}
\usepackage{subcaption}
\usetikzlibrary{arrows.meta, angles}

\def\twofigwidth{0.475\textwidth}

\def\axheight{0.4\textwidth}

\def\oneaxwidth{\textwidth}
\def\twoaxwidth{0.44\textwidth}

\definecolor{lines1}{rgb}{0.1216, 0.4667, 0.7059}
\definecolor{lines2}{rgb}{1.0000, 0.4980, 0.0549}
\definecolor{lines3}{rgb}{0.1725, 0.6275, 0.1725}
\definecolor{lines4}{rgb}{0.8392, 0.1529, 0.1569}
\definecolor{lines5}{rgb}{0.5804, 0.4039, 0.7412}
\definecolor{lines6}{rgb}{0.5490, 0.3373, 0.2941}
\definecolor{lines7}{rgb}{0.8902, 0.4667, 0.7608}
\definecolor{lines8}{rgb}{0.4980, 0.4980, 0.4980}
\definecolor{lines9}{rgb}{0.7373, 0.7412, 0.1333}
\definecolor{lines10}{rgb}{0.0902, 0.7451, 0.8118}

\newcommand\polarfeec{\href{https://gitlab.mpcdf.mpg.de/polarplasma/polarfeec.jl}{PolarFEEC.jl}}
\newif\ifshowcodeapi
\showcodeapitrue
\newcommand\codeapiurl{https://polarplasma.pages.mpcdf.de/polarfeec.jl/api/}
\colorlet{codeapicolor}{lines1!65!black}
\def\codeapifilter#1#2\relax{\ifnum\uccode`#1=`#1 type\else method\fi}
\newcommand\codeapi[3][]{\unskip\ifshowcodeapi\space\href{\codeapiurl?q=\detokenize{#2&filter=}\ifx\relax#1\relax\codeapifilter#2\relax\else#1\fi}{\textcolor{codeapicolor}{[\expandafter\nolinkurl\expandafter{\detokenize{#2#3}}]}}\fi}

\usetikzlibrary{external}
\usepackage{tikz-cd}
\usepackage{derham}
\derhamsetup{row-sep=normal,column-sep=normal} 
\AtBeginEnvironment{tikzcd}{\tikzexternaldisable}
\AtEndEnvironment{tikzcd}{\tikzexternalenable}

\tikzcdset{
  arrow style=tikz,
  diagrams={>={Latex[open,round]}}
}

\newcommand\inputtikz[1]{\tikzsetnextfilename{#1}\input{tikz/#1.tex}}

\providecommand{\diskMapAmp}{0.1} 
\providecommand{\diskMapPer}{3}   
\providecommand{\diskGridNx}{8}
\providecommand{\diskGridNy}{16}
\providecommand{\bananaKappa}{1}            
\providecommand{\bananaAngle}{60}           
\providecommand{\bananaScale}{0.894427191}  

\usepackage{natbib}

\setproofplacement{appendix}

\begin{document}
  \let\WriteBookmarks\relax
  \def\floatpagepagefraction{1}
  \def\textpagefraction{.001}

  \shorttitle{High-order FEEC on domains with a disk-like coordinate singularity}
  \shortauthors{R. Remmerswaal}

  \title[mode=title]{High-order FEEC discretization on 2D and 3D domains with a disk-like coordinate singularity}

  \author[1]{Ronald A. Remmerswaal}[orcid=0000-0003-4394-2615]
  \cormark[1]
  \ead{rwr@ipp.mpg.de}
  \credit{Conceptualization, Methodology, Formal analysis, Software, Investigation, Visualization, Writing -- original draft, Writing -- review \& editing}
  \affiliation[1]{organization={Max Planck Institute for Plasma Physics},
    addressline={Boltzmannstra\ss e 2},
    postcode={85748},
    city={Garching},
    country={Germany}}
  \cortext[1]{Corresponding author}

  \begin{abstract}[]
    Tensor-product splines mapped from a square to a disk-like domain form a discrete de Rham complex, but at the coordinate singularity the pushed-forward $0$-, $1$- and $2$-forms are discontinuous, unbounded and not square integrable.
Existing polar constructions impose either high-order regularity on the $0$-forms only, or first-order regularity on the whole complex.
For general disk-like mappings composed with a smooth diffeomorphism, we construct polar subspaces of all three tensor-product spaces that are conforming, form a cochain complex with the cohomology of the disk, and are $\funspace{C}^{\pr-m}$ at the center for $m$-forms, for any $\pr$ up to the spline degree, under a compatibility condition between the mapping and the angular basis.
The construction enriches the $2$-form space with projected polynomials and lifts them to curl-like $1$-forms by Poincar\'e-type homotopy operators, which reduce to cumulative sums of radial B-spline coefficients, so that the homotopy identities hold exactly.
The same operators yield sparse commuting projections that preserve total flux and circulation, and boundary-compatible homotopies that solve the curl--curl problem with a single mass-matrix inversion.
In 3D, the tensor product with a periodic 1D complex carries all properties over to stellarator-like solid tori.
Numerical experiments confirm optimal projection accuracy and the predicted regularity, show that the regularity conditions suppress the stiff central modes of the curl--curl operator, allowing more than ten times larger time steps than first-order polar spaces, and that compatible high-order regularity is required for convergent guiding-center orbits near the axis of a discrete vacuum field.

  \end{abstract}

  \begin{keywords}
    Finite element exterior calculus \sep
    Coordinate singularity \sep
    Polar splines \sep
    High-order regularity \sep
    Discrete de Rham complex \sep
    Commuting projections
    \MSC[2020]{65N30, 65D07, 58A12}
  \end{keywords}

  \maketitle

  \section{Introduction}\label{sec:introduction}
Structure-preserving discretizations of the de Rham complex, known as finite element exterior calculus (FEEC) \citep{Arnold_Falk_Winther_2006,arnold2010feec,arnold2018feec}, approximate potentials, forces, fluxes and densities by discrete $0$-, $1$-, $2$- and $3$-forms whose spaces are linked by the gradient, curl and divergence in the same way as their continuous counterparts.
Together with projection operators that commute with these differential operators, the discrete complex inherits the exactness of the continuous one: discrete curl-free fields are gradients and discrete divergence-free fields are curls, so that constraints such as $\divergencet\+B=0$ and charge conservation hold exactly and the spurious modes of naive discretizations of Maxwell's equations are absent.
Tensor-product B-splines provide a particularly convenient realization of this framework on mapped domains \citep{buffa2010,buffa_rivas_sangalli_vazquez2011,ratnani2012}: the de Rham complex on the unit square or cube is the tensor product of one-dimensional spline complexes, its geometric degrees of freedom are point values and edge, face and cell integrals, commonly on the Greville grid, and the commuting projections are products of one-dimensional interpolation and histopolation \citep{hiemstra2014,yaman2025}.
The high order, the smoothness and the compact support of B-splines have made them a common choice in geometric particle-in-cell methods \citep{kraus2017gempic} and in structure-preserving plasma codes \citep{holderied2021,possanner2023struphy}.

\begin{figure}
  \centering
  \captionsetup[sub]{justification=raggedright}
  \subcaptionbox{The magnetic field strength $|\+B|$ on the boundary and on two poloidal cross-sections.\label{fig:intro_w7x_3d}}
  [0.7\textwidth]{
    \tikzsetnextfilename{w7x_hero}\input{tikz/w7x_hero.tex}
  }\hfill
  \subcaptionbox{Poloidal cross-section with a disk-like coordinate singularity.\label{fig:intro_w7x_cut}}
  [0.27\textwidth]{
    \tikzsetnextfilename{w7x_hero_cut}\input{tikz/w7x_hero_cut.tex}
  }
  \caption{Magnetohydrodynamic (MHD) equilibrium of the stellarator Wendelstein 7-X, computed with GVEC \citep{hindenlang2025gvec} and discretized in straight-field-line (PEST) coordinates \citep{grimm1976pest}: the magnetic field $\+B$ is projected into the polar $2$-form space $\polar{\funspace{\Bsplinethree}}^{2,4}$ of \cref{sec:derham_hodge_3d_construction} by the commuting projection $\polar{\pi}^{2,4}$ of \cref{sec:derham_hodge_3d_commuting_projections}.
  By \cref{thm:polar_subspaces_regularity_3d}, the projected field is twice continuously differentiable at the coordinate singularity.
  The grid shows the PEST coordinate lines, which converge on the magnetic axis (dashed), and the translucent shell is the cut-away part of the boundary.
  The equilibrium is computed from the input parameters of the W7-X example in the test suite of GVEC (\texttt{test-CI/examples/w7x}, \citealp{gvec_code}), using GVEC's flexible coordinate frame \citep{hindenlang2025gvec}.
  }
  \label{fig:intro_w7x}
\end{figure}
In magnetically confined plasmas the computational toroidal domain has disk-like cross-sections as illustrated in \cref{fig:intro_w7x}, and the natural flux-aligned coordinates are polar coordinates centered on the axis \citep{bouzat2018,zoni2019}.
A tensor-product spline space on the unit square therefore reaches the physical domain through a mapping that collapses one edge of the square to the center, see \cref{fig:disk_mapping}.
The mapping is singular there, and the pushed-forward basis functions are nonconforming: physical $0$-forms are discontinuous at the center, $1$-forms grow like the inverse distance to it, and $2$-forms are not square integrable (\cref{tab:disc_basis_functions_fourier} in \cref{sec:physical_function_spaces} shows the $0$-forms).
If the singularity is ignored, the discrete operators acquire stiff modes originating at the center \citep{holderied2021,jiang2026polar}, which limit the time step of explicit schemes and distort fields and particle orbits near the axis. 
The conformity and regularity of spline spaces under such singular parameterizations have been analyzed by \citet{takacs2011,takacs2012}.

Two lines of work address the singularity by restricting the tensor-product spaces to \emph{polar} subspaces, see \cref{tab:literature_map}.
The first concerns scalar functions, that is, the $0$-form space.
\citet{toshniwal2017} impose regularity at the center through explicit constraints on the spline coefficients.
Their construction applies to any order of smoothness and is made explicit up to $\funspace{C}^2$, but the regularity is attained only if the pullbacks of polynomials of the corresponding degree are exactly representable in the logical spline space (\cref{def:pr_compatibility_mapping_logical_zeroform}), and such spaces are provided only up to $\funspace{C}^1$.
These $\funspace{C}^1$ polar splines have been used, for instance, for semi-Lagrangian guiding-center simulations on mapped disks \citep{zoni2019}.
\citet{jiang2026polar} reach high-order regularity by replacing the innermost coefficient rings with projected polynomials. 
With periodic B-splines in the angular direction, the same representability condition holds only for continuity, so that higher-order regularity is attained only in the limit of angular refinement.
The second line treats the whole de Rham complex, but with first-order regularity only: the polar $1$- and $2$-form spaces are built from the gradients and curls that the cochain property requires, either by extraction \citep{toshniwal2021}, extended to solid tori with the cohomology of the torus by \citet{patrizi2025torus} and to hierarchically refinable polar spaces by \citet{cabanas2026hierarchical}, or through commuting projections based on geometric degrees of freedom \citep{holderied2022}, as used in the hybrid code of \citet{holderied2021}, or through commuting projections that compose the tensor-product projections with local conforming projections on the spline coefficients \citep{yaman2025}.
What is missing is a polar complex that combines high-order regularity for all form degrees with sparse commuting projections.
If the polar $1$-forms near the center are only the gradients of the polar $0$-forms, then no polar $2$-form can be nonzero at the center, because the curl of a gradient vanishes (\cref{sec:homotopy_motivation}).
High-order regularity of the $2$-forms, and hence of the densities and fluxes that physical models carry there, requires polar $1$-forms with a nonzero curl at the center, and these cannot be obtained from the cochain property alone.

\begin{table}
  \centering
  \footnotesize
  \setlength{\tabcolsep}{2.5pt}
  \begin{tabular}{@{}
    >{\raggedright\arraybackslash}p{2.2cm}
    >{\raggedright\arraybackslash}p{1.2cm}
    >{\raggedright\arraybackslash}p{1.8cm}
    >{\raggedright\arraybackslash}p{2.7cm}
    >{\raggedright\arraybackslash}p{1.5cm}
    >{\raggedright\arraybackslash}p{2.4cm}
    >{\raggedright\arraybackslash}p{1.9cm}
    @{}}
    \toprule
    Reference & Forms & Regularity of $0$-forms & Mechanism & Angular basis & Geometry & Commuting projections\\
    \midrule
    \citet{toshniwal2017} & $0$ & $\funspace{C}^1, \funspace{C}^2$\textsuperscript{a} & explicit coefficient constraints (extraction) & B-splines & 2D, disk-like & --\\
    \citet{zoni2019} & $0$ & $\funspace{C}^1$ & extraction (after \citealp{toshniwal2017}) & B-splines & 2D, disk-like & --\\
    \citet{jiang2026polar} & $0$ & $\funspace{C}^{\pr}$\textsuperscript{b} & projected polynomials (extraction) & B-splines & 2D, unit disk (applied in a 3D PIC code) & --\\
    \midrule
    \citet{toshniwal2021} & $0,1,2$ & $\funspace{C}^1$ & extraction; $1$-forms from gradients of polar $0$-forms & B-splines & 2D, polar surfaces (one or two poles) & --\\
    \citet{holderied2022} & $0,1,2,3$ & $\funspace{C}^1$ & commuting projections  & B-splines & 3D torus with disk-like cross-section & yes (GDOFs)\\
    \citet{yaman2025} & $0,1,2$ & $\funspace{C}^0$, $\funspace{C}^1$ & commuting projections & B-splines & 2D, disk-like & yes (local, acting on the innermost rings)\\
    \citet{patrizi2025torus} & $0,1,2,3$ & $\funspace{C}^1$ & extraction; $1$-forms from gradients of polar $0$-forms & B-splines & 3D solid torus from a polar map & --\\
    \midrule
    This work & $0,1,2,3$ & $\funspace{C}^{\pr}$; $m$-forms are $\funspace{C}^{\pr-\min(m, 2)}$ & projected polynomials for $0$- and $2$-forms, homotopy lifting for $1$-forms (extraction and projection) & B-splines or trigonometric & 2D and 3D, disk-like composed with a diffeomorphism & yes (local, acting on innermost rings); also commute with the homotopy operators\\
    \bottomrule
  \end{tabular}
  \caption{Polar spline constructions in the literature and in this work.
  `Forms' lists the degrees $m$ of the $m$-form spaces that are made regular at the coordinate singularity, `regularity' refers to the $0$-form space at the center, and `commuting projections' indicates whether projection operators onto the polar spaces are provided that commute with the differential operators.
  Here $\pr$ denotes the regularity degree introduced in \cref{sec:zeroforms}.
  The stated regularity requires the logical mapping and the logical $0$-form space to be $\pr$-compatible (\cref{def:pr_compatibility_mapping_logical_zeroform}):
  \textsuperscript{a}\,no $2$-compatible spaces are provided, so that $\funspace{C}^2$ is not attained by the spaces considered there;
  \textsuperscript{b}\,the periodic B-splines in the angular direction are not $\pr$-compatible for $\pr \ge 1$, so that $\funspace{C}^\pr$ is then attained only in the limit of angular refinement.}
  \label{tab:literature_map}
\end{table}

This paper supplies these missing $1$-forms.
We enrich the polar $2$-form space with projected polynomial $2$-forms and lift them to curl-like polar $1$-forms with the homotopy operators of the de Rham complex, that is, the Poincar\'e operators that integrate along rays from the center \citep[Sec.~2]{gopalakrishnan2004}; their vector-valued member is the vector potential in the Poincar\'e gauge of electrodynamics \citep{brittin1982}.
The decisive observation is that on clamped radial B-splines these operators are cumulative sums of coefficients, so that the homotopy identities, which express every form as the derivative of its homotopy plus the homotopy of its derivative up to the cohomology, hold exactly on the discrete level and no projection is needed after the lift.
As a consequence, the polar subspaces form a cochain complex under the differential operators and a chain complex under the homotopy operators, the polar $1$-form space splits into gradients and homotopy lifts, exactness follows by an algebraic argument, and the $1$-form projection can be defined through the homotopy identity, which makes its commutation with the differential operators automatic and keeps all projections sparse.
Two further ingredients make the construction applicable beyond the unit disk.
The physical mapping is decomposed into a disk-like logical mapping and a smooth diffeomorphism, as also proposed in \citet[Rem.~2.2]{patrizi2025torus}; every operation at the singularity involves only the former, and the geometric complexity of the domain is pushed into the latter, which we state as a design rule in \cref{rem:pr_compatibility_sufficient}.
Regularity at the center is then guaranteed under a compatibility condition between the logical mapping and the angular basis, which the unit disk with a trigonometric angular basis satisfies for every order, whereas periodic B-splines in the angular direction satisfy it to first order only; the numerical experiments show which properties survive without compatibility.

Specifically, the paper makes the following contributions.
\begin{itemize}
  \item Polar subspaces of the 2D tensor-product $0$-, $1$- and $2$-form spaces for general disk-like mappings composed with a diffeomorphism, which are conforming, form a cochain complex with the cohomology of the disk, and whose $m$-forms are $\funspace{C}^{\pr-m}$ at the center for every regularity degree $\pr$ up to the spline degree, under $\pr$-compatibility of the mapping with the logical spaces (\cref{sec:zeroforms,sec:derham_hodge_2d_construction,sec:derham_hodge_2d_properties}).
  \item Exact discrete homotopy operators and homotopy identities on the tensor-product and polar spaces, and the resulting chain complex and splitting of the polar $1$-form space (\cref{sec:derham_hodge_2d_hodge,sec:derham_hodge_2d_properties}).
  \item Sparse commuting projections onto the polar spaces, defined on the coefficient level, which commute with the differential and with the homotopy operators and preserve the total flux and circulation (\cref{sec:derham_hodge_2d_geom_projections}).
  \item Boundary-compatible homotopy operators for homogeneous essential boundary conditions, which reduce the curl--curl problem to a single inversion of the $2$-form mass matrix (\cref{sec:boundary_homotopy}).
  \item The 3D polar complex on solid tori as the tensor product of the 2D polar complex with a periodic 1D complex, which inherits conformity, regularity along the central curve, the cohomology of the solid torus, the homotopy operators and the commuting projections (\cref{sec:high_order_derham_hodge_3d}).
  \item Numerical evidence (\cref{sec:results}) for optimal projection accuracy and for the sharpness of the regularity statements, for the suppression of the stiff modes of the curl--curl operator originating at the center, with a time-step gain of about $11$ over the first-order polar spaces, for the accurate propagation of under-resolved waves across the axis of a stellarator-like domain also without compatibility, and for the necessity of compatible high-order regularity to obtain convergent guiding-center orbits near the axis of a discrete vacuum field.
\end{itemize}

The paper is organized as follows.
\Cref{sec:notation_2d} fixes the two-step mapping, the one-dimensional spline and trigonometric complexes, their tensor products, the pushforwards and the geometric degrees of freedom; the notation is collected in \cref{app:notation}.
\Cref{sec:zeroforms} reviews the polar $0$-form constructions of the literature, introduces the splitting of the coefficients into internal and external rings and the compatibility condition, and proves conformity and regularity of the polar $0$-form space.
\Cref{sec:derham_2d} contains the 2D construction: the first-order polar complexes of the literature and the reason why their direct high-order generalization is too poor (\cref{sec:first_order_regularity_de_rham,sec:homotopy_motivation}), the homotopy operators (\cref{sec:derham_hodge_2d_hodge}), the enriched polar spaces and their bases (\cref{sec:derham_hodge_2d_construction}), their properties (\cref{sec:derham_hodge_2d_properties}), the commuting projections (\cref{sec:derham_hodge_2d_geom_projections}), and essential boundary conditions (\cref{sec:boundary_homotopy}).
\Cref{sec:high_order_derham_hodge_3d} extends the construction to 3D, \cref{sec:results} presents the numerical experiments, and \cref{sec:conclusion} concludes.
All proofs are deferred to \cref{app:deferred_proofs}, and \cref{sec:matrix_representation} gives the matrix representation of all spaces and operators.
A reader primarily interested in the implementation may skim \cref{sec:notation_2d}, read \cref{sec:derham_hodge_2d_construction,sec:matrix_representation}, and turn to \cref{sec:results}; the mathematical development is contained in \cref{sec:zeroforms,sec:derham_2d,sec:high_order_derham_hodge_3d} together with \cref{app:deferred_proofs}.

  \section{Preliminaries}\label{sec:notation_2d}
We present the 2D mapping, 1D function spaces, 2D logical function spaces as a tensor product of 1D function spaces, and 2D physical function spaces by means of the pushforwards.
The notation used throughout the paper is collected in \cref{app:notation}, together with the sections in which it is introduced.

\subsection{Coordinate transformation}\label{sec:coordtransform}
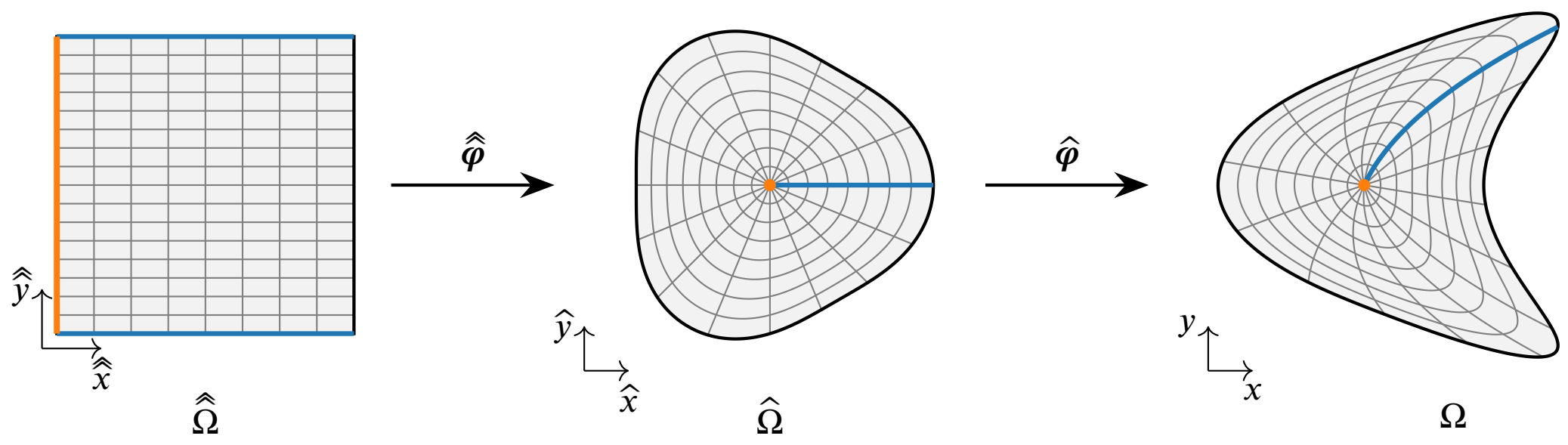
\begin{figure}
  \centering
  \tikzsetnextfilename{disk_mapping}  
  \providecommand{\diskMapAmp}{0.1} 
  \providecommand{\diskMapPer}{3}   
  \providecommand{\diskGridNx}{8}   
  \providecommand{\diskGridNy}{16}  
  \providecommand{\bananaKappa}{1} 
  \providecommand{\bananaAngle}{60}  
  \providecommand{\bananaScale}{0.894427191} 
  \begin{tikzpicture}[
      x=0.079\textwidth, y=0.079\textwidth,
      declare function={
        Rf(\r,\t) = \r*(1 + \diskMapAmp*\r*cos(360*\diskMapPer*\t));
        dx(\r,\t) = Rf(\r,\t)*cos(360*\t);
        dy(\r,\t) = Rf(\r,\t)*sin(360*\t);
        bu(\r,\t) = cos(\bananaAngle)*dx(\r,\t) - sin(\bananaAngle)*dy(\r,\t);
        bv(\r,\t) = sin(\bananaAngle)*dx(\r,\t) + cos(\bananaAngle)*dy(\r,\t);
        bx(\r,\t) = \bananaScale*(bu(\r,\t) + \bananaKappa*bv(\r,\t)^2);
        by(\r,\t) = bv(\r,\t)/\bananaScale;
      },
      grid/.style={gray, thin},
      seam/.style={lines1, line width=1.2pt},
      maparrow/.style={-{Stealth[scale=1.3]}, line width=0.8pt},
    ]
    \pgfmathsetmacro{\gridNx}{\diskGridNx}
    \pgfmathsetmacro{\gridNy}{\diskGridNy}

    \begin{scope}[shift={(0,-1)}, scale=2]
      \fill[black!5] (0,0) rectangle (1,1);
      \foreach \i in {0,...,\diskGridNx} {
        \draw[grid] (\i/\gridNx, 0) -- (\i/\gridNx, 1);
      }
      \foreach \j in {0,...,\diskGridNy} {
        \draw[grid] (0, \j/\gridNy) -- (1, \j/\gridNy);
      }
      \draw[thick] (0,0) rectangle (1,1);
      \draw[seam] (0,0) -- (1,0);
      \draw[seam] (0,1) -- (1,1);
      \draw[lines2, line width=1.5pt] (0,0) -- (0,1);
      \draw[->, thin] (-0.05,-0.05) -- (0.15,-0.05) node[below] {$\logical{x}$};
      \draw[->, thin] (-0.05,-0.05) -- (-0.05,0.15) node[left] {$\logical{y}$};
    \end{scope}
    \node at (1,-1.55) {$\logical{\Omega}$};

    \draw[maparrow] (2.25,0) -- (3.35,0) node[midway, above] {$\logicalmapping$};

    \begin{scope}[shift={(4.8,0)}]
      \fill[black!5] plot[smooth cycle, samples=200, variable=\t, domain=0:1] ({dx(1,\t)}, {dy(1,\t)});
      \foreach \i in {1,...,\diskGridNx} {
        \draw[grid] plot[smooth, samples=100, variable=\t, domain=0:1]
          ({dx(\i/\gridNx,\t)}, {dy(\i/\gridNx,\t)});
      }
      \foreach \j in {1,...,\diskGridNy} {
        \draw[grid] plot[smooth, samples=30, variable=\s, domain=0:1]
          ({dx(\s,\j/\gridNy)}, {dy(\s,\j/\gridNy)});
      }
      \draw[thick] plot[smooth cycle, samples=200, variable=\t, domain=0:1] ({dx(1,\t)}, {dy(1,\t)});
      \draw[seam] plot[smooth, samples=30, variable=\s, domain=0:1] ({dx(\s,0)}, {dy(\s,0)});
      \fill[lines2] (0,0) circle[radius=1.5pt];
      \draw[->, thin] (-1.25,-1.25) -- (-0.95,-1.25) node[below] {$\intermediatecoord{x}$};
      \draw[->, thin] (-1.25,-1.25) -- (-1.25,-0.95) node[left] {$\intermediatecoord{y}$};
      \node at (0,-1.55) {$\intermediate{\Omega}$};
    \end{scope}

    \draw[maparrow] (6.25,0) -- (7.35,0) node[midway, above] {$\intermediatemapping$};

    \begin{scope}[shift={(8.8,0)}]
      \fill[black!5] plot[smooth cycle, samples=200, variable=\t, domain=0:1] ({bx(1,\t)}, {by(1,\t)});
      \foreach \i in {1,...,\diskGridNx} {
        \draw[grid] plot[smooth, samples=100, variable=\t, domain=0:1]
          ({bx(\i/\gridNx,\t)}, {by(\i/\gridNx,\t)});
      }
      \foreach \j in {1,...,\diskGridNy} {
        \draw[grid] plot[smooth, samples=30, variable=\s, domain=0:1]
          ({bx(\s,\j/\gridNy)}, {by(\s,\j/\gridNy)});
      }
      \draw[thick] plot[smooth cycle, samples=200, variable=\t, domain=0:1] ({bx(1,\t)}, {by(1,\t)});
      \draw[seam] plot[smooth, samples=30, variable=\s, domain=0:1] ({bx(\s,0)}, {by(\s,0)});
      \fill[lines2] (0,0) circle[radius=1.5pt];
      \draw[->, thin] (-1.05,-1.25) -- (-0.75,-1.25) node[below] {$x$};
      \draw[->, thin] (-1.05,-1.25) -- (-1.05,-0.95) node[left] {$y$};
      \node at (0.6,-1.55) {$\Omega$};
    \end{scope}
  \end{tikzpicture}
  \caption{The logical domain $\logical{\Omega}$ (left) is mapped by the disk-like logical mapping $\logicalmapping$ of \cref{eqn:disk_mapping}, with $R = \logical{x} [1 + \diskMapAmp\, \logical{x}\, \cos(\diskMapPer \times 2\pi \logical{y})]$ and $(C, S) = (\cos(2\pi\logical{y}), \sin(2\pi\logical{y}))$, to the intermediate disk-like domain $\intermediate{\Omega}$ (center), which is mapped by the banana map $\intermediatemapping$ of \cref{ex:banana_mapping} with $\kappa = \bananaKappa$ and $\alpha = 2/\sqrt{5}$ to the physical domain $\Omega$ (right); the physical mapping is $\mapping = \intermediatemapping \circ \logicalmapping$.
  The logical grid has $\diskGridNx\times\diskGridNy$ subintervals.
  The left edge of the logical domain (orange) collapses to the center of the intermediate and physical domains, and the identified edges $\logical{y}=0$ and $\logical{y}=1$ (blue) form the periodic seam, which is mapped to the blue curve.}
  \label{fig:disk_mapping}
\end{figure}
The logical domain $\logical{\Omega} \defeq [0, 1]\times\mathbb{S}^1$, where $\mathbb{S}^1\defeq\mathbb{R}/\mathbb{Z}$ denotes the unit interval with its endpoints identified, is mapped to the physical domain $\Omega$ in two steps.
The two steps are kept separate on purpose.
The first step, the disk-like logical mapping, carries the coordinate singularity and is the only mapping that enters the construction of the polar spaces in \cref{sec:zeroforms,sec:derham_2d,sec:high_order_derham_hodge_3d}; in particular, the extraction matrices of the polar spaces depend on it alone (\cref{sec:derham_hodge_2d_basis}).
The second step is a smooth diffeomorphism that carries the geometric complexity of the physical domain and merely transports the polar spaces to it.
This division also lowers the requirements on the discrete spaces: the regularity of the polar spaces at the center is guaranteed under a compatibility condition between the logical mapping and the logical spaces (\cref{def:pr_compatibility_mapping_logical_zeroform}), which is easiest to satisfy when the logical mapping is as simple as possible, for example the unit disk, and all shaping is left to the diffeomorphism.
\Cref{rem:pr_compatibility_sufficient} gives an example in which the same physical mapping is compatible for one decomposition and not for another, and states this as a design rule.

First, the disk-like logical mapping $\logicalmapping$ maps $\logical{\Omega}$ to the intermediate domain $\intermediate{\Omega} \defeq \logicalmapping(\logical{\Omega})$,
  \begin{equation}\label{eqn:disk_mapping}
    \logicalmapping(\logical{\+x}) = \begin{bmatrix}
      R(\logical{x}, \logical{y}) C(\logical{y})\\
      R(\logical{x}, \logical{y}) S(\logical{y})
    \end{bmatrix},
  \end{equation}
  where $\logical{x}$ is the radial coordinate and $\logical{y}$ is the angular coordinate, $C(\logical{y})$ and $S(\logical{y})$ are angular functions, and $R(\logical{x}, \logical{y})$ is the radius function.
  The coordinates on $\intermediate{\Omega}$ are denoted by $\intermediatecoord{\+x} = (\intermediatecoord{x}, \intermediatecoord{y})$.
  See also \cref{fig:disk_mapping} for an illustration of the mappings.
  The following assumption is made on the logical mapping $\logicalmapping$.
\begin{assumption}[Disk-like mapping]\label{ass:disc_like_mapping}
  The functions $R, C, S$ are periodic in the angular coordinate, and satisfy
  \begin{equation}\label{eqn:coordinate_transformation_assumptions}
    \frac{R(\logical{x}, \logical{y})}{\logical{x}} \in \funspace{C}^\infty
    , \quad
    \lim_{\logical{x} \to 0} \frac{R(\logical{x}, \logical{y})}{\logical{x}} = R_1 > 0
    , \quad
    C S' - S C' > 0,
  \end{equation}
where $\funspace{C}^\infty$ denotes the space of infinitely differentiable functions and $R_1$ is a positive constant.
Moreover, the mapping is a smooth diffeomorphism away from the coordinate singularity at $\logical{x} = 0$.
\end{assumption}
Second, the intermediate mapping $\intermediatemapping$ maps $\intermediate{\Omega}$ to the physical domain $\Omega \defeq \intermediatemapping(\intermediate{\Omega})$.
\begin{assumption}[Diffeomorphism]\label{ass:diffeomorphism}
  The intermediate mapping $\intermediatemapping$ is a $\funspace{C}^\infty$ diffeomorphism on $\intermediate{\Omega}$, it is orientation preserving, that is, its Jacobian satisfies $\intermediate{\jac} > 0$, and it fixes the center, $\intermediatemapping(\+0) = \+0$.
\end{assumption}
The physical mapping is the composition
\begin{equation}\label{eqn:composed_mapping}
  \mapping \defeq \intermediatemapping \circ \logicalmapping
  , \qquad
  \Omega = \mapping(\logical{\Omega}),
\end{equation}
and the left edge $\logical{x} = 0$ of the logical domain is mapped to the center $\+0$ of both $\intermediate{\Omega}$ and $\Omega$.
A given physical mapping can generally be decomposed in several ways; the choice of factors is part of the discretization data, since the polar spaces constructed below depend on it.
\begin{example}[Disk]\label{ex:disk_mapping}
  A disk of radius $R_1$ is obtained by choosing $R(\logical{x}, \logical{y}) = R_1 \logical{x}$, $C(\logical{y}) = \cos(2\pi \logical{y})$, and $S(\logical{y}) = \sin(2\pi \logical{y})$.
  Note that in this case $C S' - S C' = 2 \pi$.
  With $\intermediatemapping = \mathrm{id}$, the physical domain is this disk.
\end{example}
\begin{example}[Banana-shaped domain]\label{ex:banana_mapping}
  For $\kappa \ge 0$ and $\alpha > 0$, the intermediate mapping
  \begin{equation}
    \intermediatemapping(\intermediatecoord{\+x}) = \begin{bmatrix}
      \alpha\sspar{\frac{1}{2}\rpar{\intermediatecoord{x} - \sqrt{3}\,\intermediatecoord{y}} + \frac{\kappa}{4}\rpar{\sqrt{3}\,\intermediatecoord{x} + \intermediatecoord{y}}^2}\\
      \frac{1}{2\alpha}\rpar{\sqrt{3}\,\intermediatecoord{x} + \intermediatecoord{y}}
    \end{bmatrix},
  \end{equation}
  is polynomial, has the Jacobian $\intermediate{\jac} = 1$ and a polynomial inverse, and fixes the origin.
  Hence it satisfies \cref{ass:diffeomorphism} and maps a disk-like domain to a banana-shaped domain, see \cref{fig:disk_mapping}.
\end{example}

The Jacobian matrix of the logical mapping $\logicalmapping$ is given by
\begin{equation}\label{eqn:disk_jacobian_matrix}
  \logical{\jacmat} = \begin{bmatrix}
    \pard{R}{\logical{x}} C & \pard{R}{\logical{y}} C + R C'\\ 
    \pard{R}{\logical{x}} S & \pard{R}{\logical{y}} S + R S'
  \end{bmatrix},
\end{equation}
where $C' \defeq \totd{C}{\logical{y}}$, $S' \defeq \totd{S}{\logical{y}}$, and its Jacobian is given by
\begin{equation}\label{eqn:disk_jacobian}
  \logical{\jac} 
  = \det \logical{\jacmat}
  = R \pard{R}{\logical{x}} \rpar{C S' - S C'}.
\end{equation}
Denoting the Jacobian matrix and Jacobian of the intermediate mapping $\intermediatemapping$ by $\intermediate{\jacmat}$ and $\intermediate{\jac} = \det\intermediate{\jacmat}$, functions of $\intermediatecoord{\+x}$, the chain rule gives the Jacobian matrix and Jacobian of the physical mapping $\mapping$
\begin{equation}\label{eqn:composed_jacobian}
  \jacmat = (\intermediate{\jacmat} \circ \logicalmapping)\, \logical{\jacmat}
  , \quad
  \jac = (\intermediate{\jac} \circ \logicalmapping)\, \logical{\jac}.
\end{equation}
Since $\intermediate{\jac} > 0$, this shows that the Jacobian matrix is singular at $\logical{x} = 0$: $\jac(0, \logical{y}) = \logical{\jac}(0, \logical{y}) = 0$, thus requiring a special treatment.

\subsection{De Rham complex in one dimension}
The function spaces defined on the 2D logical domain are a tensor product of 1D function spaces defined on the unit interval $[0, 1]$.
In the $\logical{x}$-direction we use clamped B-splines, and in the periodic $\logical{y}$-direction we use periodic B-splines or trigonometric functions, and we note that other choices of 1D function spaces are possible within the framework presented in this work.
The scalar B-splines have degree at least one and simple interior knots, and the periodic B-splines have uniform knots.
Consequently, the scalar basis functions are continuous and piecewise $\funspace{C}^1$ with bounded first derivatives; in particular, the periodic ones are continuous across the identified endpoints $\logical{y}=0\sim1$.
For details about the construction of B-splines we refer to \citet{deboor1978}.

The choice of the angular space involves a trade-off that recurs throughout the paper.
Periodic B-splines are compactly supported, so that mass and stiffness matrices are sparse, but they cannot represent the angular functions $\cos(2\pi\logical{y})$ and $\sin(2\pi\logical{y})$ of a disk-like mapping exactly.
Trigonometric functions represent them exactly, which is what the regularity results of \cref{sec:zeroforms,sec:derham_2d} require through the compatibility condition of \cref{def:pr_compatibility_mapping_logical_zeroform} (see \cref{ex:pr_compatibility_disc_mapping}), but they are global in the angular direction, so that the matrices couple all angular degrees of freedom.
The numerical experiments in \cref{sec:results} compare both choices: the structural properties of the polar spaces and the suppression of the stiff modes originating at the center do not depend on it, whereas the pointwise regularity at the center, which particle orbits require, does.
A periodic space that combines both advantages, compact support and the exact representation of the trigonometric functions up to a given frequency, could be obtained by extending the Tchebycheffian B-spline de Rham complexes of \citet{manni2026tbsplines} to the periodic setting; we leave this to future work.

We consider pairs of 1D function spaces which form a de Rham complex in 1D \citep{Arnold_Falk_Winther_2006,arnold2018feec}.
Throughout the paper, a de Rham complex is a sequence of function spaces connected by differential operators such that each operator maps its space into the next one and the composition of two consecutive operators vanishes; a finite-dimensional de Rham complex is also called a cochain complex.
In 1D there is a single operator, the derivative, and the complex reads
\begin{equation}
  \tikzsetnextfilename{derham_logical_1d}  \derhamcomplex[
    dimension=1,
    custom-spaces={{\funspace{\Bsplineone}}^{0,\alpha},{\bfunspace{\Bsplineone}}^{1,\alpha}},
    operators={\totd{}{\logical{x}_\alpha}}
  ]
  , \quad \alpha \in \{\xx, \yy\}.
\end{equation}
for $\logical{x}_\xx = \logical{x}$ and $\logical{x}_\yy = \logical{y}$.
That is, each pair of finite-dimensional function spaces $\funspace{\Bsplineone}^{0,\alpha}, \funspace{\Bsplineone}^{1,\alpha}$ is such that $\funspace{\Bsplineone}^{1,\alpha}$ contains the derivatives of $\funspace{\Bsplineone}^{0,\alpha}$, and the kernel of the derivative operator is contained in $\funspace{\Bsplineone}^{0,\alpha}$.
The first space $\funspace{\Bsplineone}^{0,\alpha}$ is referred to as the $0$-form space, whereas the second space $\funspace{\Bsplineone}^{1,\alpha}$ is referred to as the $1$-form space.
Besides the finite-dimensional function spaces, we also consider a pair of infinite-dimensional function spaces $\funspace{H}^{\alpha}_1$ (the space of square-integrable functions whose weak derivative is also square-integrable) and $\funspace{L}^{\alpha}_2$ (the space of square-integrable functions) which themselves form a de Rham complex in 1D.
Finally, the two de Rham complexes are connected by commuting projection operators $\pi^{0,\alpha}$ and $\pi^{1,\alpha}$ introduced in \cref{sec:geometric_dofs_1d}, which are such that the following diagram commutes
\begin{equation}\label{eqn:exact_sequence}
  \tikzsetnextfilename{derham_logical_commute_1d}  \derhamcommuting[
    dimension=1,
    custom-spaces={{\funspace{H}}^{\alpha}_1,{\funspace{L}}_2^{\alpha}},
    custom-discrete-spaces={{\funspace{\Bsplineone}}^{0,\alpha},{\bfunspace{\Bsplineone}}^{1,\alpha}},
    operators={\totd{}{\logical{x}_\alpha}},
    discrete-operators={\totd{}{\logical{x}_\alpha}},
    maps={{\pi}^{0,\alpha},{\pi}^{1,\alpha}}
  ]
  , \quad \alpha \in \{\xx, \yy\}.
\end{equation}

\subsubsection{Function spaces}\label{sec:one_dimensional_function_spaces}
\begin{figure}
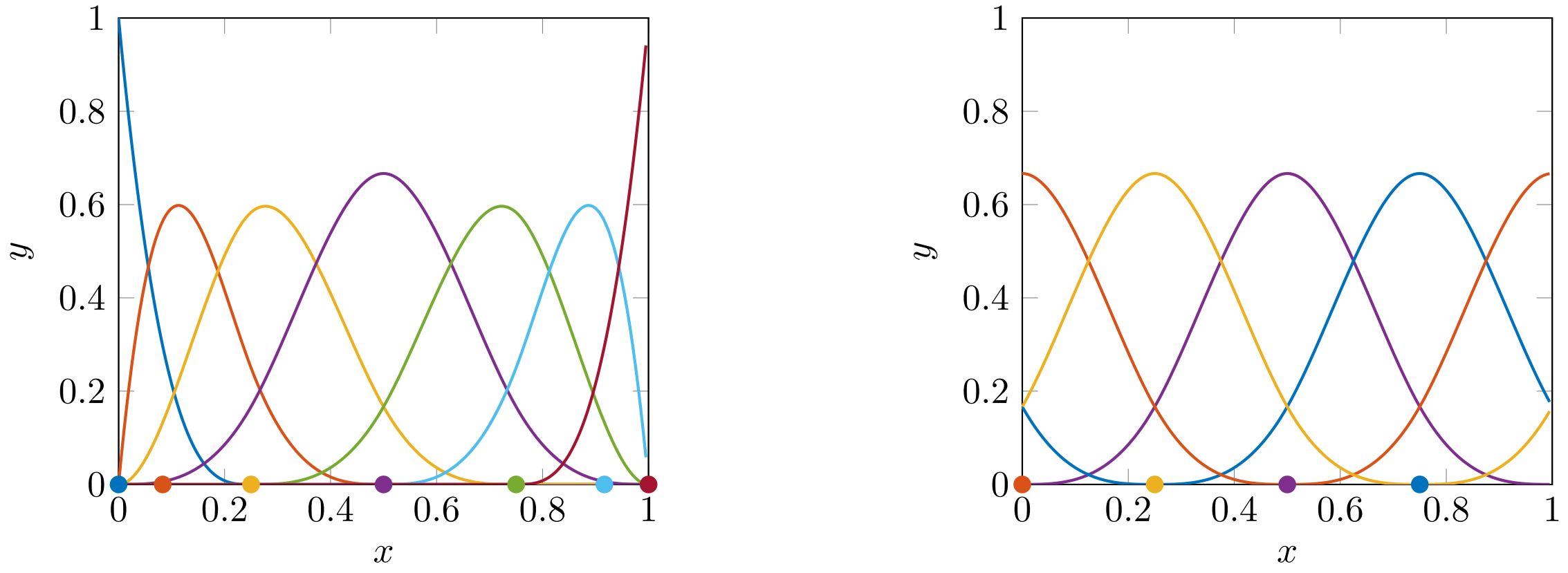

  \subcaptionbox{Clamped B-splines.\label{fig:clamped_bsplines}}
  [\twofigwidth]{
    \tikzsetnextfilename{bspline_clamped_example_r3}\input{tikz/bspline_clamped_example_r3.tex} 
  }\hfill
  \subcaptionbox{Periodic B-splines.\label{fig:periodic_bsplines}}
  [\twofigwidth]{
    \tikzsetnextfilename{bspline_periodic_example_r3}\input{tikz/bspline_periodic_example_r3.tex}
  }
  \caption{Example cubic (i.e., $r^\alpha=3$) B-splines (solid) and their Greville points (markers), using $N^\xx = N^\yy = 4$ intervals based on equidistant grid points.
  This results in $n^\xx = N^\xx + r^\xx = 7$ clamped and $n^\yy = N^\yy = 4$ periodic B-splines on the interval $[0, 1]$.}\label{fig:example_bsplines}
\end{figure}
The spaces $\funspace{\Bsplineone}^{0,\alpha}$ and $\funspace{\Bsplineone}^{1,\alpha}$ are spanned by sets of linearly independent basis functions $\Bsplineone^{0,\alpha}_i$ and $\Bsplineone^{1,\alpha}_i$, respectively
\begin{equation}
  \funspace{\Bsplineone}^{0,\alpha} \defeq \spn\{\Bsplineone^{0,\alpha}_1, \ldots, \Bsplineone^{0,\alpha}_{n^{0,\alpha}}\}
  , \quad
  \funspace{\Bsplineone}^{1,\alpha} \defeq \spn\{\Bsplineone^{1,\alpha}_1, \ldots, \Bsplineone^{1,\alpha}_{n^{1,\alpha}}\}
  , \quad
  \alpha \in \{\xx, \yy\}.
\end{equation}
The number of basis functions in each space is denoted by $n^{0,\alpha}$ and $n^{1,\alpha}$, respectively, where for the clamped B-spline space we have $n^{1,\xx} = n^{0,\xx} - 1$ and for the periodic spaces we have $n^{1,\yy} = n^{0,\yy} = n^\yy$.
When B-splines are considered, the degree of the 1D B-spline space $\funspace{\Bsplineone}^{0,\alpha}$ is denoted by $r^{0,\alpha}$, and the degree of the B-spline space $\funspace{\Bsplineone}^{1,\alpha}$ is equal to $r^{1,\alpha} = r^{0,\alpha} - 1$.
An illustration of cubic (i.e., $r^\alpha=3$) B-spline spaces is given in \cref{fig:example_bsplines}.

We define the following indexing sets for the basis functions in the 1D function spaces
\begin{equation}\label{eqn:1d_indexing_sets}
  \indexone^{0,\alpha} \defeq \{1, \ldots, n^{0,\alpha}\}
  , \quad
  \indexone^{1,\alpha} \defeq \{1, \ldots, n^{1,\alpha}\}.
\end{equation}

Some key properties of the 1D B-spline spaces are that they form a partition of unity
\begin{equation}\label{eqn:bspline_partition}
  \sum_{i\in\indexone^{0,\alpha}} \Bsplineone^{0,\alpha}_i = 1,
\end{equation}
they satisfy the following derivative property
\begin{equation}\label{eqn:bspline_derivative}
  \totd{\Bsplineone^{0,\alpha}_i}{\logical{x}_\alpha} = \Bsplineone^{1,\alpha}_{i-1} - \Bsplineone^{1,\alpha}_{i},
\end{equation}
where the $1$-form basis functions are scaled such that this identity holds as given, we define $\Bsplineone^{1,\xx}_0 = \Bsplineone^{1,\xx}_{n^{0,\xx}} = 0$ and $\Bsplineone^{1,\yy}_0 = \Bsplineone^{1,\yy}_{n^\yy}$, and the clamped B-splines satisfy the following property at the boundary $\logical{x} = 0$
\begin{equation}\label{eqn:bspline_x_atpole}
  \Bsplineone^{0,\xx}_i(0) = \begin{cases}
    1 & \text{if $i = 1$}\\
    0 & \text{if $i > 1$}\\
  \end{cases}, \quad
  \totd{\Bsplineone^{0,\xx}_i}{\logical{x}}(0) = \begin{cases}
    -\frac{1}{\Delta t^{0,\xx}_1} & \text{if $i = 1$}\\
    \frac{1}{\Delta t^{0,\xx}_1} & \text{if $i = 2$}\\
    0 & \text{otherwise}\\
  \end{cases},
\end{equation}
where $\Delta t^{0,\xx}_1 = \rpar{t^{0,\xx}_{r^{0,\xx}+2} - t^{0,\xx}_{2}}/{r^{0,\xx}}$ is defined in terms of the clamped knot sequence $\+t^{0,\xx}$.
When trigonometric functions are used in the $\logical{y}$-direction, the number of basis functions $n^\yy$ is odd, and the bases of the $0$- and $1$-form spaces are ordered by frequency,
\begin{equation}\label{eqn:trigonometric_basis}
  \begin{aligned}
    \Bsplineone^{0,\yy}_1 &= 1, &
    \Bsplineone^{0,\yy}_{2k} &= \sin(2\pi k\logical{y}), &
    \Bsplineone^{0,\yy}_{2k+1} &= \cos(2\pi k\logical{y}),\\
    \Bsplineone^{1,\yy}_1 &= 1, &
    \Bsplineone^{1,\yy}_{2k} &= 2\pi k\sin(2\pi k\logical{y}), &
    \Bsplineone^{1,\yy}_{2k+1} &= 2\pi k\cos(2\pi k\logical{y}),
  \end{aligned}
  \qquad k = 1, \ldots, \frac{n^\yy-1}{2},
\end{equation}
so that both spaces contain the frequencies up to $(n^\yy-1)/2$; we refer to each basis function as a mode.
The frequency factors in the $1$-form basis make the derivative relation free of coefficients,
\begin{equation}\label{eqn:trigonometric_derivative}
  \totd{\Bsplineone^{0,\yy}_i}{\logical{y}} = \begin{cases}
    0 & \text{if $i = 1$}\\
    \Bsplineone^{1,\yy}_{i+1} & \text{if $i$ is even}\\
    -\Bsplineone^{1,\yy}_{i-1} & \text{otherwise}\\
  \end{cases}.
\end{equation}

Throughout the paper we make the following assumption on the 1D spaces.
\begin{assumption}[Constant functions]\label{ass:constant_functions}
  Each 1D $0$-form space $\funspace{\Bsplineone}^{0,\alpha}$ contains the constant functions.
\end{assumption}
The clamped and periodic B-spline spaces satisfy this assumption by the partition of unity in \cref{eqn:bspline_partition}, and the trigonometric space by $\Bsplineone^{0,\yy}_1=1$.
We let $\mathds{1}$ denote the coefficients of the constant function $1$ in the basis $\{\Bsplineone^{0,\yy}_i\}_{i\in\indexone^{0,\yy}}$ of periodic functions.
That is
\begin{equation}\label{eqn:constant_function_coefficients}
    \mathds{1}_j^\yy \defeq 1
  \text{ for periodic B-splines, and }
  \mathds{1}_j^\yy \defeq \begin{cases}
    1 & \text{if $j = 1$}\\
    0 & \text{otherwise}
  \end{cases}
  \text{ for trigonometric functions.}
\end{equation}

\subsubsection{Differential operator}
Given a 1D B-spline function $\Bsplineone \in \funspace{\Bsplineone}^{0,\alpha}$
\begin{equation}\label{eqn:bspline_expansion_1d}
  \Bsplineone = \sum_{i\in\indexone^{0,\alpha}} \bsplineone_i \Bsplineone^{0,\alpha}_i,
\end{equation}
with coefficients $\bsplineone_i \in \mathbb{R}$, we can express its derivative as follows
\begin{equation}\label{eqn:derivative_operator_1d}
  \totd{\Bsplineone}{\logical{x}_\alpha} = \sum_{i\in \indexone^{0,\alpha}} \bsplineone_i \totd{\Bsplineone^{0,\alpha}_i}{\logical{x}_\alpha}
  = \sum_{i\in \indexone^{0,\alpha}} \bsplineone_i (\Bsplineone^{1,\alpha}_{i-1} - \Bsplineone^{1,\alpha}_{i})
  = \sum_{i\in \indexone^{1,\alpha}} (\delta_\alpha \bsplineone)_i \Bsplineone^{1,\alpha}_i,
\end{equation}
by using the derivative property from \cref{eqn:bspline_derivative}, applying summation by parts, and defining the 1D forward differencing operator
\begin{equation}\label{eqn:diff_opers_1d}
  (\delta_\alpha \bsplineone)_{i} \defeq \bsplineone_{i+1} - \bsplineone_{i}.
\end{equation}
Hence, the derivative operator $\totd{}{\logical{x}_\alpha}$ yields a 1D $1$-form in $\funspace{\Bsplineone}^{1,\alpha}$ from a 1D $0$-form in $\funspace{\Bsplineone}^{0,\alpha}$, and the coefficients of the resulting $1$-form are given by the forward differencing operator $\delta_\alpha$ applied to the coefficients of the $0$-form.

For the trigonometric basis, the coefficient differentiation operator in the periodic direction is instead defined as
\begin{equation}\label{eqn:trigonometric_coefficient_derivative}
  (\delta_\yy \bsplineone)_1\defeq0,
  \qquad
  (\delta_\yy \bsplineone)_{2k}\defeq-\bsplineone_{2k+1},
  \qquad
  (\delta_\yy \bsplineone)_{2k+1}\defeq\bsplineone_{2k},
  \quad k=1,\ldots,(n^\yy-1)/2.
\end{equation}

\subsubsection{Geometric degrees of freedom}\label{sec:geometric_dofs_1d}
Commuting projectors can most easily be constructed by defining geometric degrees of freedom (GDOFs) which are defined in terms of evaluations of $0$-forms at Greville points $\xi^{\alpha}_i \defeq \frac{t^\alpha_{i+1} + \ldots + t^\alpha_{i+r^{0,\alpha}}}{r^{0,\alpha}}$ (assuming $r^{0,\alpha} > 0$), and integrals of $1$-forms over the intervals defined by the same Greville points
\begin{equation}\label{eqn:geometric_dofs_1d}
  \gdof^{0,\alpha}_{i}(\Bsplineone) \defeq \Bsplineone(\xi^{\alpha}_i)
  , \quad
  \gdof^{1,\alpha}_{i}(\Bsplineone) \defeq \int_{\xi^{\alpha}_i}^{\xi^{\alpha}_{i+1}} \Bsplineone \diff \logical{x}_\alpha.
\end{equation}
When trigonometric functions are considered in the $\logical{y}$-direction, we use the $n^{\yy}$ distinct equidistant points $\xi_j^\yy=(j-1)/n^\yy$, $j=1,\ldots,n^\yy$, with $\xi_{n^\yy+1}^\yy=\xi_1^\yy+1$ for periodic integration.

\begin{lemma}[Unisolvence of the 1D GDOFs]\label{lem:unisolvence_1d_gdofs}
  For each $m\in\{0,1\}$ and $\alpha\in\{\xx,\yy\}$, the GDOFs in \cref{eqn:geometric_dofs_1d} are unisolvent on $\funspace{\Bsplineone}^{m,\alpha}$: an element of $\funspace{\Bsplineone}^{m,\alpha}$ is uniquely determined by its GDOFs.
\end{lemma}
\begin{statementproof}[appendix]{lem:unisolvence_1d_gdofs}
  For $m=0$, the clamped B-spline Greville collocation matrices are nonsingular by the Schoenberg--Whitney property \citep[Ch.~XIII]{deboor1978}; for uniform periodic B-splines, the collocation matrix is circulant, and its eigenvalues do not vanish \citep{schoenberg1973cardinal}; for the odd-dimensional trigonometric basis, the distinct equidistant nodes give an invertible real Fourier matrix.
  For $m=1$, let $\Bsplineone^1\in\funspace{\Bsplineone}^{1,\alpha}$ have zero integral on every Greville interval. The derivative relations in \cref{eqn:bspline_derivative,eqn:trigonometric_derivative}, together with the dimensions of the 1D spaces, show that differentiation maps onto the clamped $1$-form space and onto the zero-integral subspace of each periodic $1$-form space. The vanishing edge integrals give zero total integral in the periodic case, so $\Bsplineone^1$ has a primitive $\Bsplineone^0\in\funspace{\Bsplineone}^{0,\alpha}$. 
  Since $\Bsplineone^0$ has the same value at all Greville points, scalar unisolvence makes $\Bsplineone^0$ constant and hence $\Bsplineone^1=0$.
  Thus each square matrix of GDOFs applied to basis functions is injective and therefore invertible.
\end{statementproof}
Hence, for every $\Bsplineone$ for which the GDOFs are defined, there is a unique $\pi^{m,\alpha}\Bsplineone\in\funspace{\Bsplineone}^{m,\alpha}$ with the same GDOFs,
\begin{equation}\label{eqn:commuting_projectors_1d}
  \gdof^{m,\alpha}_i(\pi^{m,\alpha} \Bsplineone) = \gdof^{m,\alpha}_i(\Bsplineone), \quad
  \forall i \in \{1, \ldots, n^{m,\alpha}\}.
\end{equation}
This defines the projectors $\pi^{m,\alpha}$ (see \cref{sec:matrix_representation_gdof_projectors} for their matrix representation).
Since the GDOFs of a derivative are differences of nodal values, the projectors commute with differentiation, which makes the diagram in \cref{eqn:exact_sequence} commute:
\begin{equation}\label{eqn:stokes_1d}
  \gdof^{1,\alpha}_i \srpar{\pi^{1,\alpha} \totd{\Bsplineone}{\logical{x}_\alpha}} 
  =
  \gdof^{1,\alpha}_i \srpar{\totd{(\pi^{0,\alpha} \Bsplineone)}{\logical{x}_\alpha}} 
  \quad \iff \quad
  \pi^{1,\alpha} \totd{\Bsplineone}{\logical{x}_\alpha} = \totd{(\pi^{0,\alpha} \Bsplineone)}{\logical{x}_\alpha}.
\end{equation}
In the context of the de Rham complex, this means that the commuting projectors $\pi^{m,\alpha}$ are cochain maps: they are mappings between two de Rham complexes which commute with the differential operators.

\subsection{De Rham complex in two-dimensional logical space}\label{sec:logical_bspline_spaces}
In 2D, we consider three function spaces which -- together with their infinite dimensional counterparts -- form the following commuting diagram of de Rham complexes
\begin{equation}\label{eqn:logical_commuting_diagram}
  \tikzsetnextfilename{derham_logical_commute_2d}  \derhamcommuting[
    dimension=2,
    custom-spaces={\logical{\funspace{H}}_1,\logical{\funspace{H}}_{\curlt},\logical{\funspace{L}}_2},
    custom-discrete-spaces={\logical{\funspace{\Bspline}}^0,\logical{\bfunspace{\Bspline}}^1,\logical{\bfunspace{\Bspline}}^2},
    operators={\logical{\gradientt},\logical{\curlt}},
    discrete-operators={\logical{\gradientt},\logical{\curlt}},
    maps={\logical{\pi}^0,\logical{\pi}^1,\logical{\pi}^2}
  ],
\end{equation}
where each logical 2D function space $ \logical{\funspace{\Bspline}}^m$ (for $m \in \{0, 1, 2\}$) defined on the logical domain $\logical{\Omega}$ is a tensor product of the 1D function spaces defined in \cref{sec:one_dimensional_function_spaces}.
The $1$-form space $\logical{\funspace{\Bspline}}^1$ is vector-valued, and the $0$- and $2$-form spaces are scalar-valued.
The logical 2D function spaces form a cochain complex by the differential operators $\logical{\gradientt}$ and $\logical{\curlt}$, which are the gradient and curl operators in 2D logical space, where the curl is understood in the two-dimensional sense as the third component of the three-dimensional curl operator.
That is, the image of the gradient operator is contained in the $1$-form space, the image of the curl operator is contained in the $2$-form space, and the curl of the gradient is zero (this is the cochain complex property).
The corresponding infinite-dimensional spaces are $\funspace{H}_1$ (square-integrable scalar functions with square-integrable gradient), $\funspace{H}_{\curlt}$ (square-integrable vector fields with square-integrable curl), and $\funspace{L}_2$ (square-integrable functions).
For precise definitions, see \citet{Arnold_Falk_Winther_2006}.
These spaces form a de Rham complex, to which the logical 2D spaces are connected by the commuting projection operators $\logical{\pi}^m$ defined in \cref{sec:logical_geometric_dofs}.

\begin{remark}[Assumptions on the fields]
  Throughout the paper, whenever a GDOF projection is applied to a Sobolev field, we assume that its required point values (for a $0$-form), edge integrals (for a $1$-form) and face integrals (for a $2$-form) are well-defined.
  For identities involving derivatives or radial homotopies, we also assume that the resulting fields have the required GDOFs and that the cellwise fundamental theorem of calculus, Stokes' theorem, and interchange of integrals apply. 
  Smooth forms and the conforming polar forms below satisfy these conditions. 
  The vertical projection arrows in the diagrams are understood on such fields; they do not assert bounded projection on all of $\funspace{H}_1$ or $\funspace{H}_{\curlt}$, and bounded commuting projections in the sense of \citet{arnold2010feec,falk2014} are not considered in this work.
\end{remark}

\subsubsection{Function spaces}\label{sec:logical_function_spaces}
The logical 2D function spaces are defined as tensor products of the 1D function spaces from \cref{sec:one_dimensional_function_spaces}, combined degree by degree (the two-factor case of the total complex used in \cref{sec:derham_hodge_3d_construction}), resulting in
\begin{equation}\label{eqn:logical_function_spaces}
  \logical{\funspace{\Bspline}}^0 \defeq \funspace{\Bsplineone}^{0,\xx} \otimes \funspace{\Bsplineone}^{0,\yy}
  , \quad
  \logical{\funspace{\Bspline}}^1 \defeq 
  \funspace{\Bsplineone}^{1,\xx} \otimes \funspace{\Bsplineone}^{0,\yy}
  \oplus
  \funspace{\Bsplineone}^{0,\xx} \otimes \funspace{\Bsplineone}^{1,\yy}
  , \quad
  \logical{\funspace{\Bspline}}^2 \defeq \funspace{\Bsplineone}^{1,\xx} \otimes \funspace{\Bsplineone}^{1,\yy},
\end{equation}
where the de Rham complex structure of each of the 1D function spaces implies that the 2D function spaces also form a de Rham complex.
Here, the internal direct sum $\oplus$ is defined as $\funspace{A} \oplus \funspace{B} = \setdef{A + B}{A \in \funspace{A}, B \in \funspace{B}}$, and the tensor product is realized by the wedge product of differential forms, 
see also \cref{rem:vector_proxy,ex:derham_hodge_3d_tensor_product}.

We denote logical $0$-forms by $\logical{P} \in \logical{\funspace{\Bspline}}^0$, logical $1$-forms by $\logical{\+A} \in \logical{\funspace{\Bspline}}^1$, and logical $2$-forms by $\logical{B} \in \logical{\funspace{\Bspline}}^2$.
Each can be written as a linear combination of basis functions:
\begin{equation}\label{eqn:logical_zeroforms}
  \logical{P} = \sum_{(i,j)\in\indextwo^{0}} p_{ij} \logical{\Bspline}_{ij}^0
  , \quad 
  \logical{\Bspline}_{ij}^0 \defeq \Bsplineone^{0,\xx}_i \Bsplineone^{0,\yy}_j
  ,
\end{equation}
and
\begin{equation}\label{eqn:logical_oneforms}
  \logical{\+A} = 
  \sum_{(i,j)\in\indextwo^{1}_\xx} a_{\xx,ij} \logical{\+\Bspline}_{\xx,ij}^1
  + 
  \sum_{(i,j)\in\indextwo^{1}_\yy} a_{\yy,ij} \logical{\+\Bspline}_{\yy,ij}^1
  , \quad 
  \logical{\+\Bspline}_{\xx,ij}^1 \defeq \begin{bmatrix}
    \Bsplineone^{1,\xx}_i \Bsplineone^{0,\yy}_j\\
    0
  \end{bmatrix}
  , \quad
  \logical{\+\Bspline}_{\yy,ij}^1 \defeq \begin{bmatrix}
    0\\
    \Bsplineone^{0,\xx}_i \Bsplineone^{1,\yy}_j
  \end{bmatrix},
\end{equation}
as well as
\begin{equation}\label{eqn:logical_twoforms}
  \logical{B} = \sum_{(i,j)\in\indextwo^{2}} b_{ij} \logical{\Bspline}_{ij}^2
  , \quad 
  \logical{\Bspline}_{ij}^2 \defeq \Bsplineone^{1,\xx}_i \Bsplineone^{1,\yy}_j
  ,
\end{equation}
where $p_{ij}, a_{\xx,ij}, a_{\yy,ij}, b_{ij} \in \mathbb{R}$ denote the coefficients of the corresponding logical forms, and the indexing sets are given by (cf.\ \cref{eqn:1d_indexing_sets})
\begin{equation}\label{eqn:2d_indexing_sets}
  \indextwo^{0} = \indexone^{0,\xx} \times \indexone^{0,\yy}, \quad
  \indextwo^{1}_\xx = \indexone^{1,\xx} \times \indexone^{0,\yy}, \quad
  \indextwo^{1}_\yy = \indexone^{0,\xx} \times \indexone^{1,\yy}, \quad
  \indextwo^{2} = \indexone^{1,\xx} \times \indexone^{1,\yy}.
\end{equation}

\begin{remark}[Vector proxy]\label{rem:vector_proxy}
  In this paper we do not directly work with differential forms, but rather with their vector proxies.
  For example, for a $1$-form $\nu^1$ given by
  \begin{equation}
    \nu^1 = \logical{\Bspline}^1_\xx \diff \logical{x} + \logical{\Bspline}^1_\yy \diff \logical{y}
  \end{equation}
  the corresponding logical vector proxy is
  \begin{equation}
    \logical{\+\Bspline}^1 = \begin{bmatrix}
      \logical{\Bspline}^1_\xx\\
      \logical{\Bspline}^1_\yy
    \end{bmatrix}.
  \end{equation}
  Unless indicated otherwise, we refer to the (logical) vector proxies of the differential forms when we use the term (logical) `form'.
  Whenever operations on differential forms, such as the wedge product, are applied to vector proxies, they act on the associated forms.
\end{remark}

\subsubsection{Differential operators}
We extend the 1D coefficient differentiation operators $\delta_\xx,\delta_\yy$ to the 2D coefficients $p_{ij}, a_{\xx,ij}, a_{\yy,ij}, b_{ij}$ by letting $\delta_\xx$ act on the first index and $\delta_\yy$ on the second index.
They are the forward differences in \cref{eqn:diff_opers_1d} for B-splines, while $\delta_\yy$ is given by \cref{eqn:trigonometric_coefficient_derivative} for trigonometric functions.
This yields the following discrete analogues of the differential operators $\logical{\gradientt}$ and $\logical{\curlt}$, which act on the coefficients of the logical forms
\begin{equation}\label{eqn:discrete_operators}
  \gradienth p \defeq \begin{bmatrix}
    \delta_\xx p\\
    \delta_\yy p
  \end{bmatrix}, \quad
  \curlh \+a \defeq \delta_\xx a_\yy - \delta_\yy a_\xx.
\end{equation}
We note that the image of $\gradienth$ is in the kernel of $\curlh$, thereby mimicking the continuous differential operators.
Their matrix representations as Kronecker products are given in \cref{eqn:matrep_differential_2d}.

\begin{example}[Curl of a logical $1$-form]\label{ex:curl_logical_oneform}
  The curl of a logical $1$-form $\logical{\+A}$, as given by \cref{eqn:logical_oneforms}, is given by
  \begin{equation}
    \logical{\curlt} \logical{\+A} 
    = \sum_{(i,j)\in\indextwo^2} (\delta_\xx a_\yy - \delta_\yy a_\xx)_{ij} \logical{\Bspline}_{ij}^2
    = \sum_{(i,j)\in\indextwo^2} (\curlh \+a)_{ij} \logical{\Bspline}_{ij}^2.
  \end{equation}
\end{example}

\subsubsection{Geometric degrees of freedom}\label{sec:logical_geometric_dofs}
As with the function spaces and differential operators, the GDOFs in 2D are defined in terms of their 1D counterparts.
That is, the GDOFs of the logical 2D function spaces are defined in terms of point evaluations of $0$-forms at Greville nodes, line integrals of $1$-forms over Greville edges, and surface integrals of $2$-forms over Greville faces:
\begin{equation}\label{eqn:logical_geometric_dofs}
  \logical{\gdof}^{0}_{ij} \defeq \gdof^{0,\xx}_i \circ \gdof^{0,\yy}_j
  , \quad
  \logical{\gdof}^{1}_{\xx,ij} \defeq \gdof^{1,\xx}_i \circ \gdof^{0,\yy}_j
  , \quad
  \logical{\gdof}^{1}_{\yy,ij} \defeq \gdof^{0,\xx}_i \circ \gdof^{1,\yy}_j
  , \quad
  \logical{\gdof}^{2}_{ij} \defeq \gdof^{1,\xx}_i \circ \gdof^{1,\yy}_j,
\end{equation}
where the $\logical{x}$- and $\logical{y}$-components of the $1$-form GDOFs are applied to the corresponding components of the logical $1$-form.
In physical terms, the GDOFs are the quantities that the discrete forms represent: the value of a potential at a node, the circulation (or voltage) of a field along an edge, and the flux of a density through a face.
A projection that preserves the GDOFs therefore preserves Stokes' theorem on the Greville grid exactly: the circulation of a gradient along an edge is the difference of the potential at its endpoints, and the flux of a curl through a face is the circulation around its boundary.
This is the content of the commutation relations below.

\begin{example}[$1$-form GDOFs]
  The GDOFs of a logical $1$-form $\logical{\+A}$ are given by
  \begin{equation}
    \logical{\gdof}^{1}_{\xx,ij}(\logical{\+A}) = \int_{\xi^\xx_i}^{\xi^\xx_{i+1}} \logical{A}_\xx(\logical{x}, \xi_j^\yy) \diff \logical{x}
    , \quad
    \logical{\gdof}^{1}_{\yy,ij}(\logical{\+A}) = \int_{\xi^\yy_j}^{\xi^\yy_{j+1}} \logical{A}_\yy(\xi_i^\xx, \logical{y}) \diff \logical{y}.
  \end{equation}
\end{example}

The 2D logical projection operators $\logical{\pi}^m$ are then defined by imposing that they leave the GDOFs invariant
\begin{equation}\label{eqn:logical_commuting_projectors}
  \logical{\gdof}^{0}(\logical{\pi}^0 \logical{P}) = \logical{\gdof}^{0}(\logical{P})
  , \quad
  \logical{\gdof}^{1}(\logical{\pi}^1 \logical{\+A}) = \logical{\gdof}^{1}(\logical{\+A})
  , \quad
  \logical{\gdof}^{2}(\logical{\pi}^2 \logical{B}) = \logical{\gdof}^{2}(\logical{B}),
\end{equation}
which uniquely defines the projection operators $\logical{\pi}^m$ by means of the invertibility of the Vandermonde matrices associated with the GDOFs and the basis functions.
We note that, owing to the tensor-product structure of the logical function spaces and the GDOFs, the 2D Vandermonde matrices are Kronecker products of the 1D Vandermonde matrices, and thereby inherit their invertibility (see \cref{eqn:two_dimensional_gdof_vandermonde_matrices,eqn:tensor_product_gdof_projector_matrix}).
Similar to \cref{eqn:stokes_1d}, this implies commutativity of the diagram in \cref{eqn:logical_commuting_diagram}
\begin{equation}
  \logical{\pi}^1 \circ \logical{\gradientt}
  =
  \logical{\gradientt} \circ \logical{\pi}^0 
  , \quad
  \logical{\pi}^2 \circ \logical{\curlt}
  =
  \logical{\curlt} \circ \logical{\pi}^1.
\end{equation}
Hence, $\logical{\pi}^m$ are cochain maps.

\subsection{Physical function spaces}\label{sec:physical_function_spaces}
The function spaces on the physical domain $\Omega$ are obtained from the logical function spaces in two steps: by means of the pushforward operators $\logical{\pushforward}^m$ of the logical mapping $\logicalmapping$ we obtain intermediate function spaces $\intermediate{\funspace{\Bspline}}^m$ on $\intermediate{\Omega}$, which are mapped to the physical function spaces $\funspace{\Bspline}^m$ by the pushforward operators $\intermediate{\pushforward}^m$ of the intermediate mapping $\intermediatemapping$.
The pushforward operators are defined in such a way that a de Rham complex is obtained on each level, resulting in the following commuting diagram
\begin{equation}\label{eqn:physical_commuting_diagram}
  \tikzsetnextfilename{derham_logical_and_physical_2d}  \begin{derhamdiagram}[dimension=2]
    \derhamrow[name=logical-continuous,spaces={\logical{\funspace{H}}_1,\logical{\bfunspace{H}}_{\curlt},\logical{\funspace{L}}_2},operators={\logical{\gradientt},\logical{\curlt}}]
    \derhamrow[name=logical-discrete,spaces={\logical{\funspace{\Bspline}}^0,\logical{\bfunspace{\Bspline}}^1,\logical{\funspace{\Bspline}}^2},operators={\logical{\gradientt},\logical{\curlt}}]
    \derhamrow[name=intermediate-discrete,spaces={\intermediate{\funspace{\Bspline}}^0,\intermediate{\bfunspace{\Bspline}}^1,\intermediate{\funspace{\Bspline}}^2},operators={\intermediate{\gradientt},\intermediate{\curlt}}]
    \derhamrow[name=physical-discrete,spaces={\funspace{\Bspline}^0,\bfunspace{\Bspline}^1,\funspace{\Bspline}^2},operators={\gradientt,\curlt}]
    \derhamrow[name=physical-continuous,spaces={\funspace{H}_1,\bfunspace{H}_{\curlt},\funspace{L}_2},operators={\gradientt,\curlt}]
    \derhamconnect[from=logical-continuous,to=logical-discrete,degree=0,operator=\logical{\pi}^0]
    \derhamconnect[from=logical-continuous,to=logical-discrete,degree=1,operator=\logical{\pi}^1]
    \derhamconnect[from=logical-continuous,to=logical-discrete,degree=2,operator=\logical{\pi}^2]
    \derhamconnect[from=logical-discrete,to=intermediate-discrete,degree=0,operator=\logical{\pushforward}^0,style={<->}]
    \derhamconnect[from=logical-discrete,to=intermediate-discrete,degree=1,operator=\logical{\pushforward}^1,style={<->}]
    \derhamconnect[from=logical-discrete,to=intermediate-discrete,degree=2,operator=\logical{\pushforward}^2,style={<->}]
    \derhamconnect[from=intermediate-discrete,to=physical-discrete,degree=0,operator=\intermediate{\pushforward}^0,style={<->}]
    \derhamconnect[from=intermediate-discrete,to=physical-discrete,degree=1,operator=\intermediate{\pushforward}^1,style={<->}]
    \derhamconnect[from=intermediate-discrete,to=physical-discrete,degree=2,operator=\intermediate{\pushforward}^2,style={<->}]
    \derhamconnect[from=physical-continuous,to=physical-discrete,degree=0,operator=\pi^0]
    \derhamconnect[from=physical-continuous,to=physical-discrete,degree=1,operator=\pi^1]
    \derhamconnect[from=physical-continuous,to=physical-discrete,degree=2,operator=\pi^2]
  \end{derhamdiagram}
.
\end{equation}
Here, $\intermediate{\gradientt}$ and $\intermediate{\curlt}$ denote the gradient and curl with respect to the intermediate coordinates $\intermediatecoord{\+x}$.
The composed pushforward operators $\pushforward^m$, which map the logical spaces directly to the physical spaces, commute with the differential operators $\gradientt$ and $\curlt$ in physical space
\begin{equation}\label{eqn:pushforward_commutativity}
  \gradientt \circ \, \pushforward^0 = \pushforward^1 \circ \logical{\gradientt}
  , \quad
  \curlt \circ \, \pushforward^1 = \pushforward^2 \circ \logical{\curlt}.
\end{equation}
When defining the $0$-form pushforward operator simply as the composition with the inverse mapping $\mapping^{-1}$, the $1$-form pushforward operator is defined according to how the gradient of a $0$-form transforms, and the $2$-form pushforward operator is defined according to how the curl of a $1$-form transforms.
This results in
\begin{equation}\label{eqn:pushforward_operators}
  P = \pushforward^0 \logical{P} \defeq \logical{P} \circ \mapping^{-1}
  , \quad
  \+A = \pushforward^1 \logical{\+A} \defeq (\jacmat^{-\transpose} \logical{\+A}) \circ \mapping^{-1}
  , \quad
  B = \pushforward^2 \logical{B} \defeq \srpar{\frac{1}{\jac} \logical{B}} \circ \mapping^{-1},
\end{equation}
where the $1$-form pushforward is known as the covariant Piola transformation \citep[Def.~2.5]{rognes2009}.
The pushforward operators $\logical{\pushforward}^m$ and $\intermediate{\pushforward}^m$ are given by the same formulas, with $\mapping$, $\jacmat$, $\jac$ replaced by $\logicalmapping$, $\logical{\jacmat}$, $\logical{\jac}$ and by $\intermediatemapping$, $\intermediate{\jacmat}$, $\intermediate{\jac}$, respectively.
By \cref{eqn:composed_mapping,eqn:composed_jacobian}, the pushforwards compose as
\begin{equation}\label{eqn:pushforward_composition}
  \pushforward^m = \intermediate{\pushforward}^m \circ \logical{\pushforward}^m
  , \quad m \in \{0, 1, 2\}.
\end{equation}
A straightforward calculation shows that each of the three families of pushforward operators commutes with the differential operators as in \cref{eqn:pushforward_commutativity}, and therefore the finite-dimensional intermediate and physical function spaces also form de Rham complexes.
The inverses of the pushforward operators are the pullback operators, which transport physical forms back to the logical domain; in the notation of differential geometry they are the pullbacks by $\mapping$, $\logicalmapping$ and $\intermediatemapping$, respectively.
We denote them by $\pushforward^{-m}$, $\logical{\pushforward}^{-m}$ and $\intermediate{\pushforward}^{-m}$, where the superscript $-m$ indicates the inverse of the pushforward of $m$-forms, so that $\pushforward^{-m}\pushforward^{m}$ is the identity on $\logical{\funspace{\Bspline}}^m$.
Explicitly, by \cref{eqn:pushforward_operators},
\begin{equation}\label{eqn:pullback_operators}
  \pushforward^{-0} P = P \circ \mapping
  , \quad
  \pushforward^{-1} \+A = \jacmat^{\transpose} (\+A \circ \mapping)
  , \quad
  \pushforward^{-2} B = \jac \, (B \circ \mapping),
\end{equation}
and likewise for $\logical{\pushforward}^{-m}$ and $\intermediate{\pushforward}^{-m}$.

The infinite-dimensional function spaces $\funspace{H}_1, \funspace{H}_{\curlt}$ and $\funspace{L}_2$ on $\Omega$ are connected to the finite-dimensional physical function spaces $\funspace{\Bspline}^m$ by means of commuting projection operators $\pi^m$, which are defined in \cref{sec:physical_geometric_dofs}.
Similarly, $\intermediate{\funspace{H}}_1, \intermediate{\funspace{H}}_{\curlt}$ and $\intermediate{\funspace{L}}_2$ denote the corresponding spaces on $\intermediate{\Omega}$.

\subsubsection{Function spaces}
\begin{table}
  \centering
  \setlength{\tabcolsep}{1.000000000000pt}
  \begin{tabular}{@{}cccccccccccc@{}}
    \toprule
    \multicolumn{12}{c}{\makebox[0pt][c]{Cubic clamped B-splines on $N^\xx = 5$ intervals and cubic periodic B-splines on $N^\yy = 12$ intervals}}\\
    \multicolumn{1}{c}{\tikzsetnextfilename{disc_basis_tensor_1_1}\input{tikz/disc_basis_tensor_1_1.tex}\unskip}&\multicolumn{1}{c}{\tikzsetnextfilename{disc_basis_tensor_1_2}\input{tikz/disc_basis_tensor_1_2.tex}\unskip}&\multicolumn{1}{c}{\tikzsetnextfilename{disc_basis_tensor_1_3}\input{tikz/disc_basis_tensor_1_3.tex}\unskip}&\multicolumn{1}{c}{\tikzsetnextfilename{disc_basis_tensor_1_4}\input{tikz/disc_basis_tensor_1_4.tex}\unskip}&\multicolumn{1}{c}{\tikzsetnextfilename{disc_basis_tensor_1_5}\input{tikz/disc_basis_tensor_1_5.tex}\unskip}&\multicolumn{1}{c}{\tikzsetnextfilename{disc_basis_tensor_1_6}\input{tikz/disc_basis_tensor_1_6.tex}\unskip}&\multicolumn{1}{c}{\tikzsetnextfilename{disc_basis_tensor_1_7}\input{tikz/disc_basis_tensor_1_7.tex}\unskip}&\multicolumn{1}{c}{\tikzsetnextfilename{disc_basis_tensor_1_8}\input{tikz/disc_basis_tensor_1_8.tex}\unskip}&\multicolumn{1}{c}{\tikzsetnextfilename{disc_basis_tensor_1_9}\input{tikz/disc_basis_tensor_1_9.tex}\unskip}&\multicolumn{1}{c}{\tikzsetnextfilename{disc_basis_tensor_1_10}\input{tikz/disc_basis_tensor_1_10.tex}\unskip}&\multicolumn{1}{c}{\tikzsetnextfilename{disc_basis_tensor_1_11}\input{tikz/disc_basis_tensor_1_11.tex}\unskip}&\multicolumn{1}{c}{\tikzsetnextfilename{disc_basis_tensor_1_12}\input{tikz/disc_basis_tensor_1_12.tex}\unskip} \\
\bottomrule
  \end{tabular}
  \caption{The first four rings of the tensor-product physical $0$-forms $\Bspline^{0}_{ij}$, $i \in \{1, 2, 3, 4\}$, $j \in \{1, \ldots\}$, on the physical banana-shaped domain $\Omega$ of \cref{fig:disk_mapping}: the logical mapping $\logicalmapping$ and the intermediate mapping $\intermediatemapping$ coincide with those of that figure.}
  \label{tab:disc_basis_functions_fourier}
\end{table}
The intermediate and physical function spaces are defined as the pushforwards of the logical function spaces $\logical{\funspace{\Bspline}}^m$:
\begin{equation}
  \intermediate{\funspace{\Bspline}}^m \defeq \logical{\pushforward}^m \logical{\funspace{\Bspline}}^m
  , \quad
  \funspace{\Bspline}^m \defeq \intermediate{\pushforward}^m \intermediate{\funspace{\Bspline}}^m = \pushforward^m \logical{\funspace{\Bspline}}^m,
\end{equation}
resulting in intermediate basis functions $\intermediate{\Bspline}_{ij}^0 \defeq \logical{\pushforward}^0 \logical{\Bspline}_{ij}^0$, $\intermediate{\+\Bspline}_{\xx,ij}^1 \defeq \logical{\pushforward}^1 \logical{\+\Bspline}_{\xx,ij}^1$, $\intermediate{\+\Bspline}_{\yy,ij}^1 \defeq \logical{\pushforward}^1 \logical{\+\Bspline}_{\yy,ij}^1$, $\intermediate{\Bspline}_{ij}^2 \defeq \logical{\pushforward}^2 \logical{\Bspline}_{ij}^2$, and in physical basis functions defined as the pushforward of the logical basis functions
\begin{equation}
  \Bspline_{ij}^0 \defeq \pushforward^0 \logical{\Bspline}_{ij}^0
  , \quad
  \+\Bspline_{\xx,ij}^1 \defeq \pushforward^1 \logical{\+\Bspline}_{\xx,ij}^1
  , \quad
  \+\Bspline_{\yy,ij}^1 \defeq \pushforward^1 \logical{\+\Bspline}_{\yy,ij}^1
  , \quad
  \Bspline_{ij}^2 \defeq \pushforward^2 \logical{\Bspline}_{ij}^2,
\end{equation}
which by \cref{eqn:pushforward_composition} are the pushforwards $\intermediate{\pushforward}^m$ of the intermediate basis functions.
As the pushforward operators are linear, the intermediate and physical forms can be expressed in terms of the intermediate and physical basis functions, respectively, by means of the same coefficients as the logical forms.
Example physical $0$-forms are shown in \cref{tab:disc_basis_functions_fourier}, for both a purely B-spline basis and a mixed B-spline and trigonometric basis.

\begin{example}[Physical one-form]
  The physical $1$-form $\+A \in \funspace{\Bspline}^1$ is related to the logical $1$-form $\logical{\+A} \in \logical{\funspace{\Bspline}}^1$ by the pushforward operator $\pushforward^1$, and can be expressed in terms of the physical basis functions as follows (cf. \cref{eqn:logical_oneforms})
  \begin{equation}
    \+A = \pushforward^1 \logical{\+A}
    = \sum_{(i,j)\in\indextwo^{1}_\xx} a_{\xx,ij} \+\Bspline_{\xx,ij}^1
    + \sum_{(i,j)\in\indextwo^{1}_\yy} a_{\yy,ij} \+\Bspline_{\yy,ij}^1,
  \end{equation}
  where the coefficients $a_{\xx,ij}$ and $a_{\yy,ij}$ are the same as those of the logical $1$-form $\logical{\+A}$ in \cref{eqn:logical_oneforms}.
\end{example}

\subsubsection{Geometric degrees of freedom}\label{sec:physical_geometric_dofs}
The intermediate GDOFs and projection operators are defined in terms of the logical ones applied to the pullback of the intermediate forms,
\begin{equation}
  \intermediate{\gdof}^m \defeq \logical{\gdof}^m \circ \logical{\pushforward}^{-m}
  , \quad
  \intermediate{\pi}^m \defeq \logical{\pushforward}^m \circ \logical{\pi}^m \circ \logical{\pushforward}^{-m},
\end{equation}
and the physical GDOFs and projection operators are similarly defined as
\begin{equation}\label{eqn:physical_geometric_dofs}
  \gdof^m \defeq \intermediate{\gdof}^m \circ \intermediate{\pushforward}^{-m} = \logical{\gdof}^m \circ \pushforward^{-m}
  , \quad
  \pi^m \defeq \intermediate{\pushforward}^m \circ \intermediate{\pi}^m \circ \intermediate{\pushforward}^{-m} = \pushforward^m \circ \logical{\pi}^m \circ \pushforward^{-m}.
\end{equation}
It follows that the projection operators $\pi^m$ are cochain maps, and we have
\begin{equation}\label{eqn:physical_commutativity}
  \pi^1 \circ \gradientt = \gradientt \circ \, \pi^0
  , \quad
  \pi^2 \circ \curlt = \curlt \circ \, \pi^1,
\end{equation}
thereby yielding commutativity of the diagram~\eqref{eqn:physical_commuting_diagram}; the same holds for the intermediate projection operators $\intermediate{\pi}^m$.


  \section{Polar subspaces of the 2D zero-form space}\label{sec:zeroforms}
The construction of the physical $m$-form spaces outlined in \cref{sec:physical_function_spaces} is applicable for a general mapping $\mapping$, assuming that the mapping is sufficiently smooth and invertible.
The disk-like logical mapping $\logicalmapping$ considered in \cref{sec:coordtransform}, however, collapses the edge $\logical{x} = 0$ of the logical domain to the center $\+0$ and is therefore not invertible there; the same holds for $\mapping=\intermediatemapping\circ\logicalmapping$, which results in \emph{nonconforming} physical function spaces $\funspace{\Bspline}^m$.
That is, the finite-dimensional physical function spaces are not subspaces of the corresponding infinite-dimensional function spaces appearing in the physical de Rham complex found in \cref{eqn:physical_commuting_diagram}.

The nonconformity can be read off ring by ring, where the $i$-th \emph{ring} of a tensor-product space consists of the basis functions, or coefficients, with radial index $i$, that is, the $i$-th row of \cref{tab:disc_basis_functions_fourier}.
A $0$-form basis function in the first ring with a nonconstant angular factor is discontinuous at the center, because its value there depends on the direction of approach, and its gradient grows like the inverse distance to the center.
A basis function in the second ring vanishes at the center and is therefore continuous, but its gradient is in general bounded without having a limit, so that it is not differentiable.
In general, a basis function in the $i$-th ring, $i \le r^{0,\xx}+1$, behaves like the $(i-1)$-st power of the distance to the center: it is $(i-2)$ times continuously differentiable there, while its derivatives of order $i-1$ are bounded but in general discontinuous (\cref{lem:vanishing_physical_derivative_general}); basis functions in later rings vanish on a neighborhood of the center.
To obtain $0$-forms that are $\pr$ times continuously differentiable at the center, the construction of polar subspaces therefore only has to modify the first $\pr+1$ rings, and the remaining rings can be kept as they are.
The same ring-by-ring gain of regularity, shifted by the Piola factors of the pushforwards, applies to the $1$- and $2$-form spaces in \cref{sec:derham_2d}.

In this section we give an overview of techniques from the literature for the construction of so-called \emph{polar} subspaces of the physical zero-form space $\funspace{\Bspline}^0$ of regularity $\pr$, denoted by $\polar{\funspace{\Bspline}}^{0,\pr}$, and establish conditions under which they satisfy
\begin{equation}\label{eqn:polar_zeroform_space_conforming}
  \polar{\funspace{\Bspline}}^{0,\pr} \subset \funspace{\Bspline}^0 \cap \funspace{H}_1 \cap \funspace{C}^{\pr}_{0}
  .
\end{equation}
Here $\funspace{C}^{\pr}_{0}$ denotes the functions that are $\pr$-times continuously differentiable on a physical neighborhood of the center $\+0$; the same convention defines $\funspace{C}^{\infty}_{0}$.
We use $\funspace{C}^{-1}_{0}$ to denote fields bounded on a punctured physical neighborhood of $\+0$, without requiring a value or limit at the center.
Since $\intermediatemapping$ is a smooth diffeomorphism with $\intermediatemapping(\+0)=\+0$, regularity at the center can equivalently be studied on the intermediate domain $\intermediate{\Omega}$; accordingly, the proofs in this section reduce to $\intermediatemapping=\mathrm{id}$ (so that $\mapping=\logicalmapping$ and $\jacmat=\logical{\jacmat}$) where this is justified, and treat general $\intermediatemapping$ explicitly otherwise.
\begin{assumption}[Logical space degrees and sizes]\label{ass:logical_space_parameters}
Throughout the paper, the parameters of the logical $0$-form space satisfy
\begin{equation}\label{eqn:polar_zeroform_space_conditions}
  0 \le \pr \le r^{0,\xx}, \quad
  r^{0,\xx} \ge 1, \quad
  n^{0,\xx} \ge \pr+2, \quad
  n^{\yy} \ge 2\pr+1.
\end{equation}
These assumptions ensure that the logical $0$-form space has sufficiently many degrees of freedom to construct a polar subspace with $\pr$-th order regularity at the center, while $n^{0,\xx} \ge \pr+2$ ensures that the outermost radial basis function is external, which is used only when essential boundary conditions are imposed in \cref{sec:boundary_homotopy}: it makes the mean projection $\mathfrak{B}$ in \cref{eqn:boundary_homotopy_operators_2d} well defined and is used in the proof of \cref{lem:boundary_compatible_homotopy}.
\end{assumption}

In essence, two techniques are used in the literature to construct polar subspaces of the physical zero-form space $\funspace{\Bspline}^0$.
One is based on characterizing the behavior of $0$-forms at the center $\logical{x} = 0$ and deriving explicit constraints on the coefficients of the physical $0$-forms~\citep{toshniwal2017}, thereby imposing their regularity, which we briefly discuss in \cref{sec:derivative_matching}.
The other is based on ensuring that the pullback of smooth physical $0$-forms can be exactly represented in the logical $0$-form space~\citep{jiang2026polar}, which we discuss in \cref{sec:pullback_compatibility}.

\subsection{Explicit constraints on the coefficients}\label{sec:derivative_matching}
The computations in this subsection only involve the disk-like mapping; they are carried out for $\intermediatemapping=\mathrm{id}$, so that $\mapping=\logicalmapping$ and the Jacobian matrix is $\logical{\jacmat}$ from \cref{eqn:disk_jacobian_matrix}.
For a general diffeomorphism $\intermediatemapping$ the set of admissible coefficients is identical, because $P\mapsto P\circ\intermediatemapping$ leaves the coefficients unchanged and, by the chain rule with the smooth invertible $\intermediatemapping$ and $\intermediatemapping(\+0)=\+0$, preserves continuity and differentiability at the center; only the prescribed gradient limit transforms: a gradient limit $\+g$ on the intermediate domain becomes $\intermediate{\jacmat}(\+0)^{-\transpose}\+g$ on the physical domain.

\paragraph{Continuity}
We consider $P \in \funspace{\Bspline}^0$ expressed in terms of its coefficients $p_{ij}$ (cf.\ \cref{eqn:logical_zeroforms}).
The value of the $0$-form at the center $\logical{x} = 0$ is given by
\begin{equation}\label{eqn:physical_zeroform_atpole}
  \lim_{\logical{x} \rightarrow 0} P 
  = \sum_{(i,j)\in\indextwo^{0}} p_{ij} \Bsplineone^{0,\xx}_i(0) \Bsplineone^{0,\yy}_j(\logical{y})
  = \sum_{j\in\indexone^{0,\yy}} p_{1j} \Bsplineone^{0,\yy}_j(\logical{y}),
\end{equation}
by making use of \cref{eqn:bspline_x_atpole}.
Hence, for general coefficients $p_{1j}$, the value of the $0$-form at the center is a function of the logical coordinate $\logical{y}$, and thus discontinuous at the center.
Imposing continuity of $P$ at the center thus puts a constraint on the coefficients $p_{1j}$: they should be chosen such that the expression on the right-hand side of \cref{eqn:physical_zeroform_atpole} no longer depends on $\logical{y}$.
By linear independence of the 1D basis functions, this is satisfied (nontrivially) by the constraint $p_{1j} = P(\+0) \mathds{1}_j^\yy$ (where we recall \cref{eqn:constant_function_coefficients}).

Hence, the polar $0$-form space $\polar{\funspace{\Bspline}}^{0,0}$ is defined as the subspace of $\funspace{\Bspline}^0$ where the first tensor-product `ring' of $0$-forms is replaced by a single polar $0$-form $\polar{\Bspline}^{0,0}_1$.
By the constraint $p_{1j} = P(\+0) \mathds{1}_j^\yy$, the polar $0$-form $\polar{\Bspline}^{0,0}_1$ is defined as (as function of logical coordinates)
\begin{equation}
  \polar{\Bspline}^{0,0}_1 \defeq \sum_{j\in\indexone^{0,\yy}} \mathds{1}_j^\yy \Bspline^{0}_{1j} = \Bsplineone^{0,\xx}_1,
\end{equation}
up to a multiplicative constant.

\paragraph{Differentiability}
\begin{table}
  \centering
  \setlength{\tabcolsep}{1.000000000000pt}
  \begin{tabular}{@{}cccccccccccc@{}}
    \toprule
    \multicolumn{12}{c}{\makebox[0pt][c]{The first two rings replaced by degree-one polar B-splines $\polar{\Bspline}^{0,1}_l$, $l\in \{1,2,3\}$}}\\
    \multicolumn{2}{c}{\tikzsetnextfilename{disc_basis_polar_degree1_1}\input{tikz/disc_basis_polar_degree1_1.tex}\unskip}&\multicolumn{2}{c}{\tikzsetnextfilename{disc_basis_polar_degree1_2}\input{tikz/disc_basis_polar_degree1_2.tex}\unskip}&\multicolumn{2}{c}{\tikzsetnextfilename{disc_basis_polar_degree1_3}\input{tikz/disc_basis_polar_degree1_3.tex}\unskip}&&&&&& \\
    \multicolumn{1}{c}{\tikzsetnextfilename{disc_basis_tensor_3_1}\input{tikz/disc_basis_tensor_3_1.tex}\unskip}&\multicolumn{1}{c}{\tikzsetnextfilename{disc_basis_tensor_3_2}\input{tikz/disc_basis_tensor_3_2.tex}\unskip}&\multicolumn{1}{c}{\tikzsetnextfilename{disc_basis_tensor_3_3}\input{tikz/disc_basis_tensor_3_3.tex}\unskip}&\multicolumn{1}{c}{\tikzsetnextfilename{disc_basis_tensor_3_4}\input{tikz/disc_basis_tensor_3_4.tex}\unskip}&\multicolumn{1}{c}{\tikzsetnextfilename{disc_basis_tensor_3_5}\input{tikz/disc_basis_tensor_3_5.tex}\unskip}&\multicolumn{1}{c}{\tikzsetnextfilename{disc_basis_tensor_3_6}\input{tikz/disc_basis_tensor_3_6.tex}\unskip}&\multicolumn{1}{c}{\tikzsetnextfilename{disc_basis_tensor_3_7}\input{tikz/disc_basis_tensor_3_7.tex}\unskip}&\multicolumn{1}{c}{\tikzsetnextfilename{disc_basis_tensor_3_8}\input{tikz/disc_basis_tensor_3_8.tex}\unskip}&\multicolumn{1}{c}{\tikzsetnextfilename{disc_basis_tensor_3_9}\input{tikz/disc_basis_tensor_3_9.tex}\unskip}&\multicolumn{1}{c}{\tikzsetnextfilename{disc_basis_tensor_3_10}\input{tikz/disc_basis_tensor_3_10.tex}\unskip}&\multicolumn{1}{c}{\tikzsetnextfilename{disc_basis_tensor_3_11}\input{tikz/disc_basis_tensor_3_11.tex}\unskip}&\multicolumn{1}{c}{\tikzsetnextfilename{disc_basis_tensor_3_12}\input{tikz/disc_basis_tensor_3_12.tex}\unskip} \\
\bottomrule
  \end{tabular}
  \caption{Comparison of the polar $0$-form spaces on the physical banana-shaped domain of \cref{fig:disk_mapping}, which are subspaces of the tensor-product physical $0$-form space shown in the top block of \cref{tab:disc_basis_functions_fourier}. The B-spline degree is $r = 3$ and the numbers of intervals are $N^\xx = 5$ and $N^\yy = 12$. Note that the polar $0$-forms have compact support due to the truncation discussed in \cref{sec:derham_hodge_2d_construction}.}
  \label{tab:disc_basis_functions}
\end{table}
When assuming that the $0$-form is continuous at the center, we can investigate its differentiability at the center by evaluating its gradient (cf.\ \cref{eqn:pushforward_operators})
\begin{equation}\label{eqn:bspline_dx_atpole_vanishing}
  \lim_{\logical{x} \rightarrow 0} \gradientt P
  = \frac{1}{R_1 (C S' - S C') \Delta t^{0,\xx}_1} \sum_{j\in\indexone^{0,\yy}} (p_{2j} - p_{1j})
  \begin{bmatrix}
    S' & -S \\
    -C' & C
  \end{bmatrix}
  \begin{bmatrix}
    \Bsplineone^{0,\yy}_j\\
    \totd{\Bsplineone^{0,\yy}_j}{\logical{y}}
  \end{bmatrix},
\end{equation}
where we have used \lhopital, substituted \cref{eqn:bspline_x_atpole}, and used the expression for the inverse transpose of the Jacobian matrix given in \cref{eqn:disk_jacobian_matrix}, together with the assumptions on our coordinate transformation in \cref{eqn:coordinate_transformation_assumptions}, which shows that
\begin{equation}\label{eqn:toroidal_jacmat_invtranspose}
  \logical{\jacmat}^{-\transpose}
  = \frac{1}{C S' - S C'}\begin{bmatrix}
     \frac{S \pard{R}{\logical{y}}}{R \pard{R}{\logical{x}}} + \frac{S'}{\pard{R}{\logical{x}}} & -\frac{S}{R}\\ 
    -\frac{C \pard{R}{\logical{y}}}{R \pard{R}{\logical{x}}} - \frac{C'}{\pard{R}{\logical{x}}} & \frac{C}{R}
  \end{bmatrix}
  \; \implies \;
  \lim_{\logical{x} \rightarrow 0} \logical{\jacmat}^{-\transpose}
  = \lim_{\logical{x} \rightarrow 0} \frac{1}{R_1 (C S' - S C')}\begin{bmatrix}
       S' & - \frac{S}{\logical{x}}\\ 
     - C' & \frac{C}{\logical{x}}
  \end{bmatrix}
  .
\end{equation}
To obtain a $0$-form that is continuously differentiable at the center, one has to find coefficients $p_{1j} = P(\+0) \mathds{1}_j^\yy, p_{2j}$ that yield a constant value on the right-hand side of \cref{eqn:bspline_dx_atpole_vanishing}.
Letting the constant be denoted by $\+g \in \mathbb{R}^2$, we can rearrange \cref{eqn:bspline_dx_atpole_vanishing} to obtain the following linear system for the coefficients $p_{2j}$
\begin{equation}\label{eqn:toshniwal_derivative_condition_Rinv}
  \begin{bmatrix}
    C & S\\
    C' & S'
  \end{bmatrix}
  \+g
  = 
  \frac{1}{R_1 \Delta t^{0,\xx}_1} \sum_{j\in\indexone^{0,\yy}} (p_{2j} - p_{1j})
  \begin{bmatrix}
    \Bsplineone^{0,\yy}_j  \\
    \totd{\Bsplineone^{0,\yy}_j}{\logical{y}}
  \end{bmatrix}.
\end{equation}
By noting that the second row of the constraint in \cref{eqn:toshniwal_derivative_condition_Rinv} equals the derivative in $\logical{y}$ of the first row, we find only one independent constraint on the coefficients.

The approach outlined above is followed in \citet{toshniwal2017}, where the value and gradient of three linearly independent triangular Bernstein polynomials are chosen to define three sets of coefficients $\polar{\bspline}^{0,1}_{l,1j}, \polar{\bspline}^{0,1}_{l,2j}$ for $l = 1, 2, 3$ that satisfy the constraints for continuity and differentiability at the center exactly.

Importantly, satisfying the constraint in \cref{eqn:toshniwal_derivative_condition_Rinv} equates the angular functions $C(\logical{y}), S(\logical{y})$ on the left-hand side to a linear combination of the 1D basis functions $\Bsplineone^{0,\yy}_j(\logical{y})$ on the right-hand side, thus implying that differentiability imposes a relation between the mapping and the choice of 1D basis functions in the $\logical{y}$-direction.
If $C, S$ are trigonometric functions then the 1D basis functions in the $\logical{y}$-direction should be chosen as trigonometric functions, and if B-splines are preferred then the mapping should be approximated by B-splines in the $\logical{y}$-direction, as proposed in \citet{toshniwal2017}.

Whatever approach is followed, it results in three linearly independent polar $0$-forms $\polar{\Bspline}^{0,1}_l$ for $l = 1, 2, 3$ (expressed in terms of the coefficients $\polar{\bspline}^{0,1}_{l,1j}, \polar{\bspline}^{0,1}_{l,2j}$), which replace the first two tensor-product rings of $0$-forms.
An example using the shifted Zernike polynomials from \cref{eqn:shifted_zernike_basis} (following \citet{jiang2026polar}) is given in \cref{tab:disc_basis_functions}.

\paragraph{Higher-order regularity}
In principle, the approach outlined above can be extended to high-order regularity at the center by evaluating high-order derivatives of the $0$-form at the center and imposing that they are constant.
For a general disk-like mapping, however, this becomes increasingly cumbersome, and the resulting constraints on the coefficients $p_{ij}$ become increasingly complex. Hence, we do not pursue this approach further in this work.
An example of a polar $0$-form space $\polar{\funspace{\Bspline}}^{0,2}$ with second-order regularity at the center is given in \cref{tab:disc_basis_functions}.

\subsection{Pullback compatibility}\label{sec:pullback_compatibility}
The alternative approach relies on splitting the tensor-product $0$-form space into two parts, and separately ensuring that each part is sufficiently regular at the center.
The splitting is based on the following lemma, which ensures that basis functions with a sufficiently large index $i$ vanish at the center, and thus do not require any special treatment.
The lemma allows for a smooth factor $\tilde F$ and a power $\logical{x}^q$ of the radial coordinate in front of the basis function: for $0$-forms, $\tilde F=1$ and $q=0$, whereas the Piola factors of the $1$- and $2$-form pushforwards in \cref{sec:derham_2d} contribute a smooth factor and, through the inverse of the Jacobian, the power $q=-1$.

\begin{lemma}[Vanishing physical derivative]\label{lem:vanishing_physical_derivative_general}
  Let $\funspace{\Bsplineone}^\xx_{r^\xx}$ denote a clamped B-spline space of degree $r^\xx$ in the $\logical{x}$-direction with basis functions $\Bsplineone^\xx_i$ in the standard ordering (positive rescalings of the basis functions are allowed).
  Fix an integer $n \ge 0$ and assume that $\Bsplineone^\yy_j$ is $n$-times continuously differentiable in the periodic $\logical{y}$-direction.
  We also assume that $R/\logical{x}$, $C$, $S$, and $\tilde F$ extend smoothly to a closed radial-angular neighborhood of $\logical{x}=0$, with the angular seam identified periodically.
  For an integer $q$, consider the physical function
  \begin{equation}\label{eqn:vanishing_physical_derivative_general_function}
    \Bspline_{ij} \defeq \srpar{\tilde{F}(\logical{x}, \logical{y}) \, \logical{x}^q \, \Bsplineone^{\xx}_i(\logical{x}) \, \Bsplineone^{\yy}_j(\logical{y})} \circ \mapping^{-1}
    , \quad
    \Bsplineone^\xx_i \in \funspace{\Bsplineone}^\xx_{r^\xx}.
  \end{equation}
  If $i > n-q+1$, then $\Bspline_{ij}$ has a $\funspace{C}^n$ extension near the center, and all of its physical derivatives of order at most $n$ vanish there; in particular,
  \begin{equation}\label{eqn:vanishing_physical_derivative_general}
    \frac{\partial^{n} \Bspline_{ij}}{\partial x^\alpha \partial y^\beta}(\+0) = 0
    , \quad 
    i > n - q + 1,
  \end{equation}
  for $\alpha+\beta=n$.
  If $i=n-q+1$ and $n\ge1$, the same holds with $n$ replaced by $n-1$, while the physical derivatives of order $n$ are bounded on a punctured neighborhood of the center, but need not exist at the center.
\end{lemma}
\begin{statementproof}[appendix]{lem:vanishing_physical_derivative_general}
  We first assume $\intermediatemapping=\mathrm{id}$, so that $\mapping=\logicalmapping$ and $\jacmat=\logical{\jacmat}$, and comment on a general intermediate mapping at the end.
  The main observation is that every physical derivative costs at most one power of $\logical{x}$ in the pullback, with coefficients that stay bounded up to $\logical{x}=0$.
  Since $\logical{x}$ is comparable to the physical distance $|\+x|$ to the center, a $k$-th order physical derivative of $\Bspline_{ij}$ is then $\mathcal{O}(|\+x|^{m-k})$ with $m \defeq q+i-1$.

  \emph{Factorization.}
  If $i>r^\xx+1$, then $\Bsplineone^\xx_i$ vanishes on the first radial knot span $[0,h)$, so $\Bspline_{ij}$ vanishes near $\+0$ (also when $q<0$) and there is nothing to prove.
  If $i\le r^\xx+1$, the knot $0$ has multiplicity $r^\xx+2-i$ in the local knot vector of $\Bsplineone^\xx_i$, so $\Bsplineone^\xx_i$ has a zero of order $i-1$ at $\logical{x}=0$ \citep[Ch.~IX]{deboor1978}, that is, $\Bsplineone_i^\xx = \logical{x}^{i-1} b_i$ on $[0,h)$ for a polynomial $b_i$.
  Hence
  \begin{equation}\label{eqn:vanishing_physical_derivative_factorization}
    \Bspline_{ij}\circ\mapping = \logical{x}^m H_{00},
    \qquad
    H_{00} \defeq \tilde F\, b_i\, \Bsplineone_j^\yy,
  \end{equation}
  where $H_{00}$ is $\funspace{C}^n$ up to $\logical{x}=0$.
  Away from the center the mapping is a smooth diffeomorphism (\cref{ass:disc_like_mapping}), so $\Bspline_{ij}$ is $\funspace{C}^n$ on $\mapping\bigl((0,h)\times[0,1]\bigr)$, where $\Bsplineone^\xx_i$ is a polynomial, and only the behavior at $\+0$ remains to be studied.
  From now on we restrict to $0<\logical{x}<\varepsilon$ for some $\varepsilon\le h$ for which the bounds below hold.

  \emph{Derivative recursion (key step).}
  Factoring the $1/\logical{x}$ growth out of the angular entries of $\jacmat^{-\transpose}=\logical{\jacmat}^{-\transpose}$ (cf.\ \cref{eqn:toroidal_jacmat_invtranspose}) gives
  \begin{equation}\label{eqn:physical_derivative_logical_derivative}
    \pard{}{x} = \matcomp{\Gamma}_{11} \pard{}{\logical{x}} + \frac{\matcomp{\Gamma}_{12}}{\logical{x}} \pard{}{\logical{y}}
    , \quad
    \pard{}{y} = \matcomp{\Gamma}_{21} \pard{}{\logical{x}} + \frac{\matcomp{\Gamma}_{22}}{\logical{x}} \pard{}{\logical{y}}
    ,
  \end{equation}
  where, with $\tilde{R} \defeq R / \logical{x}$,
  \begin{equation}\label{eqn:gamma_matrix}
    \mat{\Gamma}
    \defeq \jacmat^{-\transpose} \begin{bmatrix}
      1 & 0\\
      0 & \logical{x}
    \end{bmatrix}
    = \frac{1}{\tilde{R} (C S' - S C')} \begin{bmatrix}
      \frac{S \pard{\tilde{R}}{\logical{y}} + \tilde{R} S'}{ {\tilde{R} + \logical{x} \pard{\tilde{R}}{\logical{x}}}} & -S\\
      \frac{- C \pard{\tilde{R}}{\logical{y}} - \tilde{R} C'}{ {\tilde{R} + \logical{x} \pard{\tilde{R}}{\logical{x}}}} & C
    \end{bmatrix}
    .
  \end{equation}
  The denominators are bounded away from zero near $\logical{x}=0$ because $\tilde R(0,\logical{y})=R_1>0$ and $CS'-SC'>0$ (\cref{eqn:coordinate_transformation_assumptions}), so the entries of $\mat{\Gamma}$ are smooth up to $\logical{x}=0$.
  Applying \cref{eqn:physical_derivative_logical_derivative} to $\logical{x}^{p}H$ gives
  \begin{equation}\label{eqn:vanishing_physical_derivative_recursion}
    \pard{}{x}\bigl(\logical{x}^{p}H\bigr)
    = \logical{x}^{p-1}\left(\matcomp{\Gamma}_{11}\Bigl(pH + \logical{x}\pard{H}{\logical{x}}\Bigr) + \matcomp{\Gamma}_{12}\pard{H}{\logical{y}}\right),
  \end{equation}
  and the same with $\matcomp{\Gamma}_{21},\matcomp{\Gamma}_{22}$ for $\pard{}{y}$.
  The two terms balance: the radial derivative lowers the power of $\logical{x}$ by at most one and has a bounded coefficient, whereas the angular derivative keeps the power but has a coefficient of order $1/\logical{x}$.
  The bracket in \cref{eqn:vanishing_physical_derivative_recursion} loses one order of regularity compared to $H$.
  Starting from \cref{eqn:vanishing_physical_derivative_factorization} and applying this $k\le n$ times shows that $(\partial_x^\alpha\partial_y^\beta\Bspline_{ij})\circ\mapping = \logical{x}^{m-k}H_{\alpha\beta}$, $\alpha+\beta=k$, where $H_{\alpha\beta}$ is $\funspace{C}^{n-k}$ and in particular bounded up to $\logical{x}=0$.
  Thus
  \begin{equation}\label{eqn:vanishing_physical_derivative_bound}
    \left|\left(\frac{\partial^k\Bspline_{ij}}{\partial x^\alpha\partial y^\beta}\right)
      \bigl(\mapping(\logical{x},\logical{y})\bigr)\right|
    \le c_{\alpha\beta}\logical{x}^{m-k},
    \qquad 0<\logical{x}<\varepsilon,\quad k=\alpha+\beta\le n,
  \end{equation}
  with constants independent of $\logical{y}$.

  \emph{Physical bound.}
  Since $|\mapping| = \logical{x}\tilde R\sqrt{C^2+S^2}$, where $\tilde R$ is bounded away from zero and $C,S$ do not vanish simultaneously (otherwise $CS'-SC'=0$), there are constants $0<c_1\le c_2$ such that $c_1\logical{x}\le|\mapping(\logical{x},\logical{y})|\le c_2\logical{x}$ for $0\le\logical{x}\le\varepsilon$, uniformly in $\logical{y}$.
  Moreover, for some $\delta>0$ every $\+x$ with $0<|\+x|<\delta$ has a unique preimage, and it satisfies $0<\logical{x}<\varepsilon$.
  Indeed, along each ray of constant $\logical{y}$, $|\mapping|$ grows continuously from $0$ to at least $c_1\varepsilon$, while the angle of $(C,S)$ increases strictly with $\logical{y}$ (because $CS'-SC'>0$) and winds exactly once (because $\mapping$ is injective for $\logical{x}>0$).
  Hence all $\+x$ with $0<|\+x|<c_1\varepsilon$ are attained with $0<\logical{x}<\varepsilon$.
  Preimages with $\logical{x}\ge\varepsilon$ are excluded by injectivity once $\delta$ is also below the positive minimum of $|\mapping|$ on the compact set $[\varepsilon,1]\times[0,1]$.
  Using $\logical{x}\le|\+x|/c_1$ when $m-k\ge0$ and $\logical{x}\ge|\+x|/c_2$ when $m-k<0$, \cref{eqn:vanishing_physical_derivative_bound} becomes
  \begin{equation}\label{eqn:vanishing_physical_derivative_bound_physical}
    \left|\frac{\partial^k\Bspline_{ij}}{\partial x^\alpha\partial y^\beta}(\+x)\right|
    \le c'_{\alpha\beta}|\+x|^{m-k},
    \qquad 0<|\+x|<\delta,\quad k=\alpha+\beta\le n,
  \end{equation}
  for every direction of approach to the center.

  \emph{Conclusion.}
  If $i>n-q+1$, then $m>n$, so all derivatives of order at most $n$ tend to zero at $\+0$ by \cref{eqn:vanishing_physical_derivative_bound_physical}.
  We show inductively in $k<n$ that the first derivatives of $G\defeq\partial_x^\alpha\partial_y^\beta\Bspline_{ij}$, $\alpha+\beta=k$, extend continuously to $\+0$ with value zero, assuming that $G$ itself does (for $k=0$ this was just shown).
  By \cref{eqn:vanishing_physical_derivative_bound_physical}, the difference quotients $|G(t,0)-G(\+0)|/|t|$ and $|G(0,t)-G(\+0)|/|t|$ are bounded by $c'_{\alpha\beta}|t|^{m-k-1}$, which tends to zero as $t\to0$ because $m-k-1>n-k-1\ge0$.
  Hence $\partial_xG(\+0)=\partial_yG(\+0)=0$, and these derivatives are continuous at $\+0$ because, by \cref{eqn:vanishing_physical_derivative_bound_physical} of order $k+1$, they tend to zero as $\+x\to\+0$.
  This gives the $\funspace{C}^n$ extension and \cref{eqn:vanishing_physical_derivative_general}.

  If $i=n-q+1$ and $n\ge1$, then $i>(n-1)-q+1$, and the case just proven with $n-1$ in place of $n$ gives the $\funspace{C}^{n-1}$ extension.
  Since $m=n$, \cref{eqn:vanishing_physical_derivative_bound_physical} only gives bounded derivatives of order $n$ on the punctured neighborhood.
  These need not have a value at $\+0$: for the disk mapping (\cref{ex:disk_mapping}) with $R_1=1$, linear radial B-splines on a single interval ($\Bsplineone^\xx_2=\logical{x}$), the constant trigonometric function $\Bsplineone^\yy_1=1$, $n=1$, $q=0$ and $\tilde F=1$, we get $\Bspline_{21}(\+x)=|\+x|$, which is not differentiable at $\+0$.

  \emph{General intermediate mapping.}
  For a general $\intermediatemapping$, $\Bspline_{ij}=\intermediate{\Bspline}_{ij}\circ\intermediatemapping^{-1}$, where $\intermediate{\Bspline}_{ij}$ is the function \cref{eqn:vanishing_physical_derivative_general_function} obtained for $\intermediatemapping=\mathrm{id}$.
  Since $\intermediatemapping^{-1}$ is smooth and fixes the center (\cref{ass:diffeomorphism}), the chain rule expresses the physical derivatives of order $k$ of $\Bspline_{ij}$ as combinations, with smooth coefficients, of the derivatives of order at most $k$ of $\intermediate{\Bspline}_{ij}$.
  Hence the $\funspace{C}^n$ extension, the vanishing derivatives at the center and the bounded derivatives on a punctured neighborhood carry over from $\intermediate{\Bspline}_{ij}$ to $\Bspline_{ij}$.
\end{statementproof}

Motivated by this result, we split the physical (and equivalently, the intermediate and logical) $0$-form space into two subspaces, the internal part containing the first $\pr + 1$ rings of tensor-product $0$-forms and the external part containing the remaining rings of tensor-product $0$-forms
\begin{equation}\label{eqn:polar_zeroform_space_odot}
  {\funspace{\Bspline}}^{0} = {\funspace{\Bspline}}^{0,\pr}_\internal \oplus {\funspace{\Bspline}}^{0,\pr}_\external
  , \quad
  {\funspace{\Bspline}}^{0,\pr}_\internal \defeq \underset{(i,j) \in \indextwo^{0,\pr}_\internal}{\spn}\{{\Bspline}^0_{ij}\}
  , \quad
  {\funspace{\Bspline}}^{0,\pr}_\external \defeq \underset{(i,j) \in \indextwo^{0,\pr}_\external}{\spn}\{{\Bspline}^0_{ij}\}
  ,
\end{equation}
where the indexing sets are defined as
\begin{equation}
  \indextwo^{0,\pr}_\internal \defeq \{(i,j) \in \indextwo^0 : i \le \pr + 1\}
  , \quad
  \indextwo^{0,\pr}_\external \defeq \indextwo^0 \setminus \indextwo^{0,\pr}_\internal.
\end{equation}
When the angular scalar basis is $\funspace{C}^\pr$, \cref{lem:vanishing_physical_derivative_general} shows that every $0$-form in $\funspace{\Bspline}^{0,\pr}_\external$ has a $\funspace{C}^\pr$ extension near the center whose physical derivatives through order $\pr$ vanish there (take $\tilde{F}=1$ and $q=0$).

The same splitting applies to the tensor-product $1$- and $2$-form spaces and is needed from \cref{sec:derham_2d} onward,
\begin{equation}\label{eqn:polar_mform_space_odot}
  \funspace{\Bspline}^{m} = \funspace{\Bspline}^{m,\pr}_\internal \oplus \funspace{\Bspline}^{m,\pr}_\external
  , \quad m \in \{0, 1, 2\},
\end{equation}
where the internal part is spanned by the basis functions of the first rings and the external part by the remaining ones: for $m=0$ as in \cref{eqn:polar_zeroform_space_odot}, and
\begin{equation}\label{eqn:polar_oneform_twoform_space_internal}
  \funspace{\Bspline}^{1,\pr}_\internal \defeq \spn\{\+\Bspline^1_{\xx,ij} : i \le \pr\} \oplus \spn\{\+\Bspline^1_{\yy,ij} : i \le \pr+1\}
  , \quad
  \funspace{\Bspline}^{2,\pr}_\internal \defeq \spn\{\Bspline^2_{ij} : i \le \pr\}.
\end{equation}
The corresponding internal indexing sets $\indextwo^{1,\pr}_{\xx,\internal}$, $\indextwo^{1,\pr}_{\yy,\internal}$ and $\indextwo^{2,\pr}_\internal$, and the external indexing sets $\indextwo^{1,\pr}_{\xx,\external}$, $\indextwo^{1,\pr}_{\yy,\external}$ and $\indextwo^{2,\pr}_\external$ of the remaining basis functions, are listed in \cref{tab:internal_rings}.
The numbers of internal rings follow from \cref{lem:vanishing_physical_derivative_general}: the polar $m$-forms are to be $\pr-m$ times continuously differentiable at the center, and the Piola factor of the pushforward in \cref{eqn:pushforward_operators} costs one power of $\logical{x}$ ($q=-1$) for the angular $1$-form component and for $2$-forms, and none ($q=0$) for $0$-forms and for the radial $1$-form component, so that $\pr-m-q+1$ rings are internal.
By the same lemma, the external $m$-forms are $\pr-m$ times continuously differentiable at the center with vanishing derivatives, provided that the angular basis is sufficiently smooth; see \cref{thm:polar_subspaces_regularity} for the precise statement.
\begin{table}
  \centering
  \small
  \begin{tabular}{@{}cllcll@{}}
    \toprule
    $m$ & basis functions & radial factor & $q$ & internal indexing set & external indexing set\\
    \midrule
    $0$ & $\Bspline^0_{ij}$ & $\Bsplineone^{0,\xx}_i$ & $0$ & $\indextwo^{0,\pr}_\internal = \{(i,j)\in\indextwo^0 : i\le\pr+1\}$ & $\indextwo^{0,\pr}_\external = \indextwo^0\setminus\indextwo^{0,\pr}_\internal$\\
    $1$ & $\+\Bspline^1_{\xx,ij}$ & $\Bsplineone^{1,\xx}_i$ & $0$ & $\indextwo^{1,\pr}_{\xx,\internal} = \{(i,j)\in\indextwo^1_\xx : i\le\pr\}$ & $\indextwo^{1,\pr}_{\xx,\external} = \indextwo^1_\xx\setminus\indextwo^{1,\pr}_{\xx,\internal}$\\
    $1$ & $\+\Bspline^1_{\yy,ij}$ & $\Bsplineone^{0,\xx}_i$ & $-1$ & $\indextwo^{1,\pr}_{\yy,\internal} = \{(i,j)\in\indextwo^1_\yy : i\le\pr+1\}$ & $\indextwo^{1,\pr}_{\yy,\external} = \indextwo^1_\yy\setminus\indextwo^{1,\pr}_{\yy,\internal}$\\
    $2$ & $\Bspline^2_{ij}$ & $\Bsplineone^{1,\xx}_i$ & $-1$ & $\indextwo^{2,\pr}_\internal = \{(i,j)\in\indextwo^2 : i\le\pr\}$ & $\indextwo^{2,\pr}_\external = \indextwo^2\setminus\indextwo^{2,\pr}_\internal$\\
    \bottomrule
  \end{tabular}
  \caption{Splitting of the tensor-product $m$-form spaces into internal and external parts, \cref{eqn:polar_mform_space_odot}: for each family of basis functions, its radial factor, the power $q$ of $\logical{x}$ contributed by the pushforward, and the internal and external indexing sets, with $\pr-m-q+1$ internal rings.}
  \label{tab:internal_rings}
\end{table}

The remainder of the polar $0$-form space is defined by projecting $n^{0,\pr}$ linearly independent smooth physical $0$-forms onto the physical $0$-form space $\funspace{\Bspline}^0$ using the commuting projection operators $\pi^0$ defined in \cref{eqn:physical_geometric_dofs}, where we note that owing to the continuity of these physical $0$-forms, the projection operator is well-defined.
Instead of physical polynomials, we use polynomials in the intermediate coordinates pushed forward to the physical domain.
Let $\intermediate{\funspace{T}}^\pr$ denote the space of polynomials of total degree at most $\pr$ in the intermediate coordinates $\intermediatecoord{\+x}$ on $\intermediate{\Omega}$, and let
\begin{equation}
  \pulledpolys{0}{\pr} \defeq \intermediate{\pushforward}^0 \intermediate{\funspace{T}}^\pr = \{\intermediate{T}\circ\intermediatemapping^{-1} : \intermediate{T}\in\intermediate{\funspace{T}}^\pr\}
\end{equation}
be the corresponding physical pushed-forward polynomials; both spaces have dimension given by $n^{0,\pr} = {(\pr+1)(\pr+2)}/{2}$.
Since $\intermediatemapping$ is smooth (\cref{ass:diffeomorphism}), the functions in $\pulledpolys{0}{\pr}$ are smooth on $\Omega$.
This results in
\begin{equation}\label{eqn:polar_zeroform_space}
  \polar{\funspace{\Bspline}}^{0,\pr} \defeq \pi^0 \pulledpolys{0}{\pr} \oplus \funspace{\Bspline}^{0,\pr}_\external.
\end{equation}

The following theorem shows that the polar $0$-form space $\polar{\funspace{\Bspline}}^{0,\pr}$ is conforming with the infinite-dimensional function space $\funspace{H}_1$.
\begin{theorem}[Conformity]\label{thm:polar_zeroform_space_conforming}
  The polar $0$-form space $\polar{\funspace{\Bspline}}^{0,\pr}$ defined in \cref{eqn:polar_zeroform_space} is conforming for any $0\le\pr\le r^{0,\xx}$
  \begin{equation}
    \polar{\funspace{\Bspline}}^{0,\pr} \subset \funspace{H}_1.
  \end{equation}
\end{theorem}
\begin{statementproof}[appendix]{thm:polar_zeroform_space_conforming}
  Let $P\in\polar{\funspace{\Bspline}}^{0,\pr}$, where we first consider $\intermediatemapping=\mathrm{id}$.
  By \cref{eqn:polar_zeroform_space}, write $P=Q+P_\external$ with $Q=\pi^0T$ for some $T\in\pulledpolys{0}{\pr}$ and $P_\external\in\funspace{\Bspline}^{0,\pr}_\external$.
  Every basis function in $P_\external$ has radial index $i\ge\pr+2\ge2$, so its value does not contribute to the first coefficient ring.

  Write $Q=\sum_{(i,j)\in\indextwo^0}q_{ij}\Bspline^0_{ij}$.
  Preservation of the $0$-form GDOFs by $\pi^0$ gives
  \begin{equation}
    \sum_{(i,j)\in\indextwo^0}
      q_{ij}\Bsplineone^{0,\xx}_i(\xi^\xx_{i'})
             \Bsplineone^{0,\yy}_j(\xi^\yy_{j'})
    =T\bigl(\mapping(\xi^\xx_{i'},\xi^\yy_{j'})\bigr),
    \qquad (i',j')\in\indextwo^0.
  \end{equation}
  Its matrix is the Kronecker product of the radial and angular scalar evaluation matrices and is invertible by the 1D unisolvence established in \cref{lem:unisolvence_1d_gdofs}.
  At the first radial node $\xi^\xx_1=0$, the endpoint values in \cref{eqn:bspline_x_atpole} and $\mapping(0,\logical{y})=\+0$ reduce this system to
  \begin{equation}\label{eqn:polar_zeroform_space_conforming_proof}
    \sum_{j\in\indexone^{0,\yy}}
      q_{1j}\Bsplineone^{0,\yy}_j(\xi^\yy_{j'})
    =T(\+0),
    \qquad j'\in\indexone^{0,\yy}.
  \end{equation}
  The coefficients $q_{1j}=T(\+0)\mathds{1}_j^\yy$ solve this system by \cref{eqn:constant_function_coefficients}; angular unisolvence makes them its unique solution.
  Since $P_\external$ has no first ring, the first ring of $P$ therefore represents the constant $c=T(\+0)$.

  Let $\logical{P}=P\circ\mapping$.
  On the first radial knot span, each $\Bsplineone_i^{0,\xx}$ is a polynomial, $\Bsplineone_1^{0,\xx}(0)=1$, and $\Bsplineone_i^{0,\xx}(0)=0$ for $i\ge2$.
  The constant first ring consequently gives
  \begin{equation}
    \logical{P}(\logical{x},\logical{y})
      =c+\logical{x}a(\logical{x},\logical{y}),
  \end{equation}
  where $a$ and its first piecewise derivatives are bounded uniformly in $\logical{y}$ near $\logical{x}=0$.
  For $P_\external$ alone the same representation holds with $c=0$.
  By \cref{eqn:physical_derivative_logical_derivative}, at points away from knot interfaces and the center,
  \begin{equation}
    (\gradientt P)\circ\mapping
      =\mat{\Gamma}
        \begin{bmatrix}
          a+\logical{x}\pard{a}{\logical{x}}\\
          \pard{a}{\logical{y}}
        \end{bmatrix}.
  \end{equation}
  The matrix $\mat{\Gamma}$ is bounded near the center by \cref{eqn:gamma_matrix}.
  Thus the piecewise physical gradient is bounded near $\+0$, while the representation $\logical{P}=c+\logical{x}a$ shows that $P$ extends continuously to the center with value $c$.
  Away from the center, the mapping is a smooth diffeomorphism and the basis functions are continuous and piecewise $\funspace{C}^1$, so $P$ and its piecewise gradient are bounded on the remaining physical knot cells.
  Therefore $P\in\funspace{H}_1$.

  Finally, let $\intermediatemapping\neq\mathrm{id}$.
  The GDOFs of $\pi^0(\intermediate{T}\circ\intermediatemapping^{-1})$ are the values $\intermediate{T}\bigl(\logicalmapping(\xi^\xx_{i'},\xi^\yy_{j'})\bigr)$, which do not depend on $\intermediatemapping$.
  Hence $P\circ\intermediatemapping$ is the function with the same coefficients in the polar $0$-form space constructed for $\intermediatemapping=\mathrm{id}$, and therefore belongs to $\intermediate{\funspace{H}}_1$ by the above.
  By the chain rule, $\gradientt P=\bigl(\intermediate{\jacmat}^{-\transpose}\intermediate{\gradientt}(P\circ\intermediatemapping)\bigr)\circ\intermediatemapping^{-1}$, which also holds for weak gradients.
  Since $\intermediate{\jacmat}^{-1}$ and $\intermediate{\jac}$ are bounded (\cref{ass:diffeomorphism}), the change of variables $\diff\+x=\intermediate{\jac}\diff\intermediatecoord{\+x}$ shows that $P$ and $\gradientt P$ are square integrable because $P\circ\intermediatemapping$ and its gradient are.
  Hence $P\in\funspace{H}_1$.
\end{statementproof}

Similar to the differentiability constraint given by \cref{eqn:toshniwal_derivative_condition_Rinv}, high-order regularity imposes a relation between the logical mapping and the choice of 1D basis functions: the pullback of the polynomials $\intermediate{\funspace{T}}^\pr$ should be exactly representable in the logical $0$-form space.
\begin{definition}[$\pr$-compatibility of the mapping and logical $0$-form space]\label{def:pr_compatibility_mapping_logical_zeroform}
  We say that the logical mapping $\logicalmapping$ and the logical $0$-form space $\logical{\funspace{\Bspline}}^0$ are \emph{$\pr$-compatible} if the pullback of the polynomials in the intermediate coordinates can be exactly represented in the logical $0$-form space
  \begin{equation}\label{eqn:pr_compatibility_zeroform}
    \logical{\pushforward}^{-0} \intermediate{\funspace{T}}^\pr \subset \logical{\funspace{\Bspline}}^0.
  \end{equation}
  Equivalently, $\pushforward^{-0}\pulledpolys{0}{\pr}\subset\logical{\funspace{\Bspline}}^0$, since $\pushforward^{-0}(\intermediate{T}\circ\intermediatemapping^{-1}) = \intermediate{T}\circ\logicalmapping = \logical{\pushforward}^{-0}\intermediate{T}$.
\end{definition}
Compatibility depends only on the logical mapping $\logicalmapping$ and the logical $0$-form space, not on the diffeomorphism $\intermediatemapping$.
The physical polar spaces do depend on the decomposition $\mapping=\intermediatemapping\circ\logicalmapping$ of a given physical mapping.
By \cref{ass:constant_functions}, every logical mapping and logical $0$-form space are $0$-compatible, so this definition only restricts the mapping and the spaces for $\pr\ge1$.

\begin{example}[$\pr$-compatibility of the disk mapping]\label{ex:pr_compatibility_disc_mapping}
  We consider the disk mapping $\logicalmapping$ defined in \cref{ex:disk_mapping}.
  The pullback of a polynomial $\intermediate{T} \in \intermediate{\funspace{T}}^\pr$ of degree $a + b \le \pr$ is given by
  \begin{equation}
    \intermediate{T} = \intermediatecoord{x}^a \intermediatecoord{y}^b
    \quad \implies \quad
    \logical{\pushforward}^{-0} \intermediate{T}
    = \intermediate{T} \circ \logicalmapping
    = (R_1 \logical{x})^{a+b} \cos^a (2\pi\logical{y}) \sin^b (2\pi\logical{y}),
  \end{equation}
  thus showing that a B-spline basis in the $\logical{x}$-direction of degree at least $r^{0,\xx} \ge \pr$ is required (because the clamped B-spline space of degree $r^{0,\xx}$ can represent univariate polynomials up to degree $r^{0,\xx}$ exactly), and that a trigonometric basis in the $\logical{y}$-direction is required containing the first $2\pr+1$ modes.
\end{example}

Provided that the logical mapping and the logical $0$-form space are $\pr$-compatible, we can show the desired regularity of the polar $0$-form space $\polar{\funspace{\Bspline}}^{0,\pr}$.
For simple-knot periodic angular B-splines, the angular regularity used below holds when $r^{0,\yy}\ge\pr+1$; for trigonometric functions it holds automatically.
\begin{theorem}[Regularity]\label{thm:polar_zeroform_space_regularity}
  Let $0\le\pr\le r^{0,\xx}$ and assume that the logical mapping $\logicalmapping$ and the logical $0$-form space $\logical{\funspace{\Bspline}}^0$ are $\pr$-compatible.
  If $\pr<r^{0,\xx}$, assume also that the angular scalar basis is $\funspace{C}^\pr$.
  Then the polar $0$-forms from $\polar{\funspace{\Bspline}}^{0,\pr}$ defined in \cref{eqn:polar_zeroform_space} satisfy
  \begin{equation}
    \polar{\funspace{\Bspline}}^{0,\pr} \subset \funspace{C}^\pr_0.
  \end{equation}

  Moreover, if $\pr = r^{0,\xx}$ then the polar $0$-forms are smooth near the center
  \begin{equation}
    \polar{\funspace{\Bspline}}^{0,r^{0,\xx}} \subset \funspace{C}^\infty_0.
  \end{equation}
\end{theorem}
\begin{statementproof}[appendix]{thm:polar_zeroform_space_regularity}
  Let $P\in\polar{\funspace{\Bspline}}^{0,\pr}$ and write $P=\pi^0T+P_\external$, where $T\in\pulledpolys{0}{\pr}$ and $P_\external\in\funspace{\Bspline}^{0,\pr}_\external$.
  Write $T=\intermediate{T}\circ\intermediatemapping^{-1}$ with $\intermediate{T}\in\intermediate{\funspace{T}}^\pr$.
  Then $\pushforward^{-0}T=T\circ\intermediatemapping\circ\logicalmapping=\intermediate{T}\circ\logicalmapping=\logical{\pushforward}^{-0}\intermediate{T}\in\logical{\funspace{\Bspline}}^0$ by $\pr$-compatibility, so $T=\pushforward^{0}\pushforward^{-0}T\in\funspace{\Bspline}^0$.
  Since $\pi^0$ is a projector onto $\funspace{\Bspline}^0$, $\pi^0T=T$, which is smooth on a neighborhood of $\Omega$ because $\intermediate{T}$ is a polynomial and $\intermediatemapping^{-1}$ is smooth (\cref{ass:diffeomorphism}).

  If $\pr<r^{0,\xx}$, every basis function in $P_\external$ has radial index $i\ge\pr+2$.
  By \cref{lem:vanishing_physical_derivative_general} with $n=\pr$, $q=0$, and $\tilde F=1$, each such function has a $\funspace{C}^\pr$ extension near $\+0$, and all its physical derivatives through order $\pr$ vanish there.
  Thus $P=T+P_\external$ belongs to $\funspace{C}^\pr_0$, with its derivatives at $\+0$ through order $\pr$ supplied by the smooth function $T$.

  If $\pr=r^{0,\xx}$, every external basis function has $i\ge r^{0,\xx}+2$ and therefore vanishes on the first radial knot span $[0,h)\times\mathbb{S}^1$.
  Its image under $\logicalmapping$ contains an intermediate neighborhood of $\+0$ (see the proof of \cref{lem:vanishing_physical_derivative_general}), whose image under the diffeomorphism $\intermediatemapping$ is a physical neighborhood of $\intermediatemapping(\+0)=\+0$.
  Hence $P_\external$ is zero on a physical neighborhood of $\+0$, where $P=T$ is smooth; in particular, $P\in\funspace{C}^\infty_0$.
  In either case, the vanishing of the external component through order $\pr$ also shows that the sum in \cref{eqn:polar_zeroform_space} is direct: if $T=\pi^0T\in\pulledpolys{0}{\pr}$ also lies in $\funspace{\Bspline}^{0,\pr}_\external$, then all derivatives of $T$ through order $\pr$ vanish at $\+0$; by the chain rule and $\intermediatemapping(\+0)=\+0$, so do those of the polynomial $\intermediate{T}=T\circ\intermediatemapping$ at $\+0$, which has degree at most $\pr$ and is therefore zero.
\end{statementproof}

\begin{remark}[How to decompose the physical mapping]\label{rem:pr_compatibility_sufficient}
  We consider two different mapping decompositions resulting in the same composed mapping $\mapping$.


  \paragraph{Based on the disk mapping.}
  Let $\logicalmapping_A$ be the disk mapping of \cref{ex:disk_mapping}, and let $\intermediatemapping_A$ be given by
  \begin{equation}\label{eqn:pr_compatibility_sufficient_deformation}
    \intermediatemapping_A(\intermediatecoord{\+x}) = \sspar{1 + \frac{0.1}{R_1^3} \intermediatecoord{x} (\intermediatecoord{x}^2 - 3 \intermediatecoord{y}^2)}\intermediatecoord{\+x},
  \end{equation}
  where the composed mapping is given by
  \begin{equation}
    \mapping_A
    = \intermediatemapping_A\circ\logicalmapping_A
    = R_1 \logical{x} \sspar{1 + 0.1 \logical{x}^3 \cos(6\pi \logical{y})} \begin{bmatrix}
      \cos(2\pi \logical{y})\\
      \sin(2\pi \logical{y})
    \end{bmatrix},
  \end{equation}
  by making use of $\cos(\theta) (\cos^2(\theta) - 3 \sin^2(\theta)) = \cos(3 \theta)$.
  Then $\logicalmapping_A$ and the logical $0$-form space are $\pr$-compatible under the conditions of \cref{ex:pr_compatibility_disc_mapping}.

  \paragraph{Using $\intermediatemapping = \mathrm{id}$.}
  We let $\intermediatemapping_B = \mathrm{id}$ and $\logicalmapping_B = \mapping_A$, so that $\mapping_B = \intermediatemapping_B\circ\logicalmapping_B = \mapping_A$.
  It follows that
  \begin{equation}
    \logicalmapping_B(\logical{\+x}) = R_1 \logical{x} \sspar{1 + 0.1 \, \logical{x}^3\cos(6\pi\logical{y})} \begin{bmatrix}
      \cos(2\pi\logical{y})\\
      \sin(2\pi\logical{y})
    \end{bmatrix}.
  \end{equation}
  Then the pullback of a polynomial $\intermediate{T}\in\intermediate{\funspace{T}}^\pr$ of degree $a+b\le\pr$ is given by
  \begin{equation}
    \intermediate{T} = \intermediatecoord{x}^a \intermediatecoord{y}^b
    \quad \implies \quad
    \logical{\pushforward}_B^{-0} \intermediate{T}
    = \intermediate{T} \circ \logicalmapping_B
    = \srpar{R_1 \logical{x} \sspar{1 + 0.1 \, \logical{x}^3\cos(6\pi\logical{y})}}^{a+b} \cos^a (2\pi\logical{y}) \sin^b (2\pi\logical{y}),
  \end{equation}
  from which it follows that $\pr$-compatibility requires a B-spline basis in the $\logical{x}$-direction of degree at least $r^{0,\xx}\ge4\pr$ and a trigonometric basis in the $\logical{y}$-direction containing the first $8\pr+1$ modes.

  Hence, the first decomposition is $\pr$-compatible provided $r^{0,\xx} \ge \pr$ whereas the second decomposition is $\pr$-compatible only if $r^{0,\xx} \ge 4\pr$, despite resulting in the same physical mapping $\mapping_A = \mapping_B$.
  It suggests the following \emph{design rule}: choose the decomposition $\mapping=\intermediatemapping\circ\logicalmapping$ such that the logical mapping $\logicalmapping$ is as simple as possible, ideally the disk mapping of \cref{ex:disk_mapping}, and transfer all geometric complexity to the intermediate mapping $\intermediatemapping$, which does not enter the compatibility condition.
  This minimizes the requirements on the logical spaces; the guiding-center experiment in \cref{sec:results_vacuum} follows this rule.
\end{remark}

\begin{remark}[A lack of $\pr$-compatibility]
  When the logical mapping and the logical $0$-form space are not $\pr$-compatible, we find that \cref{thm:polar_zeroform_space_regularity} no longer guarantees $\pr$-times continuous differentiability of the polar $0$-forms at the center.
  The polar subspace is still conforming by \cref{thm:polar_zeroform_space_conforming}, which only relies on the constant function being represented by the logical $0$-form space.
  Moreover, the results presented in \citet{jiang2026polar} and also in this work in \cref{sec:results_spectrum,sec:results_wave3d} suggest that the resulting polar $0$-form space still suppresses the stiff modes originating at the center in the discrete solution of PDEs, even when the logical mapping and the logical $0$-form space are not $\pr$-compatible.

  When regularity is desired despite a lack of $\pr$-compatibility, one could consider an approach in the spirit of the extended finite element method (XFEM) \citep{fries2010xfem}, where the tensor-product $0$-forms are enriched with the smooth physical $0$-forms (cf.\ \cref{eqn:polar_zeroform_space}):
  \begin{equation}
    \funspace{\Bspline}^{0,\pr}_{\mathrm{XFEM}} \defeq \pulledpolys{0}{\pr} \oplus \funspace{\Bspline}^{0,\pr}_\external.
  \end{equation}
  The downside of such an approach is that the basis functions are generally no longer compactly supported, and their span is not a subspace of the tensor-product $0$-form space $\funspace{\Bspline}^0$.
\end{remark}

  \newcommand\przero{0}
\section{Polar subspaces of the 2D de Rham complex}\label{sec:derham_2d}
In the previous section we have discussed the construction of high-order polar subspaces for $0$-forms, resulting in conforming subspaces in the sense of \cref{eqn:polar_zeroform_space_conforming}.
The natural question is how to extend this construction to the $1$- and $2$-form spaces, resulting in polar subspaces denoted by $\polar{\funspace{\Bspline}}^{1,\pr}$ and $\polar{\funspace{\Bspline}}^{2,\pr}$, respectively, on which high-order regularity conditions are satisfied.
The construction is based on three steps: firstly, we add gradients of polar $0$-forms to the polar $1$-form space (generalizing the work of \citet{toshniwal2021} to $\pr > 1$), which inherit the regularity of the polar $0$-form.
Secondly, we enrich the polar $2$-form space analogously to how the polar $0$-form space was enriched with projected polynomials.
Thirdly, we propose a novel construction to lift those polar $2$-forms using so-called homotopy operators to the polar $1$-form space.
Altogether, the polar subspaces will be
\begin{equation}\label{eqn:polar_subspaces_overview}
  \polar{\funspace{\Bspline}}^{0,\pr} = \pi^0 \pulledpolys{0}{\pr} \oplus \funspace{\Bspline}^{0,\pr}_\external
  , \;\;
  \polar{\funspace{\Bspline}}^{1,\pr} = \gradientt \pi^0 \pulledpolys{0}{\pr} \oplus \curlt^\dagger \reduced{\pi}^{2,1} \pulledpolys{2}{\pr-2} \oplus \funspace{\Bspline}^{1,\pr}_\external
  , \;\;
  \polar{\funspace{\Bspline}}^{2,\pr} = \reduced{\pi}^{2,1} \pulledpolys{2}{\pr-2} \oplus \funspace{\Bspline}^{2,\pr}_\external,
\end{equation}
where the $0$-form space is that of \cref{eqn:polar_zeroform_space}, the $1$- and $2$-form spaces are defined in \cref{eqn:polar_oneform_space,eqn:polar_twoform_space}, the pulled-back polynomial $2$-forms $\pulledpolys{2}{\pr-2}$ and the conforming $2$-form projector $\reduced{\pi}^{2,1}$ are introduced in \cref{sec:homotopy_motivation,sec:derham_hodge_2d_construction}, and the homotopy operator $\curlt^\dagger$ in \cref{sec:derham_hodge_2d_hodge}.
The spaces are such that following two diagrams individually commute (we leave out the logical and intermediate spaces for clarity)
\begin{equation}\label{eqn:derham_hodge_2d_diagram}
    \begin{derhamdiagram}[dimension=2]
    \derhamrow[name=polar,spaces={\polar{\funspace{\Bspline}}^{0,\pr},\polar{\funspace{\Bspline}}^{1,\pr},\polar{\funspace{\Bspline}}^{2,\pr}},operators={\gradientt,\curlt}]
    \derhamrow[name=continuous,spaces={\funspace{H}_1,\funspace{H}_{\curlt},\funspace{L}_2},operators={\gradientt,\curlt}]
    \derhamconnect[from=continuous,to=polar,degree=0,operator=\polar{\pi}^{0,\pr}]
    \derhamconnect[from=continuous,to=polar,degree=1,operator=\polar{\pi}^{1,\pr}]
    \derhamconnect[from=continuous,to=polar,degree=2,operator=\polar{\pi}^{2,\pr}]
  \end{derhamdiagram}
  \qquad
  \begin{derhamdiagram}[dimension=2]
    \derhamrow[name=polar,kind=dual,spaces={\polar{\funspace{\Bspline}}^{0,\pr},\polar{\funspace{\Bspline}}^{1,\pr},\polar{\funspace{\Bspline}}^{2,\pr}},operators={\gradientt^\dagger,\curlt^\dagger}]
    \derhamrow[name=continuous,kind=dual,spaces={\funspace{H}_1,\funspace{H}_{\curlt},\funspace{L}_2},operators={\gradientt^\dagger,\curlt^\dagger}]
    \derhamconnect[from=continuous,to=polar,degree=0,operator=\polar{\pi}^{0,\pr}]
    \derhamconnect[from=continuous,to=polar,degree=1,operator=\polar{\pi}^{1,\pr}]
    \derhamconnect[from=continuous,to=polar,degree=2,operator=\polar{\pi}^{2,\pr}]
  \end{derhamdiagram},
\end{equation}
where commuting polar projection operators $\polar{\pi}^m$ are introduced, which, contrary to the projection operators $\pi^m$, project into \emph{conforming} finite-dimensional subspaces of the infinite-dimensional function spaces.

The backward arrows denote the homotopy operators $\gradientt^\dagger$ and $\curlt^\dagger$, which are introduced in \cref{sec:derham_hodge_2d_hodge}: they integrate along rays from the center and, on the clamped radial B-splines, reduce to cumulative sums of coefficients, so that the homotopy identities, which express every form as the derivative of its homotopy plus the homotopy of its derivative (up to the value at the center for $0$-forms), hold exactly on the discrete spaces.
This second complex in the opposite direction is what the construction buys beyond the cochain property.
The homotopy identities show that every curl-free polar $1$-form is a gradient and every polar $2$-form is a curl (\cref{cor:polar_subspaces_cohomology}), they split the polar $1$-form space into gradients of polar $0$-forms and lifts of polar $2$-forms (\cref{thm:polar_subspaces_homotopy_splitting}), and they define the polar $1$-form projection in terms of the $0$- and $2$-form projections (\cref{eqn:polar_oneform_projection_operator}), which makes its commutation with the differential operators automatic and keeps all projections sparse.
The commutativity of the homotopy operators with the polar projection operators is stated precisely in \cref{cor:infinite_dimensional_polar_projection_operators}.

We start by giving an overview of approaches from the literature for imposing first-order regularity conditions on the $1$- and $2$-form spaces in \cref{sec:first_order_regularity_de_rham}, resulting in conforming subspaces of the $1$- and $2$-form spaces.
In \cref{sec:homotopy_motivation} we consider a naive generalization of the work of \citet{toshniwal2021}, and we show that the requirement of obtaining a cochain complex alone does not yield sufficiently rich high-order polar $1$- and $2$-form spaces.
Then we propose a novel generalization of the construction of high-order polar subspaces for $0$-forms to the $1$- and $2$-form spaces in \cref{sec:derham_hodge_2d_construction}, which is based on the homotopy identity for differential forms as discussed in \cref{sec:derham_hodge_2d_hodge}.
The key properties of the polar subspaces are discussed in \cref{sec:derham_hodge_2d_properties}, and we propose commuting projection operators in \cref{sec:derham_hodge_2d_geom_projections}.
In \cref{sec:boundary_homotopy} we briefly discuss the inclusion of boundary conditions in the construction of polar subspaces, and show how the resulting boundary-compatible homotopy operators can be used to efficiently solve the curl--curl problem.


\subsection{First-order regularity conditions}\label{sec:first_order_regularity_de_rham}
Several approaches have been proposed in the literature for the construction of polar subspaces for $1$- and $2$-forms, which are based on first-order regularity conditions at the center of the domain.
This subsection is a review of these approaches; the only addition is \cref{lem:conformity_conditions_tensor_product_forms}, which states the conformity conditions underlying them in the notation of this paper, as needed in the remainder of this section.

In \citet{toshniwal2021} the first-order polar $1$- and $2$-form spaces are constructed explicitly by letting the part of the space that contributes at the center be the gradient of the polar $0$-forms and the curl of the polar $1$-forms, respectively.
That is, the $1$- and $2$-form spaces are given by 
\begin{equation}\label{eqn:polar_firstorder_spaces}
  \reduced{\funspace{\Bspline}}^{1,1}
  \defeq
  \gradientt \pi^0\pulledpolys{0}{1} \oplus \funspace{\Bspline}^{1,1}_\external
  , \quad
  \reduced{\funspace{\Bspline}}^{2,1}
  \defeq
  \funspace{\Bspline}^{2,1}_\external
\end{equation}
where $\funspace{\Bspline}^{m,\pr}_\external$ denotes the external part of the tensor-product $m$-form space defined in \cref{eqn:polar_zeroform_space_odot,eqn:polar_mform_space_odot}, here with $\pr=1$.
The resulting polar subspaces yield a cochain complex by construction with the same cohomology as the continuous 2D de Rham complex on a disk-like domain \citep[Thm.~5.10]{toshniwal2021}.

A similar approach is followed in \citet{holderied2022}, where the first-order polar $1$- and $2$-form spaces are constructed by means of commuting projection operators, which are based on the geometric degrees of freedom.
Alternatively, in \citet{yaman2025} the polar $1$- and $2$-form spaces are constructed by means of commuting projection operators that compose the tensor-product projections, which are based on the geometric degrees of freedom, with local conforming projections acting on the spline coefficients, which can be cheaply computed.

Imposing first-order regularity conditions is sufficient for the construction of conforming subspaces of the $1$- and $2$-form spaces, as we show in \cref{lem:conformity_conditions_tensor_product_forms,cor:conforming_first_order_spaces}.
\begin{lemma}[Necessary and sufficient conditions for conformity of tensor-product $1$- and $2$-forms]\label{lem:conformity_conditions_tensor_product_forms}
  A $1$-form $\+A\in\funspace{\Bspline}^{1}$ belongs to $\funspace{H}_{\curlt}$ if and only if
  \begin{equation}\label{eqn:conformity_conditions_tensor_product_oneforms}
    a_{\yy,1j}=0,
    \qquad
    a_{\yy,2j}=(\delta_\yy a_\xx)_{1j},
    \qquad j\in\{1,\ldots,n^\yy\}.
  \end{equation}
  Equivalently, there is a $0$-form $P\in\funspace{\Bspline}^0$ whose first ring represents a constant, $p_{1j}=P(\+0)\mathds{1}_j^\yy$, and whose gradient has the same coefficients as $\+A$ on the following rings:
  \begin{equation}\label{eqn:conformity_conditions_tensor_product_oneforms_gradient}
    a_{\xx,1j}=(\delta_\xx p)_{1j},
    \qquad
    a_{\yy,1j}=(\delta_\yy p)_{1j}=0,
    \qquad
    a_{\yy,2j}=(\delta_\yy p)_{2j}.
  \end{equation}

  A $2$-form $B\in\funspace{\Bspline}^{2}$ belongs to $\funspace{L}_2$ if and only if
  \begin{equation}\label{eqn:conformity_conditions_tensor_product_twoforms}
    b_{1j}=0,
    \qquad j\in\{1,\ldots,n^\yy\}.
  \end{equation}
\end{lemma}
\begin{statementproof}[appendix]{lem:conformity_conditions_tensor_product_forms}
  We first assume $\intermediatemapping=\mathrm{id}$, so that $\Omega=\intermediate{\Omega}$, $\mapping=\logicalmapping$, $\jacmat=\logical{\jacmat}$, $\jac=\logical{\jac}$ and $\pushforward^m=\logical{\pushforward}^m$, and explain at the end why the result carries over to a general intermediate mapping.
  Away from the center the mapping is a smooth diffeomorphism and the basis functions are bounded, so every tensor-product form and its (piecewise) derivative are square integrable there, and only the center has to be examined.
  Near the center two effects compete with the boundedness of the logical forms: the Jacobian vanishes like $\logical{x}$, and the covariant Piola map amplifies angular components by $1/\logical{x}$.
  We first quantify both effects, then derive the $2$-form condition, the $\funspace{L}_2$ condition for $1$-forms, and finally the curl condition, which reduces to the $2$-form condition by commutativity of the pushforwards.

  \medskip
  \noindent\emph{Behavior at the center.}
  With $\tilde R\defeq R/\logical{x}$, \cref{eqn:disk_jacobian} reads $\jac=\logical{x}\,\tilde R\,(\tilde R+\logical{x}\partial_{\logical{x}}\tilde R)(CS'-SC')$, so $\jac/\logical{x}\to R_1^2(CS'-SC')>0$ as $\logical{x}\to0$, uniformly in $\logical{y}$ (\cref{ass:disc_like_mapping}).
  The matrix $\mat{\Gamma}$ from \cref{eqn:gamma_matrix} is smooth up to $\logical{x}=0$ and has $\det\mat{\Gamma}=\logical{x}/\jac$.
  Hence there are $\varepsilon>0$, smaller than the first positive radial knot, and positive constants $c_1,c_2,c_3$ such that
  \begin{equation}\label{eqn:conformity_pole_metric_bounds}
    c_1\logical{x}\le\jac\le c_2\logical{x},
    \qquad
    \|\mat{\Gamma}\|+\|\mat{\Gamma}^{-1}\|\le c_3,
    \qquad
    0<\logical{x}<\varepsilon,\ \logical{y}\in[0,1].
  \end{equation}
  On $[0,\varepsilon)$ the radial basis functions are polynomials.
  Their values at $\logical{x}=0$ are $\Bsplineone^{0,\xx}_i(0)=\delta_{i1}$ by \cref{eqn:bspline_x_atpole}, and $\Bsplineone^{1,\xx}_i(0)=\delta_{i1}/\Delta t^{0,\xx}_1$ by combining \cref{eqn:bspline_derivative} (with $\Bsplineone^{1,\xx}_0=0$) with the derivative values in \cref{eqn:bspline_x_atpole}.
  Since the angular basis functions are bounded, every tensor-product function $f$ satisfies
  \begin{equation}\label{eqn:conformity_pole_taylor}
    f(\logical{x},\logical{y})=f(0,\logical{y})+\mathcal{O}(\logical{x}),
    \qquad 0<\logical{x}<\varepsilon,
  \end{equation}
  uniformly in $\logical{y}$, and only the first radial index contributes to $f(0,\logical{y})$.
  Below, $\Omega_\varepsilon\defeq\mapping\bigl((0,\varepsilon)\times[0,1]\bigr)$ denotes the punctured neighborhood of the center.

  \medskip
  \noindent\emph{The $2$-form condition.}
  Let $q\defeq\logical{B}$, so that $B\circ\mapping=q/\jac$ by \cref{eqn:pushforward_operators}, and
  \begin{equation}\label{eqn:conformity_twoform_norm}
    \int_{\Omega_\varepsilon}|B|^2\diff^2 x
      =\int_0^\varepsilon\int_0^1
        \frac{|q(\logical{x},\logical{y})|^2}{\jac(\logical{x},\logical{y})}
        \diff\logical{y}\diff\logical{x}.
  \end{equation}
  By \cref{eqn:conformity_pole_taylor}, $q=q_0+\mathcal{O}(\logical{x})$ with
  \begin{equation}
    q_0(\logical{y})\defeq q(0,\logical{y})=\frac{1}{\Delta t_1^{0,\xx}}
      \sum_{j=1}^{n^\yy}b_{1j}\Bsplineone_j^{1,\yy}(\logical{y}).
  \end{equation}
  If $q_0=0$, then $|q|\le c_4\logical{x}$, and the lower bound on $\jac$ in \cref{eqn:conformity_pole_metric_bounds} bounds the integrand in \cref{eqn:conformity_twoform_norm} by $c_4^2\logical{x}/c_1$, so the integral is finite.
  If $q_0\ne0$, then $\|q(\logical{x},\cdot)\|_{L_2(0,1)}\ge\|q_0\|_{L_2(0,1)}/2$ for all sufficiently small $\logical{x}$, and the upper bound on $\jac$ bounds \cref{eqn:conformity_twoform_norm} from below by a positive multiple of $\int_0\diff\logical{x}/\logical{x}=\infty$.
  The angular basis functions are linearly independent, so $q_0=0$ is equivalent to \cref{eqn:conformity_conditions_tensor_product_twoforms}.

  \medskip
  \noindent\emph{The $1$-form condition.}
  Write $u$ and $v$ for the radial and angular components of $\logical{\+A}$ in \cref{eqn:logical_oneforms}.
  By the covariant Piola transformation in \cref{eqn:pushforward_operators} and the definition of $\mat{\Gamma}$ in \cref{eqn:gamma_matrix},
  \begin{equation}\label{eqn:tensor_product_oneform_pole_expansion}
    \+A\circ\mapping
      =\jacmat^{-\transpose}\begin{bmatrix}u\\v\end{bmatrix}
      =\mat{\Gamma}
        \begin{bmatrix}u\\v/\logical{x}\end{bmatrix}.
  \end{equation}
  The bounds in \cref{eqn:conformity_pole_metric_bounds} therefore give positive constants $c_5,c_6$ such that
  \begin{equation}
    c_5\int_0^\varepsilon\int_0^1
      \left(\logical{x}|u|^2+\frac{|v|^2}{\logical{x}}\right)
      \diff\logical{y}\diff\logical{x}
    \le\int_{\Omega_\varepsilon}|\+A|^2\diff^2 x
    \le c_6\int_0^\varepsilon\int_0^1
      \left(\logical{x}|u|^2+\frac{|v|^2}{\logical{x}}\right)
      \diff\logical{y}\diff\logical{x}.
  \end{equation}
  The term with $u$ is finite because $u$ is bounded.
  The angular component behaves as the $2$-form above: by \cref{eqn:conformity_pole_taylor}, $v=v_0+\mathcal{O}(\logical{x})$ with
  \begin{equation}\label{eqn:tensor_product_oneform_singular_part}
    v_0(\logical{y})\defeq v(0,\logical{y})=\sum_{j=1}^{n^\yy}a_{\yy,1j}\Bsplineone_j^{1,\yy}(\logical{y}),
  \end{equation}
  and the same two cases show that the integral of $|v|^2/\logical{x}$ is finite if and only if $v_0=0$.
  By linear independence of the angular basis, $\+A\in\funspace{L}_2$ if and only if $a_{\yy,1j}=0$ for all $j$.
  In that case $v/\logical{x}$ is bounded, so \cref{eqn:tensor_product_oneform_pole_expansion} shows that $\+A$ is bounded near the center.

  \medskip
  \noindent\emph{The curl condition.}
  By \cref{eqn:pushforward_commutativity} and \cref{ex:curl_logical_oneform}, the curl of $\+A$ is the physical $2$-form with coefficients $(\curlh\+a)_{ij}=(\delta_\xx a_\yy-\delta_\yy a_\xx)_{ij}$.
  By the $2$-form condition, it belongs to $\funspace{L}_2$ if and only if its first ring vanishes,
  \begin{equation}\label{eqn:tensor_product_oneform_second_ring_condition}
    0=(\curlh\+a)_{1j}=a_{\yy,2j}-a_{\yy,1j}-(\delta_\yy a_\xx)_{1j},
    \qquad j\in\{1,\ldots,n^\yy\}.
  \end{equation}
  Combined with the $1$-form condition $a_{\yy,1j}=0$, this is exactly \cref{eqn:conformity_conditions_tensor_product_oneforms}.

  \medskip
  \noindent\emph{The gradient interpretation.}
  Let $P\in\funspace{\Bspline}^0$ have the constant first ring $p_{1j}=c_0\mathds{1}_j^\yy$.
  By \cref{eqn:conformity_pole_taylor}, $P\circ\mapping=c_0+\mathcal{O}(\logical{x})$, so $P$ extends continuously to the center with $P(\+0)=c_0$.
  Moreover, $\delta_\yy\mathds{1}^\yy=0$ (by \cref{eqn:diff_opers_1d} for periodic B-splines and \cref{eqn:trigonometric_coefficient_derivative} for trigonometric functions).
  Hence $(\delta_\yy p)_{1j}=0$, and $(\delta_\yy p)_{2j}=(\delta_\yy\delta_\xx p)_{1j}$ because $p_{2j}=p_{1j}+(\delta_\xx p)_{1j}$.
  If \cref{eqn:conformity_conditions_tensor_product_oneforms_gradient} holds, then $a_{\yy,1j}=0$ and $a_{\yy,2j}=(\delta_\yy\delta_\xx p)_{1j}=(\delta_\yy a_\xx)_{1j}$, which is \cref{eqn:conformity_conditions_tensor_product_oneforms}.
  Conversely, if \cref{eqn:conformity_conditions_tensor_product_oneforms} holds, then
  \begin{equation}\label{eqn:tensor_product_oneform_potential_construction}
    p_{1j}=c_0\mathds{1}_j^\yy,
    \qquad
    p_{2j}=c_0\mathds{1}_j^\yy+a_{\xx,1j},
  \end{equation}
  with any $c_0\in\mathbb{R}$ and arbitrary remaining coefficients, defines such a $P$ satisfying \cref{eqn:conformity_conditions_tensor_product_oneforms_gradient}.

  \medskip
  \noindent\emph{General intermediate mapping.}
  The intermediate mapping $\intermediatemapping$ is a smooth diffeomorphism up to the boundary that fixes the center (\cref{ass:diffeomorphism}), so, unlike $\logicalmapping$, it introduces no singularity.
  By \cref{eqn:pushforward_composition}, every physical tensor-product form is the pushforward by $\intermediate{\pushforward}^m$ of the intermediate tensor-product form with the same coefficients, to which the case $\intermediatemapping=\mathrm{id}$ applies.
  The pushforwards $\intermediate{\pushforward}^m$ commute with the gradient and the curl, also in the weak sense, and, because $\intermediate{\jacmat}$, $\intermediate{\jac}$ and their inverses are bounded, they change $\funspace{L}_2$ norms by at most a bounded factor; see the end of the proof of \cref{thm:polar_zeroform_space_conforming} for $0$-forms.
  Hence $\intermediate{\pushforward}^0$, $\intermediate{\pushforward}^1$, $\intermediate{\pushforward}^2$ map $\intermediate{\funspace{H}}_1$, $\intermediate{\funspace{H}}_{\curlt}$, $\intermediate{\funspace{L}}_2$ onto $\funspace{H}_1$, $\funspace{H}_{\curlt}$, $\funspace{L}_2$, and the coefficient conditions derived above hold for every intermediate mapping.
  Similarly, composition with $\intermediatemapping^{\pm1}$ and multiplication by the smooth invertible factors $\intermediate{\jacmat}^{\pm\transpose}$ and $\intermediate{\jac}^{\pm1}$ preserve boundedness and $k$-fold continuous differentiability near the center, and thus the classes $\funspace{C}^k_0$, $k\ge-1$.
\end{statementproof}

\begin{corollary}[Conforming first-order spaces]\label{cor:conforming_first_order_spaces}
  The spaces defined in \cref{eqn:polar_firstorder_spaces} are conforming subspaces of $\funspace{H}_{\curlt}$ and $\funspace{L}_2$, respectively.
\end{corollary}
\begin{statementproof}[appendix]{cor:conforming_first_order_spaces}
  At $\pr=1$, the external index sets in \cref{eqn:polar_mform_space_odot} contain radial $1$-form indices $i\ge2$, angular $1$-form indices $i\ge3$, and $2$-form indices $i\ge2$.
  Hence an external $1$-form has $a_{\xx,1j}=a_{\yy,1j}=a_{\yy,2j}=0$, while an external $2$-form has $b_{1j}=0$.
  Both are conforming by \cref{lem:conformity_conditions_tensor_product_forms}.
  For $P\in\pi^0\pulledpolys{0}{1}\subset\polar{\funspace{\Bspline}}^{0,1}$, \cref{thm:polar_zeroform_space_conforming} gives $P\in\funspace{H}_1$; thus $\gradientt P\in\funspace{H}_{\curlt}$ because $\curlt\gradientt P=0$.
  The sum in \cref{eqn:polar_firstorder_spaces} is therefore contained in $\funspace{H}_{\curlt}$.
\end{statementproof}

\subsection{Tentative high-order polar cochain complex}\label{sec:homotopy_motivation}
The approach of \citet{toshniwal2021} can straightforwardly be generalized to high-order regularity conditions by using the following direct-sum decomposition of the polar $1$-form space 
\begin{equation}\label{eqn:reduced_oneform_space}
  \reduced{\funspace{\Bspline}}^{1,\pr} \defeq \gradientt \pi^0 \pulledpolys{0}{\pr} \oplus \funspace{\Bspline}^{1,\pr}_\external,
\end{equation}
where $\funspace{\Bspline}^{1,\pr}_\external$ is the external part of the tensor-product $1$-form space from \cref{eqn:polar_mform_space_odot}, whose elements have bounded physical derivatives up to order $\pr-1$ at the center.

When comparing to \cref{eqn:polar_zeroform_space} we find that the internal polar $1$-forms are exactly the gradients of the internal polar $0$-forms, resulting in the desired high-order regularity conditions at the center while preserving the cochain complex structure.
However, this construction implies that every polar $1$-form is \emph{purely} a gradient at the center.
Requiring the cochain property does not yield any polar $2$-forms that contribute at the center: since $\curlt \gradientt = 0$, we find that the internal polar $2$-form space is empty, and the polar $2$-form space is given by
\begin{equation}\label{eqn:reduced_twoform_space}
  \reduced{\funspace{\Bspline}}^{2,\pr} \defeq \funspace{\Bspline}^{2,\pr}_\external,
\end{equation}
where $\funspace{\Bspline}^{2,\pr}_\external$ is the external part of the tensor-product $2$-form space from \cref{eqn:polar_mform_space_odot}.
Both observations show that constructing polar spaces by imposing the cochain property alone does not yield sufficiently rich polar $1$- and $2$-form spaces.
In practice this means that a smooth $2$-form that does not vanish at the center, such as a current density or the normal component of a magnetic flux density, cannot be approximated with the intended regularity: its coefficients in the first $\pr$ rings would have to vanish, so that the flux through the central cells is forced into the $(\pr+1)$-st ring.

The following $2$-form projector projects onto the $2$-form space $\reduced{\funspace{\Bspline}}^{2,\pr}$ by transferring the first $\pr$ radial rings of $2$-form coefficients to the $\pr+1$st.
For $\pr = 1$, this gives the conforming projection operator onto the $2$-form space in \cref{eqn:polar_firstorder_spaces} (cf.\ \citet{yaman2025}).
Although the tentative spaces themselves are not used further, this projector is: it enters the definition of the enriched polar $2$-form space and of its compactly supported basis in \cref{sec:derham_hodge_2d_construction}, as well as the polar projection operators in \cref{sec:derham_hodge_2d_geom_projections}.
It acts on the radial index only, see \cref{eqn:reduced_twoform_projector_matrix} for its matrix representation.
\begin{lemma}[Conforming polar $2$-form projector]\label{lem:conforming_polar_twoform_projector}
  Let $\pr\ge1$.
  For $B\in\funspace{\Bspline}^2$ with coefficients $b_{ij}$ as in \cref{eqn:logical_twoforms}, define
  \begin{equation}\label{eqn:reduced_twoform_coefficient_projector}
    \reduced{\Pi}^{2,\pr} B
    \defeq
    \sum_{j=1}^{n^\yy} \sum_{i'=1}^{\pr+1} b_{i'j} \Bspline_{(\pr+1) j}^2
    +
    \sum_{i=\pr+2,j=1}^{n^{1,\xx},n^\yy} b_{ij} \Bspline_{ij}^2
    , \quad
    \reduced{\pi}^{2,\pr} \defeq \reduced{\Pi}^{2,\pr} \circ \pi^2
  \end{equation}
  Then $\reduced{\Pi}^{2,\pr}:\funspace{\Bspline}^2\to\reduced{\funspace{\Bspline}}^{2,\pr}$ and $\reduced{\pi}^{2,\pr}:\funspace{L}_2\to\reduced{\funspace{\Bspline}}^{2,\pr}$ are projectors onto the conforming $2$-form space in \cref{eqn:reduced_twoform_space}.
\end{lemma}
\begin{statementproof}[appendix]{lem:conforming_polar_twoform_projector}
  The coefficients of $\reduced{\Pi}^{2,\pr}B$ are
  \begin{equation}
    (\reduced{\Pi}^{2,\pr}B)_{ij}=0\quad(1\le i\le\pr),
    \qquad
    (\reduced{\Pi}^{2,\pr}B)_{(\pr+1)j}=\sum_{i'=1}^{\pr+1}b_{i'j},
    \qquad
    (\reduced{\Pi}^{2,\pr}B)_{ij}=b_{ij}\quad(i\ge\pr+2).
  \end{equation}
  Thus its image lies in $\funspace{\Bspline}^{2,\pr}_\external$ by \cref{eqn:polar_mform_space_odot}.
  If $B$ belongs to that space, then $b_{ij}=0$ for $i\le\pr$, and the displayed formula gives $\reduced{\Pi}^{2,\pr}B=B$.
  Hence $\reduced{\Pi}^{2,\pr}$ projects onto $\reduced{\funspace{\Bspline}}^{2,\pr}=\funspace{\Bspline}^{2,\pr}_\external$.
  Since $\pr\ge1$, every form in this image has $b_{1j}=0$ and belongs to $\funspace{L}_2$ by \cref{lem:conformity_conditions_tensor_product_forms}.
  The face-integral GDOFs make $\pi^2$ well-defined on $\funspace{L}_2$, and unisolvence gives $\pi^2B=B$ for every $B\in\funspace{\Bspline}^{2,\pr}_\external$.
  Consequently $\reduced{\pi}^{2,\pr}=\reduced{\Pi}^{2,\pr}\circ\pi^2$ has the same image and fixes every member of it, so it too is a projector.
\end{statementproof}

The tentative polar $m$-form spaces $\reduced{\funspace{\Bspline}}^{m,\pr}$ in \cref{eqn:reduced_oneform_space,eqn:reduced_twoform_space} are conforming subspaces of $\funspace{H}_{\curlt}$ and $\funspace{L}_2$, respectively, and they satisfy the cochain property by construction.
Moreover, they can be shown to have the desired regularity conditions at the center and the same cohomology as the continuous 2D de Rham complex on $\Omega$, which is that of the disk.
However, the spaces are not sufficiently rich: the polar $1$-form space $\reduced{\funspace{\Bspline}}^{1,\pr}$ is purely a gradient at the center, and the polar $2$-form space $\reduced{\funspace{\Bspline}}^{2,\pr}$ does not contribute at the center.
We would like to enrich the polar $2$-form space $\reduced{\funspace{\Bspline}}^{2,\pr}$ with projected polynomial $2$-forms that contribute at the center (analogously to the polar $0$-form space $\polar{\funspace{\Bspline}}^{0,\pr}$ in \cref{eqn:polar_zeroform_space}).
However, the question arises how to lift these projected polynomial $2$-forms to a curl-like part of the polar $1$-form space $\polar{\funspace{\Bspline}}^{1,\pr}$, while preserving the cochain property, the desired regularity conditions at the center, and the same cohomology as the continuous 2D de Rham complex on $\Omega$ (cf.\ \cref{cor:polar_subspaces_cohomology}).

The main tool to achieve this are homotopy operators, which is the topic of the following section.

\subsection{Homotopy operators}\label{sec:derham_hodge_2d_hodge}

We first discuss the homotopy operator of the 1D de Rham complex and its B-spline coefficient representation in \cref{sec:homotopy_1d}.
In \cref{sec:homotopy_2d} we show that the homotopy operators of the 2D polar de Rham complex reduce to this 1D operator along radial lines, so that their coefficient representation and homotopy identities follow from the 1D results.

\subsubsection{The 1D de Rham complex}\label{sec:homotopy_1d}
For a sufficiently smooth $0$-form $P$ and $1$-form $A$ on the interval $[0,1]$, the fundamental theorem of calculus and its derivative read
\begin{subequations}
  \label{eqn:homotopy_identities_1d}
  \begin{align+}
    P &= \mathfrak{P}_\xx P + \partial_\xx^\dagger \partial_\xx P,\label{eqn:homotopy_identity_0form_1d}\\
    A &= \partial_\xx \partial_\xx^\dagger A,\label{eqn:homotopy_identity_1form_1d}
  \end{align+}
\end{subequations}
where $\partial_\xx^\dagger$ denotes the anti-derivative operator
\begin{equation}
  (\partial_\xx^\dagger A)(\xx)
  \defeq
  \int_0^\xx A(s) \diff s
  ,\label{eqn:homotopy_operator_1d}
\end{equation}
and $\mathfrak{P}_\xx$ denotes the projection operator onto the space of constant $0$-forms:
\begin{equation}
  \mathfrak{P}_\xx P \defeq P(0).
\end{equation}
These are the homotopy identities of the 1D de Rham complex.

On the clamped B-spline spaces $\funspace{\Bsplineone}^{0,\xx}$ and $\funspace{\Bsplineone}^{1,\xx}$, the homotopy operator is a cumulative sum of coefficients:
\begin{equation}\label{eqn:discrete_radial_homotopy_operator}
  \partial_\xx^\dagger A
  = \partial_\xx^\dagger\sum_{i\in\indexone^{1,\xx}} a_i\Bsplineone^{1,\xx}_i
  =\sum_{i\in\indexone^{0,\xx}}(\delta_\xx^\dagger a)_i\Bsplineone^{0,\xx}_i,
  \qquad
  (\delta_\xx^\dagger a)_i\defeq\sum_{i'=1}^{i-1}a_{i'}.
\end{equation}
Indeed, the right-hand side vanishes at $0$ by \cref{eqn:bspline_x_atpole} because $(\delta_\xx^\dagger a)_1=0$, and by \cref{eqn:derivative_operator_1d} its derivative has the coefficients $\delta_\xx\delta_\xx^\dagger a$, which equal $a$ by the second relation below.
In particular, $\partial_\xx^\dagger$ maps $\funspace{\Bsplineone}^{1,\xx}$ into $\funspace{\Bsplineone}^{0,\xx}$, and $\delta_\xx^\dagger$ is represented by the cumulative-sum matrix in \cref{eqn:radial_cumulative_sum_matrix}.
For coefficients $a$ of a 1D $1$-form and $p$ of a 1D $0$-form, telescoping the sums gives
\begin{equation}\label{eqn:discrete_radial_homotopy_relations}
  (\delta_\xx^\dagger\delta_\xx p)_i=p_i-p_1
  , \quad
  \delta_\xx\delta_\xx^\dagger a=a
  .
\end{equation}
These are the coefficient versions of the homotopy identities in \cref{eqn:homotopy_identities_1d}.
In the first relation, $p_1$ is the value at $0$ by \cref{eqn:bspline_x_atpole}, and this constant has the coefficient $p_1$ for every $i$ by the partition of unity in \cref{eqn:bspline_partition}; hence $p_1$ represents $\mathfrak{P}_\xx$.
Consequently, the 1D homotopy identities hold exactly on the clamped B-spline spaces.
The anti-derivative of a clamped B-spline $1$-form is thus again a B-spline, with coefficients given by a cumulative sum.
This classical property is what the construction below rests on: in \cref{sec:homotopy_2d} we show that the Poincar\'e operators of the disk-like domain reduce to this radial anti-derivative, because the rays from the center are the images of the radial coordinate lines.

\subsubsection{The 2D polar de Rham complex}\label{sec:homotopy_2d}
The homotopy operators for the 2D polar de Rham complex are defined in a similar way, and they satisfy similar homotopy identities.
For sufficiently smooth differential forms $P$, $\+A$, and $B$ (which in particular belong to $\funspace{H}_1$, $\funspace{H}_{\curlt}$, and $\funspace{L}_2$, respectively), the homotopy identities are (cf.\ the 3D identities in \citealp[Prop.~2.1]{gopalakrishnan2004}; their coefficient versions are established in \cref{cor:discrete_homotopy_identities})
\begin{subequations}
  \label{eqn:homotopy_identities}
  \begin{align+}
    P &= \mathfrak{P} P + \gradientt^\dagger \gradientt P,\label{eqn:homotopy_identity_0form}\\
    \+A &= \gradientt \gradientt^\dagger \+A + \curlt^\dagger \curlt \+A,\label{eqn:homotopy_identity_1form}\\
    B &= \curlt \curlt^\dagger B,\label{eqn:homotopy_identity_2form}
  \end{align+}
\end{subequations}
where $\gradientt^\dagger$ and $\curlt^\dagger$ denote the homotopy operators for the gradient and curl operators, respectively.
We first define them on the intermediate domain as the following Poincar\'e representatives along straight rays from the center:
  \begin{equation}
    (\intermediate{\gradientt}^\dagger \intermediate{\+A})(\intermediatecoord{\+x})
    \defeq
    \int_0^1 \intermediate{\+A}(t\intermediatecoord{\+x}) \cdot \intermediatecoord{\+x} \diff t
    , \quad
    (\intermediate{\curlt}^\dagger \intermediate{B})(\intermediatecoord{\+x})
    \defeq
    \begin{bmatrix}
      - \intermediatecoord{y}\\
      \intermediatecoord{x}
    \end{bmatrix}
    \int_0^1 t \intermediate{B}(t \intermediatecoord{\+x}) \diff t
    .\label{eqn:homotopy_operators}
  \end{equation}
The segment from $\+0$ to $\intermediatecoord{\+x}\in\intermediate{\Omega}$ lies in $\intermediate{\Omega}$, since by \cref{eqn:disk_mapping,eqn:coordinate_transformation_assumptions} each logical radial line of constant $\logical{y}$ is mapped onto such a segment.
In the language of electrodynamics, $\intermediate{\gradientt}^\dagger\intermediate{\+A}$ is the line integral of $\intermediate{\+A}$ along the ray from the center, and $\intermediate{\curlt}^\dagger\intermediate{B}$ is the vector potential of $\intermediate{B}$ in the Poincar\'e (or multipolar) gauge $\intermediatecoord{\+x}\cdot\intermediate{\+A}=0$ \citep{brittin1982}; in the language of finite element exterior calculus, these are the Poincar\'e operators of the de Rham complex on a star-shaped domain \citep[eqs.~(2.2)--(2.4)]{gopalakrishnan2004}, cf.\ the Koszul operators in \citet[Sec.~3.2]{Arnold_Falk_Winther_2006} and their regularized versions, averaged over centers in a ball, in \citet[Sec.~3.1]{costabel_mcintosh2010}.
On the physical domain, the homotopy operators are obtained by means of the intermediate pullback and pushforward operators in \cref{eqn:pushforward_operators}:
\begin{equation}\label{eqn:homotopy_operators_transport}
  \gradientt^\dagger
  \defeq
  \intermediate{\pushforward}^0\circ\intermediate{\gradientt}^\dagger\circ\intermediate{\pushforward}^{-1}
  , \quad
  \curlt^\dagger
  \defeq
  \intermediate{\pushforward}^1\circ\intermediate{\curlt}^\dagger\circ\intermediate{\pushforward}^{-2}
  ,
\end{equation}
and for $\intermediatemapping=\mathrm{id}$ they are the Poincar\'e operators \cref{eqn:homotopy_operators} on $\Omega$.
The projection operator $\mathfrak{P}$ maps onto the space of constant $0$-forms:
\begin{equation}
  \mathfrak{P} P \defeq P(\+0).
\end{equation}
For tensor-product inputs, including nonconforming ones, \cref{eqn:homotopy_operators_transport} is read on $\Omega\setminus\{\+0\}$, with the lower endpoints taken as limits along each image of a ray.
The scalar output $\gradientt^\dagger\+A$ extends by zero to the center.
The vector output $\curlt^\dagger B$ also does so when $B$ is conforming, but may have no center limit for a nonconforming $B$.

The homotopy identities \cref{eqn:homotopy_identities} hold for the intermediate operators on $\intermediate{\Omega}$, which is star-shaped with respect to $\+0$ (cf.\ \citealp[Sec.~2]{gopalakrishnan2004}).
They carry over to $\Omega$ because $\intermediate{\pushforward}^m$ commutes with the gradient and curl by the chain rule and, since $\intermediatemapping(\+0)=\+0$, preserves center values of $0$-forms.
A direct computation shows that $\intermediate{\gradientt}^\dagger\circ\intermediate{\curlt}^\dagger=0$, since $[-\intermediatecoord{y},\intermediatecoord{x}]^\transpose\cdot\intermediatecoord{\+x}=0$; hence \cref{eqn:homotopy_operators_transport} gives
\begin{equation}
  \gradientt^\dagger \circ \curlt^\dagger = 0,
\end{equation}
which is the chain complex property in the backward direction to the de Rham complex (which relies on $\curlt \circ \gradientt = 0$).

The homotopy operators integrate along straight rays from the center in the intermediate domain, that is, along the images of these rays under $\intermediatemapping$ in the physical domain.
In logical coordinates these rays are the radial lines of constant $\logical{y}$, so the 2D homotopy operators reduce to the 1D homotopy operator $\partial_\xx^\dagger$ from \cref{sec:homotopy_1d}, applied in the radial variable $\logical{x}$ for each fixed $\logical{y}$.
Accordingly, $\delta_\xx^\dagger$ from \cref{eqn:discrete_radial_homotopy_operator} acts on the radial index of 2D coefficients, as $\delta_\xx$ does.
\begin{lemma}[Homotopy operators in logical coordinates]\label{lem:pseudo_inverse_curl_logical}
  Let $\logical{\+A}\in\logical{\funspace{\Bspline}}^{1}$ and $\logical{B}\in\logical{\funspace{\Bspline}}^{2}$.
  Away from the center, the homotopy operators from \cref{eqn:homotopy_operators}, transformed to logical coordinates, reduce to the 1D radial homotopy operator applied to the logical radial component of $\logical{\+A}$ and to $\logical{B}$, respectively:
  \begin{equation}\label{eqn:pseudo_inverse_gradient_logical}
    \logical{\gradientt}^\dagger \logical{\+A}
    \defeq \logical{\pushforward}^{-0} \intermediate{\gradientt}^\dagger \logical{\pushforward}^{1} \logical{\+A}
    = \partial_\xx^\dagger \logical{A}_{\xx}
    = \sum_{(i,j)\in \indextwo^{0}} (\delta_\xx^\dagger a_\xx)_{ij} \logical{\Bspline}^0_{ij},
  \end{equation}
  and
  \begin{equation}\label{eqn:pseudo_inverse_curl_logical}
    \logical{\curlt}^\dagger \logical{B}
    \defeq \logical{\pushforward}^{-1} \intermediate{\curlt}^\dagger \logical{\pushforward}^{2} \logical{B}
    = \begin{bmatrix}
      0\\
      \partial_\xx^\dagger \logical{B}
    \end{bmatrix}
    = \sum_{(i,j)\in \indextwo^{1}_\yy} (\delta_\xx^\dagger b)_{ij} \logical{\+\Bspline}^1_{\yy,ij}.
  \end{equation}
  Here $a_{\xx,ij}$ are the coefficients of the logical radial component $\logical{A}_{\xx}$ in \cref{eqn:logical_oneforms}, and $b_{ij}$ are the coefficients in \cref{eqn:logical_twoforms}.
  By \cref{eqn:homotopy_operators_transport,eqn:pushforward_composition}, the physical homotopy operators have the same representation, $\pushforward^{-0}\gradientt^\dagger\+A=\logical{\gradientt}^\dagger(\pushforward^{-1}\+A)$ and $\pushforward^{-1}\curlt^\dagger B=\logical{\curlt}^\dagger(\pushforward^{-2}B)$ for $\+A\in\funspace{\Bspline}^1$ and $B\in\funspace{\Bspline}^2$ on $\Omega\setminus\{\+0\}$.
\end{lemma}
\begin{statementproof}[appendix]{lem:pseudo_inverse_curl_logical}
  The second equalities in \cref{eqn:pseudo_inverse_gradient_logical,eqn:pseudo_inverse_curl_logical} are \cref{eqn:discrete_radial_homotopy_operator}, applied for each angular index $j$.
  It remains to show that the radial integrals $(\partial_\xx^\dagger\logical{A}_\xx)(\logical{x},\logical{y})=\int_0^{\logical{x}}\logical{A}_\xx(s,\logical{y})\diff s$ and $(\partial_\xx^\dagger\logical{B})(\logical{x},\logical{y})=\int_0^{\logical{x}}\logical{B}(s,\logical{y})\diff s$ are the pullbacks of the homotopy operators in \cref{eqn:homotopy_operators} applied to $\intermediate{\+A}\defeq\logical{\pushforward}^{1}\logical{\+A}$ and $\intermediate{B}\defeq\logical{\pushforward}^{2}\logical{B}$.
  By the covariant Piola transformation in \cref{eqn:pushforward_operators}, the pullback satisfies
  \begin{equation}
    \logical{A}_{\xx}(s,\logical{y})
    = \partial_s\logicalmapping(s,\logical{y}) \cdot \intermediate{\+A}(\logicalmapping(s,\logical{y}))
    = \partial_s R(s,\logical{y})
      \begin{bmatrix}C(\logical{y})\\S(\logical{y})\end{bmatrix}
      \cdot \intermediate{\+A}(\logicalmapping(s,\logical{y})).
  \end{equation}
  For fixed $\logical{y}$, the mapping assumptions give $R(0,\logical{y})=0$ and $\partial_sR(0,\logical{y})=R_1>0$.
  Away from the center, the Jacobian $\logical{\jac}$ is nonzero because $\logicalmapping$ is a diffeomorphism there (\cref{ass:disc_like_mapping}); hence $R$ and $\partial_sR$ remain positive along the ray.
  Thus $s\mapsto R(s,\logical{y})/R(\logical{x},\logical{y})$ increases from $0$ to $1$.
  We now let $\intermediatecoord{\+x} = \logicalmapping(\logical{x},\logical{y})$ and make the substitution $t R(\logical{x},\logical{y}) = R(s,\logical{y})$, resulting in
  \begin{equation}
    \logicalmapping(s,\logical{y}) = t\intermediatecoord{\+x},
    \quad
    \partial_s R(s,\logical{y})\diff s
    = R(\logical{x},\logical{y})\diff t,
    \quad
    \intermediatecoord{\+x}
    = R(\logical{x},\logical{y})
      \begin{bmatrix}C(\logical{y})\\S(\logical{y})\end{bmatrix}.
  \end{equation}
  Hence the scalar pushforward of $\logical{P}\defeq\partial_\xx^\dagger\logical{A}_\xx$ satisfies
  \begin{align}
    (\logical{\pushforward}^0\logical{P})(\intermediatecoord{\+x})
    = \int_0^{\logical{x}} \logical{A}_{\xx}(s,\logical{y}) \diff s 
    = \int_0^1
      R(\logical{x},\logical{y})
      \begin{bmatrix}C(\logical{y})\\S(\logical{y})\end{bmatrix}
      \cdot \intermediate{\+A}(t\intermediatecoord{\+x}) \diff t 
    = \int_0^1 \intermediatecoord{\+x} \cdot \intermediate{\+A}(t\intermediatecoord{\+x}) \diff t
    = (\intermediate{\gradientt}^\dagger \intermediate{\+A})(\intermediatecoord{\+x}),
  \end{align}
  which establishes the first equality in \cref{eqn:pseudo_inverse_gradient_logical}.

  We now consider the curl homotopy operator.  By the expression for the Jacobian in \cref{eqn:disk_jacobian}, the substitution $t R(\logical{x}, \logical{y}) = R(s,\logical{y})$ shows that the nonzero component of the logical vector in \cref{eqn:pseudo_inverse_curl_logical} is
  \begin{equation}
    \int_0^{\logical{x}} \logical{B}(s,\logical{y}) \diff s
    = 
    \int_0^{\logical{x}} \logical{\jac}(s,\logical{y}) \intermediate{B}(\logicalmapping(s, \logical{y})) \diff s
    = 
    R^2 \rpar{C S' - S C'} \int_0^1 t \intermediate{B}(t \intermediatecoord{\+x}) \diff t
    .
  \end{equation}
  Using the expression for the inverse transpose of the Jacobian matrix in \cref{eqn:toroidal_jacmat_invtranspose}, the covariant Piola transformation gives
  \begin{equation}
    (\logical{\pushforward}^1\logical{\curlt}^\dagger\logical{B})(\intermediatecoord{\+x})
    = \logical{\jacmat}^{-\transpose}\begin{bmatrix}0\\\int_0^{\logical{x}} \logical{B}(s,\logical{y}) \diff s\end{bmatrix}
    = \begin{bmatrix}
      -R(\logical{x}, \logical{y}) S(\logical{y})\\ 
      R(\logical{x}, \logical{y}) C(\logical{y})
    \end{bmatrix}
    \int_0^1 t \intermediate{B}(t \intermediatecoord{\+x}) \diff t
    =
    \begin{bmatrix}
      - \intermediatecoord{y}\\
      \intermediatecoord{x}
    \end{bmatrix}
    \int_0^1 t \intermediate{B}(t \intermediatecoord{\+x}) \diff t
    .
  \end{equation}
  The last expression is $(\intermediate{\curlt}^\dagger \intermediate{B})(\intermediatecoord{\+x})$ by \cref{eqn:homotopy_operators}, which establishes the first equality in \cref{eqn:pseudo_inverse_curl_logical}.
\end{statementproof}

Provided with the coefficient representation of the homotopy operators, we can now show that the homotopy identities~\eqref{eqn:homotopy_identities} hold exactly for continuous $0$-forms and bounded $1$-forms.
\begin{corollary}[Discrete homotopy identities]\label{cor:discrete_homotopy_identities}
  Let $P\in\funspace{\Bspline}^0$ be continuous at the center, that is, $p_{1j}=P(\+0)\mathds{1}_j^\yy$, let $\+A\in\funspace{\Bspline}^1$ be bounded at the center, that is, $a_{\yy,1j}=0$ (see the proof of \cref{lem:conformity_conditions_tensor_product_forms}), and let $B\in\funspace{\Bspline}^2$.
  Then the homotopy identities \eqref{eqn:homotopy_identities} hold exactly on $\Omega\setminus\{\+0\}$.
  In terms of the coefficients, they read
  \begin{equation}\label{eqn:discrete_homotopy_identities}
    p_{ij}=P(\+0)\mathds{1}_j^\yy+(\delta_\xx^\dagger\delta_\xx p)_{ij},
    \qquad
    \+a=\gradienth\delta_\xx^\dagger a_\xx
      +\begin{bmatrix}0\\\delta_\xx^\dagger\curlh\+a\end{bmatrix},
    \qquad
    b=\delta_\xx\delta_\xx^\dagger b.
  \end{equation}
\end{corollary}
\begin{statementproof}[appendix]{cor:discrete_homotopy_identities}
  We first show that \cref{eqn:homotopy_identities} is equivalent to \cref{eqn:discrete_homotopy_identities}.
  By \cref{eqn:pushforward_commutativity,eqn:discrete_operators}, the gradient and curl act on the coefficients as $\gradienth$ and $\curlh$, and $\gradientt^\dagger$ and $\curlt^\dagger$ act as $\delta_\xx^\dagger$ by \cref{lem:pseudo_inverse_curl_logical}.
  The constant $0$-form $\mathfrak{P}P=P(\+0)$ has the coefficients $P(\+0)\mathds{1}_j^\yy$ for every radial index $i$, by the radial partition of unity in \cref{eqn:bspline_partition} and the definition of $\mathds{1}^\yy$ in \cref{eqn:constant_function_coefficients}.
  The basis functions are linearly independent, so each identity in \cref{eqn:homotopy_identities} is equivalent to the corresponding identity in \cref{eqn:discrete_homotopy_identities}.

  Since $\delta_\xx$ and $\delta_\xx^\dagger$ act on the radial index only, the relations in \cref{eqn:discrete_radial_homotopy_relations} hold for 2D coefficients for each angular index $j$; in particular, $(\delta_\xx^\dagger\delta_\xx w)_{ij}=w_{ij}-w_{1j}$.
  The only additional ingredient is
  \begin{equation}\label{eqn:radial_angular_commutation}
    \delta_\yy\delta_\xx^\dagger v=\delta_\xx^\dagger\delta_\yy v,
  \end{equation}
  which holds because $\delta_\xx^\dagger$ acts only on the radial index and $\delta_\yy$ only on the angular index, for periodic B-splines as well as for trigonometric functions.
  This commutation is the only point where the two directions interact, and it produces the cancellation in the $1$-form identity.

  \emph{The $0$-form identity.}
  The first relation in \cref{eqn:discrete_radial_homotopy_relations} gives $(\delta_\xx^\dagger\delta_\xx p)_{ij}=p_{ij}-p_{1j}$, and $p_{1j}=P(\+0)\mathds{1}_j^\yy$ by continuity at the center.

  \emph{The $1$-form identity.}
  By \cref{eqn:discrete_operators}, $\gradienth\delta_\xx^\dagger a_\xx$ has the radial component $\delta_\xx\delta_\xx^\dagger a_\xx$ and the angular component $\delta_\yy\delta_\xx^\dagger a_\xx$.
  The radial component equals $a_\xx$ by the second relation in \cref{eqn:discrete_radial_homotopy_relations}.
  For the angular coefficients, the first relation and $a_{\yy,1j}=0$ give $\delta_\xx^\dagger\delta_\xx a_\yy=a_\yy$, and therefore
  \begin{equation}
    \delta_\xx^\dagger\curlh\+a
    =\delta_\xx^\dagger\delta_\xx a_\yy-\delta_\xx^\dagger\delta_\yy a_\xx
    =a_\yy-\delta_\yy\delta_\xx^\dagger a_\xx,
  \end{equation}
  where the last step uses \cref{eqn:radial_angular_commutation}.
  This is the angular component of the $1$-form identity in \cref{eqn:discrete_homotopy_identities}.

  \emph{The $2$-form identity.}
  This is the second relation in \cref{eqn:discrete_radial_homotopy_relations} with $v=b$.

  No condition on $B$ is needed because $\delta_\xx\delta_\xx^\dagger$ is the identity.
  In contrast, $\delta_\xx^\dagger\delta_\xx$ loses the first radial coefficients $w_{1j}$.
  For a $0$-form these are restored by $\mathfrak{P}P$, which requires continuity at the center.
  For the angular component of a $1$-form there is no such term, which is why $a_{\yy,1j}=0$ is required.
\end{statementproof}

\subsection{Enriched high-order polar subspaces}\label{sec:derham_hodge_2d_construction}
We proceed by enriching the polar $2$-form space, and subsequently use the homotopy operator $\curlt^\dagger$ to lift these $2$-forms to a curl-like part of the polar $1$-form space.
We use the convention $\intermediate{\funspace{T}}^k=\{0\}$, and hence $\pulledpolys{0}{k}=\pulledpolys{2}{k}=\{0\}$, for $k<0$.

The polar $2$-form space is defined as (cf.\ \cref{eqn:reduced_twoform_space})
\begin{equation}\label{eqn:polar_twoform_space}
  \polar{\funspace{\Bspline}}^{2,\pr} \defeq \reduced{\pi}^{2,1} \pulledpolys{2}{\pr-2} \oplus \funspace{\Bspline}^{2,\pr}_\external,
\end{equation}
where the projection operator $\reduced{\pi}^{2,1}$ projects onto the conforming $2$-form space $\reduced{\funspace{\Bspline}}^{2,1} = \funspace{\Bspline}^{2,1}_\external$, as discussed in \cref{lem:conforming_polar_twoform_projector}.
The projection operator $\reduced{\pi}^{2,1}$ is used to ensure that the polar $2$-form space is conforming, which, unlike the $0$-form projector $\pi^0$ used in \cref{eqn:polar_zeroform_space}, is not guaranteed by the $2$-form projection operator $\pi^2$.

The enriched part of the polar $2$-form space is obtained by projecting the pushed-forward polynomial $2$-forms of degree at most $\pr-2$,
\begin{equation}
  \pulledpolys{2}{\pr-2}
  \defeq
  \intermediate{\pushforward}^2\intermediate{\funspace{T}}^{\pr-2}
  =
  \{(\intermediate{T}/\intermediate{\jac})\circ\intermediatemapping^{-1} : \intermediate{T}\in\intermediate{\funspace{T}}^{\pr-2}\},
\end{equation}
which are smooth on $\Omega$ by \cref{ass:diffeomorphism}; for $\intermediatemapping=\mathrm{id}$ they are the physical polynomials.
Compared to the degree $\pr$ polar $0$-form space (cf.\ \cref{eqn:polar_zeroform_space}), the polar $2$-form space is reduced by two degrees.
This is because the two-form tensor-product space has one degree less in the radial direction (see \cref{sec:logical_function_spaces}), and the pullback $\pushforward^{-2}\pulledpolys{2}{\pr-2}=\logical{\pushforward}^{-2}\intermediate{\funspace{T}}^{\pr-2}$ (\cref{eqn:pushforward_composition}) involves multiplication by the Jacobian $\logical{\jac}$ of the logical mapping, which is proportional to $\logical{x}$ at the center as a consequence of \cref{eqn:coordinate_transformation_assumptions}.
The space $\pulledpolys{2}{\pr-2}$ has dimension $n^{2,\pr}=(\pr-1)\pr/2$; for $\pr=1$ it is zero, consistently with \citet{toshniwal2021,holderied2022,yaman2025}.

\begin{table}
  \centering
  \setlength{\tabcolsep}{1.000000000000pt}
  \def\allformsarrowcolor{black!60}
  \colorlet{allformstwobg}{lines8!12}
  \def\allformsscalaropacity{0.65}
  \begin{tabular}{@{}cccc@{}}
    \cmidrule(l){3-4}
    &&\multicolumn{2}{c}{\makebox[0pt][c]{$\polar{\Bspline}^{2,3}_{l} = \curlt \polar{\+\Bspline}^{1,3}_{\mathrm{C},l}$}}\\
    \tikzsetnextfilename{allforms_zero_1}\input{tikz/allforms_zero_1.tex}\unskip&\tikzsetnextfilename{allforms_zero_2}\input{tikz/allforms_zero_2.tex}\unskip&\cellcolor{allformstwobg}\tikzsetnextfilename{allforms_two_1}\input{tikz/allforms_two_1.tex}\unskip&\cellcolor{allformstwobg}\tikzsetnextfilename{allforms_two_2}\input{tikz/allforms_two_2.tex}\unskip \\
    \tikzsetnextfilename{allforms_zero_3}\input{tikz/allforms_zero_3.tex}\unskip&\tikzsetnextfilename{allforms_zero_4}\input{tikz/allforms_zero_4.tex}\unskip&\tikzsetnextfilename{allforms_zero_5}\input{tikz/allforms_zero_5.tex}\unskip&\cellcolor{allformstwobg}\tikzsetnextfilename{allforms_two_3}\input{tikz/allforms_two_3.tex}\unskip \\
    \tikzsetnextfilename{allforms_zero_6}\input{tikz/allforms_zero_6.tex}\unskip&\tikzsetnextfilename{allforms_zero_7}\input{tikz/allforms_zero_7.tex}\unskip&\tikzsetnextfilename{allforms_zero_8}\input{tikz/allforms_zero_8.tex}\unskip&\tikzsetnextfilename{allforms_zero_9}\input{tikz/allforms_zero_9.tex}\unskip \\
    \multicolumn{4}{c}{\makebox[0pt][c]{$\polar{\+\Bspline}^{1,3}_{\mathrm{G},l} = \gradientt \polar{\Bspline}^{0,3}_{l+1}$}}\\
    \bottomrule
  \end{tabular}
  \caption{Polar $m$-forms on the physical banana-shaped domain of \cref{fig:disk_mapping}: the logical mapping $\logicalmapping$ and the intermediate mapping $\intermediatemapping$ coincide with those of that figure. The $0$-forms are shown in the bottom-left triangle (except for the constant $0$-form $\polar{\Bspline}^{0,\pr}_{1}$), and correspond exactly to the lowest block of polar $0$-forms in \cref{tab:disc_basis_functions}. The $2$-forms are shown in the shaded top-right triangle. The corresponding $1$-forms are shown by the arrows.}
  \label{tab:allforms}
\end{table}

The homotopy identities~\eqref{eqn:homotopy_identities} motivate the following direct-sum decomposition of the $1$-form space (cf.\ \cref{eqn:reduced_oneform_space})
\begin{equation}\label{eqn:polar_oneform_space}
  \polar{\funspace{\Bspline}}^{1,\pr}
  \defeq \gradientt\pi^0\pulledpolys{0}{\pr}
  \oplus \curlt^\dagger\reduced{\pi}^{2,1}\pulledpolys{2}{\pr-2}
  \oplus \funspace{\Bspline}^{1,\pr}_\external.
\end{equation}
The resulting polar basis functions for $\pr=3$ are shown in \cref{tab:allforms}, and the internal parts in \cref{eqn:polar_zeroform_space,eqn:polar_oneform_space,eqn:polar_twoform_space} use the following numbers of basis functions:
\begin{equation}\label{eqn:polar_subspaces_overview_dimensions}
  n^{0,\pr} = \frac{(\pr+1)(\pr+2)}{2},\quad
  n^{1,\pr}_\mathrm{G} = \frac{(\pr + 3)\pr}{2}, \quad
  n^{1,\pr}_\mathrm{C} = \frac{(\pr-1)\pr}{2}, \quad
  n^{2,\pr} = \frac{(\pr-1)\pr}{2},
\end{equation}
Thus $n^{1,\pr}=n^{1,\pr}_\mathrm{G}+n^{1,\pr}_\mathrm{C}=(\pr+1)\pr$. 
These counts give the internal dimensions whenever the basis functions are linearly independent, as proved in \cref{lem:polar_basis_pr_compatible} under the assumption of $\pr$-compatibility; \cref{tab:polar_dimensions} lists them for $\pr\le6$.
\begin{table}
  \centering
  \begin{tabular}{@{}lccccccc@{}}
    \toprule
    $\pr$ & $0$ & $1$ & $2$ & $3$ & $4$ & $5$ & $6$\\
    \midrule
    $n^{0,\pr}$ & $1$ & $3$ & $6$ & $10$ & $15$ & $21$ & $28$\\
    $n^{1,\pr}_\mathrm{G}$ & $0$ & $2$ & $5$ & $9$ & $14$ & $20$ & $27$\\
    $n^{1,\pr}_\mathrm{C}$ & $0$ & $0$ & $1$ & $3$ & $6$ & $10$ & $15$\\
    $n^{2,\pr}$ & $0$ & $0$ & $1$ & $3$ & $6$ & $10$ & $15$\\
    \bottomrule
  \end{tabular}
  \caption{Numbers of internal polar basis functions, \cref{eqn:polar_subspaces_overview_dimensions}, for $\pr=0,\ldots,6$.}
  \label{tab:polar_dimensions}
\end{table}

\subsubsection{Choice of basis functions}\label{sec:derham_hodge_2d_basis}
We now choose a polynomial basis on the intermediate domain, retain the internal coefficients of the full $0$-form projection, and apply a compactly supported correction to the projected $2$-forms.

Let $\{\intermediate{Z}_l\}_{l=1}^{n^{0,\pr}}$ denote the Zernike basis \citep{zernike1934} of $\intermediate{\funspace{T}}^\pr$ in the intermediate coordinates $\intermediatecoord{\+x}$, in the standard (OSA/ANSI) ordering by total degree \citep[Sec.~2.1.2]{niu2022}, indexed from $l=1$ so that $\intermediate{Z}_1=1$.
Only the grading by total degree and the center-value property below are used in the construction, so any graded polynomial basis would serve; we use the Zernike basis, the standard orthogonal polynomial basis on the disk, for definiteness, and shift it so that all but the first basis function vanish at the center:
\begin{equation}\label{eqn:shifted_zernike_basis}
  \intermediate{T}_1\defeq 1,
  \quad
  \intermediate{T}_l\defeq \intermediate{Z}_l-\intermediate{Z}_l(\+0),
  \quad l=2,\ldots,n^{0,\pr},
  \quad
  \intermediate{\funspace{T}}^\pr=\spn\{\intermediate{T}_l\}_{l=1}^{n^{0,\pr}}.
\end{equation}
This invertible change of basis preserves the polynomial space, and by the ordering by total degree, $\{\intermediate{T}_l\}_{l=1}^{n^{2,\pr}}$ is a basis of $\intermediate{\funspace{T}}^{\pr-2}$.
It gives the useful center-value property
\begin{equation}\label{eqn:polynomial_basis_centre_values}
  \intermediate{T}_1(\+0)=1
  , \quad
  \intermediate{T}_l(\+0)=0
  , \quad 
  l = 2, \ldots, n^{0,\pr}
  ,
\end{equation}
which we exploit in the coefficient-based representation of the homotopy operators in \cref{thm:polar_subspaces_discrete_homotopy}.
The corresponding physical basis functions are the pushed-forward polynomial $0$- and $2$-forms
\begin{equation}
  T^0_l\defeq\intermediate{\pushforward}^0\intermediate{T}_l=\intermediate{T}_l\circ\intermediatemapping^{-1},
  \qquad
  T^2_l\defeq\intermediate{\pushforward}^2\intermediate{T}_l=(\intermediate{T}_l/\intermediate{\jac})\circ\intermediatemapping^{-1},
\end{equation}
which form bases of $\pulledpolys{0}{\pr}$ ($l\le n^{0,\pr}$) and $\pulledpolys{2}{\pr-2}$ ($l\le n^{2,\pr}$), since $\intermediate{\pushforward}^0$ and $\intermediate{\pushforward}^2$ are injective.
As $\intermediatemapping(\+0)=\+0$, the center values are preserved, $T^0_l(\+0)=\intermediate{T}_l(\+0)$.
The full pullbacks $\pushforward^{-m}T^m_l=\logical{\pushforward}^{-m}\intermediate{T}_l$ (\cref{eqn:pushforward_composition}) do not depend on $\intermediatemapping$; hence neither do their GDOFs $\gdof^m(T^m_l)=\logical{\gdof}^m(\logical{\pushforward}^{-m}\intermediate{T}_l)$ (\cref{eqn:physical_geometric_dofs}) nor the coefficients of the polar basis functions below, which are obtained from these GDOFs by coefficient operations.

With this basis for $\pulledpolys{0}{\pr}$, we can construct the basis for the polar $0$-form space by applying the projection operator $\pi^0$ to each of the pulled-back shifted Zernike polynomials, resulting in
\begin{equation}
  \polar{\Bspline}^{0,\pr}_{l} \defeq \chi^{0,\pr}_\internal \pi^0 T^0_l, \quad l \in \{1, \hdots, n^{0,\pr}\},
\end{equation}
where $\chi^{m,\pr}_\internal$ denotes the truncation operator that sets the external coefficients of a polar $m$-form to zero (i.e., the coefficients corresponding to $(i,j)\in\indextwo^{m,\pr}_\external$).
Analogously, $\chi^{m,\pr}_\external$ denotes the truncation operator that sets the internal coefficients of a polar $m$-form to zero.
The truncated basis functions are compactly supported, with coefficients on the internal rings only, which keeps the extraction matrices of \cref{sec:matrix_representation} sparse.
Note that truncation changes the basis, not the space, and does not imply that any approximation is made: the removed external contributions are contained in the external part of the polar $0$-form space, given by $\funspace{\Bspline}^{0,\pr}_\external$.
It follows that the truncated coefficient representation in terms of the $0$-form basis functions $\Bspline^0_{ij}$ is given by
\begin{equation}\label{eqn:polar_zeroform_basis}
  \polar{\Bspline}^{0,\pr}_{l} = \sum_{(i,j)\in\indextwo^{0,\pr}_\internal} \polar{\bspline}^{0,\pr}_{l,ij} \Bspline^0_{ij}
  , \quad
  l \in \{1, \hdots, n^{0,\pr}\}
  ,
\end{equation}
where $\polar{\bspline}^{0,\pr}_{l,ij}$ denote the coefficients of the $l$-th polar $0$-form basis function $\polar{\Bspline}^{0,\pr}_{l}$.

Computing the gradient of the polar $0$-form basis functions in \cref{eqn:polar_zeroform_basis} yields a basis for the gradient-like polar $1$-form space, which is given by
\begin{equation}\label{eqn:polar_oneform_gradient_basis}
  \polar{\+\Bspline}^{1,\pr}_{\mathrm{G},l} 
  \defeq 
  \gradientt \polar{\Bspline}^{0,\pr}_{l+1}
  =
  \gradientt (\chi^{0,\pr}_\internal \pi^0 T^0_{l+1})
  , \quad
  l \in \{1, \hdots, n^{1,\pr}_\mathrm{G}\}
  .
\end{equation}
We leave out the first polar $0$-form basis function because $1\in\funspace{\Bspline}^0$ by \cref{ass:constant_functions}, so that $\pi^0T^0_1=1$.
Indeed, $1-\polar{\Bspline}^{0,\pr}_1$ belongs to $\funspace{\Bspline}^{0,\pr}_\external$ by truncation, so $\gradientt\polar{\Bspline}^{0,\pr}_1=-\gradientt (1-\polar{\Bspline}^{0,\pr}_1)$ already belongs to $\funspace{\Bspline}^{1,\pr}_\external$ by the external index ranges in \cref{eqn:polar_zeroform_space_odot,eqn:polar_mform_space_odot}.
We note that although the resulting gradient-like polar $1$-form basis functions are compactly supported, their span is not a subspace of $\funspace{\Bspline}^{1,\pr}_\internal$, as a consequence of the derivative relation~\eqref{eqn:bspline_derivative} of the 1D B-spline space used in the $\logical{x}$-direction.

The compactly supported polar $2$-form basis functions are obtained from each conforming projected pushed-forward polynomial by subtracting its reduced coefficient projection $\reduced{\Pi}^{2,\pr}$, defined in \cref{eqn:reduced_twoform_coefficient_projector}:
\begin{equation}\label{eqn:polar_twoform_curl_basis}
  \polar{\Bspline}^{2,\pr}_{l}
  \defeq
  (I^2-\reduced{\Pi}^{2,\pr})\reduced{\pi}^{2,1}T^2_l,
  \qquad l\in\{1,\ldots,n^{2,\pr}\}.
\end{equation}
For a $2$-form $B$ with coefficients $b_{ij}$, the compactly supported construction has the explicit coefficients
\begin{equation}\label{eqn:polar_twoform_compact_coefficients}
  \bigl((I^2-\reduced{\Pi}^{2,\pr})B\bigr)_{ij}
  =
  \begin{cases}
    b_{ij}, & i\le\pr,\\
    -\displaystyle\sum_{i'=1}^{\pr}b_{i'j}, & i=\pr+1,\\
    0, & i\ge\pr+2.
  \end{cases}
\end{equation}
Thus the first $\pr$ rings are retained, ring $\pr+1$ contains their negative radial sum, and all later rings vanish. In particular, every compactly supported $2$-form basis function has zero radial coefficient sum for each $j$.
Because $\reduced{\Pi}^{2,\pr}\reduced{\pi}^{2,1}T^2_l$ belongs to $\funspace{\Bspline}^{2,\pr}_\external$, these compactly supported basis functions and the external basis span the space in \cref{eqn:polar_twoform_space}.

The curl-like polar $1$-form basis functions are their curl homotopy lifts:
\begin{equation}\label{eqn:polar_oneform_curl_basis_definition}
  \polar{\+\Bspline}^{1,\pr}_{\mathrm{C},l}
  \defeq
  \curlt^\dagger\polar{\Bspline}^{2,\pr}_l
  =
  \chi^{1,\pr}_\internal\curlt^\dagger(\reduced{\pi}^{2,1}T^2_l),
  \qquad l\in\{1,\ldots,n^{1,\pr}_\mathrm{C}\}.
\end{equation}
The equality follows from the cumulative-sum formula in \cref{eqn:pseudo_inverse_curl_logical}: the coefficients of the lift vanish from angular $1$-form ring $\pr+2$ onward because of the compensating $2$-form ring in \cref{eqn:polar_twoform_compact_coefficients}.
Consequently the lift is compactly supported and has the coefficient representation
\begin{equation}\label{eqn:polar_oneform_curl_basis}
  \polar{\+\Bspline}^{1,\pr}_{\mathrm{C},l} =
  \sum_{(i,j)\in\indextwo^{1,\pr}_{\xx,\internal}} \polar{\bspline}^{1,\pr}_{\mathrm{C},l,\xx,ij} \+\Bspline^1_{\xx,ij}
  +
  \sum_{(i,j)\in\indextwo^{1,\pr}_{\yy,\internal}} \polar{\bspline}^{1,\pr}_{\mathrm{C},l,\yy,ij} \+\Bspline^1_{\yy,ij},
  \qquad l\in\{1,\ldots,n^{1,\pr}_\mathrm{C}\}.
\end{equation}
By \cref{cor:discrete_homotopy_identities}, $\curlt\polar{\+\Bspline}^{1,\pr}_{\mathrm{C},l}=\polar{\Bspline}^{2,\pr}_l$.
Since $\curlt^\dagger\reduced{\Pi}^{2,\pr}\reduced{\pi}^{2,1}T^2_l$ belongs to $\funspace{\Bspline}^{1,\pr}_\external$ by \cref{eqn:pseudo_inverse_curl_logical}, replacing the unmodified lifts by these compactly supported lifts leaves the $1$-form space in \cref{eqn:polar_oneform_space} unchanged.

The basis functions for the polar $m$-form spaces with $\pr = 3$ are visualized in \cref{tab:allforms}, and the corresponding extraction matrices are given in \cref{eqn:polar_extraction_matrices}.
In terms of these basis functions, the polar $m$-forms $P \in \polar{\funspace{\Bspline}}^{0,\pr}, \+A \in \polar{\funspace{\Bspline}}^{1,\pr}, \+B \in \polar{\funspace{\Bspline}}^{2,\pr}$ can be expressed as
\begin{subequations}
  \label{eqn:polar_forms_basis_expansion}
  \begin{align+}
    P &= \sum_{l=1}^{n^{0,\pr}} p_l \polar{\Bspline}^{0,\pr}_{l} + \sum_{(i,j)\in\indextwo^{0,\pr}_\external} p_{ij} \Bspline^0_{ij}\label{eqn:polar_zeroform_basis_expansion}\\
    \+A &= \sum_{l=1}^{n^{1,\pr}_\mathrm{G}} a_{\mathrm{G},l} \polar{\+\Bspline}^{1,\pr}_{\mathrm{G},l} + \sum_{l=1}^{n^{1,\pr}_\mathrm{C}} a_{\mathrm{C},l} \polar{\+\Bspline}^{1,\pr}_{\mathrm{C},l} + \sum_{(i,j)\in\indextwo^{1,\pr}_{\xx,\external}} a_{\xx,ij} \+\Bspline^1_{\xx,ij} + \sum_{(i,j)\in\indextwo^{1,\pr}_{\yy,\external}} a_{\yy,ij} \+\Bspline^1_{\yy,ij}\label{eqn:polar_oneform_basis_expansion}\\
    B &= \sum_{l=1}^{n^{2,\pr}} b_l \polar{\Bspline}^{2,\pr}_{l} + \sum_{(i,j)\in\indextwo^{2,\pr}_\external} b_{ij} \Bspline^2_{ij}.\label{eqn:polar_twoform_basis_expansion}
  \end{align+}
\end{subequations}
Linear independence of these families, which makes them bases of the polar spaces, is established under $\pr$-compatibility in \cref{lem:polar_basis_pr_compatible} and is assumed in \cref{sec:derham_hodge_2d_geom_projections}, where it is needed for the construction of projection operators.

\subsection{Properties of the polar subspaces}\label{sec:derham_hodge_2d_properties}
In what follows we discuss the structural properties of the polar subspaces.
We first show that they form a cochain complex of conforming spaces, that the homotopy operators form a chain complex on them, and that they have the same cohomology as the continuous de Rham complex on $\Omega$, which is that of the disk.
We then introduce a compatibility condition on the mapping and the logical function spaces, under which the polar families form bases and the polar subspaces satisfy exact regularity conditions at the center of the domain.

We first verify the cochain property.
\begin{theorem}[Cochain complex]\label{thm:polar_subspaces_commutativity_conformity}
  Let $\pr\ge0$.
  The polar subspaces form a cochain complex (the forward direction of \cref{eqn:derham_hodge_2d_diagram}):
  \begin{equation}
    \gradientt\polar{\funspace{\Bspline}}^{0,\pr}
      \subset\polar{\funspace{\Bspline}}^{1,\pr},
    \qquad
    \curlt\polar{\funspace{\Bspline}}^{1,\pr}
      \subset\polar{\funspace{\Bspline}}^{2,\pr},
    \qquad
    \curlt\gradientt=0.
  \end{equation}
  At $\pr=0$, this is an algebraic coefficient complex; conformity for $\pr\ge1$ is established in \cref{thm:polar_subspaces_conformity}.
\end{theorem}
\begin{statementproof}[appendix]{thm:polar_subspaces_commutativity_conformity}
  The gradient of $\pi^0\pulledpolys{0}{\pr}$ belongs to the polar $1$-form space by \cref{eqn:polar_oneform_space}, and its curl is zero.
  The curl of $\curlt^\dagger\reduced{\pi}^{2,1}\pulledpolys{2}{\pr-2}$ is $\reduced{\pi}^{2,1}\pulledpolys{2}{\pr-2}$ by \cref{cor:discrete_homotopy_identities}, so it belongs to \cref{eqn:polar_twoform_space}.

  For an external $0$-form with coefficients $p_{ij}=0$ for $i\le\pr+1$, the radial differences $\delta_\xx p$ vanish for $i\le\pr$, and the angular derivatives $\delta_\yy p$ vanish for $i\le\pr+1$.
  Hence its gradient belongs to $\funspace{\Bspline}^{1,\pr}_\external$.
  Likewise, an external $1$-form has $a_{\xx,ij}=0$ for $i\le\pr$ and $a_{\yy,ij}=0$ for $i\le\pr+1$.
  Both terms of $\curlh\+a=\delta_\xx a_\yy-\delta_\yy a_\xx$ therefore vanish for $i\le\pr$, so its curl belongs to $\funspace{\Bspline}^{2,\pr}_\external$.
  These arguments apply to both choices of angular basis because $\delta_\yy$ preserves the radial index.
  Finally, $\curlt\gradientt=0$ follows from the commutation of $\delta_\xx$ and $\delta_\yy$.
\end{statementproof}

Next, we show that the polar $m$-form spaces are conforming.
\begin{theorem}[Conformity]\label{thm:polar_subspaces_conformity}
  We assume $\pr \ge 1$.
  Then the polar subspaces are conforming subspaces of the infinite-dimensional function spaces, i.e.,
  \begin{equation}
    \polar{\funspace{\Bspline}}^{0,\pr} \subset \funspace{H}_1, \quad
    \polar{\funspace{\Bspline}}^{1,\pr} \subset \funspace{H}_{\curlt}, \quad
    \polar{\funspace{\Bspline}}^{2,\pr} \subset \funspace{L}_2.
  \end{equation}
\end{theorem}
\begin{statementproof}[appendix]{thm:polar_subspaces_conformity}
  The polar $0$-form space belongs to $\funspace{H}_1$ by \cref{thm:polar_zeroform_space_conforming}.

  \paragraph{Conformity of the polar $2$-form space}
  For $B\in\funspace{\Bspline}^{2,\pr}_\external$, the index range in \cref{eqn:polar_mform_space_odot} gives $b_{1j}=0$ because $\pr\ge1$.
  The conforming projected pulled-back polynomial $\reduced{\pi}^{2,1}T$, $T\in\pulledpolys{2}{\pr-2}$, has $b_{1j}=0$ by \cref{lem:conforming_polar_twoform_projector}.
  Thus every polar $2$-form belongs to $\funspace{L}_2$ by \cref{lem:conformity_conditions_tensor_product_forms}.

  \paragraph{Conformity of the polar $1$-form space}
  If $P\in\pi^0\pulledpolys{0}{\pr}$, then $P\in\funspace{H}_1$ by the $0$-form result.
  Hence $\gradientt P$ is square-integrable and its weak curl vanishes, so the gradient-like part belongs to $\funspace{H}_{\curlt}$.

  For the curl-like part, take $B\in\reduced{\pi}^{2,1}\pulledpolys{2}{\pr-2}$.
  Its lift $\curlt^\dagger B$ has, by \cref{eqn:pseudo_inverse_curl_logical},
  \begin{equation}
    a_{\xx,ij}=0,\qquad a_{\yy,1j}=0,\qquad
    a_{\yy,2j}=b_{1j}=0.
  \end{equation}
  These are the $1$-form conditions in \cref{lem:conformity_conditions_tensor_product_forms}, so the lift belongs to $\funspace{H}_{\curlt}$.

  Finally, an external $1$-form has $a_{\xx,1j}=a_{\yy,1j}=a_{\yy,2j}=0$ because $\pr\ge1$ and the index ranges in \cref{eqn:polar_mform_space_odot} exclude these coefficients.
  It therefore satisfies the same conformity conditions.
  The full polar $1$-form space belongs to $\funspace{H}_{\curlt}$ by linearity.
\end{statementproof}

Using the coefficient-based representation of the polar $m$-forms from \cref{eqn:polar_forms_basis_expansion}, we can now derive the coefficient-based representation of the homotopy operators.
\begin{lemma}[Discrete homotopy operators]\label{thm:polar_subspaces_discrete_homotopy}
  Let $\pr \ge 0$.
  We further assume that the basis of $\intermediate{\funspace{T}}^\pr$ satisfies \cref{eqn:polynomial_basis_centre_values}.
  Let $P$, $\+A$, and $B$ have the coefficient expansions in \cref{eqn:polar_forms_basis_expansion}.
  Then the homotopy operators have the coefficient representations
  \begin{subequations}\label{eqn:polar_homotopy_operators}
    \begin{align}
      \mathfrak{P} P
      &=p_1,
      \label{eqn:polar_evaluation_at_centre}\\
      \gradientt^\dagger\+A
      &=\sum_{l=1}^{n^{1,\pr}_{\mathrm{G}}}
        a_{\mathrm{G},l}\polar{\Bspline}^{0,\pr}_{l+1}
      +
      \sum_{(i,j)\in\indextwo^{0,\pr}_\external}
      \srpar{\sum_{i'=\pr+1}^{i-1}a_{\xx,i'j}}
      \Bspline^0_{ij},
      \label{eqn:polar_gradient_homotopy_coefficients}\\
      \curlt^\dagger B
      &=\sum_{l=1}^{n^{2,\pr}}b_l
        \polar{\+\Bspline}^{1,\pr}_{\mathrm{C},l}
      +\sum_{(i,j)\in\indextwo^{1,\pr}_{\yy,\external}}
        \srpar{\sum_{i'=\pr+1}^{i-1}b_{i'j}}
        \+\Bspline^1_{\yy,ij}.
      \label{eqn:polar_curl_homotopy_coefficients}
    \end{align}
  \end{subequations}
  The identities in \cref{eqn:homotopy_identities} hold for every $P\in\polar{\funspace{\Bspline}}^{0,\pr}$, $\+A\in\polar{\funspace{\Bspline}}^{1,\pr}$, and $B\in\polar{\funspace{\Bspline}}^{2,\pr}$.
\end{lemma}
\begin{statementproof}[appendix]{thm:polar_subspaces_discrete_homotopy}
  The $0$-form projection operator based on GDOFs reproduces values at the center (see the proof of \cref{thm:polar_zeroform_space_conforming}), and truncation does not remove the first radial ring.
  Consequently, $\polar{\Bspline}^{0,\pr}_l(\+0)=T^0_l(\+0)=\intermediate{T}_l(\+0)$.
  Together with \cref{eqn:polynomial_basis_centre_values} and the fact that the external $0$-form basis functions vanish at the center, this proves \cref{eqn:polar_evaluation_at_centre}.

  We next derive the two homotopy formulas.  By
  \cref{eqn:polar_oneform_gradient_basis,eqn:discrete_radial_homotopy_relations},
  \begin{equation}
    \gradientt^\dagger
    \polar{\+\Bspline}^{1,\pr}_{\mathrm{G},l}
    =\gradientt^\dagger\gradientt
      \polar{\Bspline}^{0,\pr}_{l+1}
    =\polar{\Bspline}^{0,\pr}_{l+1},
  \end{equation}
  where the last equality follows from
  $\polar{\Bspline}^{0,\pr}_{l+1}(\+0)=\intermediate{T}_{l+1}(\+0)=0$.
  Every curl-like polar basis function has vanishing logical $\xx$-component,
  as do the external basis functions $\+\Bspline^1_{\yy,ij}$.  Hence
  \cref{eqn:pseudo_inverse_gradient_logical} gives
  \begin{equation}
    \gradientt^\dagger
      \polar{\+\Bspline}^{1,\pr}_{\mathrm{C},l}=0,
    \quad
    \gradientt^\dagger\+\Bspline^1_{\yy,ij}=0.
  \end{equation}
  The same coefficient formula shows that the external $\xx$-component of
  $\+A$ contributes only to the external $0$-form coefficients, with
  \begin{equation}
    (\delta_\xx^\dagger a_\xx)_{ij}
    =\sum_{i'=\pr+1}^{i-1}a_{\xx,i'j},
    \quad (i,j)\in\indextwo^{0,\pr}_\external,
  \end{equation}
  because its coefficients vanish for $i'\leq\pr$.  Applying
  these identities term by term to
  \cref{eqn:polar_oneform_basis_expansion} yields
  \cref{eqn:polar_gradient_homotopy_coefficients}.

  For each compactly supported internal $2$-form basis function, \cref{eqn:polar_oneform_curl_basis_definition} directly gives
  $\curlt^\dagger\polar{\Bspline}^{2,\pr}_l
   =\polar{\+\Bspline}^{1,\pr}_{\mathrm{C},l}$.
  On the external part, \cref{eqn:pseudo_inverse_curl_logical} gives
  \begin{equation}
    (\delta_\xx^\dagger b)_{ij}
    =\sum_{i'=\pr+1}^{i-1}b_{i'j},
    \quad (i,j)\in\indextwo^{1,\pr}_{\yy,\external},
  \end{equation}
  because an external $2$-form has zero coefficients for $i'\leq\pr$.
  This proves \cref{eqn:polar_curl_homotopy_coefficients}.

  It remains to verify the hypotheses of
  \cref{cor:discrete_homotopy_identities}.  Every polar $0$-form is continuous
  at the center: its full tensor-product coefficients on the first radial ring satisfy
  $p_{1j}=P(\+0)\mathds{1}_j^\yy$ by the center-value property of $\pi^0$,
  while the external basis functions have a vanishing first radial ring.
  Similarly, the $\yy$-coefficients in the
  first radial ring of every polar $1$-form vanish.  This holds for a
  gradient-like basis function because it is the gradient of a continuous
  $0$-form, for a curl-like basis function by the radial integration in
  \cref{eqn:pseudo_inverse_curl_logical}, and for the external basis functions
  by their index ranges in \cref{eqn:polar_mform_space_odot}.  Thus the
  assumptions of \cref{cor:discrete_homotopy_identities} are satisfied, and
  that corollary proves all three identities in
  \cref{eqn:homotopy_identities} on the polar subspaces.
\end{statementproof}

Finally, we show that the polar subspaces additionally yield a chain complex with the homotopy operators, and that the polar $1$-form space can be split into a direct sum of the image of the gradient operator under the polar $0$-form space and the image of the curl homotopy operator under the polar $2$-form space.
\begin{theorem}[Homotopy chain complex]\label{thm:polar_subspaces_homotopy_splitting}
  The homotopy operators form a chain complex on the polar subspaces (the backward direction of \cref{eqn:derham_hodge_2d_diagram}), i.e.,
  \begin{equation}\label{eqn:polar_subspaces_homotopy_preservation}
    \gradientt^\dagger\polar{\funspace{\Bspline}}^{1,\pr}
    \subset\polar{\funspace{\Bspline}}^{0,\pr},\quad
    \curlt^\dagger\polar{\funspace{\Bspline}}^{2,\pr}
    \subset\polar{\funspace{\Bspline}}^{1,\pr},\quad
    \gradientt^\dagger\curlt^\dagger=0.
  \end{equation}
  Moreover, the polar $1$-form space has the algebraic direct-sum decomposition
  \begin{equation}\label{eqn:polar_oneform_homotopy_splitting}
    \polar{\funspace{\Bspline}}^{1,\pr}
    =\gradientt\polar{\funspace{\Bspline}}^{0,\pr}
    \oplus
    \curlt^\dagger\polar{\funspace{\Bspline}}^{2,\pr}.
  \end{equation}
\end{theorem}
\begin{statementproof}[appendix]{thm:polar_subspaces_homotopy_splitting}
  The explicit formulas in \cref{eqn:polar_homotopy_operators} give the two inclusions in \cref{eqn:polar_subspaces_homotopy_preservation}.
  By \cref{eqn:pseudo_inverse_gradient_logical,eqn:pseudo_inverse_curl_logical}, $\curlt^\dagger B$ has zero logical radial component and $\gradientt^\dagger$ integrates only that component.
  Hence $\gradientt^\dagger\curlt^\dagger=0$, proving the chain-complex property.
  The cochain property in \cref{thm:polar_subspaces_commutativity_conformity} and the second inclusion show that the right-hand side of \cref{eqn:polar_oneform_homotopy_splitting} is contained in the polar $1$-form space.

  Conversely, for $\+A\in\polar{\funspace{\Bspline}}^{1,\pr}$, the homotopy identity from \cref{thm:polar_subspaces_discrete_homotopy} gives
  \begin{equation}
    \+A=\gradientt(\gradientt^\dagger\+A)
      +\curlt^\dagger(\curlt\+A).
  \end{equation}
  Here $\gradientt^\dagger\+A$ belongs to the polar $0$-form space by \cref{eqn:polar_subspaces_homotopy_preservation}, and $\curlt\+A$ belongs to the polar $2$-form space by the cochain property.
  This proves the reverse inclusion.
  The sum is direct: if $\gradientt P=\curlt^\dagger B$, then applying $\curlt$ and using $\curlt\gradientt=0$ and $\curlt\curlt^\dagger B=B$ gives $B=0$, and hence $\gradientt P=0$.
  Thus the two field components are unique.
  Since $\Omega$ is connected, $P$ is unique up to a constant; the displayed choice $P=\gradientt^\dagger\+A$ fixes $P(\+0)=0$.
\end{statementproof}

\begin{corollary}[Cohomology]\label{cor:polar_subspaces_cohomology}
  The cochain complex in \cref{thm:polar_subspaces_commutativity_conformity} has the same cohomology as the continuous de Rham complex on $\Omega$, which is that of the disk:
  \begin{equation}
    \dim(\kernel\gradientt)=1,
    \qquad
    \dim(\kernel\curlt/\image\gradientt)=0,
    \qquad
    \dim(\polar{\funspace{\Bspline}}^{2,\pr}/\image\curlt)=0.
  \end{equation}
\end{corollary}
\begin{statementproof}[appendix]{cor:polar_subspaces_cohomology}
  Since $\Omega$ is connected, every polar $0$-form with zero gradient is constant.
  Conversely, $1\in\funspace{\Bspline}^0$ by \cref{ass:constant_functions}, so $1=\pi^0 1\in\pi^0\pulledpolys{0}{\pr}\subseteq\polar{\funspace{\Bspline}}^{0,\pr}$ by \cref{eqn:polar_zeroform_space}; hence $\kernel\gradientt=\spn\{1\}$.
  By \cref{eqn:polar_oneform_homotopy_splitting}, any polar $1$-form has the form $\+A=\gradientt P+\curlt^\dagger B$ with $B$ in the polar $2$-form space.
  Applying $\curlt$ gives $\curlt\+A=B$ by the polar homotopy identity in \cref{thm:polar_subspaces_discrete_homotopy}, so $\kernel\curlt=\image\gradientt$.
  The same identity and \cref{eqn:polar_subspaces_homotopy_preservation} show that every polar $2$-form is the curl of its polar lift $\curlt^\dagger B$.
\end{statementproof}

The results above hold for every disk-like logical mapping and every intermediate mapping.
Two further properties, sufficient conditions for linear independence of the families in \cref{eqn:polar_forms_basis_expansion} and exact regularity at the center, rely on the projected pulled-back polynomials being the pulled-back polynomials themselves, that is, $\pi^0T=T$ for $T\in\pulledpolys{0}{\pr}$ and $\reduced{\pi}^{2,1}T=T$ for $T\in\pulledpolys{2}{\pr-2}$.
As for the polar $0$-form space in \cref{thm:polar_zeroform_space_regularity}, this requires a compatibility condition on the logical mapping and the logical function spaces.

\begin{definition}[$\pr$-compatibility of the logical mapping and logical spaces]\label{def:pr_compatibility_mapping_logical_spaces}
  We say that the logical mapping $\logicalmapping$ and the logical $m$-form spaces $\logical{\funspace{\Bspline}}^m$ are \emph{$\pr$-compatible} if the logical $0$-form space is $\pr$-compatible in the sense of \cref{def:pr_compatibility_mapping_logical_zeroform}, and if the pullback of the intermediate polynomial $2$-forms can be exactly represented in the logical $2$-form space:
  \begin{equation}\label{eqn:pr_compatibility_twoform}
    \logical{\pushforward}^{-2}\intermediate{\funspace{T}}^{\pr-2}
    \subset \logical{\funspace{\Bspline}}^2,
    \quad\text{equivalently}\quad
    \pushforward^{-2}\pulledpolys{2}{\pr-2}
    \subset \logical{\funspace{\Bspline}}^2,
  \end{equation}
  where the equivalence follows from $\pulledpolys{2}{\pr-2}=\intermediate{\pushforward}^2\intermediate{\funspace{T}}^{\pr-2}$ and $\pushforward^{-2}=\logical{\pushforward}^{-2}\circ\intermediate{\pushforward}^{-2}$ (\cref{eqn:pushforward_composition}).
  The condition does not involve the intermediate mapping $\intermediatemapping$.
\end{definition}
\begin{example}[$\pr$-compatibility of the disk mapping]
  This is a continuation of \cref{ex:pr_compatibility_disc_mapping}, which establishes the condition on the logical $0$-form space.
  For the disk mapping as logical mapping $\logicalmapping$, the Jacobian is $\logical{\jac}=2\pi R_1^2\logical{x}$.
  Hence, the pullback of a monomial $2$-form $\intermediate{T}\in\intermediate{\funspace{T}}^{\pr-2}$ of total degree $a+b\leq\pr-2$ is
  \begin{equation}
    \intermediate{T} = \intermediatecoord{x}^a \intermediatecoord{y}^b
    \quad \implies \quad
    \logical{\pushforward}^{-2}\intermediate{T}
    =\logical{\jac}\,(\intermediate{T}\circ\logicalmapping)
    =2\pi R_1^2\logical{x}(R_1\logical{x})^{a+b}
      \cos^a(2\pi\logical{y})\sin^b(2\pi\logical{y}).
  \end{equation}
  Therefore, for $\pr\geq2$, the logical $2$-form space must contain polynomials of degree $\pr-1$ in $\logical{x}$ and the first $2(\pr-2)+1$ trigonometric functions in $\logical{y}$.
  Together with the $0$-form requirements in \cref{ex:pr_compatibility_disc_mapping}, this establishes $\pr$-compatibility of the disk mapping and logical spaces, for every intermediate mapping $\intermediatemapping$.
\end{example}

The assumption of $\pr$-compatibility is sufficient\footnote{In our numerical experiments we find that the projected polynomials always form a basis, regardless of $\pr$-compatibility.} to guarantee that the families in \cref{eqn:polar_forms_basis_expansion} are linearly independent and hence form bases of the polar spaces.
\begin{lemma}[Polar bases under $\pr$-compatibility]\label{lem:polar_basis_pr_compatible}
  Suppose that the logical mapping and logical spaces are $\pr$-compatible in the sense of \cref{def:pr_compatibility_mapping_logical_spaces}.
  Then the families in \cref{eqn:polar_forms_basis_expansion} are linearly independent and hence form bases of the polar spaces.
\end{lemma}
\begin{statementproof}[appendix]{lem:polar_basis_pr_compatible}
  We first assume $\intermediatemapping=\mathrm{id}$, and comment on a general intermediate mapping at the end.
  Then $\mapping=\logicalmapping$, $\jac=\logical{\jac}$, and $\pulledpolys{0}{\pr}$ and $\pulledpolys{2}{\pr-2}$ consist of the polynomials of total degree at most $\pr$ and $\pr-2$ in $\+x$, with $T^0_l=T^2_l=\intermediate{T}_l$.
  Compatibility gives $\pi^0P=P$ for $P\in\pulledpolys{0}{\pr}$ and $\reduced{\pi}^{2,1}B=B$ for $B\in\pulledpolys{2}{\pr-2}$: the latter is a bounded, hence conforming, tensor-product $2$-form, so its first coefficient ring vanishes by \cref{lem:conformity_conditions_tensor_product_forms}. The polynomial bases and the indexed external tensor-product bases are each linearly independent. We must also verify that their spans intersect only at zero; this does not follow from independence within each family.

  For this purpose, observe that every polynomial $p$ of total degree at most $d$ satisfying $p\circ\mapping=\mathcal{O}(\logical{x}^{d+1})$ uniformly in $\logical{y}$ is zero.
  Indeed, by \cref{eqn:coordinate_transformation_assumptions}, $R/\logical{x}$ is continuous with value $R_1>0$ at $\logical{x}=0$, and $C^2+S^2>0$ because $CS'-SC'>0$; by periodicity in $\logical{y}$ there is $c>0$ with $\abs{\mapping(\logical{x},\logical{y})}=R\sqrt{C^2+S^2}\ge c\,\logical{x}$ near the center.
  Every $\+x$ near $\+0$ equals $\mapping(\logical{x},\logical{y})$ for some $(\logical{x},\logical{y})$ with $\logical{x}\le\abs{\+x}/c$ (see the proof of \cref{lem:vanishing_physical_derivative_general}), so $p(\+x)=\mathcal{O}(\abs{\+x}^{d+1})$.
  Along each ray $\+x=s\+\omega$, $\abs{\+\omega}=1$, the function $s\mapsto p(s\+\omega)$ is then a polynomial of degree at most $d$ that is $\mathcal{O}(s^{d+1})$, hence zero, and therefore $p=0$.

  If $P\in\pulledpolys{0}{\pr}$ belongs to $\funspace{\Bspline}^{0,\pr}_\external$, its first $\pr+1$ radial coefficient rings vanish. Hence $P\circ\mapping=\mathcal{O}(\logical{x}^{\pr+1})$ and $P=0$. Likewise, for $\pr\ge2$, let $B\in\pulledpolys{2}{\pr-2}$ belong to $\funspace{\Bspline}^{2,\pr}_\external$; then its first $\pr$ radial coefficient rings vanish. Thus $\jac(B\circ\mapping)=\mathcal{O}(\logical{x}^{\pr})$. Since $\jac/\logical{x}$ has a positive limit, uniformly in $\logical{y}$ (\cref{eqn:conformity_pole_metric_bounds}), $B\circ\mapping=\mathcal{O}(\logical{x}^{\pr-1})$, so the polynomial $B$ of degree at most $\pr-2$ is zero. This proves the direct sums in degrees zero and two.

  In degree one, consider $\gradientt P+\curlt^\dagger B+\+A_\external=\+0$, where $P$ is a combination of $T^0_l$ with $l\ge2$, $B\in\pulledpolys{2}{\pr-2}$, and $\+A_\external\in\funspace{\Bspline}^{1,\pr}_\external$. Taking the curl gives $B+\curlt\+A_\external=0$ by \cref{cor:discrete_homotopy_identities}. Since $\curlt\+A_\external$ is external, degree-two directness gives $B=0$. The logical radial component of $\gradientt P=-\+A_\external$ then gives $\partial_{\logical{x}}(P\circ\mapping)=\mathcal{O}(\logical{x}^{\pr})$ by the external index range in \cref{eqn:polar_mform_space_odot}. As $P(\+0)=0$ by \cref{eqn:polynomial_basis_centre_values}, integration gives $P\circ\mapping=\mathcal{O}(\logical{x}^{\pr+1})$. The observation above yields $P=0$, and then $\+A_\external=\+0$. Thus the unmodified $1$-form family is independent as well.

  Finally, each compactly supported basis function differs from its unmodified polynomial or homotopy counterpart by an external form: in degree zero by $\chi^{0,\pr}_\external$, in degree two by $\reduced{\Pi}^{2,\pr}$, and in degree one by the gradient or curl homotopy of those external corrections. These are triangular changes of basis with identity diagonal when the external basis functions are retained. Therefore the compactly supported families in \cref{eqn:polar_forms_basis_expansion} are linearly independent.

  \emph{General intermediate mapping.}
  The families are linearly independent if and only if their tensor-product coefficient vectors are.
  These vectors do not depend on $\intermediatemapping$ (\cref{sec:derham_hodge_2d_basis}), and neither does $\pr$-compatibility, which involves only $\logicalmapping$; hence the result carries over from $\intermediatemapping=\mathrm{id}$.
\end{statementproof}

Under the same assumption, the polar subspaces satisfy the desired regularity conditions at the center.
\begin{theorem}[Regularity]\label{thm:polar_subspaces_regularity}
  Let $1\le\pr\le r^{0,\xx}$ and assume that the logical mapping $\logicalmapping$ and the logical $m$-form spaces $\logical{\funspace{\Bspline}}^m$ are $\pr$-compatible in the sense of \cref{def:pr_compatibility_mapping_logical_spaces}.
  If $\pr<r^{0,\xx}$, assume also that the angular scalar basis is $\funspace{C}^\pr$ and the angular $1$-form basis is $\funspace{C}^{\pr-1}$.
  Then the following regularity at the center is guaranteed:
  \begin{equation}
    \polar{\funspace{\Bspline}}^{0,\pr}\subset\funspace{C}^\pr_0,
    \qquad
    \polar{\funspace{\Bspline}}^{1,\pr}\subset\funspace{C}^{\pr-1}_0,
    \qquad
    \polar{\funspace{\Bspline}}^{2,\pr}\subset\funspace{C}^{\pr-2}_0.
  \end{equation}
  For $\pr=1$, the last inclusion means boundedness on a punctured neighborhood, without a center limit.
  If $\pr=r^{0,\xx}$, all three polar spaces are smooth on a neighborhood of the center.
\end{theorem}
\begin{statementproof}[appendix]{thm:polar_subspaces_regularity}
  The $0$-form conclusion follows from \cref{thm:polar_zeroform_space_regularity}.
  For the $1$- and $2$-forms, we first assume $\intermediatemapping=\mathrm{id}$, and comment on a general intermediate mapping at the end.
  Then $\Omega=\intermediate{\Omega}$, $\mapping=\logicalmapping$, $\jac=\logical{\jac}$, $\pushforward^m=\logical{\pushforward}^m$, and $\pulledpolys{0}{\pr}$ and $\pulledpolys{2}{\pr-2}$ consist of polynomials.

  \paragraph{Internal parts}
  The polynomial $2$-forms in $\pulledpolys{2}{\pr-2}$ are smooth.
  By $\pr$-compatibility \cref{eqn:pr_compatibility_twoform}, every $B\in\pulledpolys{2}{\pr-2}$ is exactly represented in the physical tensor-product $2$-form space.
  Being bounded near the center, it has $b_{1j}=0$ by \cref{lem:conformity_conditions_tensor_product_forms}, so it belongs to the range of $\reduced{\pi}^{2,1}$ by \cref{lem:conforming_polar_twoform_projector}.
  Thus
  \begin{equation}
    \reduced{\pi}^{2,1}\pulledpolys{2}{\pr-2}=\pulledpolys{2}{\pr-2},
  \end{equation}
  where both sides are zero when $\pr=1$.
  Hence the polynomial part of the polar $2$-form space is exactly $\pulledpolys{2}{\pr-2}$ near the center, and is therefore smooth there.

  Likewise, $\pr$-compatibility gives $\pi^0\pulledpolys{0}{\pr}=\pulledpolys{0}{\pr}$, so the gradient-like polynomial part of the $1$-form space consists of gradients of polynomials and is smooth near the center.
  For the curl-like part, let $B$ be a $2$-form that is smooth on a ball around the center contained in $\Omega$.
  On that ball, \cref{eqn:homotopy_operators} gives
  \begin{equation}
    (\curlt^\dagger B)(\+x)
    =\begin{bmatrix}-y\\x\end{bmatrix}\int_0^1 tB(t\+x)\diff t,
  \end{equation}
  since the segment from $\+0$ to $\+x$ stays in the ball.
  Differentiation under the integral sign shows that $\curlt^\dagger B$ is smooth on the ball.
  Applied to $B\in\pulledpolys{2}{\pr-2}$, this shows that the curl-like polynomial part of the $1$-form space is smooth near the center.

  \paragraph{External parts}
  If $\pr=r^{0,\xx}$, the external $0$-, $1$-, and $2$-form basis functions all vanish on the first radial cell by their index ranges in \cref{eqn:polar_zeroform_space_odot,eqn:polar_mform_space_odot}.
  Thus near the center each polar form equals one of the smooth forms of the internal parts above, proving the smoothness claim.

  Suppose $\pr<r^{0,\xx}$.
  For a polar $1$-form, \cref{eqn:gamma_matrix} gives
  \begin{equation}\label{eqn:polar_oneform_space_gamma_matrix}
    \+A = \srpar{\mat{\Gamma}
    \begin{bmatrix}
      1 & 0\\
      0 & \frac{1}{\logical{x}}
    \end{bmatrix}
    \logical{\+A}} \circ \mapping^{-1}.
  \end{equation}
  Its logical radial component contributes through the first column of $\mat{\Gamma}$ with $q=0$, and its angular component through the second column with $q=-1$.
  For external $2$-forms, the corresponding smooth factor is $\logical{x}/\jac$ with $q=-1$.
  The indices and derivative orders to which \cref{lem:vanishing_physical_derivative_general} applies are
  \begin{equation}
    \begin{array}{c|c|c|c}
      \text{external contribution} & q & n & \text{radial indices}\\ \hline
      \text{$1$-form radial} & 0 & \pr-1 & i\ge\pr+1\\
      \text{$1$-form angular} & -1 & \pr-1 & i\ge\pr+2\\
      \text{$2$-form, $\pr\ge2$} & -1 & \pr-2 & i\ge\pr+1
    \end{array}
  \end{equation}
  In each row $i>n-q+1$, so the external contributions have the asserted continuous physical derivatives at the center.
  The assumed angular smoothness makes these conclusions valid across angular knot rays as well.
  When $\pr=1$, the external $2$-forms instead have $q=-1$, $n=0$, and $i\ge2=n-q+1$; \cref{eqn:vanishing_physical_derivative_bound_physical} in the proof of the same lemma, with $k=0$ and $m=i-2\ge0$, gives boundedness near the center.

  \emph{General intermediate mapping.}
  Every polar form is the pushforward by $\intermediate{\pushforward}^m$ of the polar form with the same coefficients constructed for $\intermediatemapping=\mathrm{id}$ (\cref{sec:derham_hodge_2d_basis}), and $\pr$-compatibility involves only $\logicalmapping$.
  As noted at the end of the proof of \cref{lem:conformity_conditions_tensor_product_forms}, $\intermediate{\pushforward}^m$ preserves boundedness and $k$-fold continuous differentiability near the center, and hence the classes $\funspace{C}^k_0$ and smoothness; thus the result carries over from $\intermediatemapping=\mathrm{id}$.
\end{statementproof}

\subsection{Commuting projection operators}\label{sec:derham_hodge_2d_geom_projections}
We discuss the construction of the commuting projection operators $\polar{\pi}^{m,\pr}$ for $m = 0,1,2$ that project into the polar $m$-form spaces defined in \cref{eqn:polar_zeroform_space,eqn:polar_oneform_space,eqn:polar_twoform_space}
\begin{equation}
  \polar{\pi}^{0,\pr} : \funspace{H}_1 \to \polar{\funspace{\Bspline}}^{0,\pr}, \quad
  \polar{\pi}^{1,\pr} : \funspace{H}_{\curlt} \to \polar{\funspace{\Bspline}}^{1,\pr}, \quad
  \polar{\pi}^{2,\pr} : \funspace{L}_2 \to \polar{\funspace{\Bspline}}^{2,\pr},
\end{equation}
such that they commute with the differential operators.
In principle, the approach from \citet{holderied2022} can be followed, which is based on defining the projection operators in terms of polar GDOFs of the polar $m$-form spaces, which are linear combinations of the GDOFs of the full tensor-product $m$-form spaces.
By the inherent commutative properties of the tensor-product GDOFs, this ensures that the projection operators commute with the differential operators.
However, in following such an approach, it is not clear how to define the polar GDOFs in such a way that all three projection operators are sparse and that they can be evaluated efficiently.

Instead, we follow the approach from \citet{yaman2025}, and construct the projection operators as follows:
\begin{equation}\label{eqn:polar_projection_operator_definition}
  \polar{\pi}^{m,\pr} \defeq \polar{\Pi}^{m,\pr} \circ \pi^m,
\end{equation}
where $\pi^m$ is the commuting projection operator on the full tensor-product $m$-form space defined in \cref{eqn:logical_commuting_projectors,eqn:physical_geometric_dofs} and $\polar{\Pi}^{m,\pr}$ is a commuting projection operator that acts on the coefficients of the corresponding $m$-form, and maps from the full tensor-product $m$-form space to the polar $m$-form space
\begin{equation}
  \polar{\Pi}^{m,\pr} : \funspace{\Bspline}^m \to \polar{\funspace{\Bspline}}^{m,\pr}
  .
\end{equation}
The composition with the tensor-product projection $\pi^m$ follows \citet{yaman2025}; the coefficient projections $\polar{\Pi}^{m,\pr}$ are constructed below.

\begin{assumption}[Polar bases]\label{ass:polar_bases}
  The families in \cref{eqn:polar_forms_basis_expansion} are linearly independent, and hence bases of the polar spaces $\polar{\funspace{\Bspline}}^{m,\pr}$.
\end{assumption}
The construction of the projection operators rests on this assumption: the dual basis functions in \cref{eqn:polar_dual_basis} require the internal parts of the polar basis functions to be linearly independent.
By \cref{lem:polar_basis_pr_compatible}, the assumption holds whenever the logical mapping and the logical spaces are $\pr$-compatible.
Without $\pr$-compatibility it is not guaranteed, but it held in all numerical experiments of \cref{sec:results}, and it is checked in practice by verifying that the extraction matrices of \cref{eqn:polar_extraction_matrices} have full column rank.

\subsubsection{Construction of the projection operators}
With the basis described in \cref{sec:derham_hodge_2d_construction}, we define biorthogonal (and internally supported) dual basis functions $\polar{\Bspline}^{0,\pr,\perp}_{l} \in \funspace{\Bspline}^{0,\pr}_\internal$ and $\polar{\Bspline}^{2,\pr,\perp}_{l} \in \funspace{\Bspline}^{2,\pr}_\internal$ which are expressed in terms of coefficients as
\begin{equation}\label{eqn:polar_dual_basis}
  \polar{\Bspline}^{m,\pr,\perp}_{l} = \sum_{(i,j)\in\indextwo^{m,\pr}_\internal} \polar{\bspline}^{m,\pr,\perp}_{l,ij} \Bspline^m_{ij}
  , \quad 
  l \in \{1, \hdots, n^{m,\pr}\}
  , \quad
  m \in \{0, 2\}
  ,
\end{equation} 
and satisfy the Euclidean biorthogonality conditions
\begin{equation}\label{eqn:polar_dual_basis_biorthogonality}
  \ellinprod{\polar{\Bspline}^{m,\pr,\perp}_{l}}{\polar{\Bspline}^{m,\pr}_{l'}} = \delta_{ll'}
  , \quad 
  l,l' \in \{1, \hdots, n^{m,\pr}\}
  , \quad
  m \in \{0, 2\}
  .
\end{equation}
The Euclidean inner product acting on the coefficients is defined for two $m$-forms $F, F'$ (of equal degree $m$, with $m \in \{0, 2\}$), as
\begin{equation}\label{eqn:coefficient_euclidean_inner_product}
  \ellinprod{F}{F'} \defeq \sum_{(i,j)\in\indextwo^{m}} f_{ij} f'_{ij}.
\end{equation}
\begin{remark}[Existence of internal duals]\label{rem:internal_duals}
  Under \cref{ass:polar_bases}, for which \cref{lem:polar_basis_pr_compatible} gives sufficient conditions, the internally supported duals in \cref{eqn:polar_dual_basis} exist.
  Indeed, the internal parts $\chi^{m,\pr}_\internal\polar{\Bspline}^{m,\pr}_l$ are linearly independent, since a combination of polar basis functions with vanishing internal part is external and hence zero by the direct sum in \cref{eqn:polar_forms_basis_expansion}.
  Hence the basis can be redefined, without changing the polar spaces, by orthonormalization on the internal indices with respect to the inner product in \cref{eqn:coefficient_euclidean_inner_product}.
  For such a basis, the duals are the basis functions themselves, up to truncation: $\polar{\Bspline}^{m,\pr,\perp}_{l}=\chi^{m,\pr}_\internal\polar{\Bspline}^{m,\pr}_{l}$.
  An explicit construction is given in \cref{eqn:polar_dual_coefficient_matrices}.
\end{remark}
We define the projection operators for the internal part of the tensor-product $0$- and $2$-form spaces, denoted by $\enriched{\Pi}^{m,\pr}$, as follows
\begin{equation}\label{eqn:enriched_zeroform_projection_operator}
  \enriched{\Pi}^{m,\pr} F \defeq \sum_{l=1}^{n^{m,\pr}} \ellinprod{\polar{\Bspline}^{m,\pr,\perp}_{l}}{F} \polar{\Bspline}^{m,\pr}_{l}
  , \quad m \in \{0, 2\}.
\end{equation}
These operators project onto the spans of the polar basis functions in \cref{eqn:polar_zeroform_basis,eqn:polar_twoform_curl_basis}, respectively.
\begin{remark}[Generality of the dual-basis construction]\label{rem:dual_basis_generality}
  The construction through dual bases is not restricted to internally supported duals.
  Every projection of $\funspace{\Bspline}^m$ onto $\polar{\funspace{\Bspline}}^{m,\pr}$ is of the form $F\mapsto\sum_l\ellinprod{\polar{\Bspline}^{m,\pr,\perp}_l}{F}\,\polar{\Bspline}^{m,\pr}_l$ for some family of duals in $\funspace{\Bspline}^m$ satisfying \cref{eqn:polar_dual_basis_biorthogonality}, and the duals determine the kernel of the projection.
  In matrix form (\cref{sec:matrix_representation_polar_spaces}), this reads $\polar{\mat{\Pi}}^m=\mat{E}^m(\mat{F}^m)^\transpose$ with $(\mat{F}^m)^\transpose\mat{E}^m=\mat{I}$, where the columns of $\mat{E}^m$ and $\mat{F}^m$ are the tensor-product coefficients of the polar basis functions and of their duals, respectively.
  Let $\mat{H}^m$ be the matrix of the tensor-product GDOFs evaluated at the tensor-product basis functions (\cref{eqn:two_dimensional_gdof_vandermonde_matrices}), and $\mat{M}^m$ the tensor-product mass matrix, i.e., the matrix of $\funspace{L}_2$ inner products of the tensor-product basis functions (\cref{eqn:tensor_product_mass_matrices}).
  The projections based on polar GDOFs \citep{holderied2022}, which are linear combinations $\mat{Y}$ of the tensor-product GDOFs, are then obtained with the duals $(\mat{F}^m)^\transpose=(\mat{Y}\mat{H}^m\mat{E}^m)^{-1}\mat{Y}\mat{H}^m$, and the $\funspace{L}_2$-orthogonal projection with $(\mat{F}^m)^\transpose=(\mat{E}^{m\transpose}\mat{M}^m\mat{E}^m)^{-1}\mat{E}^{m\transpose}\mat{M}^m$; both are dense in general.
  Restricting the duals to the internal coefficients is what makes the polar projections local and sparse (\cref{lem:polar_oneform_projection_operator_internal_coefficients}).
  The commutation with the differential operators in \cref{thm:full_commuting_polar_projection_operators} does not depend on this choice: its proof uses only that $\polar{\Pi}^{0,\pr}$ and $\polar{\Pi}^{2,\pr}$ are projections onto the polar spaces and the homotopy-based definition \cref{eqn:polar_oneform_projection_operator} of $\polar{\Pi}^{1,\pr}$, so that every family of duals yields commuting projections; only the commutation with the homotopy operators in \cref{thm:homotopy_commuting_polar_projection_operators} requires the center-value dual \cref{eqn:polar_constant_dual_centre_value}.
\end{remark}

For the external $0$-form part, we use coefficient truncation:
\begin{equation}\label{eqn:reduced_polar_projection_operators}
  \reduced{\Pi}^{0,\pr}P\defeq\chi^{0,\pr}_\external P.
\end{equation}
For the external $2$-form part, we use $\reduced{\Pi}^{2,\pr}$ from \cref{eqn:reduced_twoform_coefficient_projector}, with $\reduced{\Pi}^{2,0}=I^2$ when $\pr=0$.
Its radial coefficient accumulation, rather than plain external truncation, gives the conservation properties in \cref{lem:polar_projection_operators_conservation_properties}.

The two projection operators $\enriched{\Pi}^{m,\pr}$ and $\reduced{\Pi}^{m,\pr}$ are combined to define the polar $0$- and $2$-form projection operators as
\begin{equation}\label{eqn:polar_zero_two_form_projection_operator}
  \polar{\Pi}^{m,\pr} \defeq \enriched{\Pi}^{m,\pr} + \reduced{\Pi}^{m,\pr}
  , \quad
  m \in \{0, 2\}
  .
\end{equation}
The two summands annihilate each other's ranges: in degree zero by coefficient support, and in degree two by internal support of the duals and the identity $\reduced{\Pi}^{2,\pr}(I^2-\reduced{\Pi}^{2,\pr})=0$.

We recall the homotopy identity for the polar $1$-form space given by \cref{eqn:homotopy_identity_1form}, which we write as
\begin{equation}
  I^1 = \gradientt I^0 \gradientt^\dagger + \curlt^\dagger I^2 \curlt
  ,
\end{equation}
where $I^m$ denotes the identity on the full tensor-product $m$-form space; this identity is restricted to the polar $1$-form space.
This motivates the following definition of the polar $1$-form projection operator
\begin{equation}\label{eqn:polar_oneform_projection_operator}
  \polar{\Pi}^{1,\pr} \defeq \gradientt \polar{\Pi}^{0,\pr} \gradientt^\dagger + \curlt^\dagger \polar{\Pi}^{2,\pr} \curlt
  .
\end{equation}
The matrix representations of these projectors and of the resulting projection operators \cref{eqn:polar_projection_operator_definition} are given in \cref{eqn:polar_projector_matrices,eqn:polar_gdof_projector_matrix_composition}.
Defining the $1$-form projection through the homotopy identity, from the $0$- and $2$-form projections, appears to be new.
Closely related representations of commuting projections in terms of Poincar\'e maps are known on a single polynomial element, where commuting projections with $p$-independent bounds are built from Poincar\'e-type right inverses \citep[eqs.~(3.3)--(3.5)]{gopalakrishnan2004}; in projection-based interpolation, similarly, the edge contributions of a $1$-form are obtained by projecting its one-dimensional potential along each edge \citep[Sec.~5.2]{demkowicz2005}, and analogous decompositions with discrete potential operators appear in the local cochain projections of \citet[eq.~(3.3)]{falk2014} and \citet[eq.~(4.7)]{falk2015}.
For tensor-product spline spaces, commuting projections onto the $1$- and $2$-form spaces are likewise obtained by applying the $0$-form projection to antiderivatives of the input and differentiating the result \citep[eq.~(2.6)]{PintoSchnack2025}, following \citet{buffa_rivas_sangalli_vazquez2011}, and corrections of the same form extend them to multipatch spaces with non-matching interfaces \citep[eqs.~(2.10)--(2.11)]{PintoSchnack2025}.

The correction made by the polar $1$-form projection depends only on internal input coefficients, as the following lemma shows.
\begin{lemma}[Polar $1$-form projection fixes external inputs]\label{lem:polar_oneform_projection_operator_internal_coefficients}
  The projector in \cref{eqn:polar_oneform_projection_operator} satisfies
  \begin{equation}\label{eqn:polar_oneform_projection_operator_internal_coefficients}
    \polar{\Pi}^{1,\pr}=I^1+(\polar{\Pi}^{1,\pr}-I^1)\chi^{1,\pr}_\internal,
  \end{equation}
  where $\chi^{1,\pr}_\internal$ denotes the truncation onto the internal coefficients.
\end{lemma}
\begin{statementproof}[appendix]{lem:polar_oneform_projection_operator_internal_coefficients}
  We first show that $\polar{\Pi}^{1,\pr}$ fixes every external input.  Let
  $\+A_\external\in\funspace{\Bspline}^{1,\pr}_\external$.  Its radial
  coefficients satisfy (cf.\ \cref{eqn:polar_mform_space_odot})
  \begin{equation}
    a_{\xx,ij}=0\quad\text{for }i\leq\pr,
    \quad
    a_{\yy,ij}=0\quad\text{for }i\leq\pr+1.
  \end{equation}
  By \cref{eqn:pseudo_inverse_gradient_logical}, the coefficient in $0$-form
  ring $i$ depends only on radial $1$-form coefficients with indices $i'<i$.
  Hence the output coefficients of $\gradientt^\dagger$ vanish for $i\le\pr+1$, and
  $\gradientt^\dagger\+A_\external\in\funspace{\Bspline}^{0,\pr}_\external$.
  The internally supported duals defining $\enriched{\Pi}^{0,\pr}$ pair to
  zero with these external coefficients in the coefficient inner product.
  Since $\reduced{\Pi}^{0,\pr}=\chi^{0,\pr}_\external$, it follows that
  \begin{equation}\label{eqn:zeroform_projection_identity_on_external_homotopy}
    \polar{\Pi}^{0,\pr}\gradientt^\dagger\+A_\external
    =\gradientt^\dagger\+A_\external.
  \end{equation}

  For the curl coefficients, the discrete curl in \cref{eqn:discrete_operators} gives
  \begin{equation}
    b_{ij}=a_{\yy,i+1,j}-a_{\yy,ij}-(\delta_\yy a_\xx)_{ij}=0
    \quad (i\le\pr),
  \end{equation}
  by the vanishing input coefficients above.
  Thus $\curlt\+A_\external\in\funspace{\Bspline}^{2,\pr}_\external$.
  If $B\in\funspace{\Bspline}^{2,\pr}_\external$, then $b_{ij}=0$ for
  $i\leq\pr$.  Consequently, the definition of
  $\reduced{\Pi}^{2,\pr}$ in
  \cref{eqn:reduced_twoform_coefficient_projector} reduces to
  \begin{equation}
    \reduced{\Pi}^{2,\pr}B
    =\sum_{j=1}^{n^\yy}b_{(\pr+1)j}\Bspline^2_{(\pr+1)j}
     +\sum_{i=\pr+2,j=1}^{n^{1,\xx},n^\yy}b_{ij}\Bspline^2_{ij}
    =B.
  \end{equation}
  The internally supported dual basis also gives
  $\enriched{\Pi}^{2,\pr}B=0$ by coefficient support.  Therefore,
  \begin{equation}\label{eqn:twoform_projection_identity_on_external_curl}
    \polar{\Pi}^{2,\pr}\curlt\+A_\external
    =\curlt\+A_\external.
  \end{equation}

  Since $\+A_\external$ satisfies $a_{\yy,1j}=0$, the discrete $1$-form
  homotopy identity in \cref{cor:discrete_homotopy_identities}, together
  with
  \cref{eqn:zeroform_projection_identity_on_external_homotopy,eqn:twoform_projection_identity_on_external_curl},
  now yields
  \begin{align}
    \polar{\Pi}^{1,\pr}\+A_\external
    &=\gradientt\polar{\Pi}^{0,\pr}\gradientt^\dagger\+A_\external
      +\curlt^\dagger\polar{\Pi}^{2,\pr}\curlt\+A_\external\notag\\
    &=\gradientt\gradientt^\dagger\+A_\external
      +\curlt^\dagger\curlt\+A_\external
     =\+A_\external.
  \end{align}
  Since $\polar{\Pi}^{1,\pr}$ fixes every external input, inserting
  $I^1=\chi^{1,\pr}_\internal+\chi^{1,\pr}_\external$ on its right gives
  $\polar{\Pi}^{1,\pr}=\polar{\Pi}^{1,\pr}\chi^{1,\pr}_\internal+\chi^{1,\pr}_\external$.
  Rearranging yields \cref{eqn:polar_oneform_projection_operator_internal_coefficients}.
\end{statementproof}

\subsubsection{Properties of the projection operators}

\begin{theorem}[Commuting polar projection operators]\label{thm:full_commuting_polar_projection_operators}
  For $m\in\{0,1,2\}$, the operator $\polar{\Pi}^{m,\pr}$ defined in \cref{eqn:polar_zero_two_form_projection_operator,eqn:polar_oneform_projection_operator} maps the full tensor-product coefficient space $\funspace{\Bspline}^m$ onto $\polar{\funspace{\Bspline}}^{m,\pr}$ and is the identity on that polar space.
  These operators commute coefficient-wise with the differential operators in the following sense:
  \begin{equation}
    \gradientt \polar{\Pi}^{0,\pr} P = \polar{\Pi}^{1,\pr} \gradientt P, \quad
    \curlt \polar{\Pi}^{1,\pr} \+A = \polar{\Pi}^{2,\pr} \curlt \+A,
  \end{equation}
  for $P\in\funspace{\Bspline}^0$ continuous at the center ($p_{1j}=P(\+0)\mathds{1}_j^\yy$) and every $\+A\in\funspace{\Bspline}^1$.
  At $\pr=0$, these are algebraic coefficient statements; conformity for $\pr\ge1$ is established in \cref{thm:polar_subspaces_conformity}.
\end{theorem}
\begin{statementproof}[appendix]{thm:full_commuting_polar_projection_operators}
  We first prove that each of the operators is a projection operator that maps onto the corresponding polar $m$-form space, and then show that they commute with the differential operators.

  \paragraph{Projection property -- polar $0$-form space}
  Biorthogonality in \cref{eqn:polar_dual_basis_biorthogonality} shows that
  $\enriched{\Pi}^{0,\pr}$ is a projection onto the span of the internal
  polar basis functions.  Moreover, $\reduced{\Pi}^{0,\pr}$ is the
  truncation onto $\funspace{\Bspline}^{0,\pr}_\external$.  The two
  operators annihilate each other's ranges because both the primal and the
  dual basis functions are supported on internal coefficient indices.  Hence,
  $\polar{\Pi}^{0,\pr}$ maps onto the direct sum in
  \cref{eqn:polar_zeroform_space,eqn:polar_zeroform_basis_expansion} and is
  the identity on this space.

  We record one further consequence that will be used below.
  By \cref{ass:constant_functions}, $\pi^0 T^0_1=1$, so
  $\polar{\Bspline}^{0,\pr}_1$ is the internal truncation of the constant
  $0$-form.  Internal coefficient support and biorthogonality give
  \begin{equation}\label{eqn:polar_zeroform_projection_constant}
    \enriched{\Pi}^{0,\pr}1=\polar{\Bspline}^{0,\pr}_1 = \chi^{0,\pr}_\internal 1,
    \quad
    \reduced{\Pi}^{0,\pr}1=\chi^{0,\pr}_\external 1,
    \quad\implies\quad
    \polar{\Pi}^{0,\pr}1=1.
  \end{equation}

  \paragraph{Projection property -- polar $2$-form space}
  The coefficient formula in \cref{eqn:reduced_twoform_coefficient_projector} shows that $\reduced{\Pi}^{2,\pr}$ projects onto $\funspace{\Bspline}^{2,\pr}_\external$. 
  Biorthogonality likewise shows that $\enriched{\Pi}^{2,\pr}$ is a projection onto the span of the polar basis functions $\polar{\Bspline}^{2,\pr}_l$.  
  Since the dual basis is internally supported whereas the range of $\reduced{\Pi}^{2,\pr}$ is externally supported, we have
  \begin{equation}\label{eqn:internal_external_twoform_projection}
    \enriched{\Pi}^{2,\pr}\reduced{\Pi}^{2,\pr}=0.
  \end{equation}
  The reverse cancellation follows directly from the compactly supported basis definition \cref{eqn:polar_twoform_curl_basis}.  Since $\reduced{\Pi}^{2,\pr}$ is a projector by \cref{lem:conforming_polar_twoform_projector}, using \cref{eqn:polar_twoform_curl_basis} we find
  \begin{equation}
    \reduced{\Pi}^{2,\pr}\polar{\Bspline}^{2,\pr}_l
    =\reduced{\Pi}^{2,\pr}(I^2-\reduced{\Pi}^{2,\pr})
      \reduced{\pi}^{2,1}T^2_l=0.
  \end{equation}
  For $\pr=0,1$ the enriched $2$-form sum is empty, so the same cancellation is immediate.
  Thus $\reduced{\Pi}^{2,\pr}\enriched{\Pi}^{2,\pr}=0$.
  Consequently, if
  \begin{equation}
    B=B_\internal+B_\external
    \in\polar{\funspace{\Bspline}}^{2,\pr},
  \end{equation}
  where $B_\internal$ is a linear combination of the polar basis functions
  and $B_\external\in\funspace{\Bspline}^{2,\pr}_\external$, then
  \begin{equation}
    \polar{\Pi}^{2,\pr}B
    =B_\internal+\reduced{\Pi}^{2,\pr}B_\external
    =B.
  \end{equation}
  The definition in \cref{eqn:polar_zero_two_form_projection_operator} also
  shows that the range is contained in the sum of precisely these two
  spaces.  By \cref{eqn:polar_twoform_space,eqn:polar_twoform_curl_basis},
  this sum is the polar $2$-form space, which proves the claim.

  \paragraph{Projection property -- polar $1$-form space}
  The result in \cref{thm:polar_subspaces_homotopy_splitting} shows that $\polar{\Pi}^{1,\pr}$ has range contained in
  \begin{equation}
    \gradientt\polar{\funspace{\Bspline}}^{0,\pr}
    \oplus
    \curlt^\dagger\polar{\funspace{\Bspline}}^{2,\pr}
    =\polar{\funspace{\Bspline}}^{1,\pr}.
  \end{equation}
  Conversely, if
  $\+A\in\polar{\funspace{\Bspline}}^{1,\pr}$, then
  $\gradientt^\dagger\+A\in\polar{\funspace{\Bspline}}^{0,\pr}$ and
  $\curlt\+A\in\polar{\funspace{\Bspline}}^{2,\pr}$ by
  \cref{eqn:polar_subspaces_homotopy_preservation,thm:polar_subspaces_commutativity_conformity}.
  The projection properties already established in degrees zero and two,
  followed by the $1$-form homotopy identity, therefore give
  \begin{equation}
    \polar{\Pi}^{1,\pr}\+A
    =\gradientt\gradientt^\dagger\+A
    +\curlt^\dagger\curlt\+A
    =\+A.
  \end{equation}
  Thus $\polar{\Pi}^{1,\pr}$ is a projection onto the polar $1$-form space.

  \paragraph{Commutativity with the differential operators}
  Let $P$ be continuous.  Using
  $\curlt\gradientt=0$, the $0$-form homotopy identity from
  \cref{cor:discrete_homotopy_identities}, and
  \cref{eqn:polar_zeroform_projection_constant}, we obtain
  \begin{align}
    \polar{\Pi}^{1,\pr}\gradientt P
    &=\gradientt\polar{\Pi}^{0,\pr}
      \gradientt^\dagger\gradientt P
      +\curlt^\dagger\polar{\Pi}^{2,\pr}\curlt\gradientt P\notag\\
    &=\gradientt\polar{\Pi}^{0,\pr}(P-\mathfrak{P} P)
     =\gradientt\polar{\Pi}^{0,\pr}P.
  \end{align}
  Similarly, using $\curlt\gradientt=0$ and the $2$-form homotopy identity,
  \begin{align}
    \curlt\polar{\Pi}^{1,\pr}\+A
    &=\curlt\curlt^\dagger\polar{\Pi}^{2,\pr}\curlt\+A
     =\polar{\Pi}^{2,\pr}\curlt\+A.
  \end{align}

\end{statementproof}

The biorthogonality conditions in \cref{eqn:polar_dual_basis_biorthogonality} are sufficient for the projection properties and commutativity with the differential operators.  
To ensure commutativity with the gradient homotopy, we additionally require the dual of the constant $0$-form basis function to reproduce values at the center:
\begin{equation}\label{eqn:polar_constant_dual_centre_value}
  \ellinprod{\polar{\Bspline}^{0,\pr,\perp}_{1}}{P}
  =\mathfrak{P}P
  \qquad
  \text{for every $P\in\funspace{\Bspline}^0$ continuous at the center}.
\end{equation}
\begin{remark}[A dual reproducing center values]
  One explicit choice for the first member of a complete internally supported dual family is
  \begin{equation}\label{eqn:polar_constant_dual_basis}
    \polar{\Bspline}^{0,\pr,\perp}_{1}
    \defeq
    \frac{\sum_{j=1}^{n^\yy}\mathds{1}_j^\yy\Bspline^0_{1j}}
         {\sum_{j=1}^{n^\yy}(\mathds{1}_j^\yy)^2}.
  \end{equation}
  The first radial coefficient ring of $\polar{\Bspline}^{0,\pr}_{l}$ is $T^0_l(\+0)\mathds{1}_j^\yy=\intermediate{T}_l(\+0)\mathds{1}_j^\yy$. Hence \cref{eqn:polynomial_basis_centre_values} shows that this choice is biorthogonal to every internal polar $0$-form basis function. The remaining dual functions are chosen to satisfy \cref{eqn:polar_dual_basis_biorthogonality} simultaneously.
  For $P\in\funspace{\Bspline}^0$ continuous at the center, $p_{1j}=P(\+0)\mathds{1}_j^\yy$ and thus
  \begin{equation}
    \ellinprod{\polar{\Bspline}^{0,\pr,\perp}_{1}}{P}
    =
    \frac{\sum_{j=1}^{n^\yy}\mathds{1}_j^\yy p_{1j}}
         {\sum_{j=1}^{n^\yy}(\mathds{1}_j^\yy)^2}
    =
    \mathfrak{P} P,
  \end{equation}
  so this choice satisfies \cref{eqn:polar_constant_dual_centre_value}.
\end{remark}

\begin{theorem}[Homotopy-commuting polar projection operators]\label{thm:homotopy_commuting_polar_projection_operators}
  Suppose that the first polar $0$-form dual satisfies \cref{eqn:polar_constant_dual_centre_value}.
  Then, for every $\+A\in\funspace{\Bspline}^1$ and $B\in\funspace{\Bspline}^2$,
  \begin{equation}
    \gradientt^\dagger\polar{\Pi}^{1,\pr}\+A
    =\polar{\Pi}^{0,\pr}\gradientt^\dagger\+A,
    \qquad
    \curlt^\dagger\polar{\Pi}^{2,\pr}B
    =\polar{\Pi}^{1,\pr}\curlt^\dagger B.
  \end{equation}
\end{theorem}
\begin{statementproof}[appendix]{thm:homotopy_commuting_polar_projection_operators}
  The curl-homotopy identity needs no condition on the first $0$-form dual. By $\gradientt^\dagger\curlt^\dagger=0$ and the $2$-form homotopy identity,
  \begin{equation}
    \polar{\Pi}^{1,\pr}\curlt^\dagger B
    =\gradientt\polar{\Pi}^{0,\pr}
      \gradientt^\dagger\curlt^\dagger B
      +\curlt^\dagger\polar{\Pi}^{2,\pr}
       \curlt\curlt^\dagger B
    =\curlt^\dagger\polar{\Pi}^{2,\pr}B.
  \end{equation}

  Now assume \cref{eqn:polar_constant_dual_centre_value}.
  Let $P\in\funspace{\Bspline}^0$ be continuous at the center, and let $p_1$ denote the first coefficient of the polar $0$-form $\polar{\Pi}^{0,\pr}P$ in \cref{eqn:polar_zeroform_basis_expansion}. Then
  \begin{equation}\label{eqn:polar_zeroform_projection_centre_value}
    \mathfrak{P}(\polar{\Pi}^{0,\pr}P)
    =
    p_1
    =
    \ellinprod{\polar{\Bspline}^{0,\pr,\perp}_1}{P}
    =
    \mathfrak{P}P.
  \end{equation}
  The first equality is \cref{eqn:polar_evaluation_at_centre}: the basis function $\polar{\Bspline}^{0,\pr}_{l}$ has the first coefficient ring $T^0_l(\+0)\mathds{1}_j^\yy$ and therefore the center value $T^0_l(\+0)=\intermediate{T}_l(\+0)=\delta_{l1}$ by \cref{eqn:polynomial_basis_centre_values}, while the external basis functions vanish at the center.
  The second equality follows from the definition of $\polar{\Pi}^{0,\pr}$ in \cref{eqn:enriched_zeroform_projection_operator,eqn:reduced_polar_projection_operators,eqn:polar_zero_two_form_projection_operator}, and the third is \cref{eqn:polar_constant_dual_centre_value}.
  The coefficient formula in \cref{eqn:pseudo_inverse_gradient_logical}
  shows that $\gradientt^\dagger\+A$ is continuous and vanishes at the
  center for every $\+A\in\funspace{\Bspline}^1$. Hence
  \cref{eqn:polar_zeroform_projection_centre_value} shows that
  $\polar{\Pi}^{0,\pr}\gradientt^\dagger\+A$ also vanishes at the center.
  Using the $0$-form homotopy identity and
  $\gradientt^\dagger\curlt^\dagger=0$, we therefore find
  \begin{equation}
    \gradientt^\dagger\polar{\Pi}^{1,\pr}\+A
    =\gradientt^\dagger\gradientt\polar{\Pi}^{0,\pr}
      \gradientt^\dagger\+A
      +\gradientt^\dagger\curlt^\dagger
       \polar{\Pi}^{2,\pr}\curlt\+A
    =\polar{\Pi}^{0,\pr}\gradientt^\dagger\+A.
  \end{equation}
\end{statementproof}

\begin{corollary}[Commuting polar projections on admissible fields]
  \label{cor:infinite_dimensional_polar_projection_operators}
  Let $\pr\ge1$ and define $\polar{\pi}^{m,\pr}$ by \cref{eqn:polar_projection_operator_definition}.
  As in \cref{sec:logical_bspline_spaces}, the operators are applied only to fields whose GDOFs are defined.
  Each operator is a projection onto its polar space.
  For $P\in\funspace{H}_1$ and $\+A\in\funspace{H}_{\curlt}$, they commute with the differential operators:
  \begin{equation}\label{eqn:infinite_dimensional_polar_commutativity}
    \polar{\pi}^{1,\pr}\gradientt P
    =\gradientt\polar{\pi}^{0,\pr}P,
    \quad
    \polar{\pi}^{2,\pr}\curlt\+A
    =\curlt\polar{\pi}^{1,\pr}\+A.
  \end{equation}
  If, in addition, the first polar $0$-form dual satisfies \cref{eqn:polar_constant_dual_centre_value}, then for $\+A\in\funspace{H}_{\curlt}$ and $B\in\funspace{L}_2$,
  \begin{equation}\label{eqn:infinite_dimensional_polar_homotopy_commutativity}
    \gradientt^\dagger\polar{\pi}^{1,\pr}\+A
    =\polar{\pi}^{0,\pr}\gradientt^\dagger\+A,
    \quad
    \curlt^\dagger\polar{\pi}^{2,\pr}B
    =\polar{\pi}^{1,\pr}\curlt^\dagger B.
  \end{equation}
\end{corollary}
\begin{statementproof}[appendix]{cor:infinite_dimensional_polar_projection_operators}
  The proof transfers the coefficient results in \cref{thm:full_commuting_polar_projection_operators,thm:homotopy_commuting_polar_projection_operators} to fields through $\pi^m$.
  For this, $\pi^m$ must commute with the differential and homotopy operators.
  For the differential operators this is \cref{eqn:physical_commutativity}.
  For the homotopy operators, the key observation is that the GDOFs of $\gradientt^\dagger\+A$ and $\curlt^\dagger B$ are radial cumulative sums of the GDOFs of $\+A$ and $B$, see \cref{eqn:gdof_gradient_homotopy,eqn:gdof_curl_homotopy} below.
  This holds because the first clamped Greville node is $\xi^\xx_1=0$, so the logical radial integrals from the center in \cref{lem:pseudo_inverse_curl_logical} split into integrals over Greville edges.

  \emph{Projection property.}
  The polar spaces are subspaces of the tensor-product spaces, so $\pi^m$ acts as the identity on the range of $\polar{\Pi}^{m,\pr}$.
  Together with the idempotence of the coefficient projectors, this gives
  \begin{equation}
    (\polar{\pi}^{m,\pr})^2
    =\polar{\Pi}^{m,\pr}\pi^m
      \polar{\Pi}^{m,\pr}\pi^m
    =(\polar{\Pi}^{m,\pr})^2\pi^m
    =\polar{\pi}^{m,\pr}.
  \end{equation}

  \emph{Commutation with the differential operators.}
  \Cref{thm:full_commuting_polar_projection_operators} applies to $\pi^0P$ because its first coefficient ring is constant.
  Indeed, the GDOFs at the first radial Greville node are the values $P(\mapping(0,\xi^\yy_j))=P(\+0)$, so the argument in the proof of \cref{thm:polar_zeroform_space_conforming} gives $(\pi^0P)_{1j}=P(\+0)\mathds{1}_j^\yy$.
  The $1$-form part of that theorem holds for every tensor-product $1$-form, so no condition at the center is needed for $\pi^1\+A$.
  The commutation relations in \cref{eqn:physical_commutativity} then yield
  \begin{align}
    \polar{\pi}^{1,\pr}\gradientt P
    &=\polar{\Pi}^{1,\pr}\pi^1\gradientt P
     =\polar{\Pi}^{1,\pr}\gradientt\pi^0P
     =\gradientt\polar{\pi}^{0,\pr}P,\\
    \polar{\pi}^{2,\pr}\curlt\+A
    &=\polar{\Pi}^{2,\pr}\pi^2\curlt\+A
     =\polar{\Pi}^{2,\pr}\curlt\pi^1\+A
     =\curlt\polar{\pi}^{1,\pr}\+A.
  \end{align}

  \emph{GDOFs of the homotopy operators.}
  Let $\logical{A}_\xx$ and $\logical{B}$ denote the radial component of $\pushforward^{-1}\+A$ and the pullback $\pushforward^{-2}B$.
  The pullbacks of the homotopy operators are the radial integrals in \cref{lem:pseudo_inverse_curl_logical}; the change of variables in its proof does not use the discrete structure and hence applies to $\+A$ and $B$.
  By \cref{eqn:physical_geometric_dofs,eqn:logical_geometric_dofs}, the gradient homotopy therefore has the GDOFs
  \begin{equation}\label{eqn:gdof_gradient_homotopy}
    \gdof^0_{ij}(\gradientt^\dagger\+A)
    =\int_0^{\xi^\xx_i}\logical{A}_\xx(s,\xi^\yy_j)\diff s
    =\sum_{i'=1}^{i-1}\int_{\xi^\xx_{i'}}^{\xi^\xx_{i'+1}}\logical{A}_\xx(s,\xi^\yy_j)\diff s
    =\sum_{i'=1}^{i-1}\gdof^1_{\xx,i'j}(\+A)
    =\bigl(\delta_\xx^\dagger\gdof^1_\xx(\+A)\bigr)_{ij},
  \end{equation}
  where the second equality splits $[0,\xi^\xx_i]$ at the Greville nodes, using $\xi^\xx_1=0$, and the last one is \cref{eqn:discrete_radial_homotopy_operator}.
  The radial component of $\curlt^\dagger B$ vanishes, so its radial GDOFs vanish, and its angular GDOFs are
  \begin{equation}\label{eqn:gdof_curl_homotopy}
    \gdof^1_{\yy,ij}(\curlt^\dagger B)
    =\int_{\xi^\yy_j}^{\xi^\yy_{j+1}}\int_0^{\xi^\xx_i}\logical{B}(s,\logical{y})\diff s\diff\logical{y}
    =\sum_{i'=1}^{i-1}\int_{\xi^\yy_j}^{\xi^\yy_{j+1}}\int_{\xi^\xx_{i'}}^{\xi^\xx_{i'+1}}\logical{B}(s,\logical{y})\diff s\diff\logical{y}
    =\bigl(\delta_\xx^\dagger\gdof^2(B)\bigr)_{ij}.
  \end{equation}
  Here, the last equality identifies each iterated integral with the face GDOF $\gdof^2_{i'j}(B)$.

  \emph{Commutation of $\pi^m$ with the homotopy operators.}
  By \cref{eqn:pseudo_inverse_gradient_logical}, $\gradientt^\dagger\pi^1\+A$ is a discrete $0$-form.
  Applying \cref{eqn:gdof_gradient_homotopy} to $\pi^1\+A$ and to $\+A$, and using that $\pi^0$ and $\pi^1$ preserve GDOFs, gives
  \begin{equation}
    \gdof^0(\gradientt^\dagger\pi^1\+A)
    =\delta_\xx^\dagger\gdof^1_\xx(\pi^1\+A)
    =\delta_\xx^\dagger\gdof^1_\xx(\+A)
    =\gdof^0(\gradientt^\dagger\+A)
    =\gdof^0(\pi^0\gradientt^\dagger\+A).
  \end{equation}
  Scalar unisolvence gives $\gradientt^\dagger\pi^1\+A=\pi^0\gradientt^\dagger\+A$.
  Likewise, by \cref{eqn:pseudo_inverse_curl_logical}, $\curlt^\dagger\pi^2B$ is a discrete $1$-form, both $\curlt^\dagger\pi^2B$ and $\pi^1\curlt^\dagger B$ have vanishing radial GDOFs, and \cref{eqn:gdof_curl_homotopy} gives
  \begin{equation}
    \gdof^1_\yy(\curlt^\dagger\pi^2B)
    =\delta_\xx^\dagger\gdof^2(\pi^2B)
    =\delta_\xx^\dagger\gdof^2(B)
    =\gdof^1_\yy(\curlt^\dagger B)
    =\gdof^1_\yy(\pi^1\curlt^\dagger B).
  \end{equation}
  Unisolvence gives $\curlt^\dagger\pi^2B=\pi^1\curlt^\dagger B$.

  \emph{Conclusion.}
  Under \cref{eqn:polar_constant_dual_centre_value}, \cref{thm:homotopy_commuting_polar_projection_operators} applies to $\pi^1\+A$ and $\pi^2B$, and we obtain
  \begin{align}
    \gradientt^\dagger\polar{\pi}^{1,\pr}\+A
    &=\gradientt^\dagger\polar{\Pi}^{1,\pr}\pi^1\+A
     =\polar{\Pi}^{0,\pr}\gradientt^\dagger\pi^1\+A
     =\polar{\pi}^{0,\pr}\gradientt^\dagger\+A,\\
    \curlt^\dagger\polar{\pi}^{2,\pr}B
    &=\curlt^\dagger\polar{\Pi}^{2,\pr}\pi^2B
     =\polar{\Pi}^{1,\pr}\curlt^\dagger\pi^2B
     =\polar{\pi}^{1,\pr}\curlt^\dagger B.
  \end{align}
\end{statementproof}

The specific choice of the projection operator $\reduced{\Pi}^{2,\pr}$ in \cref{eqn:reduced_twoform_coefficient_projector} ensures that the polar projection operators preserve the total integral of the $2$-form, as shown in the following lemma.
\begin{lemma}[Conservation properties of the polar projection operators]\label{lem:polar_projection_operators_conservation_properties}
  Let $\pr\ge1$. The $2$-form polar projection operator $\polar{\pi}^{2,\pr}$ defined in \cref{eqn:polar_projection_operator_definition} preserves the total integral of every $2$-form $B\in\funspace{L}_2$,
  \begin{equation}
    \int_\Omega \polar{\pi}^{2,\pr} B \diff^2 x = \int_\Omega B \diff^2 x,
  \end{equation}
  whereas the $1$-form polar projection operator $\polar{\pi}^{1,\pr}$ preserves the total circulation of every $1$-form $\+A\in\funspace{H}_{\curlt}$ that is bounded at the center,
  \begin{equation}
    \int_{\partial\Omega} (\polar{\pi}^{1,\pr} \+A) \cdot \+t \diff l = \int_{\partial\Omega} \+A \cdot \+t \diff l.
  \end{equation}
  The boundary integrals are the sums of classical Greville edge integrals in the positive boundary orientation.
\end{lemma}
\begin{statementproof}[appendix]{lem:polar_projection_operators_conservation_properties}
  We first prove the result for the coefficient projection operator
  $\polar{\Pi}^{2,\pr}$.  We show that it preserves the radial sum of the
  coefficients for every angular index.  Denote the coefficients of the
  enriched and full projections by
  \begin{equation}
    \enriched{\Pi}^{2,\pr}B
    =\sum_{(i,j)\in\indextwo^2}\enriched{b}_{ij}\Bspline^2_{ij},
    \quad
    \polar{\Pi}^{2,\pr}B
    =\sum_{(i,j)\in\indextwo^2}\polar{b}_{ij}\Bspline^2_{ij}.
  \end{equation}
  For any $2$-form $F$ with coefficients $f_{ij}$, the definition in
  \cref{eqn:reduced_twoform_coefficient_projector} gives
  \begin{equation}
    \sum_i\bigl(\reduced{\Pi}^{2,\pr}F\bigr)_{ij}
    =\sum_{i=1}^{\pr+1}f_{ij}
     +\sum_{i=\pr+2}^{n^{1,\xx}}f_{ij}
    =\sum_i f_{ij}.
  \end{equation}
  The enriched image is a linear combination of the compactly supported $2$-form basis functions in \cref{eqn:polar_twoform_curl_basis}. Their coefficient formula \eqref{eqn:polar_twoform_compact_coefficients} gives, also for their linear combinations,
  \begin{equation}
    \enriched{b}_{(\pr+1)j}
      =-\sum_{i=1}^{\pr}\enriched{b}_{ij},
    \qquad
    \enriched{b}_{ij}=0\quad(i\ge\pr+2),
    \qquad\text{hence}\qquad
    \sum_i\enriched{b}_{ij}=0.
  \end{equation}
  Using \cref{eqn:polar_zero_two_form_projection_operator}, the reduced-projector sum identity above now yields
  \begin{equation}
    \sum_i\polar{b}_{ij}
    =\sum_i\enriched{b}_{ij}
     +\sum_i\bigl(\reduced{\Pi}^{2,\pr}B\bigr)_{ij}
    =\sum_i b_{ij}.
  \end{equation}

  It remains to relate this radial sum to the integral.  The integral of a
  physical $2$-form is invariant under pullback, while the normalization in
  \cref{eqn:bspline_derivative} and the clamped endpoint values imply
  \begin{equation}
    \int_0^1\Bsplineone_i^{1,\xx}\diff\logical{x}=1.
  \end{equation}
  Therefore, using the radial-sum identity above, we obtain
  \begin{equation}
    \int_\Omega\polar{\Pi}^{2,\pr}B\diff^2x
    = \sum_{(i,j)\in\indextwo^2} \polar{b}_{ij} \int_0^1\Bsplineone_j^{1,\yy}\diff\logical{y} 
    = \sum_{(i,j)\in\indextwo^2} b_{ij} \int_0^1\Bsplineone_j^{1,\yy}\diff\logical{y}
    = \int_{\logical{\Omega}}\logical{B} \diff\logical{x}\diff\logical{y}
    = \int_\Omega B\diff^2x.
  \end{equation}

  Let $\+A\in\funspace{\Bspline}^1$ be bounded at the center. Boundedness gives $\logical{A}_\yy(0,\logical{y})=0$ by the pole expansion in \cref{lem:conformity_conditions_tensor_product_forms}. Periodicity gives $\logical{A}_\xx(\logical{x},1)=\logical{A}_\xx(\logical{x},0)$. Thus the fundamental theorem of calculus in logical coordinates yields
  \begin{equation}
    \int_\Omega\curlt\+A\diff^2x
    =\int_0^1\int_0^1
      \bigl(\partial_{\logical{x}}\logical{A}_\yy
            -\partial_{\logical{y}}\logical{A}_\xx\bigr)
      \diff\logical{x}\diff\logical{y}
    =\int_0^1 \logical{A}_\yy(1,\logical{y})\diff\logical{y}
    =\int_{\partial\Omega}\+A\cdot\+t\diff l.
  \end{equation}
  The same identity holds for $\polar{\Pi}^{1,\pr}\+A$, which is bounded by \cref{thm:polar_subspaces_conformity,lem:conformity_conditions_tensor_product_forms}. Combining these identities with \cref{thm:full_commuting_polar_projection_operators} and the integral conservation above gives
  \begin{align}
    \int_{\partial\Omega}(\polar{\Pi}^{1,\pr}\+A)\cdot\+t\diff l
    &=\int_\Omega\curlt\polar{\Pi}^{1,\pr}\+A\diff^2x\notag\\
    &=\int_\Omega\polar{\Pi}^{2,\pr}\curlt\+A\diff^2x
     =\int_\Omega\curlt\+A\diff^2x
     =\int_{\partial\Omega}\+A\cdot\+t\diff l,
  \end{align}
  which proves conservation of the total circulation for
  $\polar{\Pi}^{1,\pr}$.

  It remains to pass from the coefficient projection operators to the polar
  projection operators in \cref{eqn:polar_projection_operator_definition}.
  The $2$-form GDOFs are integrals over the Greville faces, which partition
  $\Omega$, while the tangential $1$-form GDOFs on the boundary are integrals
  over the positively oriented boundary Greville edges, including the periodic wraparound interval, which partition $\partial\Omega$.
  Since the GDOF projectors preserve each of these integrals by
  \cref{eqn:logical_commuting_projectors,eqn:physical_geometric_dofs},
  they preserve their sums:
  \begin{equation}
    \int_\Omega\pi^2B\diff^2x=\int_\Omega B\diff^2x,
    \quad
    \int_{\partial\Omega}(\pi^1\+A)\cdot\+t\diff l
    =\int_{\partial\Omega}\+A\cdot\+t\diff l.
  \end{equation}
  Moreover, $\pi^1\+A$ is bounded at the center whenever $\+A$ is: the
  angular edge moments of a bounded $\+A$ at the first radial Greville node
  $\xi^\xx_1=0$ vanish, so angular unisolvence gives $a_{\yy,1j}=0$ for
  $\pi^1\+A$, and this first-ring condition makes $\pi^1\+A$ bounded at the
  center by the pole expansion in the proof of
  \cref{lem:conformity_conditions_tensor_product_forms}.  Consequently, the definition of
  $\polar{\pi}^{m,\pr}$ and the two coefficient-level conservation results
  above give
  \begin{align}
    \int_\Omega\polar{\pi}^{2,\pr}B\diff^2x
    &=\int_\Omega\polar{\Pi}^{2,\pr}(\pi^2B)\diff^2x
     =\int_\Omega\pi^2B\diff^2x
     =\int_\Omega B\diff^2x,\\
    \int_{\partial\Omega}(\polar{\pi}^{1,\pr}\+A)\cdot\+t\diff l
    &=\int_{\partial\Omega}
      \bigl(\polar{\Pi}^{1,\pr}(\pi^1\+A)\bigr)\cdot\+t\diff l
    =\int_{\partial\Omega}(\pi^1\+A)\cdot\+t\diff l
     =\int_{\partial\Omega}\+A\cdot\+t\diff l,
  \end{align}
  as claimed.
\end{statementproof}

\subsection{Essential boundary conditions and a homotopy-based curl--curl solver}\label{sec:boundary_homotopy}
The treatment of the coordinate singularity in the previous subsections is independent of the boundary conditions, and a reader interested only in the polar spaces may skip this subsection on a first reading.
We include it because the homotopy operators, once adapted to homogeneous essential boundary conditions, offer more than a technical complement: they provide explicit particular solutions of the strong curl equation and reduce the curl--curl problem to a single inversion of the $2$-form mass matrix (\cref{sec:bvp_2d}), which we also use to prepare the initial data of the 3D wave experiment in \cref{sec:results_wave3d}.
Throughout this subsection, let $\pr\ge1$.

\subsubsection{Relative spaces and boundary-compatible homotopies}
To be able to discuss the boundary conditions, we introduce the following trace operators 
\begin{equation}
  (\operatorname{tr}^0_\bc P)(\logical{y}) \defeq P\circ\mapping(1,\logical{y}),
  \quad
  (\operatorname{tr}^1_\bc\+A)(\logical{y})
  \defeq
  (\+A\circ\mapping(1,\logical{y}))
  \cdot\partial_{\logical{y}}\mapping(1,\logical{y}).
\end{equation}
The boundary conditions $P=0$ and $\+A\times\+n=\+0$ are then equivalent to $\operatorname{tr}^0_\bc P=0$ and $\operatorname{tr}^1_\bc\+A=0$, respectively.
The boundary conditions are compatible in the sense that
\begin{equation}\label{eqn:trace_gradient_identity}
  \operatorname{tr}^1_\bc\gradientt P
  =(\gradientt P\circ\mapping(1,\logical{y}))
    \cdot\partial_{\logical{y}}\mapping(1,\logical{y})
  =\partial_{\logical{y}}(P\circ\mapping(1,\logical{y}))
  =\partial_{\logical{y}}\operatorname{tr}^0_\bc P.
\end{equation}
We call the following spaces with homogeneous boundary traces the relative polar spaces
\begin{equation}\label{eqn:relative_polar_spaces}
  \polar{\funspace{\Bspline}}^{0,\pr}_\bc
  \defeq
  \setdef{P\in\polar{\funspace{\Bspline}}^{0,\pr}}
  {\operatorname{tr}^0_\bc P=0},
  \quad
  \polar{\funspace{\Bspline}}^{1,\pr}_\bc
  \defeq
  \setdef{\+A\in\polar{\funspace{\Bspline}}^{1,\pr}}
  {\operatorname{tr}^1_\bc\+A=0}.
\end{equation}
The compactly supported internal polar $0$- and $1$-form basis functions have zero boundary value and tangential trace, respectively. Thus imposing these conditions removes the tensor-product $0$-forms $\Bspline^{0}_{n^{0,\xx}j}$ and angular $1$-forms $\+\Bspline^{1}_{\yy,n^{0,\xx}j}$ of the outermost radial ring, which are external by the size condition $n^{0,\xx}\ge\pr+2$ of \cref{ass:logical_space_parameters}, while the $2$-form space remains unchanged.
The ambient gradient inclusion in \cref{thm:polar_subspaces_commutativity_conformity} and the trace-gradient identity~\eqref{eqn:trace_gradient_identity} give $\gradientt\polar{\funspace{\Bspline}}^{0,\pr}_\bc \subset \polar{\funspace{\Bspline}}^{1,\pr}_\bc$.
The ambient curl inclusion and $\curlt\gradientt=0$ then show that the relative spaces form the cochain complex
\begin{equation}
  \derhamcomplex[
    dimension=2,
    custom-spaces={\polar{\funspace{\Bspline}}^{0,\pr}_\bc,\polar{\funspace{\Bspline}}^{1,\pr}_\bc,\polar{\funspace{\Bspline}}^{2,\pr}},
    operators={\gradientt,\curlt}
  ]
  .
\end{equation}

\begin{figure}
  \centering
  \tikzsetnextfilename{bdy_homotopy}\input{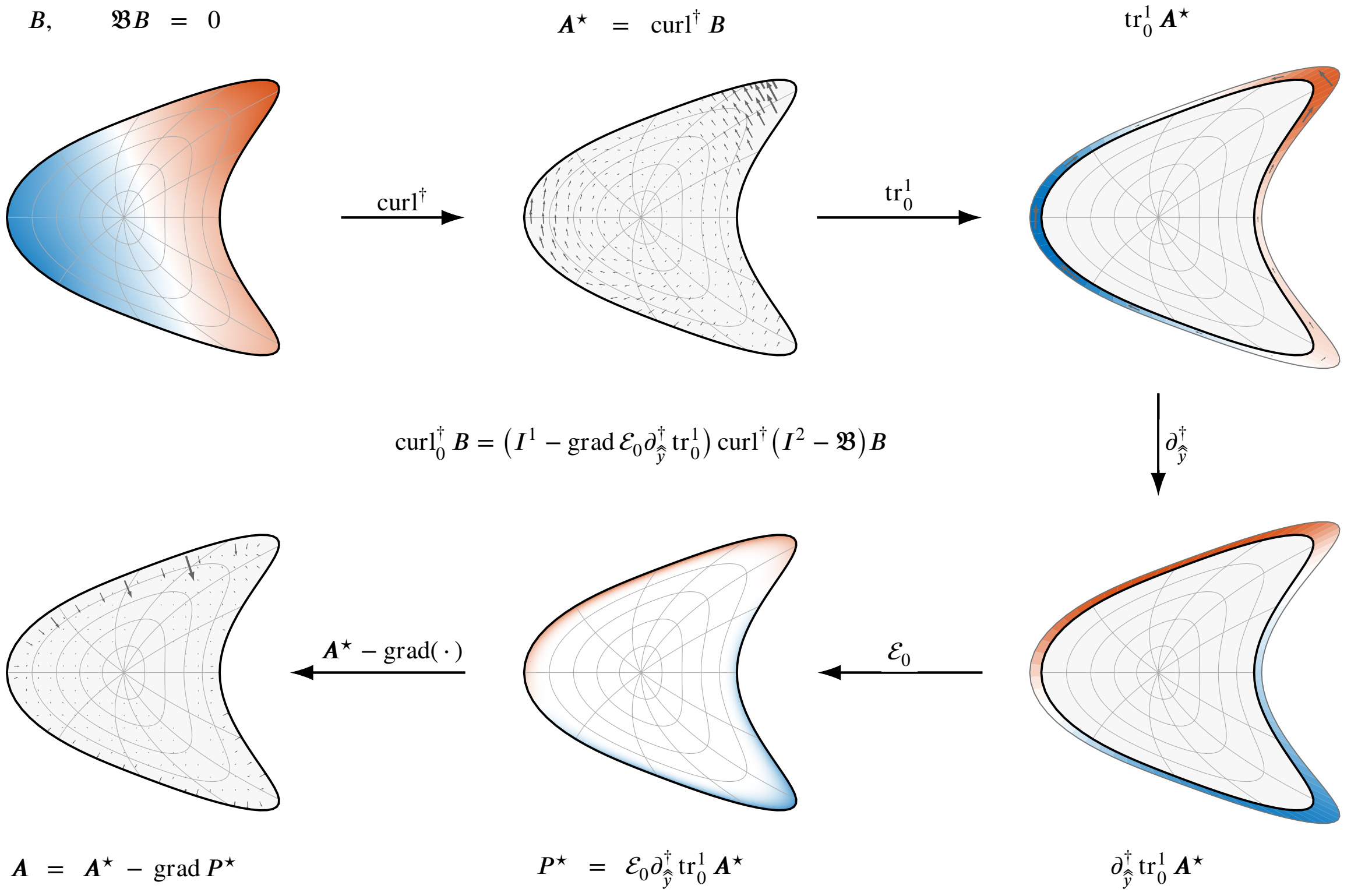}
  \caption{Illustration of the construction of the boundary-compatible homotopy operator $\curlt^\dagger_\bc$ in \cref{eqn:boundary_curl_homotopy_2d} for a zero-integral polar $2$-form $B$ on the banana-shaped domain of \cref{fig:disk_mapping}. For a general input, the $2$-form is projected into the zero-integral subspace by $I^2-\mathfrak{B}$. The tangential trace of $\curlt^\dagger B$ on the boundary is computed, followed by its periodic primitive using a 1D homotopy operator. This primitive is extended into the interior of the domain, and its gradient is subtracted from $\curlt^\dagger B$ to remove the tangential trace.}
  \label{fig:bdy_homotopy}
\end{figure}
By Stokes' theorem, the curl of a $1$-form with zero tangential trace has zero total integral. Therefore, a boundary-compatible curl homotopy can recover only the zero-integral part of an arbitrary $2$-form.
Even for a zero-integral $B$, the ordinary homotopy from \cref{sec:derham_hodge_2d_hodge} may have a nonzero tangential trace on the boundary,
\begin{equation}
  \operatorname{tr}^1_\bc\curlt^\dagger B
  \neq 0.
\end{equation}
This trace can be removed by subtracting the gradient of a $0$-form whose boundary value is a periodic primitive of the trace.
We first compute this primitive using a 1D homotopy operator and then extend it into the interior of the domain, as illustrated in \cref{fig:bdy_homotopy}.
A similar approach can be applied to the $\gradientt^\dagger$ homotopy operator. 
Let $\polar{\pi}^{2,\pr}_{\funspace{L}_2}$ denote the $\funspace{L}_2$-orthogonal projection onto $\polar{\funspace{\Bspline}}^{2,\pr}$.
This leads to the following boundary-compatible homotopy operators
\begin{subequations}\label{eqn:boundary_homotopy_operators_2d}
  \begin{align+}
    \gradientt^\dagger_\bc\+A
    &\defeq
    \srpar{
      I^0 - \operatorname{tr}^0_\bc
    }\gradientt^\dagger \srpar{I^1-\curlt^\dagger_\bc\curlt}\+A
    ,\\
    \curlt^\dagger_\bc B
    &\defeq
    \srpar{I^1 - \gradientt\mathcal{E}_\bc \partial^\dagger_{\logical{y}} \operatorname{tr}^1_\bc} \curlt^\dagger \srpar{I^2-\mathfrak{B}} B
    \label{eqn:boundary_curl_homotopy_2d}
    ,\\
    \mathfrak{B}B
    &\defeq
    \frac{\polar{\pi}^{2,\pr}_{\funspace{L}_2}1}
      {\int_\Omega \polar{\pi}^{2,\pr}_{\funspace{L}_2}1\diff^2x}
    \int_\Omega B\diff^2x
    .
  \end{align+}
\end{subequations}

The polar $2$-form space contains an external form of nonzero integral: the size condition in \cref{ass:logical_space_parameters} provides an external radial basis function, and the periodic angular $1$-form space contains a function of nonzero integral. Hence $\polar{\pi}^{2,\pr}_{\funspace{L}_2}1$ is nonzero. Orthogonality shows that its integral equals its positive squared $\funspace{L}_2$ norm. Consequently, $\mathfrak{B}$ is a well-defined, integral-preserving $\funspace{L}_2$-orthogonal projector on the polar $2$-form space, with kernel given by the zero-integral forms. This keeps the mean correction in the discrete space even when the physical constant $1$ does not belong to it.
Here and below, a constant boundary $0$-form is identified with its constant
extension to $\Omega$.  As shown in the proof of
\cref{lem:boundary_compatible_homotopy}, the trace in the definition of
$\gradientt^\dagger_\bc$ is indeed constant.  Furthermore,
$\mathcal{E}_\bc$ extends a 1D $0$-form defined on the boundary into the interior
of the domain,
\begin{equation}
  \mathcal{E}_\bc P
  \defeq
  \pushforward^0\left(
    \Bsplineone^{0,\xx}_{n^{0,\xx}}(\logical{x}) P(\logical{y})
  \right).
\end{equation}
and $\partial^\dagger_{\logical{y}}$ computes the periodic primitive of zero angular mean for a boundary $1$-form $A$ with zero integral,
\begin{equation}
  (\partial^\dagger_{\logical{y}} A)(\logical{y})
  \defeq
  \int_0^{\logical{y}}A(s)\diff s
  -\int_0^1\int_0^t A(s)\diff s\diff t.
\end{equation}
The radial homotopies use cumulative sums of coefficients, see \cref{eqn:boundary_homotopy_matrices} for the matrix representation. The boundary primitive uses cumulative sums with a mean adjustment for periodic B-splines, or direct integration of the trigonometric modes. For a fixed space, $\polar{\pi}^{2,\pr}_{\funspace{L}_2}1$ and its normalization can be computed once by solving a system with the polar $2$-form mass matrix. Subsequent applications of $\mathfrak{B}$ only require the total integral and multiplication by this fixed normalized form.

With the boundary-compatible homotopy operators, we can now state their corresponding homotopy identities, which are the boundary-compatible version of \cref{eqn:homotopy_identities}.
\begin{lemma}[Boundary-compatible homotopy identities]\label{lem:boundary_compatible_homotopy}
  The operators in \cref{eqn:boundary_homotopy_operators_2d} satisfy the following chain complex properties
  \begin{equation}
    \gradientt^\dagger_\bc
    \polar{\funspace{\Bspline}}^{1,\pr}_\bc
    \subset
    \polar{\funspace{\Bspline}}^{0,\pr}_\bc,
    \quad
    \curlt^\dagger_\bc
    \polar{\funspace{\Bspline}}^{2,\pr}
    \subset
    \polar{\funspace{\Bspline}}^{1,\pr}_\bc,   
    \quad
    \gradientt^\dagger_\bc\curlt^\dagger_\bc=0.
  \end{equation}
  Moreover, the following homotopy identities for $P \in\polar{\funspace{\Bspline}}^{0,\pr}_\bc$,
  $\+A\in\polar{\funspace{\Bspline}}^{1,\pr}_\bc$, and $B\in\polar{\funspace{\Bspline}}^{2,\pr}$ (cf.\ \cref{eqn:homotopy_identities}) are satisfied
  \begin{subequations}\label{eqn:boundary_homotopy_identities}
    \begin{align}
      P&=\gradientt^\dagger_\bc\gradientt P,
      \label{eqn:boundary_homotopy_identity_zero}\\
      \+A&=\gradientt\gradientt^\dagger_\bc\+A
      +\curlt^\dagger_\bc\curlt\+A,
      \label{eqn:boundary_homotopy_identity_one}\\
      B&=\mathfrak{B}B+\curlt\curlt^\dagger_\bc B.
      \label{eqn:boundary_homotopy_identity_two}
    \end{align}
  \end{subequations}
  In particular, the relative polar complex $\polar{\funspace{\Bspline}}^{0,\pr}_\bc\to\polar{\funspace{\Bspline}}^{1,\pr}_\bc\to\polar{\funspace{\Bspline}}^{2,\pr}$ has cohomology of dimensions $0,0,1$, as the relative de Rham complex of the disk, and $\mathfrak{B}$ projects onto the nontrivial class.
\end{lemma}
\begin{statementproof}[appendix]{lem:boundary_compatible_homotopy}
  The proof has five steps: the curl homotopy and the $2$-form identity, the gradient homotopy and the $1$-form identity, the $0$-form identity, the chain-complex property, and the cohomology.
  The main observation is that, by Stokes' theorem, a zero-integral $2$-form has a lift whose tangential trace has zero circulation along the boundary.
  This trace therefore has a periodic primitive, and subtracting the gradient of an extension of that primitive removes the trace without changing the curl.

  \emph{The curl homotopy and the $2$-form identity.}
  For $B\in\polar{\funspace{\Bspline}}^{2,\pr}$, set
  \begin{equation}
    B_0\defeq B-\mathfrak{B}B,
    \qquad
    A_\bc\defeq
    \operatorname{tr}^1_\bc\curlt^\dagger B_0.
  \end{equation}
  The definition of $\mathfrak{B}$ gives $B_0\in\polar{\funspace{\Bspline}}^{2,\pr}$ with zero integral.
  Stokes' theorem and the ordinary $2$-form homotopy identity give
  \begin{equation}
    \int_0^1 A_\bc(\logical{y})\diff\logical{y}
    =\int_0^1
      (\curlt^\dagger B_0\circ\mapping(1,\logical{y}))
      \cdot\partial_{\logical{y}}\mapping(1,\logical{y})
      \diff\logical{y}
    =\int_\Omega\curlt\curlt^\dagger B_0\diff^2x
    =\int_\Omega B_0\diff^2x
    =0.
  \end{equation}
  Consequently, $\partial^\dagger_{\logical{y}} A_\bc$ is periodic and satisfies $\partial_{\logical{y}}\partial^\dagger_{\logical{y}} A_\bc=A_\bc$.
  The trace $A_\bc$ lies in the periodic angular $1$-form space, and the derivative relations in \cref{eqn:bspline_derivative,eqn:trigonometric_derivative} show that its zero-integral primitive lies in the corresponding angular $0$-form space.
  Since $\operatorname{tr}^0_\bc\mathcal{E}_\bc=I^0$ by the definition of the extension operator, and the traces commute with the gradient,
  \begin{equation}
    \operatorname{tr}^1_\bc\curlt^\dagger_\bc B
    =A_\bc-
      \operatorname{tr}^1_\bc\gradientt
      \mathcal{E}_\bc\partial^\dagger_{\logical{y}} A_\bc
    =A_\bc-
      \partial_{\logical{y}}\partial^\dagger_{\logical{y}} A_\bc
    =0.
  \end{equation}
  Thus $\curlt^\dagger_\bc B$ satisfies the homogeneous boundary condition.
  Moreover, $\curlt^\dagger B_0$ belongs to the polar $1$-form space by \cref{thm:polar_subspaces_discrete_homotopy}, and $\mathcal{E}_\bc\partial^\dagger_{\logical{y}} A_\bc$ belongs to the polar $0$-form space because $n^{0,\xx}\ge\pr+2$ (\cref{ass:logical_space_parameters}) makes the outermost radial basis function external.
  Hence $\curlt^\dagger_\bc B\in\polar{\funspace{\Bspline}}^{1,\pr}_\bc$.
  Since the correction term is a gradient, it does not change the curl:
  \begin{equation}
    \curlt\curlt^\dagger_\bc B
    =\curlt\curlt^\dagger B_0
    =B_0
    =B-\mathfrak{B}B.
  \end{equation}
  This proves the mapping property of $\curlt^\dagger_\bc$ and \cref{eqn:boundary_homotopy_identity_two}.

  \emph{The gradient homotopy and the $1$-form identity.}
  Let $\+A\in\polar{\funspace{\Bspline}}^{1,\pr}_\bc$ and define
  \begin{equation}
    \+Z_{\+A}\defeq
    (I^1-\curlt^\dagger_\bc\curlt)\+A.
  \end{equation}
  The boundary condition and Stokes' theorem give
  \begin{equation}
    \int_\Omega\curlt\+A\diff^2x
    =\int_0^1\operatorname{tr}^1_\bc\+A\diff\logical{y}
    =0,
  \end{equation}
  so $\mathfrak{B}\curlt\+A=0$, and the $2$-form identity with $B=\curlt\+A$ gives
  \begin{equation}
    \curlt\+Z_{\+A}
    =\curlt\+A-\curlt\curlt^\dagger_\bc\curlt\+A
    =0.
  \end{equation}
  The ordinary $1$-form homotopy identity therefore reduces to $\+Z_{\+A}=\gradientt Q_{\+A}$ with $Q_{\+A}\defeq\gradientt^\dagger\+Z_{\+A}$.
  Both terms defining $\+Z_{\+A}$ have vanishing tangential trace, so trace compatibility gives
  \begin{equation}
    \partial_{\logical{y}}\operatorname{tr}^0_\bc Q_{\+A}
    =\operatorname{tr}^1_\bc\gradientt Q_{\+A}
    =\operatorname{tr}^1_\bc\+Z_{\+A}
    =0.
  \end{equation}
  Hence $\operatorname{tr}^0_\bc Q_{\+A}$ is constant, and with the convention that this constant is extended unchanged to $\Omega$, the definition in \cref{eqn:boundary_homotopy_operators_2d} reads
  \begin{equation}
    \gradientt^\dagger_\bc\+A
    =Q_{\+A}-\operatorname{tr}^0_\bc Q_{\+A}.
  \end{equation}
  This function has zero trace, and both $Q_{\+A}$ and the constant belong to the polar $0$-form space, which proves the mapping property of $\gradientt^\dagger_\bc$.
  Subtracting a constant does not change the gradient, so
  \begin{equation}
    \gradientt\gradientt^\dagger_\bc\+A
    =\gradientt Q_{\+A}
    =\+Z_{\+A}
    =\+A-\curlt^\dagger_\bc\curlt\+A,
  \end{equation}
  which is \cref{eqn:boundary_homotopy_identity_one}.

  \emph{The $0$-form identity.}
  Let $P\in\polar{\funspace{\Bspline}}^{0,\pr}_\bc$.
  Since $\curlt\gradientt P=0$, we have $\+Z_{\gradientt P}=\gradientt P$, and the ordinary $0$-form homotopy identity gives $Q_{\gradientt P}=\gradientt^\dagger\gradientt P=P-P(\+0)$.
  Because $P$ vanishes on the boundary, the trace of $Q_{\gradientt P}$ is the constant $-P(\+0)$, and subtracting it gives $\gradientt^\dagger_\bc\gradientt P=P$, which is \cref{eqn:boundary_homotopy_identity_zero}.

  \emph{The chain-complex property.}
  Let $B\in\polar{\funspace{\Bspline}}^{2,\pr}$ and $\+A\defeq\curlt^\dagger_\bc B$; we show $\gradientt^\dagger_\bc\+A=0$.
  By definition, $\gradientt^\dagger_\bc\+A$ is obtained by applying $(I^0-\operatorname{tr}^0_\bc)\gradientt^\dagger$ to $\+Z_{\+A}=\+A-\curlt^\dagger_\bc\curlt\+A$.
  It therefore suffices to show that $\+Z_{\+A}=\+0$, that is,
  \begin{equation}\label{eqn:boundary_curl_homotopy_reproduces_lift}
    \curlt^\dagger_\bc\curlt\+A=\+A.
  \end{equation}
  First, the $2$-form identity proved above gives $\curlt\+A=\curlt\curlt^\dagger_\bc B=B-\mathfrak{B}B$.
  Second, $\curlt^\dagger_\bc$ applies $I^2-\mathfrak{B}$ to its argument before anything else, see \cref{eqn:boundary_curl_homotopy_2d}.
  Since $\mathfrak{B}$ is a projector, $(I^2-\mathfrak{B})(B-\mathfrak{B}B)=(I^2-\mathfrak{B})B$, and hence
  \begin{equation}
    \curlt^\dagger_\bc\curlt\+A
    =\curlt^\dagger_\bc(B-\mathfrak{B}B)
    =\curlt^\dagger_\bc B
    =\+A.
  \end{equation}
  This proves \cref{eqn:boundary_curl_homotopy_reproduces_lift}, so $\+Z_{\+A}=\+0$ and $\gradientt^\dagger_\bc\curlt^\dagger_\bc B=0$.

  \emph{The cohomology.}
  By \cref{eqn:boundary_homotopy_identity_zero}, the gradient is injective on $\polar{\funspace{\Bspline}}^{0,\pr}_\bc$, and by \cref{eqn:boundary_homotopy_identity_one} and the mapping property of $\gradientt^\dagger_\bc$, every curl-free $\+A\in\polar{\funspace{\Bspline}}^{1,\pr}_\bc$ is the gradient of a relative polar $0$-form.
  By \cref{eqn:boundary_homotopy_identity_two}, every polar $2$-form is the sum of the curl of a relative polar $1$-form and an element of the one-dimensional image of $\mathfrak{B}$.
  This image contains no nonzero curls, since the curl of a relative $1$-form has zero integral by Stokes' theorem, whereas a nonzero element of the image of $\mathfrak{B}$ has nonzero integral, as the kernel of the projector $\mathfrak{B}$ consists of the zero-integral forms.
  Hence the second cohomology is one-dimensional and represented by the image of $\mathfrak{B}$.
\end{statementproof}


\subsubsection{Solving the curl--curl problem}\label{sec:bvp_2d}
The following is not needed for the polar construction, but shows what the boundary-compatible homotopy operators buy.
The boundary-compatible homotopy operators solve the strong curl equation
\begin{equation}
  \text{Find $\+A\in\polar{\funspace{\Bspline}}^{1,\pr}_\bc$ such that: }
  \curlt\+A = B,
\end{equation}
for any polar $2$-form $B \in \polar{\funspace{\Bspline}}^{2,\pr}$ with zero mean, i.e., $\mathfrak{B}B=0$.
A particular solution is given by $\+A = \curlt^\dagger_\bc B$.

The homotopy operators can also be used to solve the following curl--curl problem with homogeneous boundary conditions
\begin{equation}
  \text{Find $\+A\in\polar{\funspace{\Bspline}}^{1,\pr}_\bc$ such that: }
  \inprod{\curlt \+\Bspline^1}{\curlt \+A} = \inprod{\+\Bspline^1}{\+J}
  , \quad
  \forall \+\Bspline^1\in\polar{\funspace{\Bspline}}^{1,\pr}_\bc
  ,
\end{equation}
for any polar $1$-form $\+J \in \polar{\funspace{\Bspline}}^{1,\pr}_\bc$ that is weakly divergence-free
\begin{equation}\label{eqn:weak_divergence_free}
  \inprod{\gradientt \Bspline^0}{\+J} = 0
  , \quad
  \forall \Bspline^0\in\polar{\funspace{\Bspline}}^{0,\pr}_\bc
  .
\end{equation}
An additional gauge condition is required to ensure uniqueness of the solution $\+A$, but we do not discuss this here.

The homotopy identity~\eqref{eqn:boundary_homotopy_identity_one} can be used to express both the test and trial functions in terms of the boundary-compatible homotopy operators
\begin{equation}
  \+\Bspline^1 = \gradientt \Bspline^0 + \curlt^\dagger_\bc \Bspline^2
  , \quad
  \+A = \gradientt P + \curlt^\dagger_\bc B,
\end{equation}
resulting in the following equivalent problem
\begin{equation}
  \text{Find $B \in \polar{\funspace{\Bspline}}^{2,\pr}$ such that: }
  \inprod{\Bspline^2}{B - \mathfrak{B}B} 
  = 
  \inprod{\curlt^\dagger_\bc \Bspline^2}{\+J}
  , \quad
  \forall \Bspline^2\in\polar{\funspace{\Bspline}}^{2,\pr}
  ,
\end{equation}
where we have made use of \cref{eqn:weak_divergence_free,eqn:boundary_homotopy_identity_two}, as well as the $\funspace{L}_2$-orthogonality of $\mathfrak{B}$.
Dropping the $\inprod{\Bspline^2}{\mathfrak{B} B}$ term selects the zero-mean part of $B$: testing with $\Bspline^2=\polar{\pi}^{2,\pr}_{\funspace{L}_2}1$ gives $\int_\Omega B\diff^2x=0$, since $\curlt^\dagger_\bc\polar{\pi}^{2,\pr}_{\funspace{L}_2}1=\+0$.
Hence, solving the curl--curl problem in 2D requires only the inversion of the polar $2$-form mass matrix, provided that the right-hand side is already weakly divergence-free.
The particular solution can be written as
\begin{equation}
  \+A 
  = \curlt^\dagger_\bc (\curlt^\dagger_\bc)^\transpose \+J
  ,
\end{equation}
where $(\curlt^\dagger_\bc)^\transpose$ is the adjoint of $\curlt^\dagger_\bc$ with respect to the $\funspace{L}_2$ inner product, and thus involves the inversion of the polar $2$-form mass matrix, see \cref{eqn:polar_curlcurl_matrix}.

  \section{Polar subspaces of the 3D de Rham complex}\label{sec:high_order_derham_hodge_3d}
We extend the 2D construction to 3D physical domains, eventually resulting in polar subspaces for which the following two diagrams individually commute (we again omit the logical and intermediate spaces for clarity):
\begin{equation}\label{eqn:derham_hodge_3d_diagram}
  \tikzsetnextfilename{derham_physical_3d}\begin{derhamdiagram}[dimension=3]
  \derhamrow[name=polar,spaces={\polar{\funspace{\Bsplinethree}}^{0,\pr},\polar{\bfunspace{\Bsplinethree}}^{1,\pr},\polar{\bfunspace{\Bsplinethree}}^{2,\pr},\polar{\funspace{\Bsplinethree}}^{3,\pr}},operators={\gradientt,\curlt,\divergencet}]
  \derhamrow[name=continuous,spaces={\funspace{H}_1,\bfunspace{H}_{\curlt},\bfunspace{H}_{\divergencet},\funspace{L}_2},operators={\gradientt,\curlt,\divergencet}]
  \derhamconnect[from=continuous,to=polar,degree=0,operator=\polar{\pi}^{0,\pr}]
  \derhamconnect[from=continuous,to=polar,degree=1,operator=\polar{\pi}^{1,\pr}]
  \derhamconnect[from=continuous,to=polar,degree=2,operator=\polar{\pi}^{2,\pr}]
  \derhamconnect[from=continuous,to=polar,degree=3,operator=\polar{\pi}^{3,\pr}]
\end{derhamdiagram}
\qquad
\begin{derhamdiagram}[dimension=3]
  \derhamrow[name=polar,kind=dual,spaces={\polar{\funspace{\Bsplinethree}}^{0,\pr},\polar{\bfunspace{\Bsplinethree}}^{1,\pr},\polar{\bfunspace{\Bsplinethree}}^{2,\pr},\polar{\funspace{\Bsplinethree}}^{3,\pr}},operators={\gradientt^\dagger,\curlt^\dagger,\divergencet^\dagger}]
  \derhamrow[name=continuous,kind=dual,spaces={\funspace{H}_1,\bfunspace{H}_{\curlt},\bfunspace{H}_{\divergencet},\funspace{L}_2},operators={\gradientt^\dagger,\curlt^\dagger,\divergencet^\dagger}]
  \derhamconnect[from=continuous,to=polar,degree=0,operator=\polar{\pi}^{0,\pr}]
  \derhamconnect[from=continuous,to=polar,degree=1,operator=\polar{\pi}^{1,\pr}]
  \derhamconnect[from=continuous,to=polar,degree=2,operator=\polar{\pi}^{2,\pr}]
  \derhamconnect[from=continuous,to=polar,degree=3,operator=\polar{\pi}^{3,\pr}]
\end{derhamdiagram},
\end{equation}
where $\gradientt, \curlt, \divergencet$ denote the usual gradient, curl and divergence operators in 3D, and the corresponding homotopy operators are denoted by $\gradientt^\dagger, \curlt^\dagger, \divergencet^\dagger$.
The forward arrows in the first diagram form the de Rham cochain complex, whereas the backward arrows in the second diagram form the corresponding chain complex of homotopy operators.  
The latter are constructed in \cref{sec:derham_hodge_3d_homotopy}.
The newly introduced space $\bfunspace{H}_{\divergencet}$ is defined as the space of square-integrable vector fields with square-integrable divergence.
As before, zero-forms are denoted by $P$, one-forms by $\+A$, two-forms by $\+B$, and the newly introduced three-forms by $Q$.
The polar $m$-form spaces will again be subspaces of the tensor-product B-spline spaces
\begin{equation}
  \polar{\funspace{\Bsplinethree}}^{m,\pr} \subset \funspace{\Bsplinethree}^{m}, \quad m = 0, 1, 2, 3,
\end{equation}
which are defined such that the resulting $m$-forms have desirable regularity properties at the coordinate singularity and the spaces are conforming subspaces of the corresponding infinite-dimensional Hilbert spaces.
Polar spline complexes with the cohomology of the solid torus have been constructed before with first-order regularity, by extraction \citep{patrizi2025torus} and through commuting projections \citep{holderied2022}.
The contribution of this section is not the existence of such a complex, but its high-order regularity along the central curve, which is inherited from the 2D enrichment and homotopy construction through the tensor product, together with the homotopy operators and the commuting projections.

Throughout the discussion of the 3D construction, we will redefine several objects that were already defined in the 2D construction: the 2D equivalent will be denoted with a subscript $\mathrm{2D}$, and the 3D equivalent will be denoted without a subscript.
For example, $\logicalmapping_{\mathrm{2D}}$ denotes the 2D logical mapping, while $\logicalmapping$ denotes the 3D logical mapping.

\subsection{Coordinate transformation}\label{sec:coordtransform_3d}
\begin{figure}
  \centering
  \tikzsetnextfilename{coordinate_transformation_3d}\input{tikz/coordinate_transformation_3d.tex}
  \caption{The logical domain $\logical{\Omega}$ and the physical domain $\Omega=\mapping(\logical{\Omega})$ of \cref{ex:tube_mapping_3d}; the intermediate cylinder $\intermediate{\Omega}$ is not drawn.
  The orange logical face $\logical{x}=0$ collapses to the orange central curve $\Gamma$, along whose part exposed by the cut the local frame $({\color{lines1}\+T},{\color{lines3}\+N_1},{\color{lines4}\+N_2})$ is shown at the left-most poloidal plane and at three points in between.}
  \label{fig:coordinate_transformation_3d}
\end{figure}
As in 2D, the logical domain is mapped to the physical domain in two steps (see \cref{fig:coordinate_transformation_3d}).
The logical domain is $\logical{\Omega}\defeq\logical{\Omega}_\mathrm{2D} \times \mathbb{S}^1$, with $\mathbb{S}^1=\mathbb{R}/\mathbb{Z}$ as in \cref{sec:coordtransform}, with radial, poloidal and toroidal coordinates $\logical{\+x} = (\logical{x},\logical{y},\logical{z})$, and we write $\logical{\+x}_\mathrm{2D}=(\logical{x},\logical{y})$ for the cross-sectional coordinates.
The logical mapping is the product of a 2D logical mapping $\logicalmapping_\mathrm{2D}$ of the form \cref{eqn:disk_mapping}, which satisfies \cref{ass:disc_like_mapping}, with the identity in the toroidal direction,
\begin{equation}\label{eqn:logical_mapping_3d}
  \logicalmapping(\logical{\+x}_\mathrm{2D},\logical{z})
  \defeq
  \bigl(\logicalmapping_\mathrm{2D}(\logical{\+x}_\mathrm{2D}),\logical{z}\bigr)
  , \quad\text{i.e.,}\quad
  \logicalmapping=\logicalmapping_\mathrm{2D}\times\mathrm{id}.
\end{equation}
It maps $\logical{\Omega}$ onto the intermediate domain, the straight periodic cylinder
\begin{equation}
  \intermediate{\Omega}\defeq\logicalmapping(\logical{\Omega})=\intermediate{\Omega}_\mathrm{2D}\times\mathbb{S}^1
  , \quad
  \intermediate{\Omega}_\mathrm{2D}\defeq\logicalmapping_\mathrm{2D}(\logical{\Omega}_\mathrm{2D}),
\end{equation}
with coordinates $(\intermediatecoord{\+x}_\mathrm{2D},\intermediatecoord{z})$, where $\intermediatecoord{\+x}_\mathrm{2D}=(\intermediatecoord{x},\intermediatecoord{y})$ and the periodic coordinate $\intermediatecoord{z}=\logical{z}$.
The logical face $\logical{x}=0$ collapses to the axis $\{\+0\}\times\mathbb{S}^1$.

The physical mapping is the composition
\begin{equation}\label{eqn:torus_mapping_3d}
  \mapping\defeq\intermediatemapping\circ\logicalmapping
  , \qquad
  \Omega=\mapping(\logical{\Omega}),
\end{equation}
and the logical face $\logical{x}=0$ collapses to the central curve
\begin{equation}\label{eqn:central_curve_3d}
  \Gamma\defeq\intermediatemapping(\{\+0\}\times\mathbb{S}^1),
\end{equation}
a smooth closed curve without self-intersections in the interior of $\Omega$, which we orient by increasing $\intermediatecoord{z}$.
As in 2D, a given physical mapping can be decomposed in several ways, and the polar spaces constructed below depend on the decomposition.
\begin{assumption}[Embedded solid torus]\label{ass:diffeomorphism_3d}
  The intermediate mapping $\intermediatemapping$ is defined on $U\times\mathbb{S}^1$ for an open neighborhood $U\subset\mathbb{R}^2$ of the closure of $\intermediate{\Omega}_\mathrm{2D}$, that is, it is periodic in $\intermediatecoord{z}$.
  It is a $\funspace{C}^\infty$ diffeomorphism from $U\times\mathbb{S}^1$ onto its image in $\mathbb{R}^3$, and it is orientation preserving, that is, its Jacobian satisfies $\intermediate{\jac}>0$.
\end{assumption}

The Jacobian matrix and Jacobian of the logical mapping are
\begin{equation}\label{eqn:logical_jacobian_3d}
  \logical{\jacmat}
  =
  \begin{bmatrix}
    \logical{\jacmat}_\mathrm{2D} & \+0\\
    \+0^\transpose & 1
  \end{bmatrix}
  , \quad
  \logical{\jac}=\logical{\jac}_\mathrm{2D},
\end{equation}
where $\logical{\jacmat}_\mathrm{2D}$ and $\logical{\jac}_\mathrm{2D}$ are given by \cref{eqn:disk_jacobian_matrix,eqn:disk_jacobian} for $\logicalmapping_\mathrm{2D}$ and depend on $\logical{\+x}_\mathrm{2D}$ only.
Denoting the Jacobian matrix and Jacobian of $\intermediatemapping$ by $\intermediate{\jacmat}$ and $\intermediate{\jac}=\det\intermediate{\jacmat}$, functions of $(\intermediatecoord{\+x}_\mathrm{2D},\intermediatecoord{z})$, the chain rule gives the Jacobian matrix and Jacobian of the physical mapping (cf.\ \cref{eqn:composed_jacobian})
\begin{equation}\label{eqn:composed_jacobian_3d}
  \jacmat=(\intermediate{\jacmat}\circ\logicalmapping)\,\logical{\jacmat}
  , \quad
  \jac=(\intermediate{\jac}\circ\logicalmapping)\,\logical{\jac}_\mathrm{2D}.
\end{equation}
Hence $\jac$ vanishes on the face $\logical{x}=0$ only, and all singular behavior stems from $\logicalmapping_\mathrm{2D}$.
Indeed, $\logicalmapping_\mathrm{2D}$ maps $\{\logical{x}>0\}$ diffeomorphically onto $\intermediate{\Omega}_\mathrm{2D}\setminus\{\+0\}$ by \cref{ass:disc_like_mapping}, and $\intermediatemapping$ maps $(\intermediate{\Omega}_\mathrm{2D}\setminus\{\+0\})\times\mathbb{S}^1$ diffeomorphically onto $\Omega\setminus\Gamma$ by \cref{ass:diffeomorphism_3d}.
Thus $\mapping$ is a smooth diffeomorphism from $\{\logical{x}>0\}$ onto $\Omega\setminus\Gamma$, with periodic coordinates identified, and its inverse
\begin{equation}\label{eqn:inverse_3d_via_inverse_2d}
  \mapping^{-1}=(\logicalmapping_\mathrm{2D}^{-1}\times\mathrm{id})\circ\intermediatemapping^{-1}
\end{equation}
is singular at $\Gamma$ only through $\logicalmapping_\mathrm{2D}^{-1}$.

The pushforward operators of the physical mapping are the 3D Piola transformations (cf.\ \cref{eqn:pushforward_operators})
\begin{equation}\label{eqn:pushforward_operators_3d}
  \pushforward^0 \logical{P} \defeq \logical{P} \circ \mapping^{-1}
  , \quad
  \pushforward^1 \logical{\+A} \defeq (\jacmat^{-\transpose} \logical{\+A}) \circ \mapping^{-1}
  , \quad
  \pushforward^2 \logical{\+B} \defeq \srpar{\frac{1}{\jac} \jacmat \logical{\+B}} \circ \mapping^{-1}
  , \quad 
  \pushforward^3 \logical{Q} \defeq \srpar{\frac{1}{\jac} \logical{Q}} \circ \mapping^{-1}.
\end{equation}
The pushforward operators $\logical{\pushforward}^m$ of $\logicalmapping$ and $\intermediate{\pushforward}^m$ of $\intermediatemapping$ are given by the same formulas, with $\mapping$, $\jacmat$, $\jac$ replaced by $\logicalmapping$, $\logical{\jacmat}$, $\logical{\jac}$ and by $\intermediatemapping$, $\intermediate{\jacmat}$, $\intermediate{\jac}$, respectively, and the pullbacks are denoted by $\pushforward^{-m}$, $\logical{\pushforward}^{-m}$ and $\intermediate{\pushforward}^{-m}$.
By \cref{eqn:torus_mapping_3d,eqn:composed_jacobian_3d}, the pushforwards compose as
\begin{equation}\label{eqn:pushforward_composition_3d}
  \pushforward^m=\intermediate{\pushforward}^m\circ\logical{\pushforward}^m
  , \quad m\in\{0,1,2,3\}.
\end{equation}
By the chain rule, each of the three families commutes with the 3D gradient, curl and divergence, where for $\pushforward^m$ and $\logical{\pushforward}^m$ this holds away from the face $\logical{x}=0$.

\begin{example}[Tube domains]\label{ex:tube_mapping_3d}
  Let $\+\gamma:\mathbb{S}^1\to\mathbb{R}^3$ parameterize a closed curve, let the columns of $\mat{N}(\intermediatecoord{z})\in\mathbb{R}^{3\times2}$ be orthonormal vectors normal to $\+\gamma'(\intermediatecoord{z})$ such that $[\+\gamma'\ \mat{N}]$ is positively oriented, and let $\intermediatemapping_\mathrm{2D}(\cdot;\intermediatecoord{z})$ be cross-section maps that satisfy \cref{ass:diffeomorphism} on $\intermediate{\Omega}_\mathrm{2D}$ and depend smoothly and periodically on $\intermediatecoord{z}$.
  The tube mapping
  \begin{equation}\label{eqn:tube_mapping_3d}
    \intermediatemapping(\intermediatecoord{\+x})
    \defeq
    \+\gamma(\intermediatecoord{z})
    +
    \mat{N}(\intermediatecoord{z})\,\intermediatemapping_\mathrm{2D}(\intermediatecoord{\+x}_\mathrm{2D};\intermediatecoord{z})
  \end{equation}
  maps each cross-section $\intermediate{\Omega}_\mathrm{2D}\times\{\intermediatecoord{z}\}$ to the 2D domain $\intermediatemapping_\mathrm{2D}(\intermediate{\Omega}_\mathrm{2D};\intermediatecoord{z})$ in the normal plane of the curve at $\+\gamma(\intermediatecoord{z})$, so that $\Gamma=\+\gamma(\mathbb{S}^1)$.
  The cross-sectional coordinates of $\+x=\intermediatemapping(\intermediatecoord{\+x})$ are denoted by
  \begin{equation}
    \+\mu\defeq\mat{N}(\intermediatecoord{z})^\transpose(\+x-\+\gamma(\intermediatecoord{z}))=\intermediatemapping_\mathrm{2D}(\intermediatecoord{\+x}_\mathrm{2D};\intermediatecoord{z}).
  \end{equation}

  The domain in \cref{fig:coordinate_transformation_3d} is of this type.
  Its logical mapping $\logicalmapping_\mathrm{2D}$ is the three-petal disk-like mapping of \cref{fig:disk_mapping}, and its central curve is a helix with four periods,
  \begin{equation}
    \+\gamma(\intermediatecoord{z})
    =\begin{bmatrix}
      \rho\cos(2\pi\intermediatecoord{z})\\
      \rho\sin(2\pi\intermediatecoord{z})\\
      0.12\sin(8\pi\intermediatecoord{z})
    \end{bmatrix}
    , \quad
    \rho=1+0.12\cos(8\pi\intermediatecoord{z})
    ,
  \end{equation}
  and the columns $\+N_1$, $\+N_2$ of $\mat{N}$ are obtained by rotating $\widetilde{\+N}_1=\+T\times\+e_z/\abs{\+T\times\+e_z}$, with the vertical unit vector $\+e_z$ (the tangent of $\+\gamma$, denoted by $\+T$, is never vertical), and $\widetilde{\+N}_2=\+T\times\widetilde{\+N}_1$ about $\+T$ by the angle $4\pi\intermediatecoord{z}$, that is, $\+N_1=\cos(4\pi\intermediatecoord{z})\widetilde{\+N}_1+\sin(4\pi\intermediatecoord{z})\widetilde{\+N}_2$ and $\+N_2=\cos(4\pi\intermediatecoord{z})\widetilde{\+N}_2-\sin(4\pi\intermediatecoord{z})\widetilde{\+N}_1$.
  The cross-section map is $\intermediatemapping_\mathrm{2D}(\intermediatecoord{\+x}_\mathrm{2D};\intermediatecoord{z})=0.22\,\intermediatemapping_{\kappa(\intermediatecoord{z})}(\intermediatecoord{\+x}_\mathrm{2D})$, where $\intermediatemapping_\kappa$ denotes the banana map of \cref{ex:banana_mapping} with $\alpha=2/\sqrt{5}$ and bending $\kappa$, and $\kappa(\intermediatecoord{z})=(1+\sin(2\pi\intermediatecoord{z}))/2\in[0,1]$.

  The tube mapping satisfies \cref{ass:diffeomorphism_3d} whenever it is injective and its Jacobian is positive, which is the case if the cross-sections are small compared to the reach of $\Gamma$, that is, the largest distance from $\Gamma$ within which every point has a unique nearest point on $\Gamma$ \citep[Def.~4.1]{federer1959curvature}; for the domain of \cref{fig:coordinate_transformation_3d} this was verified numerically.
\end{example}

\begin{remark}[Practical applicability]\label{rem:practical_applicability_3d}
  In practice, only the logical mapping $\logicalmapping$ and the physical mapping $\mapping$ are needed: the intermediate mapping itself is not explicitly needed. 
  The extraction matrices of the polar spaces depend on $\logicalmapping$ alone (\cref{sec:derham_hodge_2d_basis}), while the geometric degrees of freedom and the mass matrices only involve $\mapping$ and its Jacobian.
  Equilibrium codes such as GVEC \citep{hindenlang2025gvec} and DESC \citep{dudt2020desc} can therefore be used directly, by letting the logical mapping be the straight cylinder over the unit disk of \cref{ex:disk_mapping} and the physical mapping be the equilibrium mapping.
  See also the introductory \cref{fig:intro_w7x}.
\end{remark}

\subsection{High-order polar subspaces}\label{sec:derham_hodge_3d_construction}
We construct the 3D spaces from the already available 2D polar spaces and a newly introduced 1D periodic cochain complex in the $\logical{z}$-direction formed by the function spaces $\funspace{\Bsplineone}^{m,\zz}$, for $m\in\{0, 1\}$, together with the differential operator $\totd{}{\logical{z}}$ (cf.\ \cref{eqn:exact_sequence}).
The periodic spaces are constructed as in \cref{sec:one_dimensional_function_spaces}, so $\funspace{\Bsplineone}^{0,\zz}\subset\funspace{H}_1^\zz$ and $\funspace{\Bsplineone}^{1,\zz}\subset\funspace{L}_2^\zz$.
The main observation is that the logical de Rham complex is the tensor product of the polar sequence on a disk-like cross-section and the periodic de Rham complex in the toroidal direction.
This lets us reuse the 2D polar spaces, their conformity and regularity properties, their homotopy operators and their commuting projection operators without introducing a separate and more expensive 3D construction.

We make use of the standard tensor product of two cochain complexes of differential forms \citep[Sec.~3]{arnold2015curvilinear} (cf.\ \cref{eqn:logical_function_spaces})
\begin{equation}\label{eqn:tensor_product_complex_spaces_3d}
  \logicalpolar{\funspace{\Bsplinethree}}^{m,\pr}
    \defeq\bigoplus_{p+q=m}
    \logicalpolar{\funspace{\Bspline}}^{p,\pr} \otimes \funspace{\Bsplineone}^{q,\zz}
  , \quad
  m \in \{0,1,2,3\},
\end{equation}
where $\logicalpolar{\funspace{\Bspline}}^{p,\pr} \subset \logical{\funspace{\Bspline}}^{p}$ are the logical 2D polar spaces from \cref{sec:derham_hodge_2d_construction}, constructed with the logical mapping $\logicalmapping_\mathrm{2D}$, and the corresponding products of elements are given in \cref{eqn:tensor_product_complex_spaces_3d_explicit}.
A similar procedure is used for computing the tensor product of the 2D and 1D differential operators, resulting in the usual 3D gradient, curl and divergence operators.
The intermediate and physical polar spaces are defined by means of the pushforward operators in \cref{eqn:pushforward_operators_3d}:
\begin{equation}\label{eqn:polar_spaces_3d_pushforward}
  \intermediatepolar{\funspace{\Bsplinethree}}^{m,\pr}\defeq\logical{\pushforward}^m \logicalpolar{\funspace{\Bsplinethree}}^{m,\pr}
  , \quad
  {\polar{\funspace{\Bsplinethree}}}^{m,\pr}\defeq\intermediate{\pushforward}^m\intermediatepolar{\funspace{\Bsplinethree}}^{m,\pr} = \pushforward^m \logicalpolar{\funspace{\Bsplinethree}}^{m,\pr},
\end{equation}
where the last equality is \cref{eqn:pushforward_composition_3d}; as in 2D, the forms on all three levels have the same coefficients.

The tensor product in \cref{eqn:tensor_product_complex_spaces_3d} combines a 2D $p$-form proxy with a periodic $q$-form proxy according to the wedge product, resulting in a 3D $(p+q)$-form proxy
\begin{equation}
  \logicalpolar{\funspace{\Bspline}}^{p,\pr} \otimes \funspace{\Bsplineone}^{q,\zz}
  \defeq
  \spn\setdef{
    \logicalpolar{\Bspline}^{p,\pr} \wedge \Bsplineone^q
  }{
    \logicalpolar{\Bspline}^{p,\pr} \in \logicalpolar{\funspace{\Bspline}}^{p,\pr}, \Bsplineone^q \in \funspace{\Bsplineone}^{q,\zz}
  }
  , \quad
  p \in \{0,1,2\}, \quad q \in \{0,1\}.
\end{equation}
For $\logical{V}^p\in\logical{\funspace{\Bspline}}^p$ and $U^q\in\funspace{\Bsplineone}^{q,\zz}$, we define
\begin{equation}\label{eqn:tensor_product_complex_spaces_3d_explicit}
  \logical{V}^p\wedge U^q
  \defeq\mat{\Lambda}^{p,q}\logical{V}^p U^q,
  \qquad p\in\{0,1,2\},\quad q\in\{0,1\},
\end{equation}
where the following constant matrices place the 2D proxy in the corresponding 3D components, including the signs from the wedge product:
\begin{equation}\label{eqn:tensor_product_proxy_matrices_3d}
    \mat{\Lambda}^{0,0} = \mat{\Lambda}^{2,1} = 1, \quad
    \mat{\Lambda}^{0,1} = \mat{\Lambda}^{2,0} = \begin{bmatrix}0\\0\\1\end{bmatrix}, \quad
    \mat{\Lambda}^{1,0} = \begin{bmatrix}1&0\\0&1\\0&0\end{bmatrix}, \quad
    \mat{\Lambda}^{1,1} = \begin{bmatrix}0&1\\-1&0\\0&0\end{bmatrix},
\end{equation}
see also \cref{ex:derham_hodge_3d_tensor_product}.
The columns of each matrix are orthonormal, and their images for fixed total degree $p+q$ are orthogonal and span the corresponding proxy space.
Thus $(\mat{\Lambda}^{p,q})^\transpose$ extracts the 2D proxy from its 3D component.
We write $\wedge$ for these products of elements, since \cref{eqn:tensor_product_complex_spaces_3d_explicit} is the vector proxy of the wedge product of the forms associated with $\logical{V}^p$ and $U^q$, pulled back to the product domain (see \cref{ex:derham_hodge_3d_tensor_product}).
For linear operators $A$ and $B$ acting on the 2D and periodic spaces, respectively, $A\otimes B$ denotes the operator defined by $(A\otimes B)(\logical{V}^p\wedge U^q)\defeq(A\logical{V}^p)\wedge(BU^q)$.

\begin{example}[The $(1,1)$ tensor product]\label{ex:derham_hodge_3d_tensor_product}
  We let $\nu^1 = \logical{\Bspline}^1_\xx \diff \logical{x} + \logical{\Bspline}^1_\yy \diff \logical{y}$ be the $1$-form associated with the vector proxy $\logical{\+\Bspline}^1$ (see also \cref{rem:vector_proxy}), and $\mu^1 = \Bsplineone^1 \diff \logical{z}$ be the $1$-form associated with the vector proxy $\Bsplineone^1$.
  Then the $(1,1)$ product $\nu^1 \wedge \mu^1$ is given by
  \begin{equation}
    \nu^1 \wedge \mu^1
    = 
    \logical{\Bspline}^1_\xx \Bsplineone^1 \diff \logical{x} \wedge \diff \logical{z}
    + 
    \logical{\Bspline}^1_\yy \Bsplineone^1 \diff \logical{y} \wedge \diff \logical{z}
    = 
    \logical{\Bspline}^1_\yy \Bsplineone^1 \diff \logical{y} \wedge \diff \logical{z}
    -
    \logical{\Bspline}^1_\xx \Bsplineone^1 \diff \logical{z} \wedge \diff \logical{x}
    ,
  \end{equation}
  by making use of the antisymmetry of the wedge product: $\diff \logical{x} \wedge \diff \logical{z} = -\diff \logical{z} \wedge \diff \logical{x}$.
  The vector proxy of the resulting $2$-form $\omega^2 = \logical{\Bsplinethree}^2_\xx \diff \logical{y} \wedge \diff \logical{z} + \logical{\Bsplinethree}^2_\yy \diff \logical{z} \wedge \diff \logical{x} + \logical{\Bsplinethree}^2_\zz \diff \logical{x} \wedge \diff \logical{y}$ is given exactly by the expression for $\logical{\+\Bspline}^1 \wedge \Bsplineone^1$ in \cref{eqn:tensor_product_complex_spaces_3d_explicit}.
\end{example}

The resulting polar $m$-form spaces are subspaces of the corresponding tensor-product spaces defined by
\begin{equation}
  \logical{\funspace{\Bsplinethree}}^{m}
    \defeq\bigoplus_{p+q=m}
    \logical{\funspace{\Bspline}}^{p} \otimes \funspace{\Bsplineone}^{q,\zz}
  , \quad
  \intermediate{\funspace{\Bsplinethree}}^m \defeq \logical{\pushforward}^m \logical{\funspace{\Bsplinethree}}^{m}
  , \quad
  \funspace{\Bsplinethree}^m \defeq \intermediate{\pushforward}^m \intermediate{\funspace{\Bsplinethree}}^{m} = \pushforward^m \logical{\funspace{\Bsplinethree}}^{m}
  , \quad
  m \in \{0,1,2,3\}.
\end{equation}


\subsubsection{Choice of basis functions}\label{sec:derham_hodge_3d_basis}
We extend the discussion found in \cref{sec:derham_hodge_2d_basis} to the 3D case.
We use the tensor-product structure of the 3D polar spaces from \cref{eqn:tensor_product_complex_spaces_3d_explicit} and take the products of the 2D polar basis functions with the 1D periodic basis functions.
In particular, the 2D $2$-form factors use the compactly supported basis in \cref{eqn:polar_twoform_curl_basis}, and the curl-like 2D $1$-form factors use its homotopy lift.
As these products are defined in logical coordinates, we define the 3D polar basis functions in logical coordinates and then push them forward to physical coordinates using the pushforward operators from \cref{eqn:pushforward_operators_3d}.
The matrix representations of the bases, which are the extraction matrices of the 3D polar spaces, are Kronecker products of the 2D polar extraction matrices and the 1D periodic extraction matrices, see \cref{eqn:matrep_extraction_3d}.
We write out the $1$-form case, which contains both bidegrees $(1,0)$ and $(0,1)$; the other degrees follow the same pattern and are summarized in \cref{tab:polar_basis_functions_3d}.

\paragraph{Polar $1$-form basis functions}
There are two types of 3D polar $1$-form basis functions, which we denote by $\logicalpolar{\+\Bsplinethree}^{(1,0),\pr}$ and $\logicalpolar{\+\Bsplinethree}^{(0,1),\pr}$, respectively.
They are defined as
\begin{equation}
  \logicalpolar{\+\Bsplinethree}^{(1,0),\pr}_{l,k}
  \defeq\mat{\Lambda}^{1,0}\logicalpolar{\+\Bspline}^{1,\pr}_l\Bsplineone^{0,\zz}_k
  , \quad
  l \in \{1, \ldots, n^{1,\pr}_\mathrm{2D}\}, \quad
  k \in \indexone^{0,\zz},
\end{equation}
and
\begin{equation}
  \logicalpolar{\+\Bsplinethree}^{(0,1),\pr}_{l,k}
  \defeq\mat{\Lambda}^{0,1}\logicalpolar{\Bspline}^{0,\pr}_l\Bsplineone^{1,\zz}_k
  , \quad
  l \in \{1, \ldots, n^{0,\pr}_\mathrm{2D}\}, \quad
  k \in \indexone^{1,\zz},
\end{equation}
where we now linearly index the internal 2D polar $1$-form basis functions defined in \cref{eqn:polar_oneform_gradient_basis,eqn:polar_oneform_curl_basis_definition}, with $n^{1,\pr}_\mathrm{2D} = n^{1,\pr}_{\mathrm{2D},\mathrm{G}} + n^{1,\pr}_{\mathrm{2D},\mathrm{C}}$.
The remaining external basis functions are defined as the usual tensor-product basis functions
\begin{equation}
  \logical{\+\Bsplinethree}^1_{\alpha,ijk}
  \defeq\mat{\Lambda}^{1,0}\logical{\+\Bspline}^1_{\alpha,ij}\Bsplineone^{0,\zz}_k,
  \quad\alpha\in\{\xx,\yy\},
  \qquad
  \logical{\+\Bsplinethree}^1_{\zz,ijk}
  \defeq\mat{\Lambda}^{0,1}\logical{\Bspline}^0_{ij}\Bsplineone^{1,\zz}_k.
\end{equation}
It follows that any $1$-form $\+A \in \polar{\funspace{\Bsplinethree}}^{1,\pr}$ can be expressed as
\begin{equation}
  \+A = \sum_{k\in\indexone^{0,\zz}} \sum_{l=1}^{n^{1,\pr}_\mathrm{2D}} a^{(1,0)}_{l,k} \polar{\+\Bsplinethree}^{(1,0),\pr}_{l,k}
  + 
  \sum_{k\in\indexone^{1,\zz}} \sum_{l=1}^{n^{0,\pr}_\mathrm{2D}} a^{(0,1)}_{l,k} \polar{\+\Bsplinethree}^{(0,1),\pr}_{l,k}
  + 
  \sum_{\alpha\in\{\xx,\yy,\zz\}}\sum_{(i,j,k) \in \indexthree^1_{\alpha,\external}} a_{\alpha,ijk} \+\Bsplinethree^1_{\alpha,ijk}
  ,
\end{equation}
where $\indexthree^1_{\alpha,\external} \defeq \indextwo^1_{\alpha,\external} \times \indexone^{0,\zz}$ for $\alpha\in\{\xx,\yy\}$, and $\indexthree^1_{\zz,\external} \defeq \indextwo^0_\external \times \indexone^{1,\zz}$.

\paragraph{The other degrees}
The polar $0$-, $2$- and $3$-form basis functions are obtained in the same way.
For each bidegree $(p,q)$ with $p+q=m$, the internal basis functions are $\mat{\Lambda}^{p,q}$ applied to the products of the internal 2D polar basis functions of degree $p$ with the toroidal basis functions of degree $q$, the external basis functions are the corresponding products of the external 2D tensor-product basis functions, and a polar $m$-form is the sum of the resulting internal and external expansions over the bidegrees $p+q=m$; for the bidegree $(1,1)$, $\mat{\Lambda}^{1,1}$ swaps the two components of the 2D $1$-form and changes one sign, see \cref{ex:derham_hodge_3d_tensor_product}.
\Cref{tab:polar_basis_functions_3d} lists the factors for all six bidegrees, where $n^{0,\pr}_\mathrm{2D}=n^{0,\pr}$ and $n^{2,\pr}_\mathrm{2D}=n^{2,\pr}$ denote the numbers of internal 2D polar $0$- and $2$-form basis functions from \cref{eqn:polar_subspaces_overview_dimensions}.
In the matrix representation, the coefficients are ordered by bidegree blocks, see \cref{sec:matrix_representation_3d}.

\begin{table}
  \centering
  \footnotesize
  \setlength{\tabcolsep}{2.5pt}
  \begin{tabular}{@{}cclllll@{}}
    \toprule
    $(p,q)$ & $m$ & internal 2D factor & number & external 2D factor & toroidal factor & 3D components\\
    \midrule
    $(0,0)$ & $0$ & $\logicalpolar{\Bspline}^{0,\pr}_{l}$ & $n^{0,\pr}_\mathrm{2D}$ & $\logical{\Bspline}^0_{ij}$, $(i,j)\in\indextwo^{0,\pr}_\external$ & $\Bsplineone^{0,\zz}_k$ & scalar\\
    $(1,0)$ & $1$ & $\logicalpolar{\+\Bspline}^{1,\pr}_{\mathrm{G},l}$, $\logicalpolar{\+\Bspline}^{1,\pr}_{\mathrm{C},l}$ & $n^{1,\pr}_\mathrm{2D}$ & $\logical{\+\Bspline}^1_{\alpha,ij}$, $(i,j)\in\indextwo^{1,\pr}_{\alpha,\external}$ & $\Bsplineone^{0,\zz}_k$ & $\xx$, $\yy$\\
    $(0,1)$ & $1$ & $\logicalpolar{\Bspline}^{0,\pr}_{l}$ & $n^{0,\pr}_\mathrm{2D}$ & $\logical{\Bspline}^0_{ij}$, $(i,j)\in\indextwo^{0,\pr}_\external$ & $\Bsplineone^{1,\zz}_k$ & $\zz$\\
    $(2,0)$ & $2$ & $\logicalpolar{\Bspline}^{2,\pr}_{l}$ & $n^{2,\pr}_\mathrm{2D}$ & $\logical{\Bspline}^2_{ij}$, $(i,j)\in\indextwo^{2,\pr}_\external$ & $\Bsplineone^{0,\zz}_k$ & $\zz$\\
    $(1,1)$ & $2$ & $\logicalpolar{\+\Bspline}^{1,\pr}_{\mathrm{G},l}$, $\logicalpolar{\+\Bspline}^{1,\pr}_{\mathrm{C},l}$ & $n^{1,\pr}_\mathrm{2D}$ & $\logical{\+\Bspline}^1_{\alpha,ij}$, $(i,j)\in\indextwo^{1,\pr}_{\alpha,\external}$ & $\Bsplineone^{1,\zz}_k$ & $\xx$, $\yy$\\
    $(2,1)$ & $3$ & $\logicalpolar{\Bspline}^{2,\pr}_{l}$ & $n^{2,\pr}_\mathrm{2D}$ & $\logical{\Bspline}^2_{ij}$, $(i,j)\in\indextwo^{2,\pr}_\external$ & $\Bsplineone^{1,\zz}_k$ & scalar\\
    \bottomrule
  \end{tabular}
  \caption{The 3D polar basis functions by bidegree $(p,q)$, generalizing the $1$-form case written out in the text: the internal 2D factor from \cref{eqn:polar_zeroform_basis,eqn:polar_oneform_gradient_basis,eqn:polar_oneform_curl_basis_definition,eqn:polar_twoform_curl_basis}, its number, the external 2D factor with its index set from \cref{eqn:polar_zeroform_space_odot,eqn:polar_mform_space_odot}, the toroidal factor, and the components of the 3D vector proxy in which $\mat{\Lambda}^{p,q}$ from \cref{eqn:tensor_product_proxy_matrices_3d} places the product.
  }
  \label{tab:polar_basis_functions_3d}
\end{table}

\subsection{Properties of the polar spaces}\label{sec:derham_hodge_3d_properties}
The construction of the 3D polar subspaces as a total complex of the 2D polar subspaces and the 1D periodic complex automatically yields several properties, which we formalize in the following theorem.
\begin{theorem}[Cochain complex and cohomology of the 3D polar subspaces]\label{thm:polar_3d_subspaces_commutativity_conformity}
  Let $\pr\ge0$.
  The polar spaces defined in \cref{eqn:tensor_product_complex_spaces_3d,eqn:polar_spaces_3d_pushforward} form a cochain complex
  \begin{equation}
    \gradientt \polar{\funspace{\Bsplinethree}}^{0,\pr} \subset \polar{\funspace{\Bsplinethree}}^{1,\pr}
    , \quad
    \curlt \polar{\funspace{\Bsplinethree}}^{1,\pr} \subset \polar{\funspace{\Bsplinethree}}^{2,\pr}
    , \quad
    \divergencet \polar{\funspace{\Bsplinethree}}^{2,\pr} \subset \polar{\funspace{\Bsplinethree}}^{3,\pr}.
  \end{equation}
  Their cohomology has dimensions
  \begin{equation}
    \dim (\kernel \gradientt) =
    \dim (\kernel \curlt / \image \gradientt) = 1
    , \quad
    \dim (\kernel \divergencet / \image \curlt) =
    \dim (\polar{\funspace{\Bsplinethree}}^{3,\pr} / \image \divergencet) = 0,
  \end{equation}
  agreeing with the continuous de Rham complex on the solid torus $\Omega$.
  At $\pr=0$, these are algebraic coefficient statements; for $\pr\ge1$, conformity in the physical Hilbert complex is established in \cref{thm:polar_subspaces_conformity_3d}.
\end{theorem}
\begin{statementproof}[appendix]{thm:polar_3d_subspaces_commutativity_conformity}
  The logical 2D polar spaces form a cochain complex by \cref{thm:polar_subspaces_commutativity_conformity} and the commutation of the 2D pushforwards with the differential operators, as do the 1D periodic spaces by the derivative relations in \cref{eqn:bspline_derivative,eqn:trigonometric_derivative}.
  Hence their tensor product \cref{eqn:tensor_product_complex_spaces_3d} is a cochain complex \citep[Sec.~3]{arnold2015curvilinear}.
  In logical coordinates, let $\logical{\mathrm d}_\mathrm{2D}$ denote the 2D gradient in degree zero and curl in degree one, and let $\mathrm d_\zz$ denote the periodic derivative $\partial_\zz$ in degree zero. Both vanish in their respective top degrees.
  For $\logical{\Bspline}^p\in\logicalpolar{\funspace{\Bspline}}^{p,\pr}$ and $\Bsplineone^q\in\funspace{\Bsplineone}^{q,\zz}$, the tensor-product differential is the exterior derivative of the product form \citep[eq.~(4)]{arnold2015curvilinear},
  \begin{equation}\label{eqn:logical_3d_tensor_product_differential}
    \logical{\mathrm d}(\logical{\Bspline}^p\wedge\Bsplineone^q)
    =\logical{\mathrm d}_\mathrm{2D}\logical{\Bspline}^p\wedge\Bsplineone^q
     +(-1)^p\logical{\Bspline}^p\wedge\mathrm d_\zz\Bsplineone^q.
  \end{equation}
  Under the vector-proxy identifications in \cref{eqn:tensor_product_complex_spaces_3d_explicit}, $\logical{\mathrm d}$ is the logical 3D gradient, curl, or divergence in the corresponding degree.
  By \cref{eqn:polar_spaces_3d_pushforward}, the physical polar forms are the pushforwards $\pushforward^m$ of the logical ones with the same coefficients, and $\pushforward^m=\intermediate{\pushforward}^m\circ\logical{\pushforward}^m$ commutes with the differential operators on $\Omega\setminus\Gamma$, where $\mapping$ is a diffeomorphism (\cref{eqn:inverse_3d_via_inverse_2d}); hence the physical gradient, curl and divergence act on the coefficients as the logical ones, for every intermediate mapping.

  For the cohomology we use the algebraic K\"unneth formula \citep[Thm.~3.6.3]{weibel1994}: both factor complexes consist of finite-dimensional real vector spaces, so the cohomology of their tensor product is the graded tensor product of the factor cohomologies.
  Let $\mathfrak{H}^m$, $\mathfrak{H}_\mathrm{2D}^p$, and $\mathfrak{H}_\zz^q$ denote the cohomology dimensions of the 3D polar spaces, the 2D polar spaces, and the 1D periodic complex, respectively. Taking dimensions gives
  \begin{equation}
    \mathfrak{H}^m
    =\sum_{p+q=m}
    \mathfrak{H}_\mathrm{2D}^p\mathfrak{H}_\zz^q,
    \qquad m\in\{0,1,2,3\},
  \end{equation}
  where the sum is over the degrees of the two factor complexes.
  By \cref{cor:polar_subspaces_cohomology}, $\mathfrak{H}_\mathrm{2D}^0=1$ and $\mathfrak{H}_\mathrm{2D}^1=\mathfrak{H}_\mathrm{2D}^2=0$, with the constant $0$-form representing the remaining class.
  For the periodic B-spline or trigonometric pair in the toroidal direction,
  \begin{equation}
    \kernel\partial_\zz=\spn\{1\},
    \qquad
    \partial_\zz\funspace{\Bsplineone}^{0,\zz}
    =\setdef{\Bsplineone^1\in\funspace{\Bsplineone}^{1,\zz}}
      {\int_0^1\Bsplineone^1\diff\logical{z}=0}.
  \end{equation}
  Indeed, a periodic scalar with zero derivative is constant, and every periodic derivative has zero integral. Conversely, the primitive of a zero-integral periodic $1$-form is periodic and belongs to the scalar space by the derivative relations in \cref{eqn:bspline_derivative,eqn:trigonometric_derivative}. Since the constant $1$-form has unit integral, it represents the remaining degree-one class. Thus $\mathfrak{H}_\zz^0=\mathfrak{H}_\zz^1=1$, and the dimension formula yields
  \begin{equation}
    \mathfrak{H}^0=\mathfrak{H}^1=1,
    \qquad
    \mathfrak{H}^2=\mathfrak{H}^3=0.
  \end{equation}
\end{statementproof}

Conformity and regularity of the 3D polar subspaces follow almost directly from their 2D counterparts.
By construction, the logical polar forms are sums of products of logical 2D polar forms and periodic forms, and the logical mapping $\logicalmapping=\logicalmapping_\mathrm{2D}\times\mathrm{id}$ acts on the cross-section only.
Hence, on the intermediate cylinder, the polar forms remain sums of products, now of 2D polar forms on $\intermediate{\Omega}_\mathrm{2D}$, which are conforming and regular at the center by the 2D results, and of periodic forms in $\intermediatecoord{z}$; such products inherit the conformity and regularity of their factors.
The only remaining step is the intermediate mapping $\intermediatemapping$, a smooth diffeomorphism that introduces no singularity and therefore preserves both properties.
We make both steps precise in \cref{lem:polar_3d_product_structure} below, and then use it to show that the 3D polar subspaces are conforming subspaces of the infinite-dimensional function spaces, and moreover that they satisfy the desired regularity conditions on the central curve $\Gamma$ of the domain.

Let
\begin{equation}
  \intermediatepolar{\funspace{\Bspline}}^{p,\pr}_\mathrm{2D}
  \defeq\logical{\pushforward}^p_\mathrm{2D}\logicalpolar{\funspace{\Bspline}}^{p,\pr}
  , \quad p\in\{0,1,2\},
\end{equation}
denote the intermediate 2D polar spaces
For a 2D $p$-form proxy $\intermediate{V}^p$ on $\intermediate{\Omega}_\mathrm{2D}$ and $U^q\in\funspace{\Bsplineone}^{q,\zz}$, we define $\intermediate{V}^p\wedge U^q\defeq\mat{\Lambda}^{p,q}\intermediate{V}^p U^q$ as a function of $(\intermediatecoord{\+x}_\mathrm{2D},\intermediatecoord{z})$, analogously to \cref{eqn:tensor_product_complex_spaces_3d_explicit}.
We write $\funspace{C}^k_\Gamma$, for $k\ge0$, for scalar functions and vector fields that are $k$ times continuously differentiable on a physical neighborhood of the entire central curve $\Gamma$.
\begin{lemma}[Product structure of the 3D polar spaces]\label{lem:polar_3d_product_structure}
  For $p\in\{0,1,2\}$ and $q\in\{0,1\}$, let $\logical{\Bspline}^p\in\logical{\funspace{\Bspline}}^{p}$ and $\Bsplineone^q\in\funspace{\Bsplineone}^{q,\zz}$, and interpret their product in the vector-proxy sense of \cref{eqn:tensor_product_complex_spaces_3d_explicit}.
  Then, elementwise on $\intermediate{\Omega}\setminus(\{\+0\}\times\mathbb{S}^1)$
  \begin{equation}\label{eqn:logical_tensor_product_intermediate_product_3d}
    \logical{\pushforward}^{p+q}
    \left(\logical{\Bspline}^p\wedge\Bsplineone^q\right)
    =
    \left(\logical{\pushforward}^p_\mathrm{2D}\logical{\Bspline}^p\right)\wedge\Bsplineone^q.
  \end{equation}
  Consequently, the intermediate polar spaces from \cref{eqn:polar_spaces_3d_pushforward} are sums of products of intermediate 2D polar forms and periodic forms,
  \begin{equation}\label{eqn:intermediate_polar_spaces_3d}
    \intermediatepolar{\funspace{\Bsplinethree}}^{m,\pr}
    =\bigoplus_{p+q=m}
      \intermediatepolar{\funspace{\Bspline}}^{p,\pr}_\mathrm{2D}\otimes\funspace{\Bsplineone}^{q,\zz}
    , \quad
    m\in\{0,1,2,3\}.
  \end{equation}
  The pushforwards $\intermediate{\pushforward}^0,\ldots,\intermediate{\pushforward}^3$ commute with the weak gradient, curl and divergence, and map $\intermediate{\funspace{H}}_1$, $\intermediate{\funspace{H}}_{\curlt}$, $\intermediate{\funspace{H}}_{\divergencet}$, $\intermediate{\funspace{L}}_2$ one-to-one onto $\funspace{H}_1$, $\funspace{H}_{\curlt}$, $\funspace{H}_{\divergencet}$, $\funspace{L}_2$.
  Moreover, $\intermediate{\pushforward}^m$ maps forms that are $\funspace{C}^k$, respectively bounded, on a (punctured) neighborhood of the axis $\{\+0\}\times\mathbb{S}^1$ to forms in $\funspace{C}^k_\Gamma$, respectively bounded on a punctured neighborhood of $\Gamma$.
\end{lemma}
\begin{statementproof}[appendix]{lem:polar_3d_product_structure}
  \emph{Pushforward of products.}
  We consider $\logical{x}>0$, where $\logicalmapping$ is a diffeomorphism, using elementwise evaluations for piecewise-defined factors.
  By \cref{eqn:logical_mapping_3d}, a product $f(\logical{\+x}_\mathrm{2D})\,g(\logical{z})$ satisfies $(fg)\circ\logicalmapping^{-1}=(f\circ\logicalmapping_\mathrm{2D}^{-1})(\intermediatecoord{\+x}_\mathrm{2D})\,g(\intermediatecoord{z})$, and by \cref{eqn:logical_jacobian_3d},
  \begin{equation}\label{eqn:logical_jacobian_inverse_transpose_3d}
    \logical{\jacmat}^{-\transpose}
    =
    \begin{bmatrix}
      \logical{\jacmat}_\mathrm{2D}^{-\transpose} & \+0\\
      \+0^\transpose & 1
    \end{bmatrix}
    , \quad
    \frac{1}{\logical{\jac}}\logical{\jacmat}
    =
    \begin{bmatrix}
      \logical{\jac}_\mathrm{2D}^{-1}\logical{\jacmat}_\mathrm{2D} & \+0\\
      \+0^\transpose & \logical{\jac}_\mathrm{2D}^{-1}
    \end{bmatrix}.
  \end{equation}
  We verify \cref{eqn:logical_tensor_product_intermediate_product_3d} using the matrices $\mat{\Lambda}^{p,q}$ from \cref{eqn:tensor_product_proxy_matrices_3d} and the 2D pushforwards $\logical{\pushforward}^p_\mathrm{2D}$, which are given by \cref{eqn:pushforward_operators} with $\logicalmapping_\mathrm{2D}$, $\logical{\jacmat}_\mathrm{2D}$, $\logical{\jac}_\mathrm{2D}$.
  It suffices to show that, before composition with $\logicalmapping^{-1}$, the factor of $\logical{\pushforward}^{p+q}$ applied to $\mat{\Lambda}^{p,q}\logical{\Bspline}^p\Bsplineone^q$ equals $\mat{\Lambda}^{p,q}$ times the factor of $\logical{\pushforward}^p_\mathrm{2D}$ applied to $\logical{\Bspline}^p$, times $\Bsplineone^q$; the factors are $1$, $\logical{\jacmat}_\mathrm{2D}^{-\transpose}$ and $\logical{\jac}_\mathrm{2D}^{-1}$ in 2D degrees $0$, $1$ and $2$.
  For $(p,q)=(0,0)$ this is trivial, and for $(0,1)$ the last column of $\logical{\jacmat}^{-\transpose}$ is the last unit vector $\mat{\Lambda}^{0,1}$.
  For $(1,0)$, $\logical{\jacmat}^{-\transpose}\mat{\Lambda}^{1,0}\logical{\+\Bspline}^1=\mat{\Lambda}^{1,0}\logical{\jacmat}_\mathrm{2D}^{-\transpose}\logical{\+\Bspline}^1$.
  For $(2,0)$, the last column of $\logical{\jac}^{-1}\logical{\jacmat}$ is $\logical{\jac}_\mathrm{2D}^{-1}\mat{\Lambda}^{2,0}$, and for $(2,1)$ the factor is $\logical{\jac}^{-1}=\logical{\jac}_\mathrm{2D}^{-1}$.
  For $(1,1)$, $\mat{\Lambda}^{1,1}=\begin{bmatrix}\mat{Q}\\\+0^\transpose\end{bmatrix}$ with $\mat{Q}\defeq\begin{bmatrix}0&1\\-1&0\end{bmatrix}$ by \cref{eqn:tensor_product_proxy_matrices_3d}, and
  \begin{equation}
    \frac{1}{\logical{\jac}}\logical{\jacmat}\mat{\Lambda}^{1,1}\logical{\+\Bspline}^1
    =\begin{bmatrix}\logical{\jac}_\mathrm{2D}^{-1}\logical{\jacmat}_\mathrm{2D}\mat{Q}\logical{\+\Bspline}^1\\0\end{bmatrix}
    =\mat{\Lambda}^{1,1}\logical{\jacmat}_\mathrm{2D}^{-\transpose}\logical{\+\Bspline}^1,
  \end{equation}
  by the identity
  \begin{equation}
    \jacmat_\mathrm{2D} \mat{Q} = (\det\jacmat_\mathrm{2D})\, \mat{Q} \jacmat_\mathrm{2D}^{-\transpose},
  \end{equation}
  which holds for every invertible $2\times2$ matrix $\jacmat_\mathrm{2D}$, since $\jacmat_\mathrm{2D}\mat{Q}\jacmat_\mathrm{2D}^\transpose=(\det\jacmat_\mathrm{2D})\mat{Q}$ by direct computation.
  Composing with $\logicalmapping^{-1}$ in all six cases gives \cref{eqn:logical_tensor_product_intermediate_product_3d}, and linearity gives it for finite sums; with \cref{eqn:tensor_product_complex_spaces_3d,eqn:polar_spaces_3d_pushforward}, this proves \cref{eqn:intermediate_polar_spaces_3d}.

  \emph{Intermediate mapping.}
  The intermediate mapping $\intermediatemapping$ is a smooth diffeomorphism up to the boundary that maps the axis $\{\+0\}\times\mathbb{S}^1$ onto $\Gamma$ (\cref{ass:diffeomorphism_3d}), so, unlike $\logicalmapping$, it introduces no singularity.
  Its pushforwards $\intermediate{\pushforward}^m$ are the vector proxies of the pullback of differential forms by $\intermediatemapping^{-1}$, so they commute with the gradient, curl and divergence \citep[Sec.~2.2]{Arnold_Falk_Winther_2006}, and, because $\intermediate{\jacmat}$, $\intermediate{\jac}$ and their inverses are bounded, they change $\funspace{L}_2$ norms by at most a bounded factor.
  Hence $\intermediate{\pushforward}^0,\ldots,\intermediate{\pushforward}^3$ map $\intermediate{\funspace{H}}_1$, $\intermediate{\funspace{H}}_{\curlt}$, $\intermediate{\funspace{H}}_{\divergencet}$, $\intermediate{\funspace{L}}_2$ one-to-one onto $\funspace{H}_1$, $\funspace{H}_{\curlt}$, $\funspace{H}_{\divergencet}$, $\funspace{L}_2$, with inverses given by the pushforwards of $\intermediatemapping^{-1}$; for $\intermediate{\pushforward}^1$ and $\intermediate{\pushforward}^2$, the covariant and contravariant Piola transformations, see \citet[Sec.~2]{rognes2009}.
  Similarly, composition with $\intermediatemapping^{\pm1}$ and multiplication by the smooth invertible factors $\intermediate{\jacmat}^{\pm\transpose}$, $\intermediate{\jacmat}^{\pm1}$ and $\intermediate{\jac}^{\pm1}$ preserve boundedness and $k$-fold continuous differentiability; hence $\intermediate{\pushforward}^m$ maps forms that are $\funspace{C}^k$, respectively bounded, on a (punctured) neighborhood of the axis to forms in $\funspace{C}^k_\Gamma$, respectively bounded on a punctured neighborhood of $\Gamma$.
\end{statementproof}

\begin{theorem}[Conformity of the 3D polar subspaces]\label{thm:polar_subspaces_conformity_3d}
  We assume $\pr \ge 1$.
  Then the 3D polar subspaces are conforming subspaces of the infinite-dimensional function spaces, i.e.,
  \begin{equation}
    \polar{\funspace{\Bsplinethree}}^{0,\pr} \subset \funspace{H}_1, \quad
    \polar{\funspace{\Bsplinethree}}^{1,\pr} \subset \funspace{H}_{\curlt}, \quad
    \polar{\funspace{\Bsplinethree}}^{2,\pr} \subset \funspace{H}_{\divergencet}, \quad
    \polar{\funspace{\Bsplinethree}}^{3,\pr} \subset \funspace{L}_2.
  \end{equation}
\end{theorem}
\begin{statementproof}[appendix]{thm:polar_subspaces_conformity_3d}
  By \cref{eqn:polar_spaces_3d_pushforward} and \cref{lem:polar_3d_product_structure}, every physical polar form is the pushforward by $\intermediate{\pushforward}^m$ of the intermediate polar form with the same coefficients, and $\intermediate{\pushforward}^m$ maps the intermediate Hilbert spaces onto the physical ones.
  Hence it suffices to prove the inclusions for the intermediate polar spaces, which by \cref{eqn:intermediate_polar_spaces_3d} consist of finite sums of products $\intermediate{V}^p\wedge U^q$ with $\intermediate{V}^p\in\intermediatepolar{\funspace{\Bspline}}^{p,\pr}_\mathrm{2D}$ and $U^q\in\funspace{\Bsplineone}^{q,\zz}$.

  By \cref{thm:polar_subspaces_conformity}, applied to the 2D construction with the logical mapping $\logicalmapping_\mathrm{2D}$ and the identity as 2D intermediate mapping, the intermediate 2D polar spaces $\intermediatepolar{\funspace{\Bspline}}^{p,\pr}_\mathrm{2D}$ belong to $\funspace{H}_1$, $\funspace{H}_{\curlt}$, and $\funspace{L}_2$ on $\intermediate{\Omega}_\mathrm{2D}$ in degrees zero, one, and two, respectively.
  The periodic scalar and $1$-form spaces belong to $\funspace{H}_1^\zz$ and $\funspace{L}_2^\zz$, respectively, and the periodic derivative maps the scalar space into the $1$-form space.

  Fubini's theorem shows that the products of these factors, with the vector-proxy identifications $\intermediate{V}^p\wedge U^q=\mat{\Lambda}^{p,q}\intermediate{V}^pU^q$, are square-integrable on $\intermediate{\Omega}$.
  Integration by parts in each factor, against test functions on $\intermediate{\Omega}$ restricted to fixed $\intermediatecoord{z}$ or fixed $\intermediatecoord{\+x}_\mathrm{2D}$, gives the tensor-product differential rule in \cref{eqn:logical_3d_tensor_product_differential}, now using the weak derivatives of the intermediate cross-sectional factors.
  Its terms contain the 2D weak gradient or curl when $p<2$, and the toroidal weak derivative when $q=0$.
  Every term is again a product of square-integrable factors and hence belongs to $\intermediate{\funspace{L}}_2$.
  Consequently, the finite sums of these products belong to $\intermediate{\funspace{H}}_1$, $\intermediate{\funspace{H}}_{\curlt}$, $\intermediate{\funspace{H}}_{\divergencet}$, and $\intermediate{\funspace{L}}_2$ in degrees zero through three, which proves the four inclusions.
\end{statementproof}

\begin{theorem}[Regularity of the 3D polar subspaces]
  \label{thm:polar_subspaces_regularity_3d}
  We assume $1\le\pr\le r^{0,\xx}$ and that the 2D logical mapping $\logicalmapping_\mathrm{2D}$ and the 2D logical $m$-form spaces $\logical{\funspace{\Bspline}}^m$ are $\pr$-compatible in the sense of \cref{def:pr_compatibility_mapping_logical_spaces}.
  For $q\in\{0,1\}$, assume that the toroidal spaces satisfy $\funspace{\Bsplineone}^{q,\zz}\subset\funspace{C}^{\pr-q}$ and, if $\pr<r^{0,\xx}$, that the angular scalar basis is $\funspace{C}^{\pr}$ and the angular $1$-form basis is $\funspace{C}^{\pr-1}$, as in \cref{thm:polar_subspaces_regularity}.
  Then, for $\pr\ge2$, the following regularity near $\Gamma$ is guaranteed:
  \begin{equation}\label{eqn:polar_subspaces_regularity_3d}
    \polar{\funspace{\Bsplinethree}}^{0,\pr} \subset \funspace{C}_\Gamma^{\pr}, \quad
    \polar{\funspace{\Bsplinethree}}^{1,\pr} \subset \funspace{C}_\Gamma^{\pr-1}, \quad
    \polar{\funspace{\Bsplinethree}}^{2,\pr} \subset \funspace{C}_\Gamma^{\pr-2}, \quad
    \polar{\funspace{\Bsplinethree}}^{3,\pr} \subset \funspace{C}_\Gamma^{\pr-2}.
  \end{equation}
  For $\pr=1$, the first two inclusions hold, while the polar $2$- and $3$-form proxies are bounded on a punctured neighborhood of $\Gamma$.
\end{theorem}
\begin{statementproof}[appendix]{thm:polar_subspaces_regularity_3d}
  By \cref{eqn:polar_spaces_3d_pushforward}, every physical polar form is $\intermediate{\pushforward}^m$ of the intermediate polar form with the same coefficients, $\pr$-compatibility involves only $\logicalmapping_\mathrm{2D}$, and, by \cref{lem:polar_3d_product_structure}, $\intermediate{\pushforward}^m$ maps forms that are $\funspace{C}^k$, respectively bounded, near the axis $\{\+0\}\times\mathbb{S}^1$ to forms that are $\funspace{C}^k$, respectively bounded, near $\Gamma$.
  Hence it suffices to prove the inclusions for the intermediate polar spaces \cref{eqn:intermediate_polar_spaces_3d}, with $\Gamma$ replaced by the axis.

  Fix $p\in\{0,1,2\}$ with $p\le\pr$ and $q\in\{0,1\}$.
  By \cref{thm:polar_subspaces_regularity}, applied to the 2D construction with the logical mapping $\logicalmapping_\mathrm{2D}$ and the identity as 2D intermediate mapping, and since $\intermediatepolar{\funspace{\Bspline}}^{p,\pr}_\mathrm{2D}$ is finite-dimensional, there is a ball $\mathbb{B}\subset\intermediate{\Omega}_\mathrm{2D}$ around $\+0$ on which every intermediate 2D polar $p$-form $\intermediate{V}^p\in\intermediatepolar{\funspace{\Bspline}}^{p,\pr}_\mathrm{2D}$ is $\funspace{C}^{\pr-p}$.
  The periodic factor satisfies $U^q\in\funspace{C}^{\pr-q}(\mathbb{S}^1)$.
  Since $\intermediate{V}^p\wedge U^q=\mat{\Lambda}^{p,q}\intermediate{V}^p(\intermediatecoord{\+x}_\mathrm{2D})\,U^q(\intermediatecoord{z})$ with a constant matrix $\mat{\Lambda}^{p,q}$, the product rule gives $\funspace{C}^{\min\{\pr-p,\pr-q\}}$ regularity on the neighborhood $\mathbb{B}\times\mathbb{S}^1$ of the axis, that is,
  \begin{equation}\label{eqn:regularity_bidegree_3d}
    \logical{\pushforward}^{p+q}
    \left(
      \logicalpolar{\funspace{\Bspline}}^{p,\pr}
      \otimes\funspace{\Bsplineone}^{q,\zz}
    \right)
    =
    \intermediatepolar{\funspace{\Bspline}}^{p,\pr}_\mathrm{2D}
      \otimes\funspace{\Bsplineone}^{q,\zz}
    \subset
    \funspace{C}^{\min\{\pr-p,\pr-q\}}(\mathbb{B}\times\mathbb{S}^1)
    =\funspace{C}^{\pr-\max\{p,q\}}(\mathbb{B}\times\mathbb{S}^1),
  \end{equation}
  where the equality is \cref{eqn:logical_tensor_product_intermediate_product_3d} from \cref{lem:polar_3d_product_structure}, which holds off the axis and, since the right-hand side is continuous on $\mathbb{B}\times\mathbb{S}^1$, extends pointwise to the axis.
  For $\pr\ge2$, this covers all six bidegrees.
  In degree two, the $(2,0)$ part has regularity $\pr-2$ and the $(1,1)$ part has regularity $\pr-1$; degree three contains only the $(2,1)$ part, also with regularity $\pr-2$.
  Taking finite sums in the total-complex decomposition \cref{eqn:intermediate_polar_spaces_3d} proves \cref{eqn:polar_subspaces_regularity_3d} for the intermediate polar spaces, and hence for the physical ones.

  For $\pr=1$, the same argument proves the first two inclusions and gives $\funspace{C}^0$ regularity near the axis for the $(1,1)$ part.
  The intermediate 2D polar $2$-form proxies are bounded on $\mathbb{B}\setminus\{\+0\}$ for a ball $\mathbb{B}$ around $\+0$, by \cref{thm:polar_subspaces_regularity}.
  Multiplication by the constant matrices $\mat{\Lambda}^{2,q}$ and the continuous, hence bounded, periodic factors preserves these bounds on $(\mathbb{B}\setminus\{\+0\})\times\mathbb{S}^1$ for the $(2,0)$ and $(2,1)$ parts.
  Together with the continuous $(1,1)$ part, this proves the remaining endpoint claim, again transported to $\Omega$ by $\intermediate{\pushforward}^m$.
\end{statementproof}


\subsection{Homotopy operators}\label{sec:derham_hodge_3d_homotopy}
Analogously to the total-complex construction of the 3D polar subspaces given in \cref{eqn:tensor_product_complex_spaces_3d}, as well as the construction of the 3D differential operators from the 2D and 1D differential operators, we construct the 3D homotopy operators from the 2D and 1D homotopy operators, so that the polar complex of the solid torus is again a complex in both directions.
This completes the theory in the same way as in 2D: the homotopy identities make the cohomology of \cref{thm:polar_3d_subspaces_commutativity_conformity} explicit, with scalar and vector potentials given by formulas rather than by solving a problem; the two nonzero cohomology projections have a physical meaning, the value at a point of the axis and the circulation along the axis, which in plasma physics is the loop voltage; and the 3D polar projections commute with the homotopy operators and preserve this circulation (\cref{sec:derham_hodge_3d_commuting_projections}).
Unlike in 2D, the homotopy operators are not needed to define the 3D polar spaces, and the details of their construction are deferred to the proofs in the appendix.

The logical 3D homotopy operators can be expressed in terms of radial integrals, as well as a toroidal integral stemming from the homotopy operator of the 1D periodic complex.
Radial integration collapses each cross-section to its center, after which the zero-mean part is integrated along the central circle.
Consequently, a curl-free polar $1$-form has a scalar potential precisely when its circulation along the central curve $\Gamma$ vanishes, while every divergence-free polar $2$-form has a vector potential.
We normalize the toroidal primitive at the first toroidal Greville point $\xi^\zz_1$ from \cref{eqn:geometric_dofs_1d}, using periodic extensions in the toroidal integrals.
\begin{definition}[Logical 3D homotopy operators]\label{def:logical_3d_homotopy_vector_proxies}
  Let $\logical{\gradientt}^\dagger_\mathrm{2D}$ and
  $\logical{\curlt}^\dagger_\mathrm{2D}$ denote the logical 2D homotopies
  from \cref{lem:pseudo_inverse_curl_logical} for the logical mapping $\logicalmapping_\mathrm{2D}$, acting on
  $(\logical{x},\logical{y})$ at each fixed $\logical{z}$.
  Using the matrices from \cref{eqn:tensor_product_proxy_matrices_3d},
  define the logical 3D homotopies by
  \begin{subequations}\label{eqn:logical_3d_homotopy_vector_proxies}
    \begin{align+}
      \logical{\gradientt}^\dagger\logical{\+A}
      &\defeq
          \logical{\gradientt}^\dagger_\mathrm{2D}
            (\mat{\Lambda}^{1,0})^\transpose\logical{\+A}
          +\int_{\xi^\zz_1}^{\logical{z}}
            \left(
              \logical{A}_\zz(0,\logical{y},s)
              -\int_0^1
                \logical{A}_\zz(0,\logical{y},\zeta)\diff\zeta
            \right)\diff s,
        \label{eqn:logical_3d_gradient_homotopy_vector_proxy}\\
      \logical{\curlt}^\dagger\logical{\+B}
      &\defeq
        \mat{\Lambda}^{1,0}\logical{\curlt}^\dagger_\mathrm{2D}
          (\mat{\Lambda}^{2,0})^\transpose\logical{\+B}
        +\mat{\Lambda}^{0,1}\logical{\gradientt}^\dagger_\mathrm{2D}
          (\mat{\Lambda}^{1,1})^\transpose\logical{\+B},
        \label{eqn:logical_3d_curl_homotopy_vector_proxy}\\
      \logical{\divergencet}^\dagger\logical{Q}
      &\defeq
        \mat{\Lambda}^{1,1}\logical{\curlt}^\dagger_\mathrm{2D}
          \logical{Q}.
        \label{eqn:logical_3d_divergence_homotopy_vector_proxy}
    \end{align+}
  \end{subequations}
  The two nonzero logical cohomology projections are
  \begin{equation}\label{eqn:logical_3d_cohomology_projector_vector_proxies}
    \logical{\mathfrak{P}}\logical{P}
    =\logical{P}(0,\logical{y},\xi^\zz_1),
    \qquad
    \logical{\mathfrak{A}}\logical{\+A}
    =\mat{\Lambda}^{0,1}\int_0^1
      \logical{A}_\zz(0,\logical{y},s)\diff s.
  \end{equation}
  The center traces $\logical{A}_\zz(0,\logical{y},\logical{z})$ and
  $\logical{P}(0,\logical{y},\logical{z})$ are the restrictions to $\logical{x}=0$ and are independent of
  $\logical{y}$ for polar $1$- and $0$-forms, respectively.
\end{definition}
Their matrix representations are given in \cref{eqn:matrep_homotopy_3d}.

With the logical 3D homotopy operators, we define the physical 3D homotopy operators by composition with the 3D Piola pushforwards and pullbacks
\begin{subequations}\label{eqn:physical_3d_tensor_product_homotopy}
  \begin{equation}
    \gradientt^\dagger
    \defeq\pushforward^0\circ\logical{\gradientt}^\dagger
      \circ\pushforward^{-1}, \quad
    \curlt^\dagger
    \defeq\pushforward^1\circ\logical{\curlt}^\dagger
      \circ\pushforward^{-2}, \quad
    \divergencet^\dagger
    \defeq\pushforward^2\circ\logical{\divergencet}^\dagger
      \circ\pushforward^{-3}.
  \end{equation}
\end{subequations}
The associated physical cohomology projections are
\begin{equation}\label{eqn:physical_3d_cohomology_projectors}
  \mathfrak{P}
  \defeq\pushforward^0\circ\logical{\mathfrak{P}}
    \circ \pushforward^{-0},
  \qquad
  \mathfrak{A}
  \defeq\pushforward^1\circ\logical{\mathfrak{A}}
    \circ \pushforward^{-1}.
\end{equation}
By \cref{eqn:pushforward_composition_3d}, the physical operators are the transports by $\intermediatemapping$ of the intermediate operators on the cylinder, for example $\gradientt^\dagger=\intermediate{\pushforward}^0\circ\intermediate{\gradientt}^\dagger\circ\intermediate{\pushforward}^{-1}$ with $\intermediate{\gradientt}^\dagger\defeq\logical{\pushforward}^0\circ\logical{\gradientt}^\dagger\circ\logical{\pushforward}^{-1}$, and analogously for the other operators.
By \cref{lem:pseudo_inverse_curl_logical}, the intermediate operators combine the 2D Poincar\'e operators \cref{eqn:homotopy_operators} along straight rays in each cross-section $\intermediate{\Omega}_\mathrm{2D}\times\{\intermediatecoord{z}\}$ with a toroidal integral along the axis; the physical operators therefore integrate along the images of these rays under $\intermediatemapping$ and along $\Gamma$.
In particular, $\mathfrak{P}P$ is the constant equal to the value of $P$ at $\intermediatemapping(\+0,\xi^\zz_1)\in\Gamma$, and $\mathfrak{A}\+A$ is the circulation of $\+A$ along $\Gamma$ times the pushforward $\intermediate{\pushforward}^1\mat{\Lambda}^{0,1}$ of the constant toroidal $1$-form $\diff\intermediatecoord{z}$.

The following theorem shows that the physical homotopy operators satisfy the desired homotopy identities, and that they form a chain complex.
\begin{theorem}[Physical homotopy identities]
  \label{thm:physical_3d_homotopy_identities}
We assume $\pr\ge1$.
The physical homotopy operators preserve the physical polar spaces and
satisfy
\begin{equation}
  \gradientt^\dagger
  \polar{\funspace{\Bsplinethree}}^{1,\pr}
  \subset\polar{\funspace{\Bsplinethree}}^{0,\pr},
  \qquad
  \curlt^\dagger
  \polar{\funspace{\Bsplinethree}}^{2,\pr}
  \subset\polar{\funspace{\Bsplinethree}}^{1,\pr},
  \qquad
  \divergencet^\dagger
  \polar{\funspace{\Bsplinethree}}^{3,\pr}
  \subset\polar{\funspace{\Bsplinethree}}^{2,\pr}.
\end{equation}
For $P\in\polar{\funspace{\Bsplinethree}}^{0,\pr}$, $\+A\in\polar{\funspace{\Bsplinethree}}^{1,\pr}$, $\+B\in\polar{\funspace{\Bsplinethree}}^{2,\pr}$, and $Q\in\polar{\funspace{\Bsplinethree}}^{3,\pr}$, they obey the homotopy identities
\begin{subequations}\label{eqn:physical_3d_homotopy_identities_by_degree}
  \begin{align}
    P&=\mathfrak{P}P+\gradientt^\dagger\gradientt P,\\
    \+A&=\mathfrak{A}\+A
      +\gradientt\gradientt^\dagger\+A
      +\curlt^\dagger\curlt\+A,\\
    \+B&=\curlt\curlt^\dagger\+B
      +\divergencet^\dagger\divergencet\+B,\\
    Q&=\divergencet\divergencet^\dagger Q.
  \end{align}
\end{subequations}
Moreover,
\begin{equation}\label{eqn:physical_3d_homotopy_chain_property}
  \gradientt^\dagger\curlt^\dagger=0,
  \qquad
  \curlt^\dagger\divergencet^\dagger=0,
\end{equation}
so the physical homotopy operators form a chain complex.
\end{theorem}
\begin{statementproof}[appendix]{thm:physical_3d_homotopy_identities}
  For a periodic $0$-form proxy
  $P\in\funspace{\Bsplineone}^{0,\zz}$ and $1$-form proxy $A\in\funspace{\Bsplineone}^{1,\zz}$, define the constant periodic forms
  \begin{equation}\label{eqn:logical_1d_cohomology_projectors_3d}
    \mathfrak{P}_\zz P\defeq P(\xi^\zz_1),
    \qquad
    \mathfrak{A}_\zz A\defeq\int_0^1A(s)\diff s,
  \end{equation}
  and define the normalized periodic primitive by
  \begin{equation}\label{eqn:logical_1d_periodic_homotopy_3d}
    (\partial_\zz^\dagger A)(\logical{z})
    \defeq
    \int_{\xi^\zz_1}^{\logical{z}}
      \left(A(s)-\mathfrak{A}_\zz A\right)\diff s.
  \end{equation}
  Its integrand has zero mean, so the primitive is periodic and vanishes at $\logical{z}=\xi^\zz_1$.
  By the derivative relations in \cref{eqn:bspline_derivative,eqn:trigonometric_derivative}, every zero-integral periodic $1$-form has a primitive in $\funspace{\Bsplineone}^{0,\zz}$.
  Since constants belong to the scalar space, the normalized primitive also belongs to that space.
  The fundamental theorem of calculus gives the one-dimensional homotopy identities (cf.\ \cref{eqn:homotopy_identities_1d})
  \begin{subequations}\label{eqn:logical_1d_contraction_identity_3d}
    \begin{align}
      P &=\mathfrak{P}_\zz P
        +\partial_\zz^\dagger\partial_\zz P,\\
      A&=\mathfrak{A}_\zz A
        +\partial_\zz\partial_\zz^\dagger A.
    \end{align}
  \end{subequations}
  For later use, on periodic scalars in $\funspace{H}_1^\zz$ and $1$-form proxies in $\funspace{L}_2^\zz$, the same formulas define operators compatible with the 1D projectors:
  \begin{equation}\label{eqn:logical_1d_homotopy_commutativity_3d}
    \partial_\zz^\dagger\pi^{1,\zz}
    =\pi^{0,\zz}\partial_\zz^\dagger,
    \qquad
    \mathfrak{P}_\zz\pi^{0,\zz}=\mathfrak{P}_\zz,
    \qquad
    \mathfrak{A}_\zz\pi^{1,\zz}=\mathfrak{A}_\zz.
  \end{equation}
  The scalar identity follows from preservation of the value GDOF at $\xi^\zz_1$.
  The integral GDOFs cover a complete period, so their sum proves the mean-preservation identity.
  The two sides of the homotopy identity then have the same weak derivative by \cref{eqn:stokes_1d} and preservation of constants, and both vanish at $\xi^\zz_1$.

  We use the tensor-product differential from \cref{eqn:logical_3d_tensor_product_differential}.
  Let $\logical{\mathrm d}_\mathrm{2D}^\dagger$ collectively denote the logical 2D
  homotopy operators from \cref{lem:pseudo_inverse_curl_logical}, and set
  \begin{equation}\label{eqn:logical_2d_cohomology_projector_3d}
    \mathfrak{P}_\mathrm{2D}P=P(\+0),
    \qquad
    \mathfrak{P}_\mathrm{2D}^0\defeq\mathfrak{P}_\mathrm{2D},
    \quad
    \mathfrak{P}_\mathrm{2D}^1
    =\mathfrak{P}_\mathrm{2D}^2\defeq0.
  \end{equation}
  The 2D homotopy identities from
  \cref{thm:polar_subspaces_discrete_homotopy} become
  \begin{equation}\label{eqn:logical_2d_contraction_identity_3d}
    \logical{\mathrm d}_\mathrm{2D}\logical{\mathrm d}_\mathrm{2D}^\dagger
    +\logical{\mathrm d}_\mathrm{2D}^\dagger\logical{\mathrm d}_\mathrm{2D}
    =\mathrm I-\mathfrak{P}_\mathrm{2D}.
  \end{equation}

  For $\logical{V}^p\in\logicalpolar{\funspace{\Bspline}}^{p,\pr}$ and $U^q\in\funspace{\Bsplineone}^{q,\zz}$, the standard ordered tensor-product homotopy
  \citep[Chap.~V, Sec.~9, Prop.~9.1, Eq.~(9.3)]{maclane1995} is
  \begin{equation}\label{eqn:logical_3d_tensor_product_homotopy}
    \logical{\mathrm d}^\dagger(\logical{V}^p\wedge U^q)
    \defeq
    \logical{\mathrm d}_\mathrm{2D}^\dagger \logical{V}^p\wedge U^q
    +(-1)^p\mathfrak{P}_\mathrm{2D}^p\logical{V}^p
      \wedge\partial_\zz^\dagger U^q,
  \end{equation}
  where $\partial_\zz^\dagger$ vanishes on $\funspace{\Bsplineone}^{0,\zz}$.  Expanding
  \cref{eqn:logical_3d_tensor_product_homotopy} in total degrees one, two and
  three, using $(\mat{\Lambda}^{p,q})^\transpose$ to extract the 2D
  input proxy and $\mat{\Lambda}^{p-1,q}$ to place its homotopy for
  $p\ge1$, gives exactly
  \cref{eqn:logical_3d_homotopy_vector_proxies}.
  The toroidal term occurs only in total degree one and integrates the
  center trace of $\logical{A}_\zz$ after removing its mean.

  Substitution of \cref{eqn:logical_3d_tensor_product_differential} and
  \cref{eqn:logical_3d_tensor_product_homotopy} shows that the mixed terms involving $\mathrm d_\mathrm{2D}\mathfrak{P}_\mathrm{2D}$ or $\mathfrak{P}_\mathrm{2D}\mathrm d_\mathrm{2D}$ vanish, while the remaining mixed terms cancel by their opposite graded signs.
  By \cref{eqn:logical_2d_contraction_identity_3d,eqn:logical_1d_contraction_identity_3d}, on bidegree $(p,q)$ this gives
  \begin{equation}
    \begin{aligned}
      \logical{\mathrm d}\logical{\mathrm d}^\dagger
      +\logical{\mathrm d}^\dagger\logical{\mathrm d}
      &=(\mathrm I-\mathfrak{P}_\mathrm{2D}^p)\otimes\mathrm I
        +\mathfrak{P}_\mathrm{2D}^p\otimes
        \begin{cases}
          \mathrm I-\mathfrak{P}_\zz,&q=0,\\
          \mathrm I-\mathfrak{A}_\zz,&q=1,
        \end{cases}\\
      &=\mathrm I-\mathfrak{P}_\mathrm{2D}^p\otimes
        \begin{cases}
          \mathfrak{P}_\zz,&q=0,\\
          \mathfrak{A}_\zz,&q=1.
        \end{cases}
    \end{aligned}
  \end{equation}
  The residual operator is $\logical{\mathfrak{P}}$ in total degree zero,
  $\logical{\mathfrak{A}}$ on the $(0,1)$ component in total degree one,
  and zero otherwise.  These are the logical counterparts of
  \cref{eqn:physical_3d_homotopy_identities_by_degree}.  Preservation of the
  tensor-product polar spaces follows from preservation by the factor
  homotopies and projections.

  The chain property follows from the normalized factor homotopies, which satisfy
  $\mathrm d_\mathrm{2D}^\dagger\circ
   \mathrm d_\mathrm{2D}^\dagger=0$, as well as the side conditions
  $\mathrm d_\mathrm{2D}^\dagger\mathfrak{P}_\mathrm{2D}
  =\mathfrak{P}_\mathrm{2D}\mathrm d_\mathrm{2D}^\dagger=0$ and
  $\partial_\zz^\dagger\mathfrak{A}_\zz
  =\mathfrak{P}_\zz\partial_\zz^\dagger=0$.
  The radial and periodic primitives vanish at their respective base points, and mean subtraction annihilates constant periodic $1$-forms.
  Composing \cref{eqn:logical_3d_tensor_product_homotopy} with itself
  therefore gives
  $\logical{\mathrm d}^\dagger\circ\logical{\mathrm d}^\dagger=0$.

  The preceding arguments involve only the logical spaces and hold for every intermediate mapping.
  By \cref{eqn:inverse_3d_via_inverse_2d}, $\mapping$ is a diffeomorphism from $\{\logical{x}>0\}$ onto $\Omega\setminus\Gamma$, where the pushforwards $\pushforward^m=\intermediate{\pushforward}^m\circ\logical{\pushforward}^m$ commute with the differential operators.
  Conjugating the logical identities by them, as in \cref{eqn:physical_3d_tensor_product_homotopy,eqn:physical_3d_cohomology_projectors}, gives
  \cref{eqn:physical_3d_homotopy_identities_by_degree,eqn:physical_3d_homotopy_chain_property} elementwise on $\Omega\setminus\Gamma$.
  By \cref{thm:polar_subspaces_conformity_3d}, these are identities of conforming discrete fields, hence hold almost everywhere on $\Omega$ with weak differential operators.
  Whenever all terms have continuous representatives, the identities also hold on $\Gamma$ by continuity.
  Preservation of the physical polar spaces follows from their definition \cref{eqn:polar_spaces_3d_pushforward} as the pushforwards of the logical polar spaces.
\end{statementproof}

\subsection{Commuting projection operators}\label{sec:derham_hodge_3d_commuting_projections}
We construct the 3D polar projection operators from the 2D operators in \cref{sec:derham_hodge_2d_geom_projections}, using the same tensor products as for the spaces.
The logical 2D coefficient projectors are
\begin{equation}
  \logicalpolar{\Pi}^{p,\pr}_\mathrm{2D}
  \defeq\logical{\pushforward}^{-p}_\mathrm{2D}\circ
    \intermediatepolar{\Pi}^{p,\pr}_\mathrm{2D}\circ\logical{\pushforward}^p_\mathrm{2D},
  \qquad p\in\{0,1,2\},
\end{equation}
where $\intermediatepolar{\Pi}^{p,\pr}_\mathrm{2D}$ are the 2D coefficient projectors of \cref{sec:derham_hodge_2d_geom_projections}, constructed with the logical mapping $\logicalmapping_\mathrm{2D}$ and with the identity as intermediate mapping; they project onto the 2D polar spaces $\logical{\pushforward}^p_\mathrm{2D}\logicalpolar{\funspace{\Bspline}}^{p,\pr}$ on $\intermediate{\Omega}_\mathrm{2D}$.
Their direct sum defines the logical 3D coefficient projector:
\begin{equation}\label{eqn:logical_3d_polar_projection_operator}
  \logicalpolar{\Pi}^{m,\pr}
  \defeq\bigoplus_{p+q=m}
    \logicalpolar{\Pi}^{p,\pr}_\mathrm{2D}\otimes I^{q,\zz},
  \qquad m\in\{0,1,2,3\},
\end{equation}
where $I^{q,\zz}$ is the identity on $\funspace{\Bsplineone}^{q,\zz}$.
Using the matrices from \cref{eqn:tensor_product_proxy_matrices_3d}, its action on $\logical{\Bsplinethree}\in\logical{\funspace{\Bsplinethree}}^m$ is
\begin{equation}\label{eqn:logical_3d_polar_projection_vector_proxy}
  \logicalpolar{\Pi}^{m,\pr}\logical{\Bsplinethree}
  =\sum_{p+q=m}\mat{\Lambda}^{p,q}
    \logicalpolar{\Pi}^{p,\pr}_\mathrm{2D}
    \left((\mat{\Lambda}^{p,q})^\transpose\logical{\Bsplinethree}\right).
\end{equation}
Here the 2D operator acts on $(\logical{x},\logical{y})$ for each fixed $\logical{z}$, or equivalently on each toroidal coefficient.
Thus the same matrices that construct the tensor-product proxies also extract their factors, and all coefficient operations are inherited from the 2D projectors; see \cref{eqn:matrep_extraction_3d} for the matrix representation.
The physical coefficient projectors are formally defined as
\begin{equation}\label{eqn:physical_3d_polar_projection_operator}
  \polar{\Pi}^{m,\pr}
  \defeq\pushforward^m\circ
    \logicalpolar{\Pi}^{m,\pr}\circ\pushforward^{-m}
  =\intermediate{\pushforward}^m\circ
    \intermediatepolar{\Pi}^{m,\pr}\circ\intermediate{\pushforward}^{-m}
  , \quad
  \intermediatepolar{\Pi}^{m,\pr}
  \defeq\logical{\pushforward}^m\circ
    \logicalpolar{\Pi}^{m,\pr}\circ\logical{\pushforward}^{-m},
\end{equation}
where we note that -- by virtue of being coefficient projection operators -- their action can simply be expressed in terms of the coefficients of the $m$-form, which are the same in the logical, intermediate and physical spaces.

Let $\pi^m$ denote the 3D GDOF projector onto the full tensor-product space $\funspace{\Bsplinethree}^m$.
As in \cref{sec:physical_geometric_dofs}, it is the unique discrete form matching the point values, oriented edge integrals, oriented face integrals, or cell integrals on the mapped Greville grid for $m=0,1,2,3$, respectively.
Since integrals of forms over mapped edges, faces and cells are invariant under pullback, $\pi^m=\pushforward^m\circ\logical{\pi}^m\circ\pushforward^{-m}=\intermediate{\pushforward}^m\circ\intermediate{\pi}^m\circ\intermediate{\pushforward}^{-m}$ as in \cref{eqn:physical_geometric_dofs}, where $\logical{\pi}^m$ is the logical GDOF projector on $\logical{\funspace{\Bsplinethree}}^m$ and $\intermediate{\pi}^m\defeq\logical{\pushforward}^m\circ\logical{\pi}^m\circ\logical{\pushforward}^{-m}$.
As in \cref{eqn:polar_projection_operator_definition}, we define the infinite-dimensional polar projectors by composition:
\begin{equation}\label{eqn:polar_3d_projection_operator_definition}
  \polar{\pi}^{m,\pr}\defeq\polar{\Pi}^{m,\pr}\circ\pi^m,
  \qquad m\in\{0,1,2,3\}.
\end{equation}

\begin{theorem}[Commuting 3D polar projection operators]
  \label{thm:commuting_polar_projection_operators_3d}
  Let $\pr\ge1$.
  For each $m\in\{0,1,2,3\}$, $\polar{\pi}^{m,\pr}$ defined in \cref{eqn:polar_3d_projection_operator_definition} is a projection onto $\polar{\funspace{\Bsplinethree}}^{m,\pr}$ and is the identity on that space.

  Let $P\in\funspace{H}_1$, $\+A\in\funspace{H}_{\curlt}$, $\+B\in\funspace{H}_{\divergencet}$ and $Q\in\funspace{L}_2$ be fields whose geometric degrees of freedom are defined (as in \cref{cor:infinite_dimensional_polar_projection_operators}, for example continuous fields; for $Q$ the cell integrals exist for every $Q\in\funspace{L}_2$).
  These projectors commute with the differential operators:
  \begin{equation}\label{eqn:polar_3d_projection_differential_commutativity}
      \gradientt\polar{\pi}^{0,\pr}P
      =\polar{\pi}^{1,\pr}\gradientt P, \quad
      \curlt\polar{\pi}^{1,\pr}\+A
      =\polar{\pi}^{2,\pr}\curlt\+A, \quad
      \divergencet\polar{\pi}^{2,\pr}\+B
      =\polar{\pi}^{3,\pr}\divergencet\+B.
  \end{equation}
  If the first 2D polar $0$-form dual satisfies \cref{eqn:polar_constant_dual_centre_value}, then
  \begin{equation}\label{eqn:polar_3d_projection_homotopy_commutativity}
      \gradientt^\dagger\polar{\pi}^{1,\pr}\+A
      =\polar{\pi}^{0,\pr}\gradientt^\dagger\+A, \quad
      \curlt^\dagger\polar{\pi}^{2,\pr}\+B
      =\polar{\pi}^{1,\pr}\curlt^\dagger\+B, \quad
      \divergencet^\dagger\polar{\pi}^{3,\pr}Q
      =\polar{\pi}^{2,\pr}\divergencet^\dagger Q.
  \end{equation}
\end{theorem}
\begin{statementproof}[appendix]{thm:commuting_polar_projection_operators_3d}
  We may assume $\intermediatemapping=\mathrm{id}$ on the cylinder: by \cref{lem:polar_3d_product_structure}, $\intermediate{\pushforward}^m$ maps the intermediate Hilbert spaces one-to-one onto the physical ones and commutes with the weak differential operators, and by \cref{eqn:physical_3d_polar_projection_operator,eqn:polar_spaces_3d_pushforward,eqn:physical_3d_tensor_product_homotopy,eqn:pushforward_composition_3d}, the physical projectors, polar spaces and homotopies are $\intermediate{\pushforward}^m$-conjugates of the intermediate ones; hence each statement on $\Omega$ is equivalent to the corresponding statement on $\intermediate{\Omega}$.
  Then $\pushforward^m=\logical{\pushforward}^m$, and the remaining argument works with the logical representations.
  The product structure of the Greville GDOFs, unisolvence, and interchange of integrals give the logical representation of the composition in \cref{eqn:polar_3d_projection_operator_definition}:
  \begin{equation}
    \pushforward^{-m}\circ\polar{\pi}^{m,\pr}\circ\pushforward^m
    =\bigoplus_{p+q=m}
      \left(\logicalpolar{\Pi}^{p,\pr}_\mathrm{2D}
        \circ\logical{\pi}^p_\mathrm{2D}\right)
      \otimes\pi^{q,\zz}.
  \end{equation}
  Each factor is a projection onto the corresponding 2D polar or periodic space and fixes that space, by \cref{cor:infinite_dimensional_polar_projection_operators,eqn:commuting_projectors_1d}.
  Their tensor products and direct sum therefore project onto the logical space in \cref{eqn:tensor_product_complex_spaces_3d}; conjugation by the Piola maps gives the stated physical projection property.

  Tensor products of cochain maps are cochain maps: this is the functoriality of the tensor product of complexes in \citet[Chap.~V, Sec.~9, p.~164]{maclane1995}, with chain and cochain gradings related as in \citet[Chap.~II, Sec.~2, p.~41]{maclane1995}.
  The full GDOF factors commute with their differentials by \cref{eqn:logical_commuting_diagram,eqn:stokes_1d}.
  After the full GDOF projection, the collapsed scalar nodes and toroidal edges give center traces independent of $\logical{y}$ for the $(0,0)$ and $(0,1)$ components.
  Thus the scalar factors are continuous at the cross-section center, as required by the 2D coefficient commutation result in \cref{thm:full_commuting_polar_projection_operators}.
  Applying tensor-product functoriality to the full and coefficient projectors, and then composing these cochain maps, proves logical differential commutation on the total complex.
  The Piola maps give \cref{eqn:polar_3d_projection_differential_commutativity}.

  For homotopy commutation, the 2D factor commutes with the radial homotopies by \cref{thm:homotopy_commuting_polar_projection_operators,cor:infinite_dimensional_polar_projection_operators} and their GDOF accumulation proofs.
  The 2D gradient homotopy requires \cref{eqn:polar_constant_dual_centre_value}, whereas the 2D curl homotopy does not.
  Under that condition, constant reproduction and center-value preservation in \cref{eqn:polar_zeroform_projection_constant,eqn:polar_zeroform_projection_centre_value}, together with preservation of the collapsed-node GDOFs, also give
  \begin{equation}
    \mathfrak{P}_\mathrm{2D}^0
      \left(\logicalpolar{\Pi}^{0,\pr}_\mathrm{2D}
        \circ\logical{\pi}^0_\mathrm{2D}\right)
    =\mathfrak{P}_\mathrm{2D}^0
    =\left(\logicalpolar{\Pi}^{0,\pr}_\mathrm{2D}
        \circ\logical{\pi}^0_\mathrm{2D}\right)
      \mathfrak{P}_\mathrm{2D}^0.
  \end{equation}
  In positive degrees $\mathfrak{P}_\mathrm{2D}^p=0$.
  The periodic factor commutes with the normalized primitive by \cref{eqn:logical_1d_homotopy_commutativity_3d}.
  Its GDOF argument applies also to the center trace: the primitive vanishes at $\xi^\zz_1$, and its nodal values are accumulated toroidal edge integrals with the preserved mean subtracted.
  Hence the factor projectors commute with both terms of the ordered tensor-product homotopy in \cref{eqn:logical_3d_tensor_product_homotopy}, proving \cref{eqn:polar_3d_projection_homotopy_commutativity} after Piola conjugation.
  The total gradient and curl homotopies use the 2D gradient homotopy on bidegrees $(1,0)$ and $(1,1)$, respectively, so both require the additional dual condition.
  The total divergence homotopy uses only the 2D curl homotopy on $(2,1)$ and therefore commutes without that condition.
\end{statementproof}

\begin{lemma}[Conservation properties of the 3D polar projection operators]
  \label{lem:polar_3d_projection_operators_conservation_properties}
  Let $\pr\ge1$.
  The polar projectors preserve the total integral of every $Q\in\funspace{L}_2$ and the total boundary flux of every $\+B\in\funspace{H}_{\divergencet}$ with defined geometric degrees of freedom:
  \begin{equation}\label{eqn:polar_3d_projection_conservation}
    \begin{aligned}
      \int_\Omega\polar{\pi}^{3,\pr}Q\diff^3x
      &=\int_\Omega Q\diff^3x,\\
      \int_{\partial\Omega}(\polar{\pi}^{2,\pr}\+B)\cdot\+n\diff S
      &=\int_{\partial\Omega}\+B\cdot\+n\diff S.
    \end{aligned}
  \end{equation}
  If the first 2D polar $0$-form dual satisfies \cref{eqn:polar_constant_dual_centre_value}, the $1$-form projector preserves circulation along the center curve for every continuous $\+A\in\funspace{H}_{\curlt}$:
  \begin{equation}\label{eqn:polar_3d_projection_centre_circulation_conservation}
    \int_\Gamma(\polar{\pi}^{1,\pr}\+A)\cdot\+T\diff l
    =\int_\Gamma\+A\cdot\+T\diff l,
  \end{equation}
  where $\+T$ is the unit tangent to the center curve $\Gamma$, oriented by increasing $\intermediatecoord{z}$.
\end{lemma}
\begin{statementproof}[appendix]{lem:polar_3d_projection_operators_conservation_properties}
  The proof is written with the full pushforwards and holds for every intermediate mapping.
  In top degree there is only the $(2,1)$ component.
  The coefficient projector preserves its cross-section integral by the coefficient argument in \cref{lem:polar_projection_operators_conservation_properties}; integrating also in $\logical{z}$ preserves the total integral.
  The full GDOF projector preserves the total integral because its cell moments sum to that integral.
  Composition and invariance of top-form integrals under pullback therefore give the first identity in \cref{eqn:polar_3d_projection_conservation}.
  The divergence commutation relation in \cref{thm:commuting_polar_projection_operators_3d} and the weak divergence theorem give the boundary-flux identity.

  Finally, under \cref{eqn:polar_constant_dual_centre_value}, the scalar coefficient projector preserves center values by \cref{eqn:polar_zeroform_projection_centre_value}.
  The toroidal component of $\pi^1\+A$ has a center trace independent of $\logical{y}$, as proved above, so the coefficient projector leaves that trace unchanged.
  The full projector preserves the toroidal edge integrals along $\Gamma$, which are the line integrals along the images under $\mapping$ of the logical edges at $\logical{x}=0$ in the $\logical{z}$-direction, and which sum to the circulation along $\Gamma$.
  This proves \cref{eqn:polar_3d_projection_centre_circulation_conservation}.
\end{statementproof}

  \section{Numerical results}\label{sec:results}
We illustrate the main properties of the polar subspaces with five numerical experiments of increasing complexity.
The polar subspaces and bases, as well as their differential, homotopy and projection operators, all have concise matrix representations which are detailed in \cref{sec:matrix_representation} and have been implemented in the Julia package \polarfeec.

First, we show that the 2D commuting projections converge at the optimal rates for all $\pr$, so that the regularity conditions do not cost accuracy (\cref{sec:results_projections}).
Second, we show that the regularity guaranteed by \cref{thm:polar_subspaces_regularity} is attained and sharp: the derivatives of projected $0$- and $2$-forms are continuous at the center up to exactly the guaranteed order (\cref{sec:results_regularity}).
Third, the spectrum of the 2D curl--curl operator shows that increasing the regularity suppresses the stiff eigenmodes originating at the center, which allows for larger time steps, while the smallest eigenvalues keep converging at the optimal rate, also with spaces that are not $\pr$-compatible (\cref{sec:results_spectrum}).
Fourth, we solve the 3D $1$-form wave equation on a tube around a helical axis with spaces that are not $\pr$-compatible, and find that high regularity removes the artifacts at the axis and improves the accuracy for under-resolved pulses (\cref{sec:results_wave3d}).
Finally, we trace guiding-center orbits that pass close to the axis of a discrete vacuum magnetic field on a similar domain, and show that, unlike the waves, they require $\pr$-compatible spaces with $\pr\ge4$ for the energy error to converge at the expected rate (\cref{sec:results_vacuum}).

In all experiments, $\pr=1$ yields the first-order polar spaces of \citet{toshniwal2021,holderied2022,yaman2025} (cf.\ \cref{sec:first_order_regularity_de_rham}), so that every comparison across $\pr$ includes this baseline from the literature.
We note that an extensive study of the properties of the $0$-form space can be found in \citet{jiang2026polar}.

\subsection{Accuracy of 2D commuting projections}\label{sec:results_projections}
\begin{figure}
  \centering
  \subcaptionbox{$0$-forms.\label{fig:results_projections_zeroform}}
  [0.32\textwidth]{
    \tikzsetnextfilename{result_projections_commuting_form0}
\begin{tikzpicture}
\begin{axis}[height={0.25\textwidth}, width={0.21\textwidth}, scale only axis, xmode={log}, ymode={log}, log basis x={2}, xmin={0.019048079724204754}, xmax={0.1367154434658053}, ymin={1.0e-9}, ymax={1.0}, xlabel={$h$}, ytick={{1e0, 1e-2, 1e-4, 1e-6, 1e-8}}, ylabel={relative $\funspace{L}_2$ error}, title={$\pr = r^{0,\xx} = r^{0,\yy} = {\color{lines2}\mathbf{1}}, {\color{lines3}\mathbf{2}}, {\color{lines4}\mathbf{3}}, {\color{lines5}\mathbf{4}}, {\color{lines6}\mathbf{5}}, {\color{lines7}\mathbf{6}}$}, title style={font=\footnotesize}, tick label style={font=\scriptsize}, label style={font=\footnotesize}]
    \addplot[mark={*}, mark size={1.2pt}, color={lines2}, line width={0.8pt}]
        coordinates {
            (0.125,0.11678233265820982)
            (0.08333333333333333,0.05382881746302692)
            (0.0625,0.030823398768250528)
            (0.041666666666666664,0.013885777742296863)
            (0.03125,0.007847862691676257)
            (0.020833333333333332,0.003499815693330142)
        }
        ;
    \addplot[mark={*}, mark size={1.2pt}, color={lines3}, line width={0.8pt}]
        coordinates {
            (0.125,0.03726042770573217)
            (0.08333333333333333,0.009867763762302261)
            (0.0625,0.0033080584401806857)
            (0.041666666666666664,0.0007284101506102958)
            (0.03125,0.00027085523142435523)
            (0.020833333333333332,7.268744124742729e-5)
        }
        ;
    \addplot[mark={*}, mark size={1.2pt}, color={lines4}, line width={0.8pt}]
        coordinates {
            (0.125,0.03507954267716044)
            (0.08333333333333333,0.007059097577639407)
            (0.0625,0.0015010093321320506)
            (0.041666666666666664,0.00019838790727535655)
            (0.03125,5.302938755299899e-5)
            (0.020833333333333332,9.255536178303983e-6)
        }
        ;
    \addplot[mark={*}, mark size={1.2pt}, color={lines5}, line width={0.8pt}]
        coordinates {
            (0.125,0.025937790202356325)
            (0.08333333333333333,0.005146754335257032)
            (0.0625,0.0007764338181136588)
            (0.041666666666666664,4.4513275288410156e-5)
            (0.03125,7.12780157448369e-6)
            (0.020833333333333332,6.514160700020783e-7)
        }
        ;
    \addplot[mark={*}, mark size={1.2pt}, color={lines6}, line width={0.8pt}]
        coordinates {
            (0.125,0.030341329392356942)
            (0.08333333333333333,0.004981048640372672)
            (0.0625,0.00038288458702664854)
            (0.041666666666666664,1.5810474545404782e-5)
            (0.03125,1.7654106858677024e-6)
            (0.020833333333333332,1.0809820400964214e-7)
        }
        ;
    \addplot[mark={*}, mark size={1.2pt}, color={lines7}, line width={0.8pt}]
        coordinates {
            (0.125,0.024300589181663942)
            (0.08333333333333333,0.0038570143846026884)
            (0.0625,0.0003620840230054992)
            (0.041666666666666664,6.208894184154224e-6)
            (0.03125,3.939865829055719e-7)
            (0.020833333333333332,1.2178299022026631e-8)
        }
        ;
    \addplot[mark={none}, color={gray}, line width={0.6pt}, forget plot]
        coordinates {
            (0.02083333333333333,5.349051676348061e-9)
            (0.03125,9.139356262635323e-8)
            (0.03125,5.349051676348061e-9)
            (0.02083333333333333,5.349051676348061e-9)
        }
        ;
    \node at (axis cs:0.03125, 2.2110379675028563e-8) [color=gray,right,font=\footnotesize] {$\mathcal{O}(h^{7})$};
\end{axis}
\end{tikzpicture}
  }\hfill
  \subcaptionbox{$1$-forms.\label{fig:results_projections_oneform}}
  [0.32\textwidth]{
    \tikzsetnextfilename{result_projections_commuting_form1}
\begin{tikzpicture}
\begin{axis}[height={0.25\textwidth}, width={0.21\textwidth}, scale only axis, xmode={log}, ymode={log}, log basis x={2}, xmin={0.019048079724204754}, xmax={0.1367154434658053}, ymin={1.0e-9}, ymax={1.0}, xlabel={$h$}, ytick={{1e0, 1e-2, 1e-4, 1e-6, 1e-8}}, ylabel={relative $\funspace{L}_2$ error}, title={$\pr = r^{0,\xx} = r^{0,\yy} = {\color{lines2}\mathbf{1}}, {\color{lines3}\mathbf{2}}, {\color{lines4}\mathbf{3}}, {\color{lines5}\mathbf{4}}, {\color{lines6}\mathbf{5}}, {\color{lines7}\mathbf{6}}$}, title style={font=\footnotesize}, tick label style={font=\scriptsize}, label style={font=\footnotesize}]
    \addplot[mark={*}, mark size={1.2pt}, color={lines2}, line width={0.8pt}]
        coordinates {
            (0.125,0.29987346939159115)
            (0.08333333333333333,0.2001481421202082)
            (0.0625,0.15020396024471547)
            (0.041666666666666664,0.10018692491051376)
            (0.03125,0.07515403189170719)
            (0.020833333333333332,0.05010929637020056)
        }
        ;
    \addplot[mark={*}, mark size={1.2pt}, color={lines3}, line width={0.8pt}]
        coordinates {
            (0.125,0.06417204244109108)
            (0.08333333333333333,0.024481609275124987)
            (0.0625,0.012174013362112436)
            (0.041666666666666664,0.00484574323901392)
            (0.03125,0.002622937027011912)
            (0.020833333333333332,0.0011347181619022226)
        }
        ;
    \addplot[mark={*}, mark size={1.2pt}, color={lines4}, line width={0.8pt}]
        coordinates {
            (0.125,0.05609237543910169)
            (0.08333333333333333,0.012242719023667338)
            (0.0625,0.0034172503147677238)
            (0.041666666666666664,0.0007129448228562764)
            (0.03125,0.0002619441752849061)
            (0.020833333333333332,6.995749683129377e-5)
        }
        ;
    \addplot[mark={*}, mark size={1.2pt}, color={lines5}, line width={0.8pt}]
        coordinates {
            (0.125,0.04500035284965992)
            (0.08333333333333333,0.006397149712316883)
            (0.0625,0.0013683360455578383)
            (0.041666666666666664,0.0001508618765479867)
            (0.03125,3.7758770904356505e-5)
            (0.020833333333333332,6.222816805896815e-6)
        }
        ;
    \addplot[mark={*}, mark size={1.2pt}, color={lines6}, line width={0.8pt}]
        coordinates {
            (0.125,0.050307634738664606)
            (0.08333333333333333,0.006211148814711496)
            (0.0625,0.0006713731090673659)
            (0.041666666666666664,3.891554077374017e-5)
            (0.03125,6.481238396287212e-6)
            (0.020833333333333332,6.429006874047367e-7)
        }
        ;
    \addplot[mark={*}, mark size={1.2pt}, color={lines7}, line width={0.8pt}]
        coordinates {
            (0.125,0.04477179417311974)
            (0.08333333333333333,0.0028660283586616386)
            (0.0625,0.0004091827448980252)
            (0.041666666666666664,1.1746993010509006e-5)
            (0.03125,1.244901777878538e-6)
            (0.020833333333333332,7.277630840177516e-8)
        }
        ;
    \addplot[mark={none}, color={gray}, line width={0.6pt}, forget plot]
        coordinates {
            (0.02083333333333333,3.196540286544519e-8)
            (0.03125,3.6410591701421216e-7)
            (0.03125,3.196540286544519e-8)
            (0.02083333333333333,3.196540286544519e-8)
        }
        ;
    \node at (axis cs:0.03125, 1.0788323467087761e-7) [color=gray,right,font=\footnotesize] {$\mathcal{O}(h^{6})$};
\end{axis}
\end{tikzpicture}
  }\hfill
  \subcaptionbox{$2$-forms.\label{fig:results_projections_twoform}}
  [0.32\textwidth]{
    \tikzsetnextfilename{result_projections_commuting_form2}
\begin{tikzpicture}
\begin{axis}[height={0.25\textwidth}, width={0.21\textwidth}, scale only axis, xmode={log}, ymode={log}, log basis x={2}, xmin={0.019048079724204754}, xmax={0.1367154434658053}, ymin={1.0e-9}, ymax={1.0}, xlabel={$h$}, ytick={{1e0, 1e-2, 1e-4, 1e-6, 1e-8}}, ylabel={relative $\funspace{L}_2$ error}, title={$\pr = r^{0,\xx} = r^{0,\yy} = {\color{lines2}\mathbf{1}}, {\color{lines3}\mathbf{2}}, {\color{lines4}\mathbf{3}}, {\color{lines5}\mathbf{4}}, {\color{lines6}\mathbf{5}}, {\color{lines7}\mathbf{6}}$}, title style={font=\footnotesize}, tick label style={font=\scriptsize}, label style={font=\footnotesize}]
    \addplot[mark={*}, mark size={1.2pt}, color={lines2}, line width={0.8pt}]
        coordinates {
            (0.125,0.18366037210310282)
            (0.08333333333333333,0.12317711053132785)
            (0.0625,0.09266205080975366)
            (0.041666666666666664,0.06198550821414032)
            (0.03125,0.04658714868669255)
            (0.020833333333333332,0.031143874426424425)
        }
        ;
    \addplot[mark={*}, mark size={1.2pt}, color={lines3}, line width={0.8pt}]
        coordinates {
            (0.125,0.033014099272356184)
            (0.08333333333333333,0.013739956292321198)
            (0.0625,0.007561246678271702)
            (0.041666666666666664,0.0033147577178316763)
            (0.03125,0.0018571561216217686)
            (0.020833333333333332,0.0008237550430085688)
        }
        ;
    \addplot[mark={*}, mark size={1.2pt}, color={lines4}, line width={0.8pt}]
        coordinates {
            (0.125,0.010489603458589325)
            (0.08333333333333333,0.002148503464703449)
            (0.0625,0.0007474439738772255)
            (0.041666666666666664,0.0001919654396903165)
            (0.03125,7.708467947380575e-5)
            (0.020833333333333332,2.2054834113801415e-5)
        }
        ;
    \addplot[mark={*}, mark size={1.2pt}, color={lines5}, line width={0.8pt}]
        coordinates {
            (0.125,0.006760437676382178)
            (0.08333333333333333,0.0008129916361663766)
            (0.0625,0.0001630510239720979)
            (0.041666666666666664,2.2830169000611855e-5)
            (0.03125,6.398619822995454e-6)
            (0.020833333333333332,1.1578725682162907e-6)
        }
        ;
    \addplot[mark={*}, mark size={1.2pt}, color={lines6}, line width={0.8pt}]
        coordinates {
            (0.125,0.00510231260003394)
            (0.08333333333333333,0.00043265685492334426)
            (0.0625,5.557751811724448e-5)
            (0.041666666666666664,3.910615189842183e-6)
            (0.03125,7.325059346090355e-7)
            (0.020833333333333332,8.039393922247156e-8)
        }
        ;
    \addplot[mark={*}, mark size={1.2pt}, color={lines7}, line width={0.8pt}]
        coordinates {
            (0.125,0.005112656770689156)
            (0.08333333333333333,0.00036008507573762045)
            (0.0625,2.4543798161168955e-5)
            (0.041666666666666664,8.59053963559206e-7)
            (0.03125,1.0627387969950544e-7)
            (0.020833333333333332,6.958232069518716e-9)
        }
        ;
    \addplot[mark={none}, color={gray}, line width={0.6pt}, forget plot]
        coordinates {
            (0.02083333333333333,3.0562513573167327e-9)
            (0.03125,3.481261311693593e-8)
            (0.03125,3.0562513573167327e-9)
            (0.02083333333333333,3.0562513573167327e-9)
        }
        ;
    \node at (axis cs:0.03125, 1.0314848330943976e-8) [color=gray,right,font=\footnotesize] {$\mathcal{O}(h^{6})$};
\end{axis}
\end{tikzpicture}
  }
  \caption{Relative $\funspace{L}_2$ errors of the commuting projections $\polar{\pi}^{m,\pr}$ of a smooth $0$-, $1$- and $2$-form onto the polar spaces with $\pr=r^{0,\xx}=r^{0,\yy}$, against the mesh size $h=1/N$, for B-spline degrees $1,\ldots,6$.}
  \label{fig:results_projections}
\end{figure}
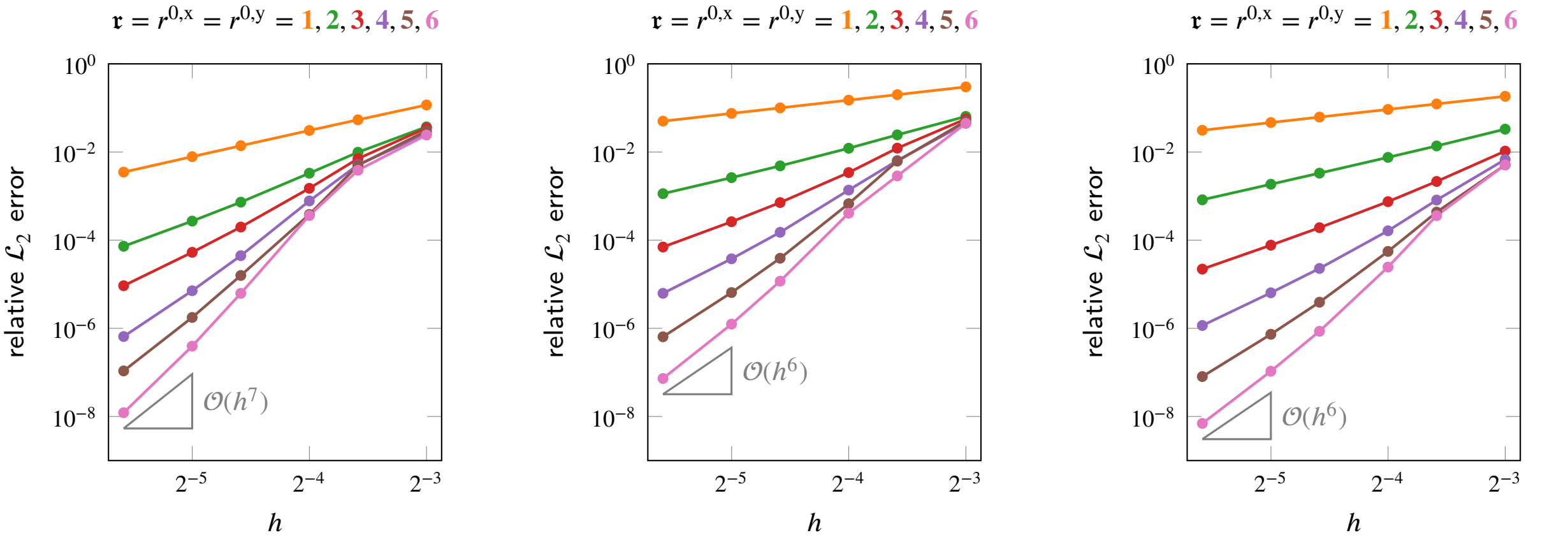
The first experiment asks whether the regularity conditions at the center cost approximation accuracy.
The commuting projections introduced in \cref{sec:derham_hodge_2d_geom_projections} commute with the differential operators by construction; we have verified this to machine precision, and the verification is part of the test suite of the Julia package.
Here we additionally verify the accuracy of the projection operators $\polar{\pi}^{m,\pr}$ from \cref{eqn:polar_projection_operator_definition}, by computing their $\funspace{L}_2$ approximation error on the polar $m$-form spaces.
We consider the composed mapping from \cref{fig:disk_mapping}.
We use the polar spaces with $\pr=r^{0,\xx} = r^{0,\yy} = 1, \ldots, 6$, built from degree-$r^{0,\alpha}$ clamped radial and periodic angular B-splines on $(N,2N)$ cells, and project the smooth $0$-, $1$- and $2$-forms
\begin{equation}\label{eqn:results_projection_forms}
  P = \sin(3x+1)\cos(2y-0.5)
  , \quad
  \+A = \begin{bmatrix}\cos(2x-y+0.3)+0.5y^2\\\sin(x+2y)-0.3x\end{bmatrix}
  , \quad
  B = \mathrm{e}^{x/2}\sin(2y+0.4)+0.3x^2.
\end{equation}
\Cref{fig:results_projections} shows the relative $\funspace{L}_2$ errors of the commuting projections.
They converge at least at the optimal rates $h^{p+1}$ for $0$-forms and $h^p$ for $1$- and $2$-forms, so the regularity conditions at the center do not degrade the accuracy.

\subsection{Regularity at the center of a 2D disk-like domain}\label{sec:results_regularity}
\begin{figure}
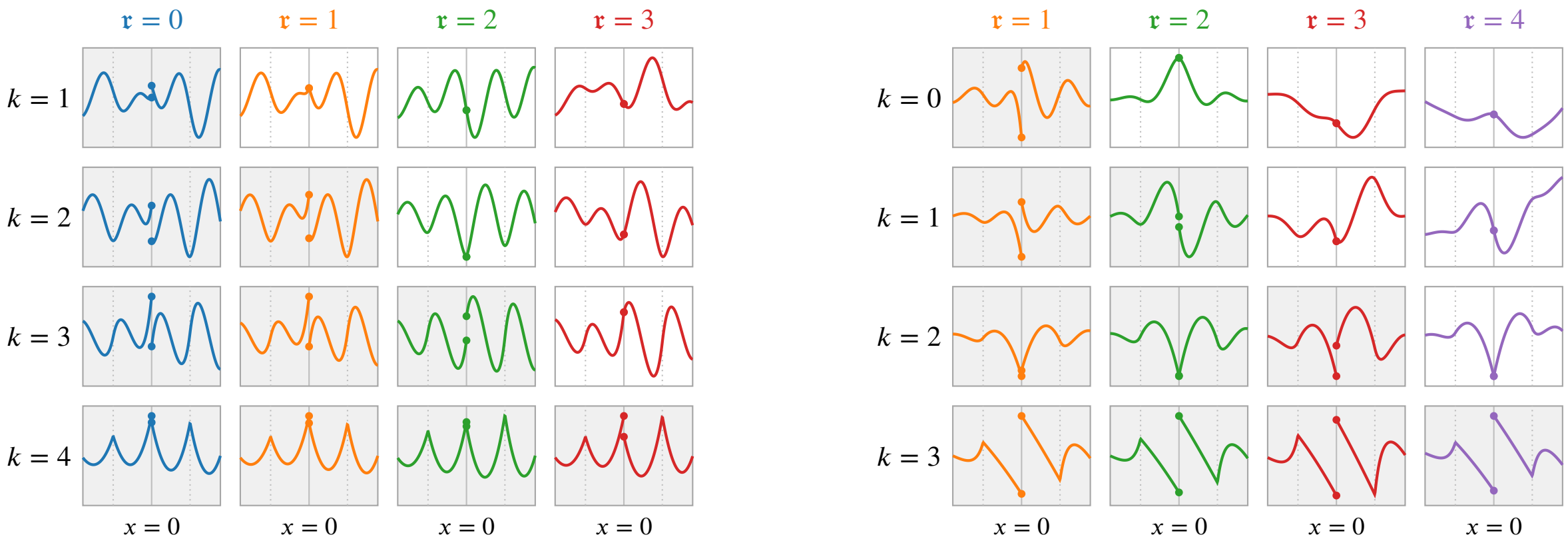

  \centering
  \subcaptionbox{$0$-forms.\label{fig:results_regularity_zeroform}}
  [\twofigwidth]{
    \tikzsetnextfilename{result_regularity_zeroform}\input{tikz/result_regularity_zeroform.tex}
  }\hfill
  \subcaptionbox{$2$-forms.\label{fig:results_regularity_twoform}}
  [\twofigwidth]{
    \tikzsetnextfilename{result_regularity_twoform}\input{tikz/result_regularity_twoform.tex}
  }
  \caption{Derivatives $\partial_x^k(P_h-P)$ along the $x$-axis, for $|x|\le0.2$, of the error of the commuting projection $P_h = \polar{\pi}^{m,\pr} P$ of the $0$-form $P$ of \cref{eqn:results_projection_forms} onto the polar $0$-form spaces (left) and $2$-form spaces (right), each normalized by its maximum over the plotted window.
  The dots mark the one-sided limits at the center, and the dotted lines the radial knots.
  A gray background indicates the derivative orders for which \cref{thm:polar_subspaces_regularity} does not guarantee continuity ($k>\pr$ for $0$-forms and $k>\pr-2$ for $2$-forms).}
  \label{fig:results_regularity}
\end{figure}
The second experiment asks whether the regularity guaranteed by \cref{thm:polar_subspaces_regularity} is attained and sharp.
We use the banana-shaped domain obtained from the unit disk of \cref{ex:disk_mapping} with $R_1=1$ as logical mapping, which is $\pr$-compatible with the logical spaces for all $\pr$ considered here, and the banana map of \cref{ex:banana_mapping} with $\kappa=1$ and $\alpha=2/\sqrt{5}$ as intermediate mapping, together with degree-$5$ clamped B-splines on $N=8$ radial cells and a trigonometric angular basis with $n^\yy=49$ modes.
We project the smooth $0$-form $P$ of \cref{eqn:results_projection_forms} with the commuting projection operators onto the polar $0$-form spaces with $\pr=0,\ldots,3$, and, as the density of a $2$-form, onto the polar $2$-form spaces with $\pr=1,\ldots,4$.

In \cref{fig:results_regularity} we show the derivative error $\partial_x^k(P_h-P)$ of the projections along the $x$-axis near the center.
Since $P$ is smooth, every discontinuity at $x=0$ stems from $P_h = \polar{\pi}^{m,\pr} P$.
The derivatives are continuous at the center precisely for $k\le\pr$ for $0$-forms and for $k\le\pr-2$ for $2$-forms.
For $\pr=1$, the polar $2$-forms are bounded but have no limit at the center.
Hence the observed regularity coincides exactly with the one guaranteed by \cref{thm:polar_subspaces_regularity}.

\subsection{Spectrum of the 2D curl--curl operator}\label{sec:results_spectrum}
The third experiment concerns the stiff modes originating at the center: the eigenmodes with the largest eigenvalues, which are concentrated at the center and irregular there (quantified below by the regularity deficit), and which limit the time step of explicit schemes.
These modes are not spurious in the sense of spectral pollution, since the discrete cohomology is that of the disk for every $\pr$ (\cref{cor:polar_subspaces_cohomology}) and, as verified below, the smallest eigenvalues converge with the correct multiplicities for every $\pr$; the question is by how much the regularity conditions reduce the stiffness, and whether they do so without $\pr$-compatibility.
We consider the eigenvalue problem of the curl--curl operator on the unit disk, with homogeneous essential boundary conditions:
\begin{equation}
  \text{Find} \quad (\+A,\lambda)\in \polar{\funspace{\Bspline}}^{1,\pr}_\bc \times \mathbb{R} 
  \quad \text{s.t.} \quad
  \inprod{\curlt\+\Bspline^1}{\curlt\+A}=\lambda\inprod{\+\Bspline^1}{\+A}
  \quad
  \forall \+\Bspline^1\in\polar{\funspace{\Bspline}}^{1,\pr}_\bc.
\end{equation}
The nonzero eigenvalues of the continuous problem are $\lambda=(j'_{m,k})^2$, where $j'_{m,k}$ denotes the $k$-th positive zero of the derivative of the Bessel function $J_m$, and the curls of the corresponding eigenfunctions are $J_m(j'_{m,k}r)\cos(m\theta)$ and $J_m(j'_{m,k}r)\sin(m\theta)$.
We use degree-$6$ clamped radial and periodic angular B-splines on $(N,2N)$ cells with $N=10,\ldots,32$, or on $(10,14)$ cells, and the polar spaces with $\pr=0,\ldots,5$.
Since periodic B-splines cannot represent $\cos(2\pi\logical{y})$ and $\sin(2\pi\logical{y})$ exactly, the spaces are merely $0$-compatible (\cref{ex:pr_compatibility_disc_mapping}), so that this experiment also tests the expectation of \cref{sec:pullback_compatibility} that the polar spaces suppress the stiff modes originating at the center without $\pr$-compatibility.
For $\pr=0$, the polar $1$-forms are not $\funspace{H}_{\curlt}$-conforming (\cref{thm:polar_subspaces_conformity}) and the results are quadrature-regularized: the Gauss quadrature, which has no node at $\logical{x}=0$, yields finite matrices, but the largest eigenvalue doubles when the number of quadrature points per cell is doubled from $7$ to $15$, whereas it is unchanged for $\pr=1$, so that $\pr=0$ serves only as a nonconforming reference and $\pr=1$ is the conforming baseline of the literature.

\begin{figure}
  \centering
  \subcaptionbox{Largest eigenvalue on $(N,2N)$ cells.\label{fig:results_spectrum_lambda_max}}
  [\twofigwidth]{
    \tikzsetnextfilename{result_spectrum_lambda_max}
\begin{tikzpicture}
\begin{axis}[height={0.29\textwidth}, width={0.33\textwidth}, scale only axis, xmode={log}, ymode={log}, xmin={9.435012960710349}, xmax={33.91622262020801}, xlabel={$N$}, ylabel={$\lambda_{\max}$}, xtick={{10,12,16,20,24,28,32}}, xticklabels={{10,12,16,20,24,28,32}}, ymin={1000.0}, ymax={1.0e9}, title={$\mathfrak{r} = {\color{lines1}\mathbf{0}}, {\color{lines2}\mathbf{1}}, {\color{lines3}\mathbf{2}}, {\color{lines4}\mathbf{3}}, {\color{lines5}\mathbf{4}}, {\color{lines6}\mathbf{5}}$}]
    \addplot[color={lines1}, mark={*}, mark size={1.7pt}, line width={0.6pt}]
        coordinates {
            (10.0,4.246470068359092e6)
            (12.0,8.792373467604112e6)
            (16.0,2.775058658271756e7)
            (20.0,6.770898931303367e7)
            (24.0,1.4035475460627675e8)
            (28.0,2.5997245880360645e8)
            (32.0,4.4344390805305696e8)
        }
        ;
    \addplot[color={lines2}, mark={*}, mark size={1.5pt}, line width={0.6pt}]
        coordinates {
            (10.0,179128.2109091091)
            (12.0,370204.3573770557)
            (16.0,1.1664599402531516e6)
            (20.0,2.843870865175447e6)
            (24.0,5.892634623342946e6)
            (28.0,1.0911903544535007e7)
            (32.0,1.860976946309372e7)
        }
        ;
    \addplot[color={lines3}, mark={*}, mark size={1.3pt}, line width={0.6pt}]
        coordinates {
            (10.0,27942.41264988669)
            (12.0,57661.35796022847)
            (16.0,181436.80315757287)
            (20.0,442081.73675458634)
            (24.0,915714.7369481185)
            (28.0,1.695374889025138e6)
            (32.0,2.8910181299446705e6)
        }
        ;
    \addplot[color={lines4}, mark={*}, mark size={1.1pt}, line width={0.6pt}]
        coordinates {
            (10.0,7887.462046412836)
            (12.0,16237.766977460726)
            (16.0,50986.789925413636)
            (20.0,124117.23274143685)
            (24.0,256964.38515379038)
            (28.0,475607.0549629623)
            (32.0,810865.9423189599)
        }
        ;
    \addplot[color={lines5}, mark={*}, mark size={0.9pt}, line width={0.6pt}]
        coordinates {
            (10.0,5854.4001422726815)
            (12.0,8490.108647972593)
            (16.0,19355.780001202533)
            (20.0,47045.27569427889)
            (24.0,97319.59155638881)
            (28.0,180037.3471163362)
            (32.0,306849.83132196614)
        }
        ;
    \addplot[color={lines6}, mark={*}, mark size={0.7pt}, line width={0.6pt}]
        coordinates {
            (10.0,5841.036049813462)
            (12.0,8485.141633668294)
            (16.0,15204.830945672118)
            (20.0,23853.340015274964)
            (24.0,44875.331595377786)
            (28.0,82942.69589315444)
            (32.0,141283.20078484068)
        }
        ;
    \addplot[mark={none}, color={gray}, line width={0.6pt}, forget plot]
        coordinates {
            (27.999999999999996,51614.563009001)
            (32.0,88052.16579961186)
            (32.0,51614.563009001)
            (27.999999999999996,51614.563009001)
        }
        ;
    \node at (axis cs:32.0, 51614.563009001) [color=gray,anchor=north east,font=\footnotesize] {$\mathcal{O}(N^{4})$};
    \addplot[mark={none}, color={gray}, line width={0.6pt}, forget plot]
        coordinates {
            (16.0,9487.96741723801)
            (20.0,14824.949089434373)
            (20.0,9487.96741723801)
            (16.0,9487.96741723801)
        }
        ;
    \node at (axis cs:20.0, 9487.96741723801) [color=gray,anchor=north east,font=\footnotesize] {$\mathcal{O}(N^{2})$};
\end{axis}
\end{tikzpicture}
  }\hfill
  \subcaptionbox{Convergence of the $32$ smallest eigenvalues on $(N,2N)$ cells: the maximum error (thick lines) and the error of each eigenvalue (thin lines).\label{fig:results_spectrum_convergence_eigenvalue}}
  [\twofigwidth]{
    \tikzsetnextfilename{result_spectrum_convergence_eigenvalue}\input{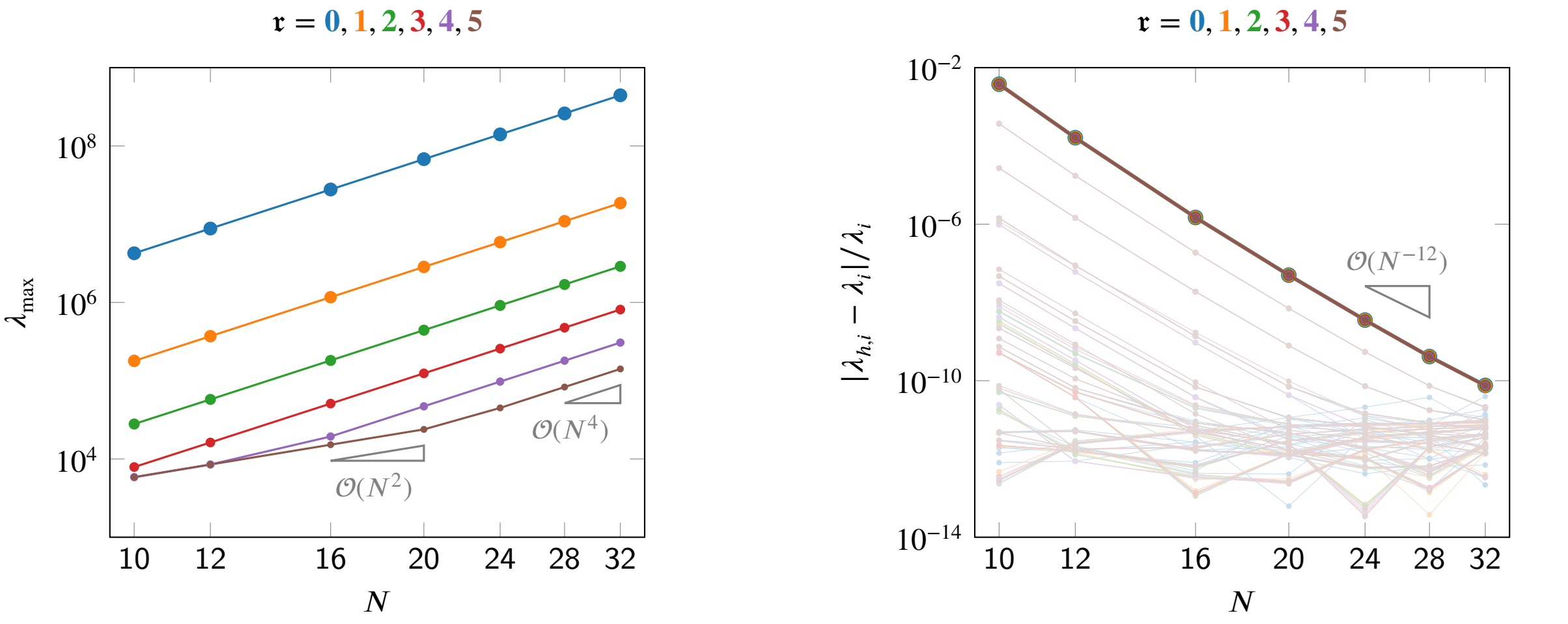}
  }
  \caption{Spectrum of the curl--curl operator on the unit disk with degree-$6$ B-splines.}
  \label{fig:results_spectrum}
\end{figure}
In \cref{fig:results_spectrum_lambda_max} we show the largest eigenvalue under mesh refinement with $(N,2N)$ cells for each regularity degree $\pr$.
For $\pr=0$, the largest eigenvalue is within about $1\%$ of that of the full tensor-product space, and it is caused by modes that oscillate at the center.
Near the center, the angular cell size is of order $N^{-2}$, so that these modes grow like $N^4$ instead of $N^2$.
For $\pr\le3$, the curves are parallel: on every mesh the largest eigenvalue decreases by a factor of about $23$ from the first-order polar spaces ($\pr=1$) to $\pr=3$, and the quadrature-regularized $\pr=0$ spaces lie another factor of about $24$ higher.
For $\pr=4,5$ and coarse meshes, the largest eigenvalue grows like $N^2$, so that it is no longer associated with the center, until the center modes dominate again from about $N=16$ and $N=20$, respectively.
At $N=32$, the largest eigenvalue decreases by a factor of about $130$ from the first-order polar spaces ($\pr=1$) to $\pr=5$; this ratio is reached from $N=24$ on and no longer changes under refinement.
Since the stable time step of an explicit time integrator scales with $\lambda_{\max}^{-1/2}$, the regularity conditions allow for a time step that is about $11$ times larger than with the first-order polar spaces (and about $56$ times larger than with the quadrature-regularized $\pr=0$ spaces).
\Cref{fig:results_spectrum_convergence_eigenvalue} shows the convergence of the $32$ smallest eigenvalues (counted with multiplicity) under mesh refinement with $(N,2N)$ cells, which matches the theoretical rate of $2r^{0,\alpha} = 12$.
Thus the regularity conditions do not affect the accuracy of the smooth eigenmodes, as desired.

\begin{figure}
  \centering
  \tikzsetnextfilename{result_spectrum_deficit}\input{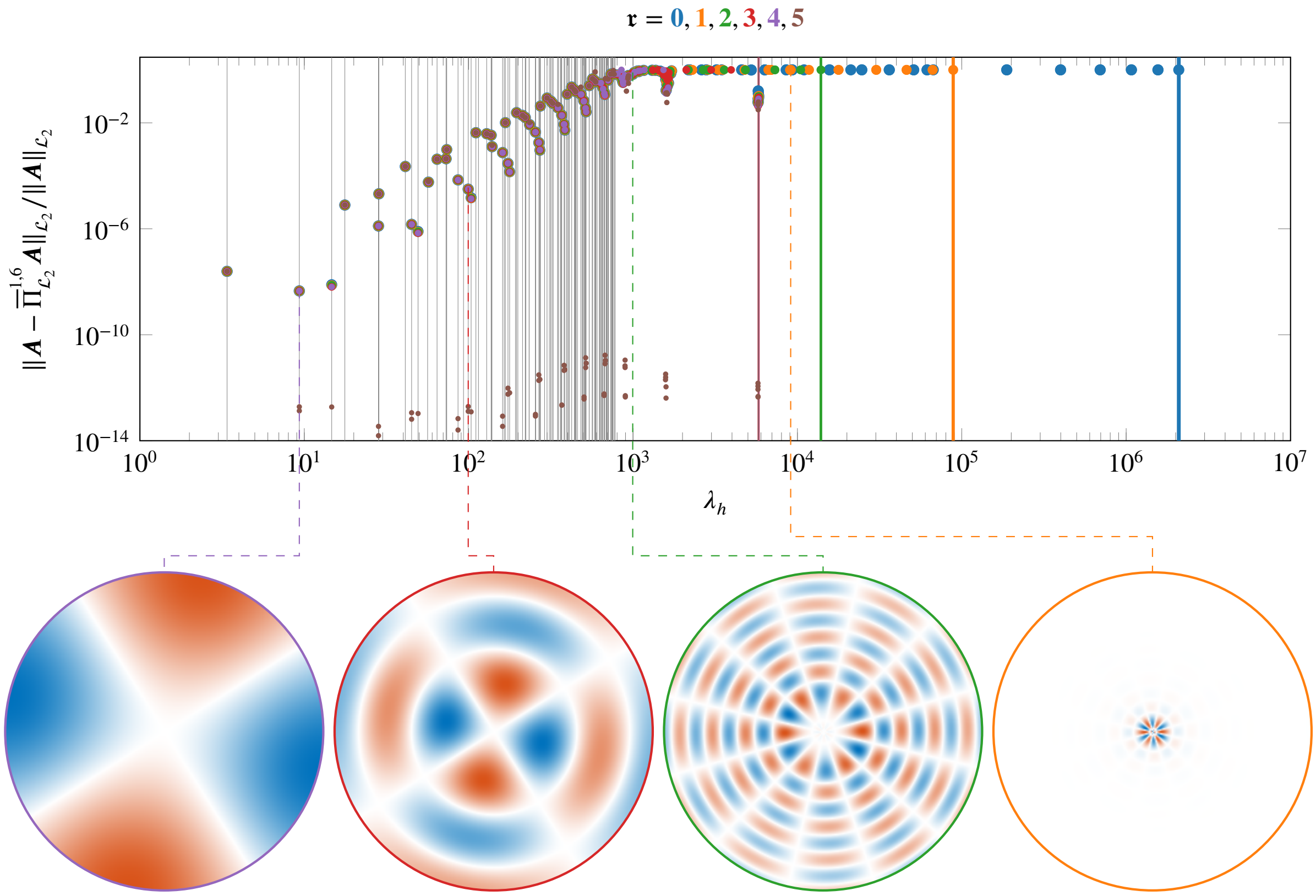}
  \caption{Eigenmodes of the curl--curl operator on the unit disk.
  The plot shows the regularity deficit $\|\+A-\polar{\Pi}^{1,6}_{\funspace{L}_2}\+A\|_{\funspace{L}_2}/\|\+A\|_{\funspace{L}_2}$, where $\polar{\Pi}^{1,6}_{\funspace{L}_2}$ denotes the $\funspace{L}_2$ projection onto the polar $1$-form space with $\pr=6$; point colors indicate $\pr$, and vertical lines mark the largest eigenvalue at each order. The four panels below show curls of modes selected at $\pr=4,3,2,1$ from left to right, with eigenvalues near $10,10^2,10^3,10^4$. Colored connectors identify each mode.}
  \label{fig:results_spectrum_eigenmodes}
\end{figure}

In \cref{fig:results_spectrum_eigenmodes} we show the spectrum for each of the regularity degrees using $(10, 14)$ cells.
The vertical scale is given by the regularity deficit: the relative error when projecting the polar $1$-form into the maximal regularity space
\begin{equation}
  \frac{\|\+A - \polar{\Pi}^{1,6}_{\funspace{L}_2}\+A\|_{\funspace{L}_2}}{\|\+A\|_{\funspace{L}_2}},
\end{equation}
where a value close to zero indicates that $\+A$ is regular near the center (according to $\polar{\funspace{\Bspline}}^{1,6}$), and a value close to one indicates that $\+A$ is irregular near the center.
Increasing $\pr$ removes the irregular, stiff eigenmodes of this coarse space while leaving the regular ones untouched.
The colored lines connect four representative modes of increasing irregularity and eigenvalue, selected at $\pr=4,3,2,1$ from left to right.
The gray vertical lines correspond to the first $N^\star$ exact eigenvalues, where $N^\star$ is the size of the polar $1$-form space minus the size of the polar $0$-form space, to account for the kernel of the curl operator (i.e., all gradients of $0$-forms have eigenvalue zero).

\subsection{Wave equation on a 3D tube mapping}\label{sec:results_wave3d}
The fourth experiment asks how the polar spaces perform in 3D on a wave problem whose spaces are not $\pr$-compatible.
We solve the $1$-form wave equation
\begin{equation}\label{eqn:results_wave3d}
  \text{Find} \quad \+A(t)\in \polar{\funspace{\Bsplinethree}}^{1,\pr}_\bc
  \quad \text{s.t.} \quad
  \inprod{\+\Bsplinethree^1}{\pard{^2\+A}{t^2}}+\inprod{\curlt\+\Bsplinethree^1}{\curlt\+A}=0
  \quad
  \forall \+\Bsplinethree^1\in\polar{\funspace{\Bsplinethree}}^{1,\pr}_\bc,
\end{equation}
on a 3D tube mapping from \cref{ex:tube_mapping_3d} around a helical axis with three field periods, with a frame that turns by one third of a turn per period and a fixed three-petal cross-section: the logical mapping is $\logicalmapping_\mathrm{2D}\times\mathrm{id}$ with $R=0.3\logical{x}(1+0.2\logical{x}\cos(6\pi\logical{y}))$ and $(C,S)=(\cos(2\pi\logical{y}),\sin(2\pi\logical{y}))$, and $\intermediatemapping_\mathrm{2D}=\mathrm{id}$.
We use degree-$4$ B-splines on $(N,2N,8N)$ radial, poloidal and toroidal cells, and the polar spaces with homogeneous essential boundary conditions, for $\pr=1,\ldots,4$.
These relative polar spaces $\polar{\funspace{\Bsplinethree}}^{m,\pr}_\bc$ consist of the $0$-, $1$- and $2$-forms with vanishing trace, tangential trace and normal trace on $\partial\Omega$, respectively; their pullbacks are the tensor products of the relative 2D spaces of \cref{sec:boundary_homotopy} with the periodic 1D spaces.
Since the three-petal shape is carried by the logical mapping and the angular basis consists of periodic B-splines, the logical spaces are only $0$-compatible with the mapping.

\begin{figure}
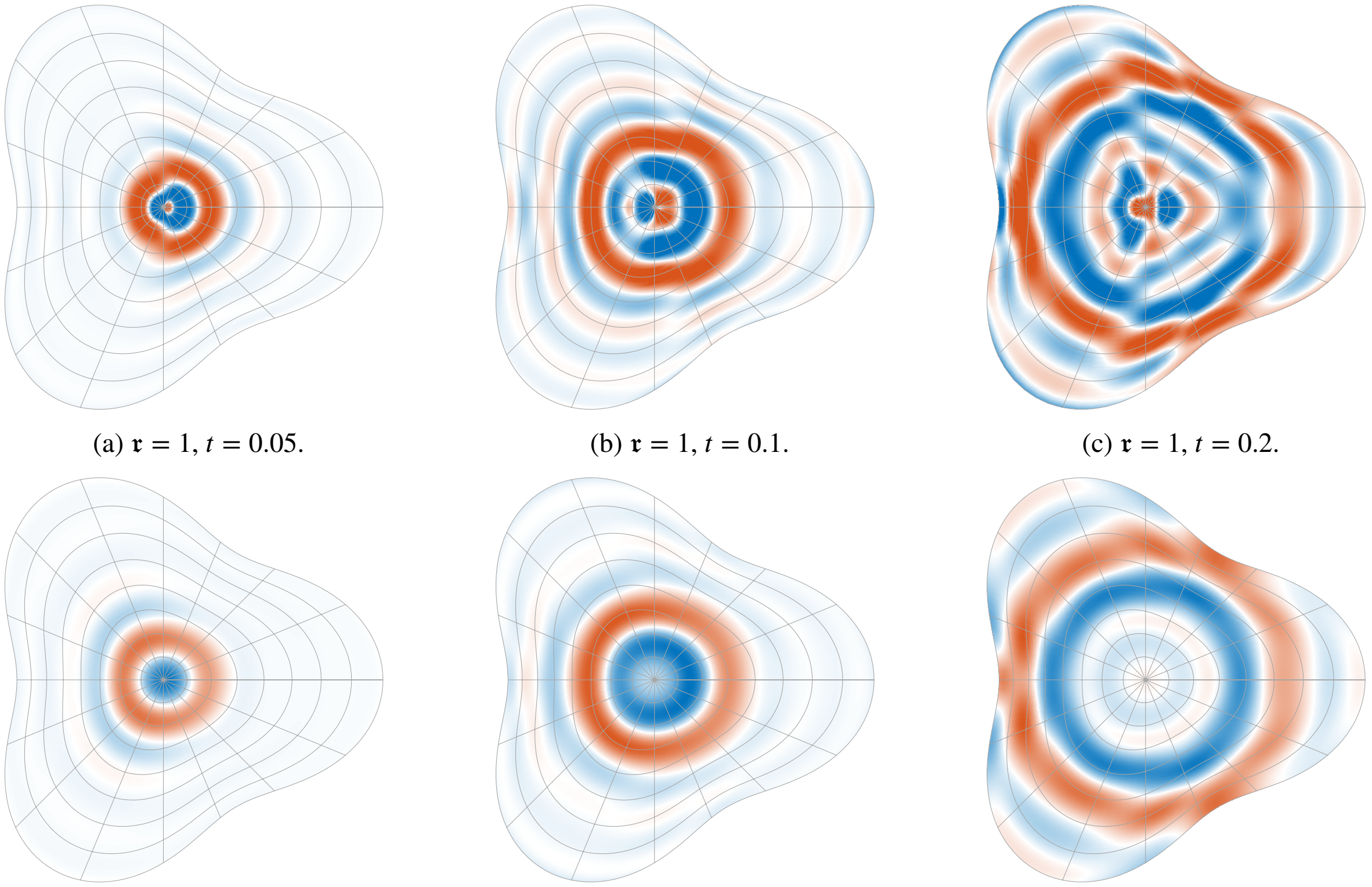

  \centering
  \subcaptionbox{$\pr=1$, $t=0.05$.\label{fig:results_wave3d_snapshot_r1_t1}}
  [0.32\textwidth]{
    \tikzsetnextfilename{result_wave3d_snapshot_r1_t1}\input{tikz/result_wave3d_snapshot_r1_t1.tex}
  }\hfill
  \subcaptionbox{$\pr=1$, $t=0.1$.\label{fig:results_wave3d_snapshot_r1_t2}}
  [0.32\textwidth]{
    \tikzsetnextfilename{result_wave3d_snapshot_r1_t2}\input{tikz/result_wave3d_snapshot_r1_t2.tex}
  }\hfill
  \subcaptionbox{$\pr=1$, $t=0.2$.\label{fig:results_wave3d_snapshot_r1_t3}}
  [0.32\textwidth]{
    \tikzsetnextfilename{result_wave3d_snapshot_r1_t3}\input{tikz/result_wave3d_snapshot_r1_t3.tex}
  }\\[1ex]
  \subcaptionbox{$\pr=4$, $t=0.05$.\label{fig:results_wave3d_snapshot_r4_t1}}
  [0.32\textwidth]{
    \tikzsetnextfilename{result_wave3d_snapshot_r4_t1}\input{tikz/result_wave3d_snapshot_r4_t1.tex}
  }\hfill
  \subcaptionbox{$\pr=4$, $t=0.1$.\label{fig:results_wave3d_snapshot_r4_t2}}
  [0.32\textwidth]{
    \tikzsetnextfilename{result_wave3d_snapshot_r4_t2}\input{tikz/result_wave3d_snapshot_r4_t2.tex}
  }\hfill
  \subcaptionbox{$\pr=4$, $t=0.2$.\label{fig:results_wave3d_snapshot_r4_t3}}
  [0.32\textwidth]{
    \tikzsetnextfilename{result_wave3d_snapshot_r4_t3}\input{tikz/result_wave3d_snapshot_r4_t3.tex}
  }
  \caption{Flux density of $\+B=\curlt\+A$ through the poloidal cut $\logical{z}=0$ for the 3D wave equation \cref{eqn:results_wave3d} on $(8,16,64)$ cells.
  The color range is symmetric and equal for both regularities at a given time, and is chosen according to the solution for $\pr=4$.}
  \label{fig:results_wave3d_snapshots}
\end{figure}
The initial $2$-form $\+B=\curlt\+A$ is a Gaussian pulse just off the magnetic axis, made discretely divergence-free by an orthogonal projection onto the range of the curl and lifted to a relative $1$-form $\+A(0)$ with $\curlt\+A(0)=\+B(0)$ and $\partial_t\+A(0)=\+0$; since gradients lie in the kernel of the curl--curl operator, $\+B(t)$ does not depend on the choice of the lift.
The semidiscrete problem is integrated in time to machine precision by a Chebyshev expansion of $\cos(t\sqrt{\cdot})$ \citep{talezer1986,kosloff1989rem}.
Snapshots of the flux density of $\+B$ through the poloidal cut $\logical{z}=0$ through the pulse, on the coarsest mesh $N=8$, are shown in \cref{fig:results_wave3d_snapshots}.
For $\pr=1$, artifacts emerge at the magnetic axis and grow in time, whereas for $\pr=4$ the wave fronts pass the axis unaffected.
\begin{figure}
  \centering
  \subcaptionbox{Against the mesh size $h=1/N$.\label{fig:results_wave3d_convergence_h}}
  [\twofigwidth]{
    \tikzsetnextfilename{result_wave3d_convergence_linf}
\begin{tikzpicture}
\begin{axis}[height={0.29\textwidth}, width={0.33\textwidth}, scale only axis, xmode={log}, ymode={log}, xmin={0.029157280985525225}, xmax={0.13397168281703667}, ymin={0.010000000000000002}, ymax={10.0}, log basis x={2}, title={$\mathfrak{r} = {\color{lines2}\mathbf{1}}, {\color{lines3}\mathbf{2}}, {\color{lines4}\mathbf{3}}, {\color{lines5}\mathbf{4}}$}, xlabel={$h$}, ylabel={$\|\bm B_h - \bm B_{\mathrm{ref}}\|_{\funspace{L}_\infty} / \|\bm B_{\mathrm{ref}}\|_{\funspace{L}_\infty}$}]
    \addplot[mark={*}, mark size={1.5pt}, color={lines2}, line width={1pt}]
        coordinates {
            (0.125,8.180180595863044)
            (0.08333333333333333,5.097874912404686)
            (0.0625,1.6661494272060835)
        }
        ;
    \addplot[mark={*}, mark size={1.5pt}, color={lines3}, line width={1pt}]
        coordinates {
            (0.125,2.4106217544290645)
            (0.08333333333333333,0.8740413781539108)
            (0.0625,0.6755304709059049)
            (0.041666666666666664,0.15989884038753546)
            (0.03125,0.024069013975821604)
        }
        ;
    \addplot[mark={*}, mark size={1.5pt}, color={lines4}, line width={1pt}]
        coordinates {
            (0.125,0.6744395811920213)
            (0.08333333333333333,0.5035103506193196)
            (0.0625,0.31532442587450676)
            (0.041666666666666664,0.09489643113130106)
            (0.03125,0.023869394190008553)
        }
        ;
    \addplot[mark={*}, mark size={1.5pt}, color={lines5}, line width={1pt}]
        coordinates {
            (0.125,0.8352246494818151)
            (0.08333333333333333,0.6980898052754995)
            (0.0625,0.46935520485136323)
            (0.041666666666666664,0.17551939448578552)
            (0.03125,0.05282659321359969)
        }
        ;
    \addplot[mark={none}, color={gray}, line width={0.6pt}, forget plot]
        coordinates {
            (0.031249999999999997,0.01884369886664526)
            (0.041666666666666664,0.05955539394890352)
            (0.041666666666666664,0.01884369886664526)
            (0.031249999999999997,0.01884369886664526)
        }
        ;
    \node at (axis cs:0.041666666666666664, 0.03349990909625824) [color=gray,right,font=\footnotesize] {$\mathcal{O}(h^p)$};
\end{axis}
\end{tikzpicture}
  }\hfill
  \subcaptionbox{Against time on the coarsest mesh.\label{fig:results_wave3d_convergence_t}}
  [\twofigwidth]{
    \tikzsetnextfilename{result_wave3d_convergence_linf_time_nr8}
\begin{tikzpicture}
\begin{axis}[height={0.29\textwidth}, width={0.33\textwidth}, scale only axis, ymode={log}, ymin={0.010000000000000002}, ymax={10.0}, xlabel={$t$}, title={$\mathfrak{r} = {\color{lines2}\mathbf{1}}, {\color{lines3}\mathbf{2}}, {\color{lines4}\mathbf{3}}, {\color{lines5}\mathbf{4}}$}, ylabel={$\|\bm B_h - \bm B_{\mathrm{ref}}\|_{\funspace{L}_\infty} / \|\bm B_{\mathrm{ref}}\|_{\funspace{L}_\infty}$}]
    \addplot[mark={none}, color={lines2}, line width={1pt}]
        coordinates {
            (0.0,1.0849821458736748)
            (0.00625,0.6501039333285057)
            (0.0125,0.22726108004081116)
            (0.01875,0.8234298264347714)
            (0.025,0.8185332380113072)
            (0.03125,0.4351418546361931)
            (0.0375,0.5378379129219178)
            (0.04375,0.780642684755575)
            (0.05,0.6269834310252403)
            (0.05625,0.37911333973923467)
            (0.0625,0.5296879549708086)
            (0.06875,0.7438339152889074)
            (0.075,0.8481228025465952)
            (0.08125,0.3528242414398438)
            (0.0875,0.6457608607643743)
            (0.09375,0.9098560196210822)
            (0.1,0.8431202564192548)
            (0.10625,0.5849741796631953)
            (0.1125,0.7037899916993519)
            (0.11875,0.8585467087014388)
            (0.125,0.7541494324904036)
            (0.13125,0.455376699968967)
            (0.1375,0.6977235948479112)
            (0.14375,0.8455379995417308)
            (0.15,0.5756209661740809)
            (0.15625,0.34817816878351315)
            (0.1625,0.7963058174238068)
            (0.16875,0.6627718069810656)
            (0.175,0.3256523307774289)
            (0.18125,0.7182335434056882)
            (0.1875,0.8795066200261351)
            (0.19375,0.5451125220205011)
            (0.2,0.3406146792240938)
        }
        ;
    \addplot[mark={none}, color={lines3}, line width={1pt}]
        coordinates {
            (0.0,0.4603533066289436)
            (0.00625,0.261171051375738)
            (0.0125,0.13696031516352733)
            (0.01875,0.3574335752169739)
            (0.025,0.2919838779779027)
            (0.03125,0.18751588991853008)
            (0.0375,0.21858918539285374)
            (0.04375,0.21141160267241774)
            (0.05,0.17357153080411303)
            (0.05625,0.13885130767862305)
            (0.0625,0.18175788738667636)
            (0.06875,0.25034256587819814)
            (0.075,0.2289943890799701)
            (0.08125,0.1254880644842869)
            (0.0875,0.0991584222874753)
            (0.09375,0.16270372791578902)
            (0.1,0.24143871867671465)
            (0.10625,0.23524158605125176)
            (0.1125,0.14734262447935326)
            (0.11875,0.07895957634572696)
            (0.125,0.13572098301488572)
            (0.13125,0.2327743030068769)
            (0.1375,0.25224068454033194)
            (0.14375,0.18734353067780352)
            (0.15,0.0667766834387999)
            (0.15625,0.08916937472160291)
            (0.1625,0.20786450024436257)
            (0.16875,0.25918227190138127)
            (0.175,0.2265648594235358)
            (0.18125,0.12162427814986611)
            (0.1875,0.069487448586807)
            (0.19375,0.15689080928836843)
            (0.2,0.23819733064917722)
        }
        ;
    \addplot[mark={none}, color={lines4}, line width={1pt}]
        coordinates {
            (0.0,0.6830454943345092)
            (0.00625,0.45289974228268837)
            (0.0125,0.14686628498326074)
            (0.01875,0.35626179300473054)
            (0.025,0.39398813994946646)
            (0.03125,0.26049810748431107)
            (0.0375,0.19061336995421804)
            (0.04375,0.19195486773800718)
            (0.05,0.18522470466219484)
            (0.05625,0.15892789552914713)
            (0.0625,0.12098354841417802)
            (0.06875,0.09952159468574026)
            (0.075,0.10125249598350268)
            (0.08125,0.0980328615131676)
            (0.0875,0.09996893337115365)
            (0.09375,0.094880065335754)
            (0.1,0.09031158788285858)
            (0.10625,0.08331010867454844)
            (0.1125,0.07750074585659535)
            (0.11875,0.0768726219461265)
            (0.125,0.07798874688304613)
            (0.13125,0.07915799235081926)
            (0.1375,0.07992887071594651)
            (0.14375,0.07696297794443599)
            (0.15,0.07207573828517544)
            (0.15625,0.06847587667153576)
            (0.1625,0.06960172177087691)
            (0.16875,0.07058599318532871)
            (0.175,0.0710737889764628)
            (0.18125,0.07236173839115324)
            (0.1875,0.06950889242578685)
            (0.19375,0.06510634693737356)
            (0.2,0.07251356733689013)
        }
        ;
    \addplot[mark={none}, color={lines5}, line width={1pt}]
        coordinates {
            (0.0,0.8012270638526946)
            (0.00625,0.5593144817172387)
            (0.0125,0.16634420819974505)
            (0.01875,0.3223687550135702)
            (0.025,0.40820551868543714)
            (0.03125,0.30762373789028347)
            (0.0375,0.19749132230907346)
            (0.04375,0.1817275909425341)
            (0.05,0.16788141724622382)
            (0.05625,0.14662595742216675)
            (0.0625,0.12608071331254464)
            (0.06875,0.11448681815950462)
            (0.075,0.1124571835097536)
            (0.08125,0.114640712751865)
            (0.0875,0.11533754746646098)
            (0.09375,0.10987407204384195)
            (0.1,0.10210918893728449)
            (0.10625,0.09825362612778088)
            (0.1125,0.09767259800303808)
            (0.11875,0.09755274425042208)
            (0.125,0.095713504798792)
            (0.13125,0.09519301966775988)
            (0.1375,0.09071514146846046)
            (0.14375,0.09046425327483731)
            (0.15,0.08980065902209616)
            (0.15625,0.0872079659043836)
            (0.1625,0.08664540260328672)
            (0.16875,0.08630890974102771)
            (0.175,0.08926762592565093)
            (0.18125,0.0872558456315927)
            (0.1875,0.08069106603540681)
            (0.19375,0.08060653840136256)
            (0.2,0.08123235398078463)
        }
        ;
    \node[anchor=south east] at (rel axis cs:0.97,0.03) {$(N^{\mathrm{x}}, N^{\mathrm{y}}, N^{\mathrm{z}}) = (8, 16, 64)$};
\end{axis}
\end{tikzpicture}
  }
  \caption{Relative maximum error of $\+B=\curlt\+A$ for the 3D wave equation \cref{eqn:results_wave3d}.}
  \label{fig:results_wave3d_convergence}
\end{figure}
\Cref{fig:results_wave3d_convergence_h} shows the relative maximum error during the final $10\%$ of the time interval under mesh refinement, compared to a reference solution using $N = 40$ with $\pr = 4$.
Note that the simulations using $\pr = 1$ and $N = 24, 32$ were skipped intentionally because the stiff eigenmodes originating at the center (cf.\ \cref{fig:results_spectrum_lambda_max}) make the time integration prohibitively expensive.
We find that asymptotically all regularity degrees achieve similar accuracy, but for poorly resolved pulses the high-regularity solution is significantly more accurate than the merely conforming solution.
In \cref{fig:results_wave3d_convergence_t} we show the evolution of the relative maximum error as function of time on the coarsest mesh.

\subsection{Guiding-center orbits through the polar axis}\label{sec:results_vacuum}
The last experiment asks when $\pr$-compatibility is needed.
The wave equation only requires the fields to be square integrable, whereas the equations of motion of a charged particle in the guiding-center approximation involve the pointwise evaluations of the magnetic field $\+B$ and its first and second derivatives.
We therefore trace guiding-center orbits that pass close to the polar axis $\Gamma$ in a discrete vacuum magnetic field, and compare the $\pr$-compatible logical spaces of \cref{thm:polar_subspaces_regularity_3d} with logical spaces that are not $\pr$-compatible.

The 2D logical mapping is now simply the unit disk of \cref{ex:disk_mapping} ($R_1=1$), and, following the design rule of \cref{rem:pr_compatibility_sufficient}, all shaping is done by the intermediate mapping: the tube mapping of \cref{ex:tube_mapping_3d} with the helical axis and rotating frame of \cref{sec:results_wave3d} and the cross-section map $\intermediatemapping_\mathrm{2D}(\intermediatecoord{\+x}_\mathrm{2D})=0.3\bigl(1+0.2\intermediatecoord{x}(\intermediatecoord{x}^2-3\intermediatecoord{y}^2)\bigr)\intermediatecoord{\+x}_\mathrm{2D}$, a smooth three-petal deformation of the disk of radius $0.3$ of the type in \cref{eqn:pr_compatibility_sufficient_deformation}.
We use degree-$4$ B-splines on $N=8$ radial and $8N$ toroidal cells, the polar spaces with homogeneous essential boundary conditions for $\pr=0,\ldots,4$, and we compare two angular spaces: a trigonometric basis with $n^\yy=9$ modes, which is $\pr$-compatible with the disk mapping for all $\pr\le4$ (\cref{ex:pr_compatibility_disc_mapping}), and periodic B-splines of degree $4$ on $2N$ cells, which are merely $0$-compatible.
As in \cref{sec:results_spectrum}, the $\pr=0$ spaces are not conforming and the inner products in \cref{eqn:results_vacuum_field} are then quadrature-regularized; we include them only as the nonconforming reference.
The discrete vacuum magnetic field is the harmonic $2$-form of the polar complex with homogeneous essential boundary conditions, normalized to unit toroidal flux:
\begin{equation}\label{eqn:results_vacuum_field}
  \text{Find} \quad \+B\in \polar{\funspace{\Bsplinethree}}^{2,\pr}_\bc
  \quad \text{s.t.} \quad
  \divergencet\+B=0,
  \quad
  \inprod{\curlt\+\Bsplinethree^1}{\+B}=0
  \quad
  \forall \+\Bsplinethree^1\in\polar{\funspace{\Bsplinethree}}^{1,\pr}_\bc,
  \quad
  \int_\Sigma\+B\cdot\+n\diff S=1,
\end{equation}
where $\Sigma$ is a poloidal cross-section.
The relative polar complex is the tensor product of the relative 2D complex, whose cohomology has the dimensions $0,0,1$ by \cref{lem:boundary_compatible_homotopy}, with the periodic 1D complex, so that the K\"unneth formula, as in the proof of \cref{thm:polar_3d_subspaces_commutativity_conformity}, gives the dimensions $0,0,1,1$ of the relative de Rham complex of the solid torus.
Hence the discrete harmonic $2$-forms form a one-dimensional space, and $\+B$ is unique.

We integrate the gyrogauge-invariant guiding-center equations of motion of \citet[(3.48a,b)]{remmerswaal2025gauge} (see also \citep{cary2009}), for the position $\+R$ and the parallel velocity $U_\shortparallel$ of a particle with mass $m$, charge $q$ and magnetic moment $M$.
Their right-hand sides involve $\+B$, its gradient and, through the curl of the vector field $\+w_0$ in \citet[(3.32)]{remmerswaal2025gauge}, which contains the gyrogauge term that makes the equations independent of the choice of the normal frame, its second derivatives, all evaluated pointwise along the orbit.
The guiding-center energy $\mathcal{E} \defeq \tfrac12mU_\shortparallel^2+M|\+B|$ is an exact invariant of these equations and serves as the error indicator below.

\begin{figure}
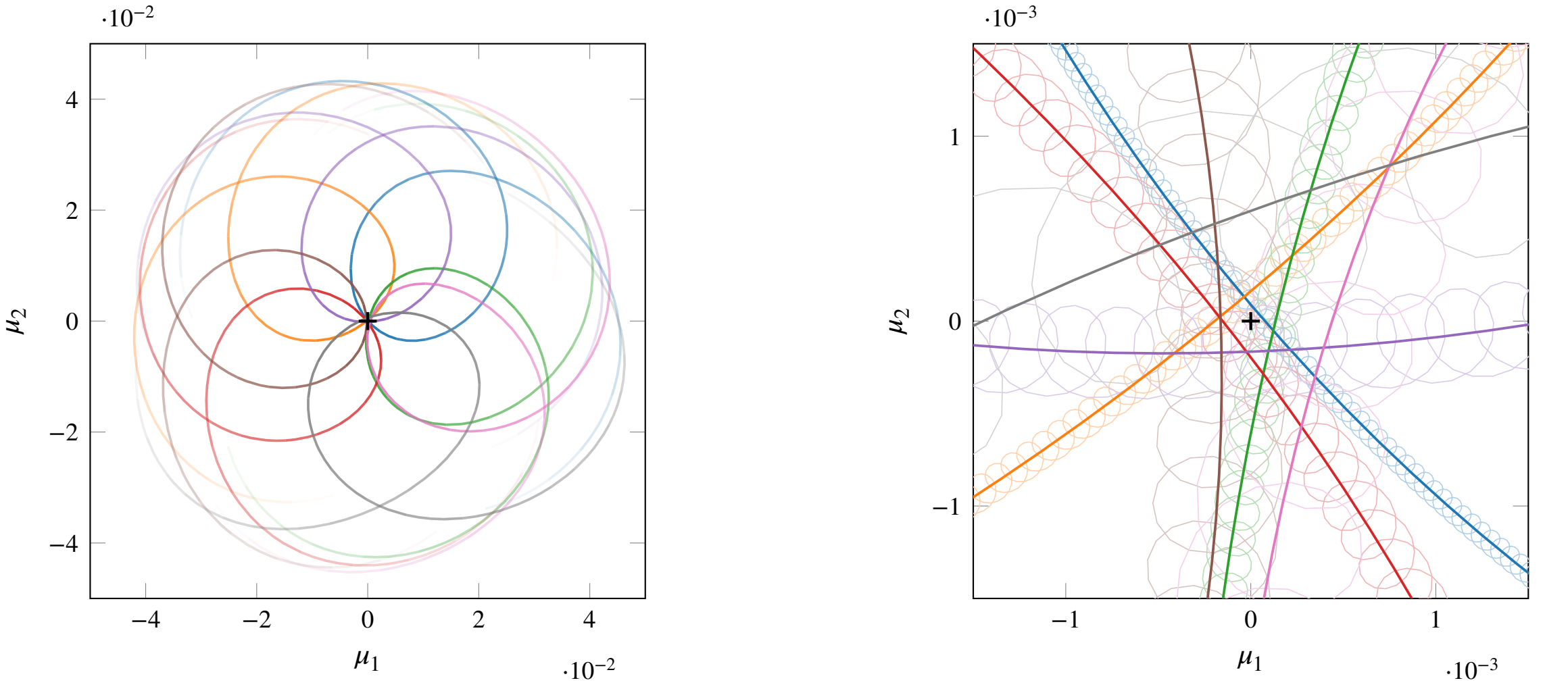

  \centering
  \subcaptionbox{Guiding-center orbits.\label{fig:results_vacuum_orbits_overview}}
  [\twofigwidth]{
    \tikzsetnextfilename{result_vacuum_orbits_overview}\input{tikz/result_vacuum_orbits_overview.tex}
  }\hfill
  \subcaptionbox{Close-up around the polar axis, with the full orbits.\label{fig:results_vacuum_orbits_closeup}}
  [\twofigwidth]{
    \tikzsetnextfilename{result_vacuum_orbits_closeup}\input{tikz/result_vacuum_orbits_closeup.tex}
  }
  \caption{Orbits near their closest passage of the polar axis ($+$) for the trigonometric angular basis and $\pr=4$, in the cross-sectional coordinates $\+\mu$ of \cref{ex:tube_mapping_3d}, one color per orbit: guiding-center orbit for $\Delta t=0.005$ (solid) and, in (b), full orbit (light).
  }
  \label{fig:results_vacuum_orbits}
\end{figure}
We use the classical fourth-order Runge--Kutta method for the integration of the equations of motion.
The particles have mass $m=1$, initial parallel velocity $U_\shortparallel=1$ and perpendicular velocity $0.3$.
We start eight such strongly passing particles near the magnetic axis, with charges $q$ from $1369$ to $137$ chosen such that $\varepsilon_B=\rho|\gradient|\+B||/|\+B|$ at the start, with $\rho$ the gyroradius, is log-uniformly distributed between $10^{-4}$ and $10^{-3}$; larger values of $\varepsilon_B$ result in orbits wherein the perpendicular drifts play a more significant role.
The starting points are chosen such that each guiding-center orbit passes $\Gamma$ at about one gyroradius, and the particles are followed for $20$ toroidal transits.
Example orbits are shown in \cref{fig:results_vacuum_orbits_overview}.
As a reference, each particle is also integrated as a full orbit with the Boris scheme~\citep{boris1970relativistic} using $32$ steps per gyration, as illustrated in \cref{fig:results_vacuum_orbits_closeup}.

\begin{figure}
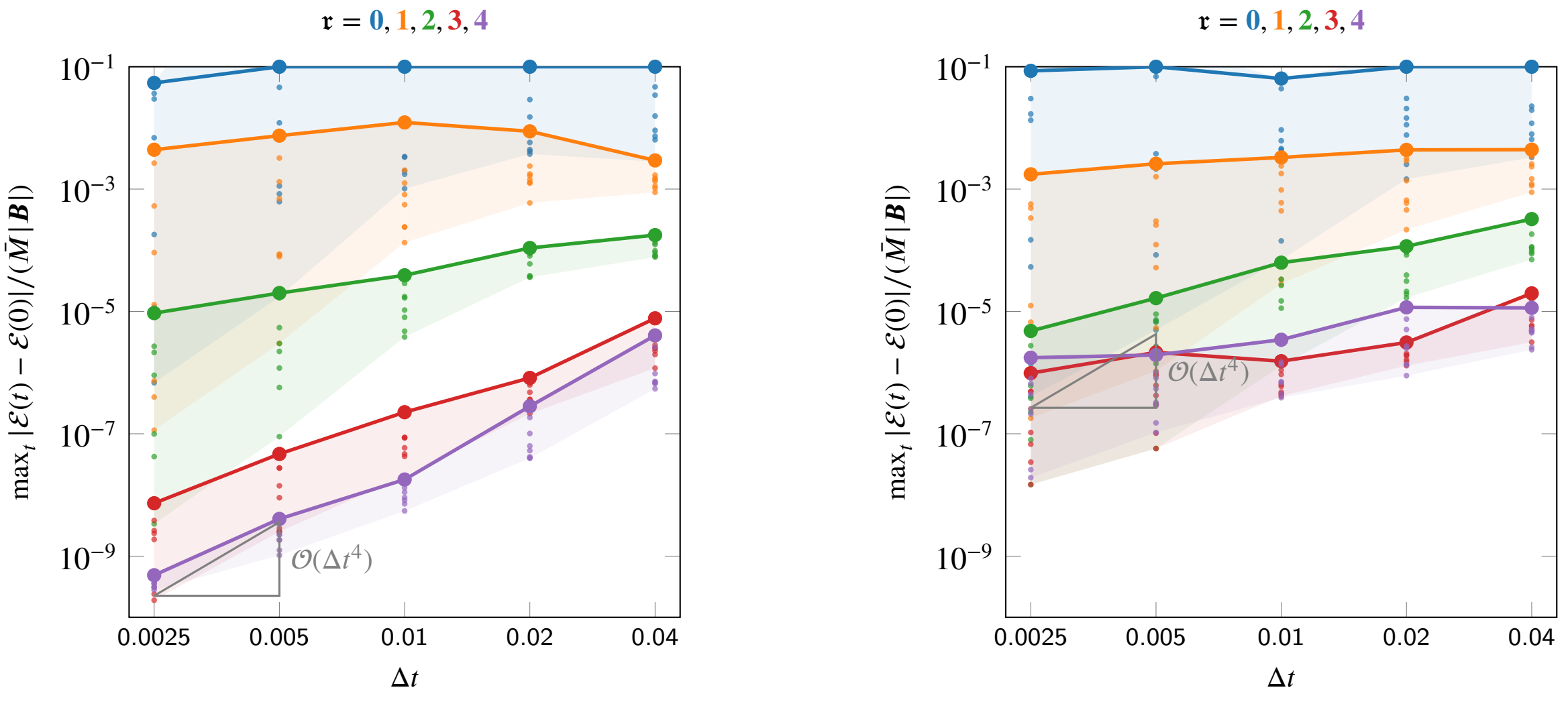

  \centering
  \subcaptionbox{Trigonometric angular basis ($\pr$-compatible).\label{fig:results_vacuum_energy_compatible}}
  [\twofigwidth]{
    \tikzsetnextfilename{result_vacuum_orbits_energy_compatible}\input{tikz/result_vacuum_orbits_energy_compatible.tex}
  }\hfill
  \subcaptionbox{B-spline angular basis (not $\pr$-compatible).\label{fig:results_vacuum_energy_bspline}}
  [\twofigwidth]{
    \tikzsetnextfilename{result_vacuum_orbits_energy_bspline}\input{tikz/result_vacuum_orbits_energy_bspline.tex}
  }
  \caption{Maximum energy error of the guiding-center orbits over $20$ toroidal transits against the time step $\Delta t$, relative to $M|\+B|$; lines show the maximum over the eight orbits, small marks the individual orbits and shaded bands their range.
  Runs that exceed the plotted range ($\pr=0$) are drawn at its upper limit.}
  \label{fig:results_vacuum_energy}
\end{figure}
Since the guiding-center energy is an invariant, it can be used as an error indicator for the time integration.
In \cref{fig:results_vacuum_energy} we show the energy error as function of time-step refinement for the eight selected orbits.
The lines indicate the maximum error over the eight orbits, and the shaded envelope covers the range of values.
With $\pr$-compatibility we find in \cref{fig:results_vacuum_energy_compatible} that the energy converges for all orbits with the expected fourth-order rate only for $\pr=4$.
This can be expected because the magnetic field is guaranteed to be twice continuously differentiable at the center only for $\pr \ge 4$ according to \cref{thm:polar_subspaces_regularity_3d}.
When $\pr$-compatibility is not satisfied we find that the error does not converge for any of the regularity degrees, as shown in \cref{fig:results_vacuum_energy_bspline}.

  \section{Conclusion}\label{sec:conclusion}
We have constructed polar subspaces of the tensor-product spaces on disk-like domains and on solid tori that are conforming, form a cochain complex with the cohomology of the continuous de Rham complex, and whose $m$-forms are $\funspace{C}^{\pr-\min(m, 2)}$ at the center of the disk, and correspondingly regular along the central curve of the torus, for every regularity degree $\pr$ up to the radial spline degree (\cref{thm:polar_subspaces_conformity,thm:polar_subspaces_regularity,thm:polar_subspaces_regularity_3d}).
The construction resolves the obstruction identified in \cref{sec:homotopy_motivation}: if the polar $1$-forms near the center are only gradients of polar $0$-forms, then no polar $2$-form can be nonzero there, so that densities and fluxes cannot be represented at the center with any regularity.
Rather than starting from the $0$-forms and differentiating, we enrich the $2$-form space with projected polynomials and integrate backwards, lifting them to curl-like $1$-forms with the Poincar\'e homotopy operators that integrate along rays from the center.
The decisive property is that on clamped radial B-splines these operators are cumulative sums of coefficients (\cref{lem:pseudo_inverse_curl_logical}), so that the homotopy identities hold exactly on the discrete level and the lifts need no further projection (\cref{cor:discrete_homotopy_identities,thm:polar_subspaces_discrete_homotopy}).
The polar spaces are thereby a complex in both directions, a cochain complex under the differential operators (\cref{thm:polar_subspaces_commutativity_conformity}) and a chain complex under the homotopy operators, which splits the polar $1$-form space into gradients and lifts and yields exactness by an algebraic argument (\cref{thm:polar_subspaces_homotopy_splitting,cor:polar_subspaces_cohomology}).

The homotopy identity also defines the polar $1$-form projection in terms of the $0$- and $2$-form projections (\cref{eqn:polar_oneform_projection_operator}).
As a consequence, all polar projections commute with the differential operators by construction, commute with the homotopy operators, act locally on the internal coefficient rings, and preserve the total flux and circulation (\cref{thm:full_commuting_polar_projection_operators,thm:homotopy_commuting_polar_projection_operators,lem:polar_projection_operators_conservation_properties}).
Adapted to homogeneous essential boundary conditions, the homotopy operators satisfy boundary-compatible homotopy identities (\cref{lem:boundary_compatible_homotopy}), provide explicit vector potentials and reduce the curl--curl problem to a single inversion of the $2$-form mass matrix (\cref{sec:bvp_2d}).
In 3D, the polar complex is the tensor product of the 2D polar complex with a periodic 1D complex, so that conformity, regularity along the central curve and the cohomology of the solid torus (\cref{thm:polar_3d_subspaces_commutativity_conformity}), the homotopy operators (\cref{thm:physical_3d_homotopy_identities}) and the commuting projections (\cref{thm:commuting_polar_projection_operators_3d}) are inherited without a separate construction.
Throughout, the physical mapping is decomposed into a disk-like logical mapping, which alone determines the extraction matrices, and a smooth diffeomorphism, which carries the geometry and enters only through the metric.
Regularity is guaranteed under $\pr$-compatibility of the logical mapping with the logical spaces (\cref{def:pr_compatibility_mapping_logical_spaces}), which suggests the design rule of \cref{rem:pr_compatibility_sufficient}: keep the logical mapping as simple as possible, ideally the unit disk, and let the diffeomorphism do the shaping, so that, for instance, numerical MHD equilibria can be used directly as physical mappings (\cref{rem:practical_applicability_3d}).
All spaces and operators have explicit matrix representations in terms of Kronecker products and sparse extraction matrices (\cref{sec:matrix_representation}), and are implemented in the open-source package \polarfeec.

The numerical experiments separate the properties that require $\pr$-compatibility from those that do not.
The commuting projections converge at the optimal rates for every $\pr$, and the regularity guaranteed by \cref{thm:polar_subspaces_regularity} is attained exactly and is sharp (\cref{sec:results_projections,sec:results_regularity}).
The structural properties hold for any disk-like mapping and any angular basis, and so does the suppression of the stiff modes originating at the center: with periodic angular B-splines, which are only $0$-compatible, the largest eigenvalue of the curl--curl operator decreases by a factor of about $130$ from the first-order polar spaces to $\pr=5$, which allows a time step about $11$ times larger, while the smooth eigenvalues keep converging at their theoretical rate (\cref{sec:results_spectrum}), and under-resolved waves cross the axis of a stellarator-like domain without the artifacts at the axis that appear for the first-order spaces (\cref{sec:results_wave3d}).
Pointwise regularity at the center, in contrast, does require compatibility: guiding-center orbits in a discrete vacuum field, whose equations of motion evaluate $\+B$ and its first two derivatives along the orbit, converge at the expected rate only with a compatible angular basis and $\pr\ge4$, as predicted by \cref{thm:polar_subspaces_regularity_3d} (\cref{sec:results_vacuum}).
Since $\pr=1$ reproduces the first-order polar complexes of the literature, these comparisons also quantify what the additional regularity buys.

Several questions remain open.
A promising route to compatibility without the global coupling of a trigonometric basis is to generalize the Tchebycheffian B-spline de Rham complexes of \citet{manni2026tbsplines} from the interval to $\mathbb{S}^1$: this yields a compactly supported periodic 1D de Rham complex that contains the trigonometric functions up to a prescribed frequency, so that the angular spaces would be $\pr$-compatible with the disk mapping while keeping the mass and stiffness matrices sparse, removing the trade-off described in \cref{sec:one_dimensional_function_spaces}.
The polar projections are oblique and are defined on fields with well-defined geometric degrees of freedom. 
Their stability and approximation properties have not been analyzed, so that the optimal rates observed in \cref{sec:results_projections} await a proof, and it remains open whether bounded commuting projections in the sense of \citet{arnold2010feec,falk2014} can be derived from them.
The boundary-compatible homotopy operators and the resulting mass-matrix solvers extend naturally to the relative complex of the solid torus, where the grad--div and curl--curl problems can be treated in the same way.
Finally, the proposed construction is a natural ingredient
for structure-preserving plasma models \citep{kraus2017gempic,holderied2021,possanner2023struphy}, where the regularity at the magnetic axis established here is precisely what field evaluation and particle pushing near the axis require.

  \printcredits

  \section*{Declaration of competing interest}
  The author declares no known competing financial interests or personal relationships that could have appeared to influence the work reported in this paper.

  \section*{Data availability}
  All results were obtained with version~0.3.0 of the open-source Julia package \polarfeec, which contains the code to reproduce them.

  \section*{Declaration of generative AI and AI-assisted technologies in the writing process}
  Generative AI tools were used for proof verification and formalization, as well as for code development.
  After using these tools, the author reviewed and edited the content as needed and takes full responsibility for the content of the publication.

  \section*{Acknowledgments}
  The author would like to thank Roman Hatzky, Martin Campos Pinto and Florian Hindenlang for fruitful discussions that have improved the quality of this work.

  \appendix
  \section{Notation}\label{app:notation}
\Cref{tab:notation,tab:notation_polar} collect the notation used throughout the paper, together with the sections, equations and statements in which it is introduced.

\begin{table}[htbp]
  \centering
  \footnotesize
  \setlength{\tabcolsep}{3pt}
  \begin{tabular}{@{}l>{\raggedright\arraybackslash}p{7cm}l@{}}
    \toprule
    Symbol & Meaning & Introduced in\\
    \midrule
    \multicolumn{3}{@{}l}{\itshape Conventions}\\
    $\logical{\square}$, $\intermediate{\square}$, $\square$ & logical, intermediate and physical version of $\square$ & \cref{sec:coordtransform,sec:physical_function_spaces}\\
    $\polar{\square}$ & polar (regular at the center) space, basis or operator & \cref{sec:zeroforms}\\
    $\reduced{\square}$, $\enriched{\square}$ & reduced (external) and enriched (internal) operator & \cref{sec:derham_hodge_2d_geom_projections}\\
    $\square_\mathrm{2D}$ & 2D object reused in the 3D construction & \cref{sec:high_order_derham_hodge_3d}\\
    \midrule
    \multicolumn{3}{@{}l}{\itshape Domains and mappings}\\
    $\logical{\Omega}$, $\intermediate{\Omega}$, $\Omega$ & logical, intermediate and physical domain & \cref{sec:coordtransform}\\
    $\logical{x}$, $\logical{y}$, $\logical{z}$ & radial, angular and toroidal logical coordinates & \cref{sec:coordtransform,sec:coordtransform_3d}\\
    $\logicalmapping$ & disk-like logical mapping & \cref{eqn:disk_mapping}\\
    $\intermediatemapping$ & smooth diffeomorphism onto the physical domain & \cref{ass:diffeomorphism}\\
    $\mapping$ & physical mapping $\intermediatemapping\circ\logicalmapping$ & \cref{eqn:composed_mapping}\\
    $R$, $C$, $S$ & radius and angular functions of $\logicalmapping$ & \cref{eqn:disk_mapping}\\
    $R_1$ & limit of $R/\logical{x}$ at the center & \cref{ass:disc_like_mapping}\\
    $\jacmat$, $\jac$ & Jacobian matrix and Jacobian of a mapping & \cref{eqn:disk_jacobian_matrix,eqn:composed_jacobian}\\
    $\Gamma$ & central curve of the 3D domain & \cref{eqn:central_curve_3d}\\
    \midrule
    \multicolumn{3}{@{}l}{\itshape One-dimensional spaces ($m\in\{0,1\}$, $\alpha\in\{\xx,\yy,\zz\}$)}\\
    $\funspace{\Bsplineone}^{m,\alpha}$ & 1D $m$-form space in direction $\alpha$ & \cref{sec:one_dimensional_function_spaces}\\
    $\Bsplineone^{m,\alpha}_i$ & its basis functions & \cref{sec:one_dimensional_function_spaces}\\
    $n^{m,\alpha}$, $r^{m,\alpha}$ & dimension and B-spline degree of $\funspace{\Bsplineone}^{m,\alpha}$ & \cref{sec:one_dimensional_function_spaces}\\
    $\indexone^{m,\alpha}$ & index set $\{1,\ldots,n^{m,\alpha}\}$ & \cref{eqn:1d_indexing_sets}\\
    $\xi^\alpha_i$ & Greville points & \cref{sec:geometric_dofs_1d}\\
    $\gdof^{m,\alpha}_i$ & 1D geometric degrees of freedom (GDOFs) & \cref{eqn:geometric_dofs_1d}\\
    $\pi^{m,\alpha}$ & 1D commuting projectors & \cref{eqn:commuting_projectors_1d}\\
    $\delta_\alpha$ & derivative acting on coefficients & \cref{eqn:diff_opers_1d}\\
    $\delta^\dagger_\xx$ & radial cumulative sum of coefficients & \cref{eqn:discrete_radial_homotopy_operator}\\
    $\mathds{1}^\yy$ & coefficients of the constant function $1$ & \cref{eqn:constant_function_coefficients}\\
    \midrule
    \multicolumn{3}{@{}l}{\itshape Tensor-product spaces ($m\in\{0,1,2\}$)}\\
    $\funspace{\Bspline}^m$ & tensor-product $m$-form space & \cref{sec:logical_function_spaces,sec:physical_function_spaces}\\
    $\Bspline^0_{ij}$ & basis of $\funspace{\Bspline}^0$; $i$ radial and $j$ angular index & \cref{eqn:logical_zeroforms}\\
    $\+\Bspline^1_{\alpha,ij}$ & basis of $\funspace{\Bspline}^1$, component $\alpha\in\{\xx,\yy\}$ & \cref{eqn:logical_oneforms}\\
    $\Bspline^2_{ij}$ & basis of $\funspace{\Bspline}^2$ & \cref{eqn:logical_twoforms}\\
    $\indextwo^m$ & index sets of the tensor-product spaces & \cref{eqn:2d_indexing_sets}\\
    $P$, $\+A$, $B$, $Q$ & $0$-, $1$-, $2$- and $3$-forms (3D) & \cref{sec:logical_function_spaces,sec:high_order_derham_hodge_3d}\\
    $p_{ij}$, $a_{\alpha,ij}$, $b_{ij}$ & their coefficients & \cref{sec:logical_function_spaces}\\
    $\gradienth$, $\curlh$ & gradient and curl acting on coefficients & \cref{eqn:discrete_operators}\\
    $\pushforward^m$ & pushforward (Piola transformation) & \cref{eqn:pushforward_operators}\\
    $\pushforward^{-m}$ & pullback, the inverse of $\pushforward^m$ & \cref{eqn:pushforward_composition}\\
    $\gdof^m$ & GDOFs of the tensor-product spaces & \cref{eqn:logical_geometric_dofs,eqn:physical_geometric_dofs}\\
    $\pi^m$ & commuting projectors onto $\funspace{\Bspline}^m$ & \cref{eqn:logical_commuting_projectors,eqn:physical_geometric_dofs}\\
    $\funspace{H}_1$, $\funspace{H}_{\curlt}$, $\funspace{H}_{\divergencet}$, $\funspace{L}_2$ & Sobolev spaces of the de Rham complex & \cref{sec:logical_bspline_spaces,sec:high_order_derham_hodge_3d}\\
    $\funspace{C}^k_0$, $\funspace{C}^k_\Gamma$ & $k$ times differentiable near the center, resp.\ $\Gamma$ & \cref{sec:zeroforms,sec:derham_hodge_3d_properties}\\
    \bottomrule
  \end{tabular}
  \caption{Notation, part I: conventions, mappings, one-dimensional and tensor-product spaces.}
  \label{tab:notation}
\end{table}

\begin{table}[htbp]
  \centering
  \footnotesize
  \setlength{\tabcolsep}{3pt}
  \begin{tabular}{@{}l>{\raggedright\arraybackslash}p{7cm}l@{}}
    \toprule
    Symbol & Meaning & Introduced in\\
    \midrule
    $\pr$ & regularity degree of the polar spaces & \cref{ass:logical_space_parameters}\\
    ring $i$ & all basis functions, or coefficients, with radial index $i$ & \cref{sec:zeroforms}\\
    $\funspace{\Bspline}^{m,\pr}_\internal$, $\funspace{\Bspline}^{m,\pr}_\external$ & internal (first rings) and external part of $\funspace{\Bspline}^m$ & \cref{eqn:polar_zeroform_space_odot,eqn:polar_mform_space_odot}\\
    $\indextwo^{m,\pr}_\internal$, $\indextwo^{m,\pr}_\external$ & internal and external index sets & \cref{eqn:polar_zeroform_space_odot,eqn:polar_mform_space_odot}\\
    $\chi^{m,\pr}_\internal$, $\chi^{m,\pr}_\external$ & truncation to the internal, resp.\ external, coefficients & \cref{sec:derham_hodge_2d_basis}\\
    $\intermediate{\funspace{T}}^k$ & polynomials of total degree at most $k$ in $\intermediatecoord{\+x}$ & \cref{sec:pullback_compatibility}\\
    $\pulledpolys{0}{k}$, $\pulledpolys{2}{k}$ & pushforwards of $\intermediate{\funspace{T}}^k$ as physical $0$- and $2$-forms & \cref{sec:pullback_compatibility,sec:derham_hodge_2d_construction}\\
    $\intermediate{T}_l$ & shifted Zernike basis of $\intermediate{\funspace{T}}^\pr$ & \cref{eqn:shifted_zernike_basis}\\
    $T^0_l$, $T^2_l$ & its pushforwards as physical $0$- and $2$-forms & \cref{sec:derham_hodge_2d_basis}\\
    $\pr$-compatible & polynomial pullbacks lie in the logical spaces & \cref{def:pr_compatibility_mapping_logical_zeroform,def:pr_compatibility_mapping_logical_spaces}\\
    $\polar{\funspace{\Bspline}}^{m,\pr}$ & polar $m$-form spaces & \cref{eqn:polar_zeroform_space,eqn:polar_oneform_space,eqn:polar_twoform_space}\\
    $\polar{\Bspline}^{0,\pr}_l$, $\polar{\Bspline}^{2,\pr}_l$ & internal polar $0$- and $2$-form basis functions & \cref{eqn:polar_zeroform_basis,eqn:polar_twoform_curl_basis}\\
    $\polar{\+\Bspline}^{1,\pr}_{\mathrm{G},l}$, $\polar{\+\Bspline}^{1,\pr}_{\mathrm{C},l}$ & gradient- and curl-like polar $1$-form basis functions & \cref{eqn:polar_oneform_gradient_basis,eqn:polar_oneform_curl_basis_definition}\\
    $n^{0,\pr}$, $n^{1,\pr}_\mathrm{G}$, $n^{1,\pr}_\mathrm{C}$, $n^{2,\pr}$ & numbers of internal polar basis functions & \cref{eqn:polar_subspaces_overview_dimensions}\\
    $\reduced{\funspace{\Bspline}}^{m,\pr}$ & tentative (gradient-only) polar spaces & \cref{sec:homotopy_motivation}\\
    $\reduced{\Pi}^{2,\pr}$, $\reduced{\pi}^{2,\pr}$ & conforming $2$-form projector; composed with $\pi^2$ & \cref{lem:conforming_polar_twoform_projector}\\
    $\gradientt^\dagger$, $\curlt^\dagger$, $\divergencet^\dagger$ & homotopy (Poincar\'e) operators & \cref{eqn:homotopy_operators,def:logical_3d_homotopy_vector_proxies}\\
    $\mathfrak{P}$ & projection onto the value at the center & \cref{sec:homotopy_2d}\\
    $\mathfrak{A}$ & projection onto the circulation along $\Gamma$ (3D) & \cref{def:logical_3d_homotopy_vector_proxies}\\
    $\mathfrak{B}$ & projection onto the mean of a $2$-form & \cref{eqn:boundary_homotopy_operators_2d}\\
    $\polar{\Bspline}^{m,\pr,\perp}_l$ & dual basis functions & \cref{eqn:polar_dual_basis}\\
    $\ellinprod{\cdot}{\cdot}$ & Euclidean inner product of coefficients & \cref{eqn:coefficient_euclidean_inner_product}\\
    $\enriched{\Pi}^{m,\pr}$, $\reduced{\Pi}^{m,\pr}$ & internal and external polar coefficient projectors & \cref{eqn:polar_zero_two_form_projection_operator}\\
    $\polar{\Pi}^{m,\pr}$ & polar coefficient projector & \cref{eqn:polar_zero_two_form_projection_operator,eqn:polar_oneform_projection_operator}\\
    $\polar{\pi}^{m,\pr}$ & polar commuting projection operator $\polar{\Pi}^{m,\pr}\circ\pi^m$ & \cref{eqn:polar_projection_operator_definition}\\
    $\operatorname{tr}^m_\bc$ & traces on the outer boundary & \cref{sec:boundary_homotopy}\\
    $\polar{\funspace{\Bspline}}^{m,\pr}_\bc$ & relative polar spaces (zero traces) & \cref{eqn:relative_polar_spaces}\\
    $\gradientt^\dagger_\bc$, $\curlt^\dagger_\bc$ & boundary-compatible homotopy operators & \cref{eqn:boundary_homotopy_operators_2d}\\
    $\funspace{\Bsplinethree}^m$, $\polar{\funspace{\Bsplinethree}}^{m,\pr}$ & 3D tensor-product and polar spaces & \cref{sec:derham_hodge_3d_construction}\\
    $\mat{\Lambda}^{p,q}$ & placement matrices of the 3D vector proxies & \cref{eqn:tensor_product_proxy_matrices_3d}\\
    $\wedge$ & product of a 2D $p$-form and a toroidal $q$-form & \cref{eqn:tensor_product_complex_spaces_3d_explicit}\\
    \bottomrule
  \end{tabular}
  \caption{Notation, part II: the polar construction.}
  \label{tab:notation_polar}
\end{table}

  \printdeferredproofs
  \section{Matrix representation}\label{sec:matrix_representation}
\begingroup
\emergencystretch=3em
This appendix collects the matrix representations of the spaces and operators used in this paper.
All operators act on coefficient vectors.
The mapping enters only through the mass matrices and through the evaluation of the GDOFs, since the pushforwards in \cref{eqn:pushforward_operators,eqn:pushforward_operators_3d} leave the coefficients unchanged.
\ifshowcodeapi
Names in brackets, such as \codeapi{derivative}{(complex, m)}, refer to the corresponding high-level functions of version~0.3.0 of the Julia package \polarfeec, which implements the constructions of this paper.
\fi

\subsection{Conventions and one-dimensional matrices}\label{sec:matrix_representation_1d}
\paragraph{Ordering and Kronecker products}
The coefficients of a 2D tensor-product form are ordered with the radial index running fastest: the coefficient $p_{ij}$ is stored at position $i+(j-1)n^{\cdot,\xx}$, where $n^{\cdot,\xx}$ is the radial dimension of the respective component.
In 3D, the toroidal index runs slowest.
With this ordering, an operator that acts by $\mat{A}_\xx$ on the radial index and by $\mat{A}_\yy$ on the angular index is represented by the Kronecker product $\mat{A}_\yy\otimes\mat{A}_\xx$.
Indeed, if the coefficients are arranged in a matrix $\mat{P}$ with entries $p_{ij}$, then $(\mat{A}_\yy\otimes\mat{A}_\xx)\operatorname{vec}(\mat{P})=\operatorname{vec}(\mat{A}_\xx\mat{P}(\mat{A}_\yy)^\transpose)$, so that Kronecker products are applied direction by direction and are never formed explicitly.
The components of a 1D or 2D form are stacked, e.g., $\mat{a}=[\mat{a}_\xx;\mat{a}_\yy]$ for a $1$-form.
Matrices representing 1D operators carry the direction as a subscript and the form degree as a superscript, e.g., $\mat{d}_\xx$, $\mat{i}^0_\yy$ and $\mat{h}^1_\zz$, whereas 1D coefficient vectors carry both as superscripts, as in \cref{eqn:constant_function_coefficients}.
A bar denotes a matrix acting on polar coefficients, except for the polar projector matrices $\polar{\mat{\Pi}}^m$, which follow the notation of \cref{sec:derham_hodge_2d_geom_projections} and act on tensor-product coefficients.
We omit the superscript $\pr$ whenever no confusion can arise.

\paragraph{One-dimensional matrices}
For $\alpha\in\{\xx,\yy,\zz\}$ and $m\in\{0,1\}$, we use the following matrices and vectors:
\begin{itemize}
  \item $\mat{i}^{m}_{\alpha}$ is the identity on $\funspace{\Bsplineone}^{m,\alpha}$, and $\mat{e}^{m,\alpha}_i$ is its $i$-th column.
  \item $\mat{d}_\alpha\in\mathbb{R}^{n^{1,\alpha}\times n^{0,\alpha}}$ represents the coefficient derivative $\delta_\alpha$ from \cref{eqn:diff_opers_1d,eqn:trigonometric_coefficient_derivative}\codeapi{derivative_matrix}{(space)}.
  For B-splines, it is the (periodic) forward difference matrix.
  \item $\mat{\mathds{1}}^{\alpha}$ contains the coefficients of the constant $0$-form $1$, see \cref{eqn:bspline_partition,eqn:constant_function_coefficients}, and $\mat{w}^\alpha\defeq\mat{\mathds{1}}^\alpha/((\mat{\mathds{1}}^\alpha)^\transpose\mat{\mathds{1}}^\alpha)$ recovers the value of a constant, $(\mat{w}^\alpha)^\transpose(c\,\mat{\mathds{1}}^\alpha)=c$.
  \item $\matcomp{q}^{m,\alpha}_i\defeq\int_0^1\Bsplineone^{m,\alpha}_i\diff\logical{x}_\alpha$ collects the integrals of the basis functions, so that $(\mat{q}^{1,\alpha})^\transpose\mat{v}$ is the integral of the $1$-form with coefficients $\mat{v}$.
  For B-splines, $\mat{q}^{1,\alpha}$ is the vector of ones.
\end{itemize}

\paragraph{Geometric degrees of freedom}\label{sec:matrix_representation_gdof_projectors}
The GDOF matrices
\begin{equation}\label{eqn:one_dimensional_gdof_vandermonde_matrix}
  \matcomp{h}^{m}_{\alpha,ij}
  \defeq
  \gdof^{m,\alpha}_{i}(\Bsplineone^{m,\alpha}_{j}),
  \qquad m\in\{0,1\},
\end{equation}
are the collocation ($m=0$) and histopolation ($m=1$) matrices of \cref{eqn:geometric_dofs_1d}\codeapi{histopolation_matrix}{(space)}.
If $\mat{g}^{m,\alpha}(F)$ collects the GDOFs of $F$, the coefficient vector of the projection in \cref{eqn:commuting_projectors_1d} is $(\mat{h}^{m}_{\alpha})^{-1}\mat{g}^{m,\alpha}(F)$\codeapi{commuting_projection}{(space, f)}.
The GDOFs of a derivative are differences of nodal values, so \cref{eqn:stokes_1d} reads $\mat{h}^{1}_{\alpha}\mat{d}_\alpha=\mat{d}^\xi_\alpha\mat{h}^{0}_{\alpha}$, where $\mat{d}^\xi_\alpha$ is the (periodic) forward difference matrix on the Greville points.

\paragraph{One-dimensional homotopy matrices}
The radial homotopy operator in \cref{eqn:discrete_radial_homotopy_operator} is the cumulative sum
\begin{equation}\label{eqn:radial_cumulative_sum_matrix}
  \mat{d}^\dagger_\xx
  \defeq
  \begin{bmatrix}
    0      & \cdots & 0\\
    1      & \ddots & \vdots\\
    \vdots & \ddots & 0\\
    1      & \cdots & 1
  \end{bmatrix}
  \in\mathbb{R}^{n^{0,\xx}\times n^{1,\xx}},
  \qquad
  \mat{d}_\xx\mat{d}^\dagger_\xx=\mat{i}^{1}_{\xx},
  \qquad
  \mat{d}^\dagger_\xx\mat{d}_\xx=\mat{i}^{0}_{\xx}-\mat{\mathds{1}}^\xx(\mat{e}^{0,\xx}_1)^\transpose,
\end{equation}
which is the matrix form of \cref{eqn:discrete_radial_homotopy_relations}\codeapi{cumulative_sum_matrix}{(n)}.
In a periodic direction $\alpha\in\{\yy,\zz\}$, the derivative has the constants in its kernel and the $1$-forms of zero integral as its range.
Let $\mat{\mathfrak{A}}_\alpha\defeq\mat{c}^\alpha(\mat{q}^{1,\alpha})^\transpose$ replace a $1$-form by the constant $1$-form with the same integral, where $\mat{c}^\alpha$ represents the constant $1$-form with unit integral.
The periodic homotopy matrix $\mat{d}^\dagger_\alpha$ is then the unique matrix with
\begin{equation}\label{eqn:periodic_primitive_matrices}
  \mat{d}_\alpha\mat{d}^\dagger_\alpha=\mat{i}^{1}_{\alpha}-\mat{\mathfrak{A}}_\alpha,
  \qquad
  \mat{d}^\dagger_\alpha\mat{\mathfrak{A}}_\alpha=\mat{0},
  \qquad
  (\mat{\ell}^\alpha)^\transpose\mat{d}^\dagger_\alpha=\mat{0},
\end{equation}
where the normalization $\mat{\ell}^\alpha$ fixes the additive constant of the primitive: the poloidal homotopy operator $\partial^\dagger_{\logical{y}}$ in \cref{eqn:boundary_homotopy_operators_2d} has zero mean, $\mat{\ell}^\yy\defeq\mat{q}^{0,\yy}$, and the toroidal homotopy operator in \cref{eqn:logical_3d_gradient_homotopy_vector_proxy} vanishes at $\xi^\zz_1$, $\mat{\ell}^\zz\defeq(\mat{h}^{0}_{\zz})^\transpose\mat{e}^{0,\zz}_1$.
For periodic B-splines, $\mat{d}^\dagger_\alpha$ is computed by applying $\mat{i}^{1}_{\alpha}-\mat{\mathfrak{A}}_\alpha$, a cumulative sum of the coefficients, and this normalization; for trigonometric functions, it integrates each mode.
The corresponding cohomology projection in degree zero is $\mat{\mathfrak{P}}_\zz\defeq\mat{\mathds{1}}^\zz(\mat{\ell}^\zz)^\transpose$, the evaluation at $\xi^\zz_1$\codeapi{periodic_primitive_matrix}{(space; normalization)}\codeapi[function]{homotopy}{(PeriodicComplex(space), 1)}\codeapi[function]{cohomology_projection}{(PeriodicComplex(space), m)}.

\subsection{Tensor-product spaces in 2D}\label{sec:matrix_representation_tensor_2d}
\paragraph{Differential and homotopy matrices}
The coefficient operators $\gradienth$ and $\curlh$ from \cref{eqn:discrete_operators} are represented by
\begin{equation}\label{eqn:matrep_differential_2d}
  \mat{G}
  \defeq
  \begin{bmatrix}
    \mat{i}^{0}_{\yy}\otimes\mat{d}_\xx\\
    \mat{d}_\yy\otimes\mat{i}^{0}_{\xx}
  \end{bmatrix},
  \qquad
  \mat{C}
  \defeq
  \begin{bmatrix}
    -\mat{d}_\yy\otimes\mat{i}^{1}_{\xx} & \mat{i}^{1}_{\yy}\otimes\mat{d}_\xx
  \end{bmatrix},
\end{equation}
which satisfy $\mat{C}\mat{G}=-\mat{d}_\yy\otimes\mat{d}_\xx+\mat{d}_\yy\otimes\mat{d}_\xx=\mat{0}$\codeapi{gradient_matrix}{(space)}\codeapi{curl_matrix}{(space)}.
By \cref{lem:pseudo_inverse_curl_logical}, the homotopy operators and the center projection $\mathfrak{P}$ are represented by\codeapi{radial_homotopy_matrix}{(space)}\codeapi{center_projection_matrix}{(space)}
\begin{equation}\label{eqn:discrete_homotopy_matrices}
  \mat{G}^\dagger
  \defeq
  \begin{bmatrix}
    \mat{i}^{0}_{\yy}\otimes\mat{d}^\dagger_\xx & \mat{0}
  \end{bmatrix},
  \qquad
  \mat{C}^\dagger
  \defeq
  \begin{bmatrix}
    \mat{0}\\
    \mat{i}^{1}_{\yy}\otimes\mat{d}^\dagger_\xx
  \end{bmatrix},
  \qquad
  \mat{\mathfrak{P}}
  \defeq
  (\mat{\mathds{1}}^\yy\otimes\mat{\mathds{1}}^\xx)
  \bigl((\mat{w}^\yy)^\transpose\otimes(\mat{e}^{0,\xx}_1)^\transpose\bigr).
\end{equation}
The dagger denotes the homotopy operator, not a Moore--Penrose inverse.
By \cref{eqn:radial_cumulative_sum_matrix}, the homotopy identities in \cref{cor:discrete_homotopy_identities} read
\begin{equation}\label{eqn:discrete_homotopy_matrix_identities}
  \mat{p}=\mat{\mathfrak{P}}\mat{p}+\mat{G}^\dagger\mat{G}\mat{p},
  \qquad
  \mat{a}=\mat{G}\mat{G}^\dagger\mat{a}+\mat{C}^\dagger\mat{C}\mat{a},
  \qquad
  \mat{b}=\mat{C}\mat{C}^\dagger\mat{b},
\end{equation}
for coefficient vectors $\mat{p}$ of $0$-forms continuous at the center, $\mat{a}$ of $1$-forms bounded at the center, and arbitrary $\mat{b}$.

\paragraph{GDOF projections}
The 2D GDOF matrices of \cref{eqn:logical_geometric_dofs} are
\begin{equation}\label{eqn:two_dimensional_gdof_vandermonde_matrices}
  \mat{H}^0
  =\mat{h}^{0}_{\yy}\otimes\mat{h}^{0}_{\xx},
  \qquad
  \mat{H}^1
  =\begin{bmatrix}
    \mat{h}^{0}_{\yy}\otimes\mat{h}^{1}_{\xx}&\mat{0}\\
    \mat{0}&\mat{h}^{1}_{\yy}\otimes\mat{h}^{0}_{\xx}
  \end{bmatrix},
  \qquad
  \mat{H}^2
  =\mat{h}^{1}_{\yy}\otimes\mat{h}^{1}_{\xx}.
\end{equation}
If $\mat{g}^m(F)$ collects the logical GDOFs of the pullback of $F$, see \cref{eqn:physical_geometric_dofs}, the coefficient vector of $\pi^mF$ is
\begin{equation}\label{eqn:tensor_product_gdof_projector_matrix}
  \mat{f}=(\mat{H}^m)^{-1}\mat{g}^m(F),
\end{equation}
where the inverse is applied direction by direction, since $(\mat{A}_\yy\otimes\mat{A}_\xx)^{-1}=(\mat{A}_\yy)^{-1}\otimes(\mat{A}_\xx)^{-1}$\codeapi{geometric_dofs}{(space, f)}\codeapi{solve_histopolation}{(space, g)}\codeapi{commuting_projection}{(space, f)}.
The commutation relations \cref{eqn:physical_commutativity} follow from the 1D relations $\mat{h}^{1}_{\alpha}\mat{d}_\alpha=\mat{d}^\xi_\alpha\mat{h}^{0}_{\alpha}$.

\paragraph{Mass matrices}
The mass matrices $\mat{M}^m$ have the entries
\begin{equation}\label{eqn:tensor_product_mass_matrices}
  \begin{aligned}
    \matcomp{M}^0_{\overline{ij},\overline{i'j'}}
    &=\int_{\logical{\Omega}}\logical{\Bspline}^0_{ij}\logical{\Bspline}^0_{i'j'}\jac\diff^2\logical{x},\\
    \matcomp{M}^1_{\alpha\beta,\overline{ij},\overline{i'j'}}
    &=\int_{\logical{\Omega}}(\logical{\+\Bspline}^1_{\alpha,ij})^\transpose(\jacmat^\transpose\jacmat)^{-1}\logical{\+\Bspline}^1_{\beta,i'j'}\jac\diff^2\logical{x},\\
    \matcomp{M}^2_{\overline{ij},\overline{i'j'}}
    &=\int_{\logical{\Omega}}\logical{\Bspline}^2_{ij}\logical{\Bspline}^2_{i'j'}\jac^{-1}\diff^2\logical{x},
  \end{aligned}
\end{equation}
where $\overline{ij}$ denotes the linear index and $\alpha,\beta\in\{\xx,\yy\}$ index the blocks of $\mat{M}^1$\codeapi{mass_matrix}{(space)}.
Their sparsity pattern is the Kronecker product of the 1D couplings, and they are assembled by quadrature with sum factorization.
Without the metric weights, they reduce to Kronecker products of the 1D mass matrices $\mat{M}^{m}_{\alpha}$, e.g., $\mat{M}^0=\mat{M}^{0}_{\yy}\otimes\mat{M}^{0}_{\xx}$.
On the collapsed disk, a nonconforming tensor-product $1$- or $2$-form can have infinite $\funspace{L}_2$ norm by \cref{lem:conformity_conditions_tensor_product_forms}, so that the exact $\mat{M}^1$ and $\mat{M}^2$ are only formal.
Only the polar mass matrices in \cref{eqn:polar_differential_homotopy_matrices} are used; with Gauss quadrature, which has no node at $\logical{x}=0$, they can be computed from the quadrature-assembled $\mat{M}^m$.

\subsection{Polar spaces in 2D}\label{sec:matrix_representation_polar_spaces}
\paragraph{Truncations}
Let $\mat{O}^{m}_{\xx,k}\defeq\operatorname{diag}(0,\ldots,0,1,\ldots,1)\in\mathbb{R}^{n^{m,\xx}\times n^{m,\xx}}$ with $k$ zeros.
The truncations $\chi^{m,\pr}_\external$ onto the external coefficients, cf.\ \cref{eqn:polar_zeroform_space_odot,eqn:polar_mform_space_odot}, are
\begin{equation}\label{eqn:internal_external_truncation_matrices}
  \mat{X}^{0}_\external
  \defeq\mat{i}^{0}_{\yy}\otimes\mat{O}^{0}_{\xx,\pr+1},
  \qquad
  \mat{X}^{1}_\external
  \defeq\begin{bmatrix}
    \mat{i}^{0}_{\yy}\otimes\mat{O}^{1}_{\xx,\pr}&\mat{0}\\
    \mat{0}&\mat{i}^{1}_{\yy}\otimes\mat{O}^{0}_{\xx,\pr+1}
  \end{bmatrix},
  \qquad
  \mat{X}^{2}_\external
  \defeq\mat{i}^{1}_{\yy}\otimes\mat{O}^{1}_{\xx,\pr},
\end{equation}
and $\mat{X}^m_\internal\defeq\mat{I}^m-\mat{X}^m_\external$.
We write $\mat{I}^m_\external$ for the columns of the identity selected by $\mat{X}^m_\external$, so that $\mat{X}^m_\external=\mat{I}^m_\external(\mat{I}^m_\external)^\transpose$.
The reduced $2$-form projector in \cref{eqn:reduced_twoform_coefficient_projector} acts on the radial index only:
\begin{equation}\label{eqn:reduced_twoform_projector_matrix}
  \reduced{\mat{\Pi}}^{2,\pr}
  =\mat{i}^{1}_{\yy}\otimes\mat{r}^{\pr}_\xx,
  \qquad
  \mat{r}^{\pr}_\xx
  \defeq\mat{O}^{1}_{\xx,\pr+1}
  +\mat{e}^{1,\xx}_{\pr+1}\sum_{i'=1}^{\pr+1}(\mat{e}^{1,\xx}_{i'})^\transpose.
\end{equation}

\paragraph{Extraction matrices}
Each polar space is the column space of an extraction matrix $\mat{E}^m$, whose columns are the tensor-product coefficient vectors of the polar basis functions in \cref{eqn:polar_forms_basis_expansion}\codeapi{extraction}{(complex, m)}.
Let $\mat{Z}^0$ and $\mat{Z}^2$ collect the coefficient vectors \cref{eqn:tensor_product_gdof_projector_matrix} of the GDOF projections of the pulled-back shifted Zernike polynomials $T^0_l$ and $T^2_l$ of degree at most $\pr$ and $\pr-2$, respectively.
Their GDOFs are those of the logical pullbacks $\logical{\pushforward}^{-m}\intermediate{T}_l$ of the intermediate polynomials (\cref{sec:derham_hodge_2d_basis}), so $\mat{Z}^0$, $\mat{Z}^2$, and hence all extraction matrices below, depend only on the logical mapping $\logicalmapping$ and not on the intermediate mapping $\intermediatemapping$.
Then \cref{eqn:polar_zeroform_basis,eqn:polar_oneform_gradient_basis,eqn:polar_twoform_curl_basis,eqn:polar_oneform_curl_basis_definition} give
\begin{subequations}\label{eqn:polar_extraction_matrices}
  \begin{equation}
    \mat{E}^0
    =
    \begin{bmatrix}\mat{E}^0_\mathrm{G}&\mat{I}^0_\external\end{bmatrix}
    , \quad
    \mat{E}^1
    =
    \begin{bmatrix}\mat{G}\mat{E}^0_{\mathrm{G},0}&\mat{C}^\dagger\mat{E}^2_\mathrm{C}&\mat{I}^1_\external\end{bmatrix}
    , \quad
    \mat{E}^2
    =
    \begin{bmatrix}\mat{E}^2_\mathrm{C}&\mat{I}^2_\external\end{bmatrix}
    ,
  \end{equation}
  with
  \begin{equation}
    \mat{E}^0_\mathrm{G}
    =
    \mat{X}^0_\internal\mat{Z}^0
    , \quad
    \mat{E}^2_\mathrm{C}
    =
    (\mat{I}^2-\reduced{\mat{\Pi}}^{2,\pr})\reduced{\mat{\Pi}}^{2,1}\mat{Z}^2
    ,
  \end{equation}
\end{subequations}
where $\mat{E}^0_{\mathrm{G},0}$ omits the first column of $\mat{E}^0_\mathrm{G}$, and where $\reduced{\mat{\Pi}}^{2,1}\mat{Z}^2$ represents $\reduced{\pi}^{2,1}T^2_l$ by \cref{lem:conforming_polar_twoform_projector}.
The columns of $\mat{C}^\dagger\mat{E}^2_\mathrm{C}$ are compactly supported by \cref{eqn:polar_twoform_compact_coefficients}, and $\mat{C}\mat{C}^\dagger\mat{E}^2_\mathrm{C}=\mat{E}^2_\mathrm{C}$ by \cref{eqn:discrete_homotopy_matrix_identities}.
Homogeneous essential boundary conditions are imposed by deleting the columns of $\mat{I}^0_\external$ and $\mat{I}^1_\external$ that belong to the outer radial ring of the $0$-forms and of the angular $1$-form component, which gives $\mat{E}^m_\bc$\codeapi{PolarComplex}{(basis, space; boundary=:dirichlet)}.

\paragraph{Polar matrices}
The polar spaces form a cochain complex and are preserved by the homotopy operators (\cref{thm:polar_subspaces_commutativity_conformity,thm:polar_subspaces_homotopy_splitting}), and the extraction matrices have full column rank (\cref{ass:polar_bases}).
Hence there are unique polar differential and homotopy matrices with\codeapi{derivative}{(complex, m)}\codeapi{homotopy}{(complex, m)}\codeapi{mass_matrix}{(complex, m)}
\begin{equation}\label{eqn:polar_differential_homotopy_matrices}
  \mat{G}\mat{E}^0=\mat{E}^1\polar{\mat{G}},
  \quad
  \mat{C}\mat{E}^1=\mat{E}^2\polar{\mat{C}},
  \quad
  \mat{G}^\dagger\mat{E}^1=\mat{E}^0\polar{\mat{G}}^\dagger,
  \quad
  \mat{C}^\dagger\mat{E}^2=\mat{E}^1\polar{\mat{C}}^\dagger,
  \qquad
  \polar{\mat{M}}^m\defeq(\mat{E}^m)^\transpose\mat{M}^m\mat{E}^m.
\end{equation}
These depend only on the extraction matrices and are obtained by multiplying from the left with any left inverse of $\mat{E}^m$, e.g., $\polar{\mat{G}}=(\mat{E}^1)^+\mat{G}\mat{E}^0$; a sparse choice is given by the dual matrices in \cref{eqn:polar_differential_homotopy_dual_formulas} below.
With $\mat{E}^m_\bc$ in place of $\mat{E}^m$, the relations for $\polar{\mat{G}}$ and $\polar{\mat{C}}$ still hold, since the gradient of a $0$-form with zero boundary trace has zero tangential boundary trace.
By \cref{eqn:polar_differential_homotopy_matrices} and the full column rank, \cref{eqn:discrete_homotopy_matrix_identities} carries over to the polar matrices, and in particular $\polar{\mat{C}}\polar{\mat{C}}^\dagger=\mat{I}$.
The polar mass matrix has the entries $\int_\Omega\polar{\Phi}_l\cdot\polar{\Phi}_{l'}\diff^2x$ for the extracted basis functions $\polar{\Phi}_l$, which are finite by \cref{thm:polar_subspaces_conformity}.

\paragraph{Dual matrices and projection matrices}
Every polar coefficient projector in \cref{sec:derham_hodge_2d_geom_projections} is of the form\codeapi{coefficient_projector}{(complex, m)}\codeapi{duals}{(complex, m)}
\begin{equation}\label{eqn:polar_projector_matrices}
  \polar{\mat{\Pi}}^m=\mat{E}^m(\mat{F}^m)^\transpose,
  \qquad
  (\mat{F}^m)^\transpose\mat{E}^m=\mat{I}.
\end{equation}
For $m\in\{0,2\}$, the internal duals \cref{eqn:polar_dual_basis} are collected in matrices $\mat{F}^0_\mathrm{G}$ and $\mat{F}^2_\mathrm{C}$ with internal support and $(\mat{F}^0_\mathrm{G})^\transpose\mat{E}^0_\mathrm{G}=\mat{I}$ and $(\mat{F}^2_\mathrm{C})^\transpose\mat{E}^2_\mathrm{C}=\mat{I}$.
One choice is
\begin{equation}\label{eqn:polar_dual_coefficient_matrices}
  \mat{F}^0_\mathrm{G}
  =\mat{X}^0_\internal\mat{E}^0_\mathrm{G}
    \bigl((\mat{E}^0_\mathrm{G})^\transpose\mat{X}^0_\internal\mat{E}^0_\mathrm{G}\bigr)^{-1},
\end{equation}
and analogously for $\mat{F}^2_\mathrm{C}$; this is the dual basis discussed in \cref{rem:internal_duals}.
The first column of $\mat{F}^0_\mathrm{G}$ may be replaced by $\mat{w}^\yy\otimes\mat{e}^{0,\xx}_1$, the coefficient vector of \cref{eqn:polar_constant_dual_basis}, without affecting biorthogonality, since the first ring of the $l$-th column of $\mat{E}^0_\mathrm{G}$ is $\intermediate{T}_l(\+0)\mat{\mathds{1}}^\yy$.
With
\begin{equation}
  \mat{F}^0=\begin{bmatrix}\mat{F}^0_\mathrm{G}&\mat{I}^0_\external\end{bmatrix},
  \qquad
  \mat{F}^2=\begin{bmatrix}\mat{F}^2_\mathrm{C}&(\reduced{\mat{\Pi}}^{2,\pr})^\transpose\mat{I}^2_\external\end{bmatrix},
\end{equation}
\cref{eqn:polar_projector_matrices} reproduces \cref{eqn:polar_zero_two_form_projection_operator}:
\begin{equation}\label{eqn:polar_zero_two_form_projector_matrices}
  \polar{\mat{\Pi}}^0
  =\mat{E}^0_\mathrm{G}(\mat{F}^0_\mathrm{G})^\transpose+\mat{X}^0_\external,
  \qquad
  \polar{\mat{\Pi}}^2
  =\mat{E}^2_\mathrm{C}(\mat{F}^2_\mathrm{C})^\transpose+\reduced{\mat{\Pi}}^{2,\pr}.
\end{equation}
Here biorthogonality uses $\reduced{\mat{\Pi}}^{2,\pr}\mat{E}^2_\mathrm{C}=\mat{0}$, which expresses the zero radial sum of the compactly supported $2$-form basis functions.
The $1$-form projector \cref{eqn:polar_oneform_projection_operator} is
\begin{equation}\label{eqn:polar_oneform_projector_homotopy_matrix}
  \polar{\mat{\Pi}}^1
  =\mat{G}\polar{\mat{\Pi}}^0\mat{G}^\dagger
   +\mat{C}^\dagger\polar{\mat{\Pi}}^2\mat{C},
\end{equation}
and its dual matrix is $(\mat{F}^1)^\transpose=(\mat{E}^1)^+\polar{\mat{\Pi}}^1$, where $(\mat{E}^1)^+$ is any left inverse of $\mat{E}^1$; this is well defined since the range of $\polar{\mat{\Pi}}^1$ is the column space of $\mat{E}^1$.
Since $(\mat{F}^m)^\transpose$ is a left inverse of $\mat{E}^m$, the matrices in \cref{eqn:polar_differential_homotopy_matrices} are
\begin{equation}\label{eqn:polar_differential_homotopy_dual_formulas}
  \polar{\mat{G}}=(\mat{F}^1)^\transpose\mat{G}\mat{E}^0,
  \quad
  \polar{\mat{C}}=(\mat{F}^2)^\transpose\mat{C}\mat{E}^1,
  \quad
  \polar{\mat{G}}^\dagger=(\mat{F}^0)^\transpose\mat{G}^\dagger\mat{E}^1,
  \quad
  \polar{\mat{C}}^\dagger=(\mat{F}^1)^\transpose\mat{C}^\dagger\mat{E}^2.
\end{equation}
With essential boundary conditions, $\mat{F}^m_\bc$ denotes the columns of $\mat{F}^m$ that belong to the columns kept in $\mat{E}^m_\bc$, so that $(\mat{F}^m_\bc)^\transpose\mat{E}^m_\bc=\mat{I}$.
In this form, \cref{lem:polar_oneform_projection_operator_internal_coefficients} reads $\polar{\mat{\Pi}}^1\mat{X}^1_\external=\mat{X}^1_\external$, and \cref{thm:full_commuting_polar_projection_operators} reads $\polar{\mat{\Pi}}^1\mat{G}\mat{p}=\mat{G}\polar{\mat{\Pi}}^0\mat{p}$ for $\mat{p}$ continuous at the center and $\mat{C}\polar{\mat{\Pi}}^1=\polar{\mat{\Pi}}^2\mat{C}$.
Finally, the polar projection operators \cref{eqn:polar_projection_operator_definition} are computed as\codeapi{commuting_projection}{(complex, m, f)}
\begin{equation}\label{eqn:polar_gdof_projector_matrix_composition}
  \operatorname{coeff}(\polar{\pi}^{m,\pr}F)
  =\mat{E}^m(\mat{F}^m)^\transpose(\mat{H}^m)^{-1}\mat{g}^m(F).
\end{equation}
These projectors are oblique, i.e., they are generally not orthogonal with respect to the mass matrices.

\paragraph{Boundary-compatible homotopy operators}
The traces on the outer boundary, the extension operator $\mathcal{E}_\bc$, and the boundary mean are represented by
\begin{equation}\label{eqn:boundary_trace_extension_matrices}
  \mat{T}^0_\bc
  \defeq\mat{i}^{0}_{\yy}\otimes(\mat{e}^{0,\xx}_{n^{0,\xx}})^\transpose,
  \quad
  \mat{T}^1_\bc
  \defeq\begin{bmatrix}\mat{0}&\mat{i}^{1}_{\yy}\otimes(\mat{e}^{0,\xx}_{n^{0,\xx}})^\transpose\end{bmatrix},
  \quad
  \mat{U}_\bc
  \defeq\mat{i}^{0}_{\yy}\otimes\mat{e}^{0,\xx}_{n^{0,\xx}},
  \quad
  \mat{\mu}_\bc
  \defeq\mat{w}^\yy\otimes\mat{e}^{0,\xx}_{n^{0,\xx}}.
\end{equation}
The projection $\mathfrak{B}$ in \cref{eqn:boundary_homotopy_operators_2d} has rank one.
In polar coefficients, with $\mat{\iota}\defeq(\mat{E}^2)^\transpose(\mat{q}^{1,\yy}\otimes\mat{q}^{1,\xx})$ the vector of integrals of the polar $2$-form basis functions, it reads
\begin{equation}\label{eqn:boundary_mean_projector_matrix}
  \polar{\mat{\mathfrak{B}}}
  =\frac{(\polar{\mat{M}}^2)^{-1}\mat{\iota}\,\mat{\iota}^\transpose}{\mat{\iota}^\transpose(\polar{\mat{M}}^2)^{-1}\mat{\iota}},
\end{equation}
since $(\polar{\mat{M}}^2)^{-1}\mat{\iota}$ represents the $\funspace{L}_2$-projection of $1$ onto the polar $2$-form space.
With $\polar{\mat{C}}_\bc\defeq(\mat{F}^2)^\transpose\mat{C}\mat{E}^1_\bc$, the boundary-compatible homotopy operators in \cref{eqn:boundary_homotopy_operators_2d} are
\begin{subequations}\label{eqn:boundary_homotopy_matrices}
  \begin{align+}
    \polar{\mat{C}}^\dagger_\bc
    &=(\mat{F}^1_\bc)^\transpose
      \bigl(\mat{I}^1-\mat{G}\mat{U}_\bc\mat{d}^\dagger_\yy\mat{T}^1_\bc\bigr)
      \mat{C}^\dagger\mat{E}^2
      \bigl(\mat{I}-\polar{\mat{\mathfrak{B}}}\bigr),\\
    \polar{\mat{G}}^\dagger_\bc
    &=(\mat{F}^0_\bc)^\transpose
      \bigl(\mat{I}^0-(\mat{\mathds{1}}^\yy\otimes\mat{\mathds{1}}^\xx)\mat{\mu}_\bc^\transpose\bigr)
      \mat{G}^\dagger\mat{E}^1_\bc
      \bigl(\mat{I}-\polar{\mat{C}}^\dagger_\bc\polar{\mat{C}}_\bc\bigr),
  \end{align+}
\end{subequations}
where $\mat{d}^\dagger_\yy$ is the zero-mean periodic primitive from \cref{eqn:periodic_primitive_matrices}, and $\mat{\mu}_\bc^\transpose$ returns the constant boundary value\codeapi{homotopy}{(complex, m)}\codeapi{cohomology_projection}{(complex, 2)}.

\paragraph{Linear systems}
In polar coefficients $\polar{\mat{u}}$, with $\mat{u}=\mat{E}^m\polar{\mat{u}}$, the $\funspace{L}_2$ projection onto the polar $m$-form space reads $\polar{\mat{M}}^m\polar{\mat{u}}=(\mat{E}^m)^\transpose\mat{f}$, where $\mat{f}$ is the tensor-product load vector\codeapi{l2_projection}{(complex, m, f)}.
A stiffness matrix of the form $\polar{\mat{C}}^\transpose\polar{\mat{M}}^2\polar{\mat{C}}$ or $\polar{\mat{G}}^\transpose\polar{\mat{M}}^1\polar{\mat{G}}$ is assembled from the polar matrices in \cref{eqn:polar_differential_homotopy_matrices}, with $\mat{E}^m_\bc$ in place of $\mat{E}^m$ for essential boundary conditions.
For a weakly divergence-free load vector $\polar{\mat{\ell}}=(\mat{E}^1_\bc)^\transpose\mat{f}$, the particular solution of the curl--curl problem in \cref{sec:bvp_2d} is
\begin{equation}\label{eqn:polar_curlcurl_matrix}
  \polar{\mat{a}}
  =\polar{\mat{C}}^\dagger_\bc(\polar{\mat{M}}^2)^{-1}(\polar{\mat{C}}^\dagger_\bc)^\transpose\polar{\mat{\ell}},
\end{equation}
since $(\polar{\mat{M}}^2)^{-1}(\polar{\mat{C}}^\dagger_\bc)^\transpose\polar{\mat{M}}^1_\bc$, with $\polar{\mat{M}}^1_\bc\defeq(\mat{E}^1_\bc)^\transpose\mat{M}^1\mat{E}^1_\bc$, represents the $\funspace{L}_2$-adjoint of $\curlt^\dagger_\bc$\codeapi{stiffness_solver}{(complex, 1; gauge=:homotopy)}\codeapi{homotopy_adjoint}{(complex, m, coefficients)}.
Thus only the polar $2$-form mass matrix is inverted.

\subsection{Polar spaces in 3D}\label{sec:matrix_representation_3d}
\paragraph{Bidegree blocks}
By \cref{eqn:tensor_product_complex_spaces_3d},\codeapi{PolarComplex3D}{(cross_section, toroidal, mapping)}\codeapi{TensorProductComplex}{(cross_section, PeriodicComplex(toroidal))} the coefficient vector of a 3D $m$-form consists of blocks of bidegree $(p,q)$ with $p+q=m$, ordered by decreasing $p$\codeapi{bidegree_blocks}{(complex, m)}.
Since the toroidal index runs slowest, the coefficient vector of a product $\logical{V}^p\wedge U^q$ is $\mat{u}\otimes\mat{v}$, where $\mat{v}$ and $\mat{u}$ are the coefficient vectors of $\logical{V}^p$ and $U^q$.
The vector-proxy components are obtained by a signed permutation $\mat{L}^m$ induced by the matrices $\mat{\Lambda}^{p,q}$ in \cref{eqn:tensor_product_proxy_matrices_3d}\codeapi{proxy_permutation}{(complex, m)}.
It is the identity for $m\in\{0,3\}$ and an unsigned permutation for $m=1$, whereas for $m=2$ it maps the blocks $(\mat{b}_{(2,0)},\mat{a}_{(1,1)})$ with $\mat{a}_{(1,1)}=[\mat{a}_\xx;\mat{a}_\yy]$ to the components $(\mat{a}_\yy,-\mat{a}_\xx,\mat{b}_{(2,0)})$.
Below, all matrices are given in the bidegree-block ordering; the proxy representation of an operator $\mat{A}$ from degree $m$ to $m'$ is $\mat{L}^{m'}\mat{A}(\mat{L}^m)^\transpose$.
If $\mat{A}_{(p,q)}$ is a matrix for each bidegree $(p,q)$ with $p+q=m$, mapping the $(p,q)$-block of the input to the $(p,q)$-block of the output, then their direct sum is the block-diagonal matrix
\begin{equation}\label{eqn:bidegree_direct_sum}
  \bigoplus_{p+q=m}\mat{A}_{(p,q)}
  \defeq
  \begin{bmatrix}
    \mat{A}_{(m,0)}&\mat{0}\\
    \mat{0}&\mat{A}_{(m-1,1)}
  \end{bmatrix},
\end{equation}
with the blocks in the bidegree ordering and the convention that blocks with $p\notin\{0,1,2\}$ or $q\notin\{0,1\}$ are absent.
That is, the direct sum is $\mat{A}_{(0,0)}$ for $m=0$, $\operatorname{diag}(\mat{A}_{(1,0)},\mat{A}_{(0,1)})$ for $m=1$, $\operatorname{diag}(\mat{A}_{(2,0)},\mat{A}_{(1,1)})$ for $m=2$, and $\mat{A}_{(2,1)}$ for $m=3$.
The blocks need not be square; for instance, the blocks of an extraction matrix map polar coefficients to tensor-product coefficients.

\paragraph{Extraction and projection matrices}
Following the notation of \cref{sec:high_order_derham_hodge_3d}, the 2D polar matrices from \cref{sec:matrix_representation_polar_spaces} carry the subscript $\mathrm{2D}$, e.g., $\mat{E}^p_\mathrm{2D}$, $\polar{\mat{G}}_\mathrm{2D}$ and $\polar{\mat{\mathfrak{P}}}_\mathrm{2D}\defeq(\mat{F}^0_\mathrm{2D})^\transpose\mat{\mathfrak{P}}_\mathrm{2D}\mat{E}^0_\mathrm{2D}$, while the 3D matrices carry no subscript.
We also write $\polar{\mat{K}}^0_\mathrm{2D}\defeq\polar{\mat{G}}_\mathrm{2D}$ and $\polar{\mat{K}}^1_\mathrm{2D}\defeq\polar{\mat{C}}_\mathrm{2D}$ for the 2D differentials, and $\polar{\mat{K}}^{1,\dagger}_\mathrm{2D}\defeq\polar{\mat{G}}^\dagger_\mathrm{2D}$ and $\polar{\mat{K}}^{2,\dagger}_\mathrm{2D}\defeq\polar{\mat{C}}^\dagger_\mathrm{2D}$ for the 2D homotopies.
The 3D extraction, dual, projection, and GDOF matrices are block diagonal with Kronecker blocks\codeapi{extraction}{(complex, m)}\codeapi{duals}{(complex, m)}:
\begin{equation}\label{eqn:matrep_extraction_3d}
  \begin{aligned}
    \mat{E}^m
    &=\bigoplus_{p+q=m}\mat{i}^{q}_{\zz}\otimes\mat{E}^p_\mathrm{2D},
    &
    \mat{F}^m
    &=\bigoplus_{p+q=m}\mat{i}^{q}_{\zz}\otimes\mat{F}^p_\mathrm{2D},\\
    \polar{\mat{\Pi}}^m
    &=\bigoplus_{p+q=m}\mat{i}^{q}_{\zz}\otimes\polar{\mat{\Pi}}^p_\mathrm{2D},
    &
    \mat{H}^m
    &=\bigoplus_{p+q=m}\mat{h}^{q}_{\zz}\otimes\mat{H}^p_\mathrm{2D}.
  \end{aligned}
\end{equation}
Here $\bigoplus$ is the direct sum over the bidegree blocks in \cref{eqn:bidegree_direct_sum}.
The third identity is the matrix form of \cref{eqn:logical_3d_polar_projection_operator}, and the polar GDOF projection \cref{eqn:polar_3d_projection_operator_definition} is computed as in \cref{eqn:polar_gdof_projector_matrix_composition}\codeapi{commuting_projection}{(complex, m, f)}.

\paragraph{Differential and homotopy matrices}
By \cref{eqn:logical_3d_tensor_product_differential}, the block from $(p,q)$ to $(p+1,q)$ of the 3D differential is $\mat{i}^{q}_{\zz}\otimes\polar{\mat{K}}^p_\mathrm{2D}$, and the block from $(p,q)$ to $(p,q+1)$ is $(-1)^p\mat{d}_\zz\otimes\polar{\mat{I}}^p_\mathrm{2D}$.
Explicitly, the 3D differential matrices are\codeapi{derivative}{(complex, m)}
\begin{equation}\label{eqn:matrep_differential_3d}
  \polar{\mat{G}}
  =\begin{bmatrix}
    \mat{i}^{0}_{\zz}\otimes\polar{\mat{G}}_\mathrm{2D}\\
    \mat{d}_\zz\otimes\polar{\mat{I}}^0_\mathrm{2D}
  \end{bmatrix},
  \quad
  \polar{\mat{C}}
  =\begin{bmatrix}
    \mat{i}^{0}_{\zz}\otimes\polar{\mat{C}}_\mathrm{2D}&\mat{0}\\
    -\mat{d}_\zz\otimes\polar{\mat{I}}^1_\mathrm{2D}&\mat{i}^{1}_{\zz}\otimes\polar{\mat{G}}_\mathrm{2D}
  \end{bmatrix},
  \quad
  \polar{\mat{D}}
  =\begin{bmatrix}
    \mat{d}_\zz\otimes\polar{\mat{I}}^2_\mathrm{2D}&\mat{i}^{1}_{\zz}\otimes\polar{\mat{C}}_\mathrm{2D}
  \end{bmatrix}.
\end{equation}
The mixed terms cancel in $\polar{\mat{C}}\polar{\mat{G}}$ and $\polar{\mat{D}}\polar{\mat{C}}$ because of the sign $(-1)^p$, and $\polar{\mat{C}}_\mathrm{2D}\polar{\mat{G}}_\mathrm{2D}=\mat{0}$ takes care of the remaining terms.
By \cref{eqn:logical_3d_tensor_product_homotopy}, the block from $(p,q)$ to $(p-1,q)$ of the 3D homotopy is $\mat{i}^{q}_{\zz}\otimes\polar{\mat{K}}^{p,\dagger}_\mathrm{2D}$, and the block from $(p,q)$ to $(p,q-1)$ is $(-1)^p\mat{d}^\dagger_\zz\otimes\polar{\mat{\mathfrak{P}}}^p_\mathrm{2D}$, where $\polar{\mat{\mathfrak{P}}}^0_\mathrm{2D}=\polar{\mat{\mathfrak{P}}}_\mathrm{2D}$ and $\polar{\mat{\mathfrak{P}}}^1_\mathrm{2D}=\polar{\mat{\mathfrak{P}}}^2_\mathrm{2D}=\mat{0}$.
Hence \cref{eqn:logical_3d_homotopy_vector_proxies} reads
\begin{equation}\label{eqn:matrep_homotopy_3d}
  \polar{\mat{G}}^\dagger
  =\begin{bmatrix}
    \mat{i}^{0}_{\zz}\otimes\polar{\mat{G}}^\dagger_\mathrm{2D}&\mat{d}^\dagger_\zz\otimes\polar{\mat{\mathfrak{P}}}_\mathrm{2D}
  \end{bmatrix},
  \quad
  \polar{\mat{C}}^\dagger
  =\begin{bmatrix}
    \mat{i}^{0}_{\zz}\otimes\polar{\mat{C}}^\dagger_\mathrm{2D}&\mat{0}\\
    \mat{0}&\mat{i}^{1}_{\zz}\otimes\polar{\mat{G}}^\dagger_\mathrm{2D}
  \end{bmatrix},
  \quad
  \polar{\mat{D}}^\dagger
  =\begin{bmatrix}
    \mat{0}\\
    \mat{i}^{1}_{\zz}\otimes\polar{\mat{C}}^\dagger_\mathrm{2D}
  \end{bmatrix},
\end{equation}
and the cohomology projections in \cref{eqn:logical_3d_cohomology_projector_vector_proxies} are $\mat{\mathfrak{P}}_\zz\otimes\polar{\mat{\mathfrak{P}}}_\mathrm{2D}$ in degree zero and $\mat{\mathfrak{A}}_\zz\otimes\polar{\mat{\mathfrak{P}}}_\mathrm{2D}$ on the $(0,1)$ block in degree one\codeapi[function]{homotopy}{(complex, m)}\codeapi[function]{cohomology_projection}{(complex, m)}.

\paragraph{Mass matrices}
The 3D mass matrices $\mat{M}^m$ are assembled by quadrature with the metric weights of \cref{eqn:pushforward_operators_3d}, analogously to \cref{eqn:tensor_product_mass_matrices}, and the polar mass matrices are $\polar{\mat{M}}^m=(\mat{E}^m)^\transpose\mat{M}^m\mat{E}^m$\codeapi{mass_matrix}{(complex, m)}.
\endgroup

  \bibliographystyle{cas-model2-names}
  \bibliography{bib}

\end{document}